\pdfoutput=1
\documentclass[11pt]{article}
\usepackage{microtype}
\usepackage{rotating}
\usepackage[usenames,dvipsnames]{xcolor}

\usepackage{listings}
\usepackage{graphicx}
\usepackage{tikz}
\usetikzlibrary{positioning,backgrounds,fit,calc,decorations.pathreplacing}
\usepackage{pgfplots}
\pgfplotsset{compat=1.18}
\usepgfplotslibrary{groupplots} % multi-panel solution-time breakdown figures
\usepackage{amsmath,amssymb,amsthm,amsfonts} % assumes amsmath package installed
\usepackage[linktocpage=true,colorlinks=true,linkcolor=blue,citecolor=blue,urlcolor=blue]{hyperref}
\usepackage[letterpaper,margin=1in]{geometry}
\usepackage{DejaVuSans}

\usepackage[acronym,nomain,nopostdot]{glossaries}
\newacronym{nlp}{NLP}{nonlinear program}
\newacronym{jso}{JSO}{JuliaSmoothOptimizers}
\newacronym{ams}{AMS}{algebraic modeling system}
\newacronym{ad}{AD}{automatic differentiation}
\newacronym{simd}{SIMD}{single-instruction, multiple-data}
\newacronym{gpu}{GPU}{graphics processing unit}
\newacronym{cpu}{CPU}{central processing unit}
\newacronym{acopf}{ACOPF}{alternating current optimal power flow}
\newacronym{opf}{OPF}{optimal power flow}
\newacronym{aot}{AOT}{ahead-of-time}
\newacronym{jit}{JIT}{just-in-time}
\newacronym{coo}{COO}{coordinate}
\newacronym{csc}{CSC}{compressed sparse column}
\newacronym{asl}{ASL}{AMPL Solver Library}
\newacronym{pde}{PDE}{partial differential equation}
\newacronym{hvp}{HVP}{Hessian--vector product}
\newacronym{sgm}{SGM}{shifted geometric mean}
\newacronym{api}{API}{application programming interface}
\newacronym{milp}{MILP}{mixed-integer linear program}
\newacronym{lp}{LP}{linear program}
\newacronym{dcp}{DCP}{disciplined convex program}
\newacronym{ocp}{OCP}{optimal control problem}
\glsdisablehyper % no glossary hyperlinks -> acronyms are not highlighted/colored

\newcommand{\st}{\mathop{\text{\normalfont s.t.}}}

\newcommand{\codetheta}{\ensuremath{\theta}}

\definecolor{jlkw}{HTML}{7B3294}      % keywords: purple
\definecolor{jltype}{HTML}{1F78B4}    % types: blue
\definecolor{jlmacro}{HTML}{E31A1C}   % macros: red
\definecolor{jlstr}{HTML}{A6761D}     % strings: brown
\definecolor{jlcmt}{HTML}{636363}     % comments: grey
\definecolor{jlnum}{HTML}{1B7837}     % numbers: dark green
\definecolor{jlbg}{HTML}{F7F7F7}      % background: light grey

\lstdefinelanguage{Julia}{%
  morekeywords={abstract,break,case,catch,const,continue,do,else,elseif,%
    end,export,false,for,function,global,if,import,in,let,local,%
    macro,module,mutable,nothing,quote,return,struct,true,try,%
    type,typealias,using,where,while},
  morekeywords=[2]{Int,Int64,Int32,Float64,Float32,Bool,String,Symbol,%
    Vector,Matrix,Array,Dict,Tuple,NamedTuple,Any,Nothing,Type,%
    AbstractFloat,AbstractArray,AbstractVector,AbstractMatrix,%
    UnitRange,StepRange},
  morekeywords=[3]{@add_var,@add_par,@add_obj,@add_con,@add_con!,%
    @add_expr,@register_univariate,@register_bivariate,%
    @simd,@inbounds,@kernel,@index,@sprintf,@info,@eval},
  sensitive=true,
  alsoletter={0123456789},
  alsoother={\$},
  morecomment=[l]\#,
  morecomment=[n]{\#=}{=\#},
  morestring=[s]{"}{"},
  morestring=[m]{'}{'},
}[keywords,comments,strings]

\makeatletter
\newcommand*\jlidstyle{\expandafter\jl@idstyle\the\lst@token\relax}
\def\jl@idstyle#1#2\relax{%
  \ifcat\noexpand#1\relax\else
    \ifnum`#1>47 \ifnum`#1<58 \color{jlnum}\fi\fi
  \fi}
\makeatother

\providecommand{\keywords}[1]{%
  \small\textbf{\textit{Keywords:}} #1%
}

\definecolor{j1}{rgb}{0.0,0.6056031611752245,0.9786801175696073}
\definecolor{j2}{rgb}{0.8888735002725198,0.43564919034818994,0.2781229361419438}
\definecolor{j3}{rgb}{0.2422242978521988,0.6432750931576305,0.3044486515341153}
\definecolor{j4}{rgb}{0.7644401754934356,0.4441117794687767,0.8242975359232758}
\definecolor{j5}{rgb}{0.6755439572114057,0.5556623322045815,0.09423433626639477}

\usepackage{algorithm}
\usepackage{algpseudocode}
\usepackage{parskip}
\usepackage[font=small,labelfont=bf]{caption}
\usepackage{booktabs}
\usepackage{longtable}
\usepackage[style=authoryear,backend=biber,natbib=true,doi=true,url=true,eprint=true,isbn=false,maxbibnames=999]{biblatex}
\usepackage[capitalise,noabbrev]{cleveref}
\crefformat{equation}{(#2#1#3)}
\Crefformat{equation}{(#2#1#3)}
\crefrangeformat{equation}{(#3#1#4)--(#5#2#6)}
\crefmultiformat{equation}{(#2#1#3)}{ and (#2#1#3)}{, (#2#1#3)}{ and (#2#1#3)}
\Crefformat{equation}{Equation~(#2#1#3)}
\Crefrangeformat{equation}{Equations~(#3#1#4)--(#5#2#6)}
\Crefmultiformat{equation}{Equations~(#2#1#3)}{ and (#2#1#3)}{, (#2#1#3)}{ and (#2#1#3)}

\usepackage{placeins}

\lstdefinelanguage{LLVM}{%
  morekeywords     = {define, ret, load, store, br, call, phi, getelementptr,
                      icmp, fmul, fadd, fsub, fdiv, shl, add, sub, sitofp, to, label},
  morekeywords     = [2]{ptr, double, float, i64, i32, i8, i1, void},
  morekeywords     = [3]{nonnull, noundef, nocapture, readonly, inbounds,
                      dereferenceable, signext, align, nuw, nsw},
  sensitive        = true,
  morecomment      = [l]{;},
}
\lstdefinestyle{llvm}{%
  language         = LLVM,
  texcl            = false,
  literate         = {},
  keywordstyle     = \bfseries\color{jlkw},
  keywordstyle     = [2]\color{jltype},
  keywordstyle     = [3]\color{Gray},
  identifierstyle  = {},
  stringstyle      = {},
  commentstyle     = \itshape\color{jlcmt},
  basicstyle       = \ttfamily\scriptsize,
  numbers          = left,
  numberstyle      = \tiny\ttfamily\color{Gray},
  breaklines       = true,
  breakatwhitespace= true,
}
\lstdefinestyle{jlrepl}{%
  language         = {},
  identifierstyle  = {},
  literate         = {θ}{{\codetheta}}1,
  moredelim        = [s][\bfseries\color{ForestGreen}]{julia}{>},
  morecomment      = [l]\#,
  commentstyle     = \color{jlcmt},
  texcl            = false,
}

\title{ExaModels.jl: an Algebraic~Modeling~System for Nonlinear~Programming on GPUs}
\author{Sungho Shin$^*$, Michel Schanen$^\ddag$, Fran\c{c}ois Pacaud$^{\dag,\ddag}$,\\ Alexis Montoison$^\P$, Mihai Anitescu$^{\ddag,\S}$}
\date{\small
  $^*$Department of Chemical Engineering, Massachusetts Institute of Technology\\
  $^\dag$MINES Paris - PSL\\
  $^\ddag$Mathematics and Computer Science Division, Argonne National Laboratory\\
  $^\P$AMD\\
  $^\S$Department of Statistics, University of Chicago\\
}

\newcommand{\bdRosSpeedupLinalg}{1.7}
\newcommand{\bdRosSpeedupResident}{1.7}

\newcommand{\bdOpfCase}{\texttt{case78484\_epigrids}}

\newcommand{\bdOpfTotCC}{147}

\newcommand{\bdOpfTotGC}{15.4}
\newcommand{\bdOpfTotGG}{7.48}

\newcommand{\bdOpfSpeedupLinalg}{9.6}
\newcommand{\bdOpfSpeedupResident}{2.1}
\newcommand{\bdOpfSpeedupTotal}{19.7}
\newcommand{\bdOpfSpeedupResidentCpuLs}{1.04}

\newcommand{\bdOpfEvalShareCC}{7}

\newcommand{\bdOpfEvalShareGC}{52}

\newcommand{\platAOneName}{Instinct~MI300X}

\newcommand{\platIOneName}{Intel~Data~Center~GPU~Max~1550}

\newcommand{\platMOneName}{M2~Pro}

\newcommand{\platNOneName}{B200}

\newcommand{\platNSixName}{RTX~PRO~6000~Blackwell}

\newcommand{\gsLvHessLarge}{76}

\newcommand{\gsLvHessAbsSmall}{26.9}
\newcommand{\gsLvHessAbsMedium}{28.3}
\newcommand{\gsLvHessAbsLarge}{35.7}
\newcommand{\gsLvJacAbsSmall}{19.2}
\newcommand{\gsLvJacAbsMedium}{19.9}
\newcommand{\gsLvJacAbsLarge}{20.4}

\newcommand{\gsLvSpanLo}{34.8}
\newcommand{\gsLvSpanHi}{76.4}
\newcommand{\gsLvAmdHessLarge}{44.5}
\newcommand{\gsLvIntelHessLarge}{27.4}
\newcommand{\gsLvAppleHessLarge}{23.3}

\newcommand{\gsCopsHessLarge}{30}

\newcommand{\gsCopsHessAbsLarge}{92.9}

\newcommand{\gsOpfPolarHessLarge}{7.3}

\newcommand{\cntObj}{115}
\newcommand{\cntCons}{115}
\newcommand{\cntGrad}{114}
\newcommand{\cntJac}{114}
\newcommand{\cntHess}{113}

\newcommand{\cmpKindsRatioCpu}{7.9}

\newcommand{\cmpSizeSpreadGpu}{0.9}
\newcommand{\cmpKindsRatioGpu}{1.21}
\newcommand{\cmpKindsLoGpu}{1}
\newcommand{\cmpKindsHiGpu}{16}
\newcommand{\cmpGpuFixedS}{16}

\newcommand{\cmpLvSizeA}{20}

\newcommand{\cmpLvJacGpuExaA}{9.6\,\textmu s}
\newcommand{\cmpLvHessGpuExaA}{16.2\,\textmu s}
\newcommand{\cmpLvSizeB}{2{,}000}

\newcommand{\cmpLvJacGpuExaB}{20.8\,\textmu s}
\newcommand{\cmpLvHessGpuExaB}{55.9\,\textmu s}
\newcommand{\cmpLvSizeC}{200{,}000}

\newcommand{\cmpLvJacGpuExaC}{41.6\,\textmu s}
\newcommand{\cmpLvHessGpuExaC}{119.3\,\textmu s}

\newcommand{\cmpOpfPerKernelHessA}{11.0\,\textmu s}
\newcommand{\cmpOpfPerKernelHessB}{10.8\,\textmu s}
\newcommand{\cmpOpfPerKernelHessC}{21.9\,\textmu s}
\newcommand{\cmpOpfPerKernelJacA}{10.0\,\textmu s}
\newcommand{\cmpOpfPerKernelJacB}{10.6\,\textmu s}
\newcommand{\cmpOpfPerKernelJacC}{11.0\,\textmu s}

\begin{document}
\maketitle
\begin{abstract}
  Large-scale nonlinear programs almost always exhibit partially separable and repetitive structure, yet most existing algebraic modeling systems do not take advantage of it.
  A nonlinear optimization solver queries the objective, the constraints, and their derivatives at every iteration, so the speed of these evaluations bears directly on the overall solution time.
  We present ExaModels.jl, a Julia-based algebraic modeling system that exploits this structure to evaluate the objective, the constraints, and their derivatives in parallel.
  At its core is a \emph{single-instruction, multiple-data abstraction} that represents a nonlinear program as a small number of algebraic patterns, each repeated over many data points.
  Because the patterns are visible at compile time, a specialized model and derivative evaluation kernel is compiled for each pattern.
  Applying that kernel independently across the data points maps naturally onto GPU parallelism and, with sufficiently many threads, yields $O(1)$ evaluation time regardless of the number of data points.
  On the largest instances of the Luk\v{s}an--Vl\v{c}ek library, GPU execution speeds up sparse Hessian evaluation by $\gsLvHessLarge\times$ over single-threaded CPU evaluation, and by $\gsCopsHessLarge\times$ on COPS and $\gsOpfPolarHessLarge\times$ on PGLIB-OPF.
\end{abstract}


\section{Introduction}
\label{sec:intro}

Large-scale \glspl{nlp} arising in engineering, statistics, and scientific applications almost always exhibit partially separable and repetitive structure~\citep{griewankUnconstrainedOptimizationPartially1982}.
The objective and constraints are sums or stacks of scalar terms, and the same algebraic pattern appears across many data points.
To solve an \gls{nlp}, the solver evaluates the NLP functions at each iteration: the objective and constraint values, together with their first- and second-order derivatives (the gradient, the constraint Jacobian, and the Lagrangian Hessian).
The speed of these evaluations therefore governs the performance of the solution procedure.
Exploiting the partially separable, repetitive structure is the key to making them fast.
This paper presents a modeling abstraction, an \gls{ad} algorithm, and their software implementation in ExaModels.jl that exploit this structure for scalable NLP function evaluation and \gls{gpu} acceleration.

To make this concrete, consider the Goddard rocket problem~\citep{goddardMethodReachingExtreme1919} from the COPS benchmark suite~\citep{dolanBenchmarkingOptimizationSoftware2001}, discretized via the trapezoidal rule with $T$ time steps:
\begin{subequations}\label{eq:rocket}
\begin{align}
  \max_{\substack{h_t,\, v_t,\, m_t, \\ \tau_t ,\,\Delta t}} \quad
  & h_T \\
  \st \quad
  & h_t = h_{t-1} + \frac{\Delta t}{2}(v_t + v_{t-1}), & t &=1,\ldots,T, \label{eq:rocket:h}\\
  & v_t = v_{t-1} + \frac{\Delta t}{2}\!\left(\frac{\tau_t - D(h_t,v_t)}{m_t} - g(h_t) + \frac{\tau_{t-1} - D(h_{t-1},v_{t-1})}{m_{t-1}} - g(h_{t-1})\right)\!, & t &=1,\ldots,T, \label{eq:rocket:v}\\
  & m_t = m_{t-1} - \frac{\Delta t}{2c}(\tau_t + \tau_{t-1}), & t &=1,\ldots,T, \label{eq:rocket:m}\\
  & h_0 = h^0,\; v_0 = 0,\; m_0 = m^0,\; m_T = m^f,
\end{align}
\end{subequations}
where $h_t$, $v_t$, $m_t$, and $\tau_t$ are the altitude, velocity, mass, and thrust at time $t$, and $\Delta t$ is the time step, itself a decision variable.
The aerodynamic drag is $D(h,v) = D_c v^2 e^{-h_c(h - h^0)/h^0}$, and the gravitational acceleration is $g(h) = g_0 (h^0/h)^2$.
The given constants $h^0$, $m^0$, $m^f$, $c$, $D_c$, $h_c$, and $g_0$ are the initial altitude, the initial mass, the final mass, the exhaust velocity, the two drag parameters, and the gravitational acceleration at the initial altitude.
The problem aims to maximize the final altitude of a vertically ascending rocket subject to aerodynamic drag and gravity.
This is a classical optimal control problem representative of the discretize-then-optimize approach widely used in trajectory optimization~\citep{bettsPracticalMethodsOptimal2010,bieglerNonlinearProgrammingConcepts2010}.

The NLP function evaluations of this problem can be naturally represented as \emph{iterators}.
For example, consider the altitude dynamics \cref{eq:rocket:h}.
The algebraic expression $h_t - h_{t-1} - \tfrac{\Delta t}{2}(v_t + v_{t-1}) = 0$ has the same functional form for every $t = 1, \ldots, T$, and only the variable indices change.
The same holds for the velocity dynamics \cref{eq:rocket:v} and the mass dynamics \cref{eq:rocket:m}.
The repetition applies not only to the functions themselves but also to their derivatives.
The partial derivatives of each pattern with respect to its variables again share a single functional form across $t$, with only the variable indices changing.

In ExaModels.jl, the modeling system developed in this paper, the constraint is written in the same form, as an algebraic pattern paired with an iterator:
\begin{lstlisting}
@add_con(core, h[t] - h[t-1] - 0.5 * dt[1] * (v[t] + v[t-1]) for t = 1:T)
\end{lstlisting}
where \texttt{core} is the model object under construction, and \texttt{h}, \texttt{v}, and \texttt{dt} are variable blocks declared beforehand.
Because the data iterator is part of the modeling syntax, the model object receives the algebraic pattern and the data points separately.
The repetitive structure is therefore visible without any analysis: each constraint type defines a single algebraic pattern, and the time indices $t = 1, \ldots, T$ are the data points over which that pattern is evaluated.
Evaluating the constraint function or its derivatives then reduces to the construction of a \texttt{for} loop over $t$ for each pattern, and these loops are embarrassingly parallel across data points.
A consequence is that when the discretization is refined, the problem size grows proportionally in $T$, but the number of distinct algebraic patterns remains fixed at five: the objective, the three dynamics constraints \cref{eq:rocket:h}--\cref{eq:rocket:m}, and one singleton pattern that collects the four initial and final conditions.
In particular, with sufficiently many parallel threads, the evaluation time of the NLP functions remains constant as $T$ grows.
By contrast, the constraints could be written out one by one, calling \texttt{@add\_con} once per $t$ inside an ordinary loop.
Each call would create a separate scalar constraint, and unless a pattern-matching procedure were applied internally, the structured, repetitive constraints would be unrolled into a flat list.
The benefits described above would then become much more difficult to realize.

This partially separable, repetitive algebraic structure is not specific to the Goddard rocket problem in \cref{eq:rocket}.
Every problem category of the COPS benchmark, optimal control, PDE-constrained optimization, parameter estimation with collocation, and mesh and shape optimization, carries a small, fixed number of algebraic patterns, with the problem size scaling in the discretization resolution~\citep{dolanBenchmarkingOptimizationSoftware2001}.
The \gls{acopf} problem~\citep{carpentierContributionLetudeDispatching1962}, one of the most widely studied \glspl{nlp} in power systems, can be expressed with only 15 distinct algebraic patterns, regardless of whether the network contains hundreds or tens of thousands of buses~\citep{shinAcceleratingOptimalPower2024}.
More broadly, large-scale \glspl{nlp} across application domains, such as \gls{pde}-constrained optimization~\citep{bieglerLargeScalePDEConstrainedOptimization2003}, optimal control, network design, and maximum likelihood estimation, routinely involve millions of variables and constraints.
Their mathematical formulations still fit on at most a few pages, because the same algebraic expressions are repeated over spatial nodes, time steps, scenarios, or data samples.
\emph{The problem size is driven by the number of data points, not the number of distinct algebraic patterns}.

A modeling system must therefore let the user write the model in a form that preserves the pattern--data separation.
When this separation is maintained, each algebraic pattern needs to be analyzed only once, and the specialized model and derivative evaluation kernels compiled for it can then be applied across all data points.
The sparsity analysis is likewise performed per pattern: the derivative sparsity is determined once for the pattern and replicated across its data points.
Furthermore, because the same instruction, the pattern's kernel, is applied independently to many data points, the evaluation maps naturally onto the \gls{simd} paradigm of parallel computing.
With a sufficient number of threads, the NLP functions can in principle be evaluated in $O(1)$ time regardless of the number of data points, effectively removing the \gls{ad} computation from the overall \gls{nlp} solution time.
Such parallel evaluation also allows the NLP functions to be evaluated directly on \glspl{gpu}.
Moreover, with the recent progress in nonlinear optimization solvers on \glspl{gpu}~\citep{shinMadNLPjl2025,pacaudCondensedspaceMethodsNonlinear2024}, this enables a fully \gls{gpu}-resident solution framework in which the \gls{ad}, the linear algebra, and the solver all remain in device memory, eliminating host--device data transfer overhead.
This observation motivates the following question:
\begin{quote}
  \emph{What modeling abstraction preserves the pattern--data separation of large-scale \glspl{nlp}, and how can it be harnessed for efficient, parallel evaluation of the NLP functions?}
\end{quote}

In answer to this question, we introduce a \emph{\gls{simd} abstraction} for \glspl{nlp}, which represents the problem as a small number of algebraic patterns, each paired with an iterator over its data points, together with an \gls{ad} scheme specialized to this structured algebraic information.
The abstraction is slightly more restrictive than general-purpose modeling interfaces: the user is \emph{required} to express the model in a structured form, where each objective term and constraint is specified as an algebraic pattern paired with a data iterator.
This is a deliberate design choice: by asking the user to expose the parallelizable structure of their problem, the modeling system can exploit it for efficient derivative evaluation.
Each pattern is represented by a \emph{parameterized expression tree}: a tree of operations that encodes the algebraic pattern and is parameterized by the problem data supplied at evaluation time.
Because each pattern is fixed and only the data varies, a single derivative evaluation kernel per pattern suffices, executed across all of its data points in a \gls{simd} fashion.
This \gls{simd} execution model corresponds directly to that of \glspl{gpu}, which is what enables \gls{gpu} acceleration for large-scale \glspl{nlp}.

Furthermore, the parameterized expression tree, the data structure realizing this abstraction, can also be exploited in sparse \gls{ad}, for which we implement a coloring-free, element-wise sparse reverse-mode algorithm.
The sparsity pattern is analyzed at the level of the individual parameterized expression tree at model-construction time, and the sparse Jacobian and Lagrangian Hessian are stored in a \emph{partially compressed} coordinate (COO) format, in which each derivative contribution is written element-wise to a precomputed position.
This structure removes the possibility of race conditions under parallel execution without incurring excessive fill-in in the sparse matrices.

Everything described above is implemented in ExaModels.jl, an open-source Julia package.
The implementation encodes the algebraic patterns in the type system of Julia, so the full algebraic structure is visible to the compiler and each evaluation kernel is compiled as pattern-specialized code.
This compile-time visibility supports \gls{aot} compilation without a separate code-generation step: the same code that runs under \gls{jit} compilation is compiled into standalone binaries or shared libraries, with no overhead in the iterative evaluation.
The pattern-specialized expression tree is furthermore a bits type, so the same kernels run on \glspl{gpu}.
Through KernelAbstractions.jl~\citep{churavyKernelAbstractionsjl2025}, a single kernel definition dispatches to NVIDIA, AMD, Intel, and Apple \glspl{gpu}, as well as to multithreaded \glspl{cpu}, providing cross-platform portability from one implementation.

The implementation additionally provides a range of features: the same model code runs in arbitrary numeric precision, including software-emulated extended precision; user-defined functions can be registered together with their derivatives; and every model supports, by default, parametric and two-stage stochastic formulations, as well as user-defined external oracles, each attached to a block of variables and constraints; batch formulations are under active development (\Cref{sec:extensions:batch}).

\subsection{Related Work}
\label{sec:intro:related}

ExaModels.jl is an \gls{ams} whose distinguishing capability lies in its \gls{ad} implementation.
The relevant prior work therefore spans two areas: algebraic modeling, and \gls{ad} together with its implementation paradigms.
At their intersection sit the \glspl{ams} that support \glspl{nlp} and ship their own \gls{ad}, with which ExaModels.jl directly competes.

\subsubsection{Algebraic Modeling Systems}
An \gls{ams} is a software system that lets the user state an optimization problem in a high-level syntax mirroring its algebraic (mathematical) form; well-known examples include AMPL~\citep{fourerModelingLanguageMathematical1990}, GAMS~\citep{bussieckGeneralAlgebraicModeling2004}, JuMP~\citep{lubinJuMP10Recent2023}, Pyomo~\citep{hartPyomoModelingSolving2011}, and CVXPY~\citep{diamondCVXPYPythonEmbeddedModeling2016}.
It sits between the user's formulation and the numerical solver, and produces from the stated problem what the solver requires.

Historically, \glspl{ams} began as standalone domain-specific languages: systems such as AMPL and GAMS define their own syntax, parse it with their own front ends, and exchange problems through their own formats.
Over time, the field has moved toward libraries embedded in high-level general-purpose languages, such as JuMP in Julia, Pyomo and CVXPY in Python, and CasADi with interfaces to both Python and MATLAB, so that models are ordinary programs with access to the host language's data structures, tooling, and ecosystem.
The two approaches can also be combined: a standalone system can serve as the backend of an embedded library, as Pyomo does by serializing models to AMPL's \texttt{.nl} format and obtaining derivatives from \gls{asl}~\citep{hartPyomoOptimizationModeling2017}.
ExaModels.jl belongs to the embedded category.
It is implemented in Julia and relies on the Julia compiler for its \gls{ad} capabilities, with no direct source code transformation.
The embedding also exposes models to the solver infrastructure, data analysis, and visualization tooling available in the Julia ecosystem.

\Glspl{ams} play a particularly important role in nonlinear programming.
For linear, quadratic, and \glspl{dcp}, what the solver requires is a fixed set of problem data that it reads once.
For a linear program, for example, it consists of the cost vector, the constraint matrix, and the bound vectors of the canonical form.
For \glspl{nlp}, by contrast, there is no such finite representation.
Instead of an explicit encoding of the constraints and objective, the \gls{ams} supplies the NLP functions as lazily evaluated callbacks, exposed to the solver as function pointers and queried at each iterate it visits.
Since these callbacks include the first- and second-order derivatives, derivative evaluation is an essential part of \gls{nlp} modeling.
Some \glspl{ams} implement their own derivative evaluation, while others delegate it to external tooling.
ExaModels.jl belongs to the former group: it is not a solver but an \gls{ams} that provides the modeling interface and the callback functions, backed by its own \gls{ad} capabilities.

\subsubsection{Automatic Differentiation and Its Implementation Paradigms}
Derivatives of computer programs can be evaluated in several different ways: by hand-written derivative evaluation code, by finite differences, by symbolic differentiation, or by \gls{ad}.
Hand-written derivatives are laborious and error-prone for large models.
Finite differences suffer from truncation error and require one perturbed evaluation per input, and symbolic differentiation suffers from expression swell~\citep{baydinAutomaticDifferentiationMachine2018}.
\Gls{ad} avoids these drawbacks by applying the chain rule to each elementary operation of the program.
It yields derivatives exact up to numerical precision at a cost proportional to the function evaluation itself, and it has become the standard tool for derivative computation in numerical computing whenever it is applicable~\citep{griewankEvaluatingDerivatives2008,naumannArtDifferentiating2012}.
A few systems use symbolic differentiation instead, such as Gravity~\citep{hijaziGravityMathematicalModeling2018} and the SymbolicAD backend of JuMP~\citep{dowsonMathOptSymbolicAD}.
The choice is justified by the fact that many algebraic expressions in classical mathematical programs are rather simple and thus do not suffer much from expression swell.
The \glspl{ams} and \gls{ad} frameworks surveyed in this section otherwise follow the \gls{ad} paradigm, as does ExaModels.jl, which ships its own \gls{ad} implementation (\Cref{sec:ad}).

\Gls{ad} runs in either \emph{forward mode} or \emph{reverse mode}, and reverse mode has long been the mode of choice in \glspl{ams} for mathematical programming.
A forward pass propagates directional derivatives to all outputs simultaneously, so forward mode scales well with the number of outputs.
A reverse pass accumulates the sensitivities of one output with respect to all inputs, so reverse mode scales well with the number of inputs~\citep{griewankEvaluatingDerivatives2008}.
As observed above, the callbacks of \glspl{nlp} with partially separable, repetitive algebraic structure consist of many small, independent, scalar-valued expressions, each depending on only a few variables.
For such scalar-valued expressions the many-output advantage of forward mode does not apply, and a single adjoint pass per expression yields all of its partial derivatives at a cost proportional to one evaluation.
ExaModels.jl uses reverse mode (\Cref{sec:ad}).

There are two paradigms for implementing \gls{ad}: \emph{operator overloading}, which redefines arithmetic on augmented types as the program runs, and \emph{source-code transformation}, which rewrites the program, or an intermediate representation of it, into one that also computes derivatives.
Existing \gls{ad} systems can mostly be classified into these two categories.
PyTorch uses operator overloading, building a representation of the computed function at every execution~\citep{paszkePyTorchImperativeStyle2019}, and ForwardDiff.jl implements forward mode by operator overloading on dual number types~\citep{revelsForwardModeAutomaticDifferentiation2016}.
JAX traces Python programs through overloaded tracer objects and differentiates the trace by program transformation~\citep{bradburyJAXComposableTransformations2018}.
Zygote.jl applies source-code transformation to Julia's intermediate representation in static single assignment form~\citep{innesDontUnrollAdjoint2018}, and Enzyme does the same at the level of the LLVM intermediate representation~\citep{mosesInsteadRewritingForeign2020}.
Among \glspl{ams} for mathematical programs, CasADi implements \gls{ad} as source-code transformation on its symbolic graph, generating new derivative expressions from captured ones~\citep{anderssonCasADiSoftwareFramework2019}.
AMPL parses models into expression graphs differentiated by \gls{asl}~\citep{gayRevisitingExpressionRepresentations2018}, and GAMS uses a similar built-in pipeline~\citep{bussieckGeneralAlgebraicModeling2004}.
JuMP captures nonlinear expressions into runtime expression graphs, which its general-purpose \gls{ad} routines then differentiate.
The capture was originally done by syntactic macros~\citep{dunningJuMPModelingLanguage2017}, and is done by operator overloading in its current interface~\citep{lubinJuMP10Recent2023}.

Operator overloading is simple to implement and inherits the full flexibility of the host language, since the program runs as ordinary code and the overloaded operators observe each operation as it executes.
That observation is also its overhead: the operators fire again at every evaluation, so the tape or graph from which the derivatives are computed is rebuilt at each execution, even when the structure of the computation has not changed.
Source-code transformation instead produces derivative code ahead of time, avoiding this overhead at the cost of a substantially more complex implementation~\citep{margossianReviewAutomaticDifferentiation2019}.
As a consequence, it has traditionally been perceived that operator overloading cannot be as efficient as source-code transformation and carries additional overhead.

However, in languages built on multiple dispatch and \gls{jit} compilation, such as Julia~\citep{bezansonJuliaFreshApproach2017}, this perception no longer holds: operator overloading can be used to compile derivative evaluation code specialized to the problem's algebraic structure, as encoded in the data types, and the resulting native code runs with no per-evaluation overhead.
The two mechanisms work together: multiple dispatch realizes operator overloading by selecting, for each operation, the method defined for the augmented types, and the \gls{jit} compiler then specializes each method for the concrete types it is invoked with, generating dedicated derivative evaluation code that is compiled once and reused.
The implementation remains simple: it depends largely on the compiler of the language and relies less on direct manipulation of the source code.
ExaModels.jl exploits this feature: the algebraic expression tree is aggressively typed, so the Julia compiler sees the full algebraic structure at compilation time and compiles the derivative evaluation code without performance overhead.
We demonstrate this through the compiler-generated LLVM intermediate representation reproduced in \Cref{app:llvm}.\footnote{Generating such low-level code is inherent to Julia's compilation, as in any compiled language, and is distinct from the separate source-level code-generation steps that some \glspl{ams} require; ExaModels.jl involves none of the latter.}

\subsubsection{Modeling Systems with Built-in Automatic Differentiation}
ExaModels.jl targets \glspl{nlp} of the structured form discussed above and supplies its own \gls{ad}, so what follows focuses on the \glspl{ams} that operate in the same space.
These are the systems that target classical mathematical programming problems, including general smooth \glspl{nlp}, and that ship with their own \gls{ad} machinery.
Five such systems are reviewed below, and \Cref{tab:ams} compares them with ExaModels.jl.
AMPL~\citep{fourerModelingLanguageMathematical1990} is a standalone language with its own parser, exposing the NLP function callbacks through the \gls{asl} via the \texttt{.nl} interchange format; GAMS~\citep{bussieckGeneralAlgebraicModeling2004} emerged around the same time and is similar in design.
JuMP~\citep{dunningJuMPModelingLanguage2017,lubinJuMP10Recent2023} is embedded in Julia; through macros and operator overloading it builds a per-scalar expression graph and runs reverse-mode \gls{ad} on the graph.
CasADi~\citep{anderssonCasADiSoftwareFramework2019} is a C++ symbolic framework whose central abstraction is the \texttt{Function}, a compiled expression graph that can be reused, embedded, or emitted as standalone C code; it supports both forward and reverse \gls{ad}, and uses graph coloring to evaluate sparse Jacobians and Hessians.
Gravity~\citep{hijaziGravityMathematicalModeling2018} introduces ``constraint templates'' to preserve the pattern structure of indexed constraints.

\begin{table}[t]
  \centering
  \caption{\Glspl{ams} with built-in \gls{ad} that target smooth \glspl{nlp}. ``Pattern-aware'' indicates whether the system preserves the repetitive structure of indexed constraints for parallel derivative evaluation.}
  \label{tab:ams}
  \footnotesize
  \setlength{\tabcolsep}{4pt}
  \begin{tabular*}{\textwidth}{@{\extracolsep{\fill}}l l l l c}
    \toprule
    System & Language & Problem classes & AD & Pattern-aware \\
    \midrule
    AMPL         & own           & LP, MILP, NLP        & ASL (built-in)            & no \\
    GAMS         & own           & LP, MILP, NLP        & built-in                  & no \\
    JuMP         & Julia         & LP, MILP, conic, NLP & reverse-mode on graph     & no \\
    CasADi       & C++/Py/MATLAB & LP, MILP, conic, NLP & source transform          & yes \\
    Gravity      & C++           & NLP (network)        & symbolic                  & yes \\
    \midrule
    ExaModels.jl & Julia         & NLP                  & reverse-mode per pattern  & yes \\
    \bottomrule
  \end{tabular*}
\end{table}

We compare these systems along two axes: how the partially separable, repetitive structure of indexed constraints is preserved and exploited, and whether the NLP functions are evaluated by compiled or interpreted code.

The first axis is structure exploitation, and most systems do not exploit the repetitive structure.
AMPL, GAMS, and JuMP build flat, per-scalar expression graphs at model construction, where every constraint instance becomes a node in a single large graph that \gls{ad} traverses at evaluation time.
JuMP is the canonical example of iterator flattening.
Indexed-constraint syntax such as \texttt{@constraint(model, [i=1:N], ...)} is lowered to a flat scalar collection before \gls{ad}, so the iteration visible to the user is no longer visible to the \gls{ad} machinery.

Two strategies preserve indexed-constraint structure at the modeling layer: vectorization and templatizing.
Vectorization expresses an indexed family as a single matrix-algebra operation, so the repetition lives in the dimensions of the matrices and vectors rather than in enumerated scalar expressions.
Templatizing declares a single algebraic expression parameterized by a data set and evaluates that one expression across all instances, so the repetition lives in the data rather than in the expression.
The two strategies differ in abstraction rather than reach.
Heterogeneous indexed structure, meaning general repeated expressions whose size or form varies per index, such as variable-size sums, can be expressed in either form.
Vectorization shifts the work to pre-constructing sparse incidence matrices, while templatizing expresses the repeated expression directly.

The representative implementation of vectorization is CasADi, which compiles each \texttt{Function} once and reuses it across calls.
CasADi supports structure preservation through \emph{vectorized} expressions of its matrix symbolic type, MX.
The slice \texttt{x[1:] - x[:-1]} stays a single matrix-level node, and both the symbolic graph and the generated C function body remain constant-size as the discretization grows; only the loop bound and workspace offsets change.
This requires the model to be written as matrix algebra.
Heterogeneous indexed structure, such as \gls{acopf} power balance over each bus's generators and incident lines, can be encoded via sparse incidence matrices and sparse matrix--vector products.
The result is mathematically equivalent to a templatized model, but uses a different abstraction: matrix algebra over pre-constructed sparsity patterns versus per-index expressions over data sets.
Machine learning frameworks such as JAX~\citep{bradburyJAXComposableTransformations2018} and PyTorch~\citep{paszkePyTorchImperativeStyle2019} are not widely used for formulating and solving classical mathematical programs, but they support the same matrix algebra syntax and admit natural vectorization.
The model is written in whole-array operations and differentiated by the framework's \gls{ad}.
They are nonetheless natural companions for \gls{gpu} evaluation of the \gls{nlp} functions: gradient evaluation is directly supported, and the sparse Jacobian and Hessian can be evaluated through coloring-based sparse \gls{ad}.

Templatizing as an evaluation strategy is a more recent idea, first explored in Gravity through its constraint templates: a constraint is declared once over an index set, and the repeated structure is exploited to reduce the evaluation time of the NLP functions, including the Jacobian and the Hessian.
Gravity is implemented as a C++ library in which the numerical accuracy of variables and parameters is specifiable by the user.
It further offers convexity detection, multithreading of subproblems, and lazy constraint generation, and the authors report function evaluation five times faster and memory use up to 60 times lower than JuMP~\citep{hijaziGravityMathematicalModeling2018}.
The indexed declaration syntax itself is common: JuMP's \texttt{@constraint(model, [i=1:N], ...)} and Pyomo's \texttt{Constraint(index\_set, rule=f)} both declare an indexed constraint family in a single statement.
To the best of our knowledge, however, neither system specializes the evaluation internally to exploit the structure this syntax exposes, and the declaration is unrolled into flat per-scalar instances.
Because this indexed style has long been standard in the \glspl{ams} of mathematical programming and operations research, such as AMPL, GAMS, JuMP, and Pyomo, the modeling syntax of ExaModels.jl remains close to that of existing modeling tools.

ExaModels.jl follows the templatizing strategy but goes beyond what Gravity has proposed.
It provides a general abstraction that can represent diverse repetitive patterns, together with an implementation built on it, namely type-stable parameterized expression trees compiled per pattern by the Julia compiler (\Cref{sec:graph}).
This design naturally leads to \gls{gpu} parallelization, which we implement through KernelAbstractions.jl and demonstrate in this paper (\Cref{sec:gpu}).
As of this writing, neither CasADi nor Gravity supports \gls{gpu}-native derivative evaluation.

The second axis is the evaluation path, compiled versus interpreted, where compilation reduces per-callback runtime overhead at the cost of longer build time.
AMPL parses models to the \texttt{.nl} format, which is evaluated by \gls{asl} with compiled operator implementations; a legacy \texttt{nlc} tool can emit standalone C.
GAMS uses a similar built-in pipeline.
CasADi provides both options: by default a \texttt{Function} is evaluated through CasADi's runtime, which walks a flat operation sequence with compiled operator implementations; alternatively, standalone C code can be emitted via codegen.
Gravity uses symbolic differentiation and dispatches compiled C++ templates.
JuMP is embedded in the \gls{jit}-compiled Julia language, but its default \gls{ad} does not compile a specialized derivative kernel per constraint type.
It is a general-purpose interpreter that walks the expression graph at each callback, so it cannot be regarded as a compiled \gls{ad} approach.
The exception is its SymbolicAD backend~\citep{dowsonMathOptSymbolicAD}, which reduces each constraint to a canonical symbolic form and computes a single symbolic derivative per distinct constraint structure.
The resulting compiled evaluator is then reused across all instances that share that structure.

ExaModels.jl is firmly in the compiled category, and it compiles aggressively.
It builds a parameterized expression tree at model-construction time and applies reverse-mode \gls{ad} directly to that tree (\Cref{sec:ad}).
Because the tree is type-stable end-to-end, each callback compiles into a single chunk of LLVM code specialized to the problem's algebraic structure (\Cref{sec:impl:specialization}).
This gives the Julia compiler full opportunity to optimize the \gls{nlp} functions, and it also makes the callbacks straightforward to convert into \gls{gpu} kernels (\Cref{sec:gpu}).

\subsection{Contributions}
\label{sec:intro:contributions}
The contributions of this paper are methodological (M), implementational (I), and experimental (E).
\begin{enumerate}
  \item[\textbf{M1.}] We introduce a \emph{\gls{simd} abstraction} for \glspl{nlp} that represents a problem as a small number of algebraic patterns, each evaluated over many data points.
The abstraction preserves the pattern--data separation throughout model construction and reduces NLP function evaluation to embarrassingly parallel loops over data points, which is the structure that \gls{gpu} execution requires.
The closest prior idea is Gravity's constraint templates~\citep{hijaziGravityMathematicalModeling2018}; our abstraction goes beyond theirs by providing a general representation that covers diverse repetitive patterns and extends naturally to parallel \gls{gpu} evaluation.
  \item[\textbf{M2.}] We develop a \emph{coloring-free sparse \gls{ad} algorithm} that exploits regular per-pattern sparsity to compute Jacobians and Hessians in a partially compressed COO form.
Unlike the coloring-based sparse \gls{ad} of standard practice~\citep{gebremedhinWhatColorYour2005}, the sparsity analysis is applied at the level of the parameterized expression tree, and the derivative values are partially compressed over each tree.
This minimizes the fill-in while eliminating the race conditions that arise when the values are accumulated in parallel, so that element-wise, coloring-free sparse \gls{ad} can be implemented for \glspl{gpu}.
To the best of our knowledge, this approach is new.
  \item[\textbf{I1.}] We present, to our knowledge, the first \gls{ams} with cross-platform \gls{gpu}-native derivative evaluation, portable across NVIDIA, AMD, Intel, and Apple devices via KernelAbstractions.jl~\citep{churavyKernelAbstractionsjl2025}; the same model code runs on all backends without modification.
  \item[\textbf{I2.}] We support \gls{aot} compilation to standalone binaries independent of the Julia runtime, with no code generation step.
Whereas systems such as CasADi and CVXPY provide separate code-generation tools for deployment, in ExaModels.jl exactly the same code serves both \gls{jit} and \gls{aot} compilation; to our knowledge, it is the first \gls{ams} with this property.
  \item[\textbf{I3.}] We support arbitrary numeric precision by exploiting the type flexibility of the Julia language.
Without modification, the same model code runs in \texttt{Float64}; in \texttt{Float32}, including on hardware that supports only single precision; and in user-defined types realizing software-emulated extended precision, such as the double-double arithmetic of DoubleFloats.jl~\citep{juliamathDoubleFloatsjl}.
Existing \glspl{ams} are typically fixed to double precision.
Among those reviewed above, Gravity allows variables and parameters to be declared over a fixed set of built-in C++ numeric types (float, double, long double, integer, and complex), though not over user-defined types, and JuMP's nonlinear interface is restricted to \texttt{Float64}.
  \item[\textbf{E1.}] We contribute ExaModels-native implementations of the Luk\v{s}an--Vl\v{c}ek~\citep{lukVsanSparsePartiallySeparable}, COPS~\citep{dolanBenchmarkingOptimizationSoftware2001}, and PGLIB-OPF benchmark libraries.
We systematically measure derivative evaluation performance on the \gls{cpu} and across multiple \gls{gpu} backends.
These benchmark implementations are a contribution in their own right: they run natively on the \gls{gpu} and can be reused in later studies for the performance benchmarking of \gls{gpu}-resident \gls{nlp} algorithms.
We report \gls{gpu} speedups of up to $\gsLvHessLarge\times$ over single-threaded \gls{cpu} evaluation.
\end{enumerate}

\paragraph{Relationship to prior work.}
The concept of the \gls{simd} abstraction first appeared in~\citet{shinAcceleratingOptimalPower2024} as part of an \gls{acopf}-specific \gls{gpu} optimization framework.
The abstraction itself, the parameterized expression tree, the coloring-free sparse \gls{ad} algorithm, and the multi-backend \gls{gpu} implementation have not previously been elaborated, and this is the role of the present paper.
Several subsequent works use ExaModels.jl as the modeling layer:~\citet{johnsonExaModelsPowerjl2025} provides a power-systems library on top of it; the \gls{gpu}-native condensed-space interior-point methods of~\citet{pacaudCondensedspaceMethodsNonlinear2024} consume its callbacks; and~\citet{pacaudGPUacceleratedNonlinearModel2024} reports a fully \gls{gpu}-resident optimization pipeline pairing ExaModels.jl with MadNLP.jl~\citep{shinMadNLPjl2025}.

\subsection{Limitations}
\label{sec:intro:limitations}
The proposed approach also carries drawbacks that follow directly from its design choices.
\begin{itemize}
  \item The approach is runtime-optimized at the cost of a compilation overhead: the NLP function callbacks are compiled specifically for each algebraic pattern, so every new expression structure must be compiled at least once before it can be evaluated.
This places a burden on the compiler and adds a one-time compilation latency ahead of the first evaluation.
The design favors runtime performance over the ability to quickly try out many formulations, which may not be ideal for rapid prototyping across structurally different models.
It is best suited, for example, to a mature model deployed for fast repeated evaluation, as in real-time decision making.
  \item The user is required to write the model in the \gls{simd}-compatible form.
It is easier in many cases to state a model through explicit loops over scalar constraints.
The pattern--data separation also requires, in some cases, the data arrays to be created before the objective and constraint patterns can be declared, which may require restructuring compared with the natural formulation in other modeling systems.
  \item The benefits are contingent on the problem being partially separable and repetitive, and a model with many distinct patterns and few data points per pattern gains little.
Components that cannot be expressed as parameterized expression trees, such as constraints with data-dependent branching or embedded simulation codes, fall outside the abstraction.
These are instead accommodated through external oracles (\Cref{sec:extensions:oracles}), which do not carry the performance guarantees of the abstraction.
  \item The templatized syntax expresses the NLP functions entirely through scalar operations.
Unless the external oracle mechanism is used, the evaluation invokes no vectorized library routines, such as BLAS kernels or sparse matrix--vector products, and therefore cannot benefit from such optimized libraries.
\end{itemize}

\subsection{Outline}
\label{sec:intro:outline}
The remainder of this paper falls into three parts, with the fundamentals developed in a top-down manner.
\Cref{sec:modeling} presents the modeling abstraction and the user- and solver-facing interfaces.
\Cref{sec:data} presents the data structures behind them: the model containers and the internal types that represent variables, parameters, objectives, and constraints, whose common building block is the \texttt{SIMDFunction}.
\Cref{sec:graph} presents the parameterized expression tree, the core idea within the \texttt{SIMDFunction}, and describes how individual nodes and leaves are put together to form the tree.
\Cref{sec:ad} explains how this tree is traversed to evaluate the derivatives with the coloring-free sparse \gls{ad} algorithm.
These fundamentals are then realized in a high-performance, portable implementation.
\Cref{sec:gpu} describes the \gls{gpu} acceleration via KernelAbstractions.jl, and \Cref{sec:impl} collects implementation-level details deferred from the main development, including the extension of the abstraction to two-stage stochastic and batch problems.
Finally, the experiments are reported in \Cref{sec:numerics}, and \Cref{sec:conclusions} concludes the paper.

\subsection{Notation}
\label{sec:intro:notation}
We use $\mathbb{R}^n$ for $n$-dimensional real space.
Code identifiers are typeset in \texttt{monospace}.
We further distinguish a mathematical quantity from its software representation: for example, $x$ denotes a vector as a mathematical object, whereas \texttt{x} refers to the array that stores it.
We also use a few standard Julia conventions throughout.
Indexing is one-based, following Julia's convention.
Variable blocks may nonetheless be declared over arbitrary integer index ranges, such as the \texttt{0:nh} blocks of the Goddard rocket model in \Cref{sec:modeling}.
The symbol \texttt{<:} denotes the subtype relation, so \texttt{Node1 <: AbstractNode} reads ``\texttt{Node1} is a subtype of \texttt{AbstractNode}.''
A trailing \texttt{!} marks a function that mutates its arguments in place (e.g., \texttt{set\_value!}).
A leading \texttt{@} marks a macro, one of Julia's metaprogramming facilities.
Rather than being called at run time, a macro receives its arguments as unevaluated syntax and rewrites them, at parse time, into an equivalent block of ordinary code that is then compiled in its place.
Writing \texttt{@add\_var}, for instance, is shorthand for the expanded sequence of statements the macro emits.
A name written with curly braces, \texttt{Type\{I, J, \ldots\}}, is a \emph{parametric type}: the entries \texttt{I}, \texttt{J}, \ldots are its type parameters, on which the compiler specializes the generated code (e.g., \texttt{DataIndexed\{I, J\}}).
An anonymous function is written \texttt{(args) -> body} (e.g., \texttt{(i, x) -> x[i+1]\^{}2}).

\section{Modeling Abstraction and Interfaces}
\label{sec:modeling}

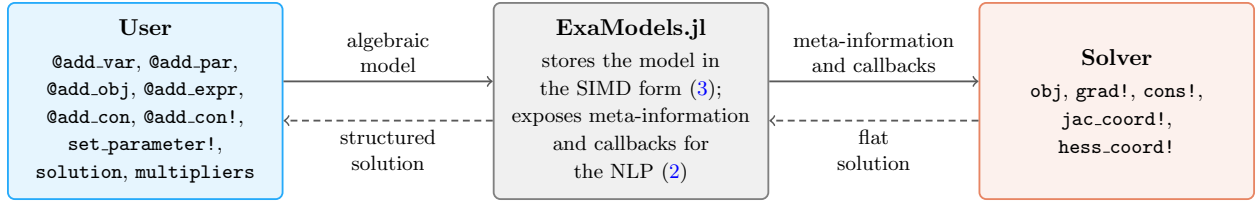
\begin{figure}[t]
  \centering
  \resizebox{\linewidth}{!}{%
\begin{tikzpicture}[
    node distance=14mm and 32mm,
    box/.style={draw=black!50, thick, rounded corners=3pt, align=center,
      inner sep=6pt, font=\small, fill=black!2,
      minimum height=30mm, minimum width=42mm},
    fwd/.style={->, thick, draw=black!65},
    bwd/.style={->, thick, densely dashed, draw=black!65},
    ]
    \node[box, draw=j1!85, fill=j1!10] (user) {\textbf{User}\\[2pt]\footnotesize \texttt{@add\_var}, \texttt{@add\_par},\\\footnotesize \texttt{@add\_obj}, \texttt{@add\_expr},\\\footnotesize \texttt{@add\_con}, \texttt{@add\_con!},\\\footnotesize \texttt{set\_parameter!},\\\footnotesize \texttt{solution}, \texttt{multipliers}};
    \node[box, right=of user, fill=black!6] (exa) {\textbf{ExaModels.jl}\\[2pt]\footnotesize stores the model in\\\footnotesize the SIMD form \cref{eq:simd};\\\footnotesize exposes meta-information\\\footnotesize and callbacks for\\\footnotesize the NLP \cref{eq:nlp}};
    \node[box, draw=j2!85, fill=j2!10, right=of exa] (solver) {\textbf{Solver}\\[2pt]\footnotesize \texttt{obj}, \texttt{grad!}, \texttt{cons!},\\\footnotesize \texttt{jac\_coord!},\\\footnotesize \texttt{hess\_coord!}};
    \draw[fwd] ([yshift=3mm]user.east) -- node[above, align=center, font=\footnotesize] {algebraic\\model} ([yshift=3mm]exa.west);
    \draw[bwd] ([yshift=-3mm]exa.west) -- node[below, align=center, font=\footnotesize] {structured\\solution} ([yshift=-3mm]user.east);
    \draw[fwd] ([yshift=3mm]exa.east) -- node[above, align=center, font=\footnotesize] {meta-information\\and callbacks} ([yshift=3mm]solver.west);
    \draw[bwd] ([yshift=-3mm]solver.west) -- node[below, align=center, font=\footnotesize] {flat\\solution} ([yshift=-3mm]exa.east);
  \end{tikzpicture}}
  \caption{ExaModels.jl between the user and the solver. Solid arrows carry the model forward from the modeling syntax to the solver callbacks; dashed arrows carry the solution back.}
  \label{fig:interfaces}
\end{figure}

ExaModels.jl sits between the user and the solver (\Cref{fig:interfaces}).
The user states a problem in mathematical terms, while the solver asks for numerical values of the NLP functions at each iteration.
An \gls{ams} must therefore expose an interface on both sides: a solver-facing callback interface for evaluating the model and a user-facing syntax for building it.

On the solver side, the interface is largely standard.
ExaModels.jl follows the convention of NLPModels.jl~\citep{jsoNLPModelsJl}, which fixes the model--solver interface as a set of meta-information accessors and callback functions for objective, constraint, gradient, Jacobian, and Hessian evaluation.
An ExaModels.jl model can therefore be passed directly to solvers for general \glspl{nlp}, such as Ipopt~\citep{wachterImplementationInteriorpointFilter2006}, Uno~\citep{vanaretImplementingUnifiedSolver2026}, and MadNLP.jl~\citep{shinMadNLPjl2025}, and to any other solver in the \gls{jso} ecosystem.
The user side is what distinguishes ExaModels.jl.
Each objective term and constraint is declared as an algebraic pattern paired with a data iterator, so a whole model is stated in a few such declarations.
The declarations also keep the pattern--data separation intact for the evaluation machinery.

\Cref{sec:modeling:nlp} presents the solver-facing interface: the canonical NLP form and the meta-information accessors and callback functions that ExaModels.jl must supply.
\Cref{sec:modeling:simd} presents the user-facing side: the \gls{simd} abstraction, which preserves the partially separable, repetitive structure of large-scale NLPs, and the modeling syntax that builds models in that form.

\subsection{Solver-Facing Interface}
\label{sec:modeling:nlp}

\subsubsection{NLP Formulation}
The solver-side \emph{mathematical abstraction} of the \gls{nlp} for ExaModels.jl is as follows:
\begin{equation}\label{eq:nlp}
  \min_{x \in \mathbb{R}^n}\; f(x;\, \theta) \quad
  \st\quad g^\flat \leq g(x;\, \theta) \leq g^\sharp,\quad x^\flat \leq x \leq x^\sharp,
\end{equation}
where $x \in \mathbb{R}^n$ is the decision vector, $f(\cdot):\mathbb{R}^n \to \mathbb{R}$ is the objective function, $g(\cdot):\mathbb{R}^n \to \mathbb{R}^m$ is the constraint function, and $x^\flat, x^\sharp \in \mathbb{R}^n$ and $g^\flat, g^\sharp \in \mathbb{R}^m$ are the variable and constraint bounds.
The vector $\theta \in \mathbb{R}^p$ collects \emph{problem parameters}: quantities that the objective and constraints depend on but that are held fixed during the solve rather than optimized over.
ExaModels.jl treats every model as such a \emph{parametric} NLP by default; a model with no parameters is simply the special case $p = 0$.
Throughout the paper, we assume that $f(\cdot)$ and $g(\cdot)$ are twice differentiable and well defined over the entire set defined by the bounds.
The form \cref{eq:nlp} is canonical: any NLP can be written this way, and it is the form in which solvers expect the problem.

As discussed in \Cref{sec:intro}, an \gls{nlp} solver does not consume a fixed set of problem data.
Instead, it queries the modeling layer for numerical evaluations~\citep{nocedalNumericalOptimization2006} at every iterate $x$ it visits during the solve.
The \emph{implementational abstraction} splits the information the solver needs from ExaModels.jl into two kinds: static \emph{meta-information} about the problem, supplied once, and per-iteration \emph{NLP function callbacks}, evaluated on demand.

\subsubsection{Meta-Information}
Before the solve begins, the solver queries static information about the problem, that is, data that is fixed throughout the solve and therefore supplied only once.
ExaModels.jl stores this data in an \texttt{NLPModels.NLPModelMeta} object and exposes it through the NLPModels.jl meta-information accessors, such as \texttt{NLPModels.get\_nvar} and \texttt{NLPModels.get\_ncon}.
The main entries are:
\begin{itemize}
  \item \texttt{nvar}, \texttt{ncon}, \texttt{nnzj}, \texttt{nnzh}: integer values for the problem dimensions, namely the numbers of variables and constraints and the numbers of Jacobian and Lagrangian Hessian nonzeros;
  \item \texttt{lvar}/\texttt{uvar}, \texttt{lcon}/\texttt{ucon}, \texttt{x0}, \texttt{y0}: vector values that pin down the instance and its starting point, namely the variable bounds $x^\flat, x^\sharp$, the constraint bounds $g^\flat, g^\sharp$, and the initial primal point and Lagrange multipliers, both zero by default;
  \item \texttt{ifix}, \texttt{ilow}, \texttt{iupp}, \texttt{irng}, \texttt{ifree} and \texttt{jfix}, \texttt{jlow}, \texttt{jupp}, \texttt{jrng}, \texttt{jfree}: index vectors partitioning the variables and constraints by bound type (fixed, lower-bounded, upper-bounded, range-bounded, and free); and
  \item \texttt{minimize}: the optimization sense, that is, whether the problem is a minimization or a maximization.
\end{itemize}

\subsubsection{Callback Functions}
At each iterate, the solver queries the numerical values of the NLP functions through five callbacks.
It supplies the iterate $x$ and, for the Lagrangian Hessian, an objective weight and the multipliers.
These callbacks constitute the \emph{solver-facing interface} of ExaModels.jl.
Except for the objective value, which is returned as a scalar, the callbacks are in-place (hence the \texttt{!} suffix) and write their results into preallocated output arrays:
\begin{itemize}
  \item \texttt{NLPModels.obj(model, x)}: the objective value $f(x) \in \mathbb{R}$;
  \item \texttt{NLPModels.grad!(model, x, g)}: the objective gradient $\nabla f(x) \in \mathbb{R}^n$, written into \texttt{g};
  \item \texttt{NLPModels.cons!(model, x, c)}: the constraint values $g(x) \in \mathbb{R}^m$, written into \texttt{c};
  \item \texttt{NLPModels.jac\_coord!(model, x, v)}: the sparse constraint Jacobian $\nabla g(x) \in \mathbb{R}^{m \times n}$, with nonzero values written into \texttt{v};
  \item \texttt{NLPModels.hess\_coord!(model, x, y, v; obj\_weight)}: the sparse $n \times n$ Lagrangian Hessian $\sigma \nabla^2 f(x) + \sum_{j=1}^m \lambda_j \nabla^2 g_j(x)$, with nonzero values written into \texttt{v},
\end{itemize}
where $\lambda \in \mathbb{R}^m$ is the vector of Lagrange multipliers, supplied in \texttt{y}, and $\sigma$ is the objective weight, supplied as the keyword argument \texttt{obj\_weight}.
Here and throughout the paper, gradients are column vectors, and the Jacobian of a vector-valued function stacks the transposed gradients of its components as rows.
These five are the most commonly used callbacks; the NLPModels.jl interface additionally specifies others, such as:
\begin{itemize}
  \item \texttt{NLPModels.jac\_structure!(model, I, J)}: the COO row and column indices of the nonzero Jacobian entries, written into \texttt{I} and \texttt{J};
  \item \texttt{NLPModels.hess\_structure!(model, I, J)}: the COO row and column indices of the nonzero Hessian entries, written into \texttt{I} and \texttt{J};
  \item \texttt{NLPModels.jprod!}: the Jacobian--vector product;
  \item \texttt{NLPModels.jtprod!}: the vector--Jacobian product;
  \item \texttt{NLPModels.hprod!}: the Hessian--vector product.
\end{itemize}
The index vectors are queried once, at the start of the solve.
Together with the value vector \texttt{v}, they form the COO representation of the sparse derivative matrices, and only \texttt{v} is recomputed at each iterate.

For the solver, the objective and constraint values are used to monitor progress (e.g., in the line search), and the first- and second-order derivatives are used to compute the step direction (e.g., the Newton step).
Because they are evaluated at every iterate, these callbacks are performance-critical; providing them efficiently is a central goal of ExaModels.jl.

This callback API was developed in NLPModels.jl by the \gls{jso} organization around the abstract type \texttt{NLPModels.AbstractNLPModel}.
It decouples the modeling layer from the solver layer: any model that implements the interface is accepted by any compatible solver, including those listed at the opening of this section.
\texttt{ExaModel} is a subtype of \texttt{NLPModels.AbstractNLPModel} and implements this callback interface.

\subsection{User-Facing Abstraction and Interface}
\label{sec:modeling:simd}

To provide the callbacks of \Cref{sec:modeling:nlp}, ExaModels.jl must store the information needed to evaluate the NLP functions.
It stores that information in the form of the \gls{simd} \gls{nlp} formulation, keeping each repetitive pattern together with its data.

\subsubsection{SIMD NLP Formulation}
The \gls{simd} \gls{nlp} formulation captures the partially separable, repetitive structure of large-scale \glspl{nlp} (\Cref{sec:intro}) at the modeling level.
Problems like the Goddard rocket \cref{eq:rocket} belong to the following class of \glspl{nlp}, in which the repetitive patterns are kept explicit:
\begin{subequations}\label{eq:simd}
\begin{align}
  \min_{x^\flat \leq x \leq x^\sharp}\;
  & \sum_{\ell=1}^{L} \sum_{i=1}^{I_\ell} f^{(\ell)}(x;\, \theta, p_i^{(\ell)}) \label{eq:simd:obj}\\
  \st\; & g^\flat \leq
  \left[ g^{(m)}(x;\, \theta, q_j^{(m)}) \right]_{j=1}^{J_m}
  + \sum_{n=1}^{N_m} \sum_{k=1}^{K_n} h^{(n)}(x;\, \theta, s_k^{(n)})
  \leq g^\sharp,
  \quad \forall\, m \in \{1, \ldots, M\}, \label{eq:simd:con}
\end{align}
\end{subequations}
where $x \in \mathbb{R}^{n}$ is the vector of decision variables with bounds $x^\flat, x^\sharp$, and $\theta \in \mathbb{R}^p$ is the problem-parameter vector of \cref{eq:nlp}.
The \emph{algebraic patterns} $f^{(\ell)} : \mathbb{R}^n \to \mathbb{R}$ and $g^{(m)} : \mathbb{R}^n \to \mathbb{R}$ are twice differentiable scalar functions of $x$, where these signatures suppress the dependence on $\theta$ and on the per-instance data.
Each augmentation pattern $h^{(n)} : \mathbb{R}^n \to \mathbb{R}^{J_m}$ is a vector-valued function whose output has a single nonzero entry.
The position of that entry within constraint block $m$ is fixed by the data $s_k^{(n)}$.
The quantities $p_i^{(\ell)}$, $q_j^{(m)}$, and $s_k^{(n)}$ are problem data that parameterize each pattern instance.
The vector $\theta$, in contrast, is shared: the same parameter vector enters every pattern and every data point.
The objective \cref{eq:simd:obj} is split into $L$ groups of objective patterns, where pattern $f^{(\ell)}$ is evaluated over $I_\ell$ data points.
The constraints \cref{eq:simd:con} consist of $M$ blocks: within block $m$, the ``base'' pattern $g^{(m)}$ defines $J_m$ constraints, and the ``augmentation'' patterns $h^{(n)}$ contribute additive terms to these constraints.

The augmentation structure arises naturally in network models, for example when arc flows are summed into nodal balance equations, and is discussed further in \Cref{sec:data:core}.
Unlike \cref{eq:nlp}, \cref{eq:simd} is not a canonical form, since not every NLP is naturally written this way.
It is, however, a form that appears throughout practice, as the examples below illustrate.
Problems with partially separable, repetitive algebraic structure fit the form \cref{eq:simd} with only a few distinct patterns; the Goddard rocket \cref{eq:rocket} is a direct instance, with the five patterns identified in \Cref{sec:intro}, and the same holds across the problem categories surveyed there.

We call the comprehensive interface built on this formulation the \emph{\gls{simd} abstraction}.
It comprises the user-facing modeling syntax and the solver-facing evaluation machinery: the user writes the model in the form of \cref{eq:simd}, and the model data is stored in a structured form that exposes the repeated structure.
This \emph{pattern--data separation} is the foundation of the abstraction.

Because the structure is partially separable, each pattern evaluation is independent across data indices, and the derivative computations inherit the same embarrassingly parallel structure.
For a given pattern $f^{(\ell)}$, the gradient contributions $\nabla_x f^{(\ell)}(x;\, \theta, p_i^{(\ell)})$ for $i = 1, \ldots, I_\ell$ can be evaluated in parallel, and so can the Jacobian and Hessian contributions.
By preserving the pattern structure at the modeling level, ExaModels.jl makes this parallelism directly available to the \gls{ad} implementation.
When the number of patterns is fixed and the problem size grows only through the number of data points, the evaluation of the NLP functions in principle attains $O(1)$ parallel time complexity, given sufficiently many threads.
The total work still grows with the number of data points; it is the parallel evaluation depth that stays constant.
How the derivative computations are carried out in this parallel form, including the accumulation of overlapping contributions without race conditions, is presented in \Cref{sec:ad}.

The structure of \cref{eq:simd} suggests how to store a model without losing its parameterized form.
Rather than recording every scalar term one by one, each pattern is kept as a \emph{bundle} of an \emph{instruction}, the algebraic expression to evaluate, and the \emph{data} over which that instruction is repeated.
Because a large-scale NLP typically comprises only a few distinct patterns, the entire model reduces to a small number of such bundles.
ExaModels.jl realizes the instruction as a \texttt{SIMDFunction}, which an \texttt{Objective} or a \texttt{Constraint} pairs with the matching data iterator.
The model is therefore stored as a collection of \texttt{SIMDFunction}s, one per pattern, together with their data.
Each callback then reduces to a loop over patterns that evaluates each pattern's \texttt{SIMDFunction} at all of its data points.
For the objective \cref{eq:simd:obj}, the value $f(x)$ is the sum of the per-data-point evaluations, and the constraint and derivative callbacks follow the same structure.
The internal construction of a \texttt{SIMDFunction} is presented in \Cref{sec:graph,sec:ad,sec:gpu}.

\subsubsection{Modeling Syntax}
\label{sec:modeling:syntax}

The syntax with which the user builds a model is shown on the Goddard rocket \cref{eq:rocket} of \Cref{sec:intro}, in four stages: creating a model, solving it, accessing the solution, and updating the data and resolving.

\paragraph{Creating a model.} The variables, objective, and constraints are declared on an \texttt{ExaCore}:
\begin{lstlisting}
using ExaModels, NLPModelsIpopt

# problem data
nh = 200   # number of time steps
h_0, v_0, m_0, g_0 = 1.0, 0.0, 1.0, 1.0
T_c, h_c, v_c, m_c = 3.5, 500.0, 620.0, 0.6
c_e = 0.5 * sqrt(g_0 * h_0)
m_f = m_c * m_0
D_c = 0.5 * v_c * (m_0 / g_0)
T_max = T_c * m_0 * g_0

# maximize final altitude
core = ExaCore(minimize = false)

# decision variables
@add_var(core, h, 0:nh; start = 1.0, lvar = 1.0)
@add_var(core, v, 0:nh; start = (i/nh * (1 - i/nh) for i = 0:nh))
@add_var(core, m, 0:nh;
    start = ((m_f - m_0) * i/nh + m_0 for i = 0:nh),
    lvar = m_f, uvar = m_0)
@add_var(core, tau, 0:nh; start = T_max/2, lvar = 0.0, uvar = T_max)
@add_var(core, dt, 1; start = 1/nh, lvar = 0.0)

# objective: final altitude
@add_obj(core, h[nh])

# altitude dynamics
@add_con(core, alt,
    -h[i] + h[i-1] + 0.5 * dt[1] * (v[i] + v[i-1]) for i = 1:nh)

# velocity dynamics
@add_con(core, vel,
    -v[i] + v[i-1] + 0.5 * dt[1] * (
        (tau[i] - D_c * v[i]^2 * exp(-h_c * (h[i] - h_0) / h_0)
            - m[i] * g_0 * (h_0 / h[i])^2) / m[i]
        + (tau[i-1] - D_c * v[i-1]^2 * exp(-h_c * (h[i-1] - h_0) / h_0)
            - m[i-1] * g_0 * (h_0 / h[i-1])^2) / m[i-1])
    for i = 1:nh)

# mass dynamics
@add_con(core, mass,
    -m[i] + m[i-1] + 0.5 * dt[1] * (-tau[i]/c_e - tau[i-1]/c_e) for i = 1:nh)

# parameter: required final mass
@add_par(core, mf, [m_f])

# boundary conditions
@add_con(core, h[0] - h_0)    # initial altitude
@add_con(core, v[0] - v_0)    # initial velocity
@add_con(core, m[0] - m_0)    # initial mass
@add_con(core, m[nh] - mf[1]) # final mass

# instantiate: ExaCore -> ExaModel
model = ExaModel(core)
\end{lstlisting}

The patterns, namely the dynamics \texttt{alt}, \texttt{vel}, and \texttt{mass}, the boundary conditions, and the objective, correspond directly to \cref{eq:rocket}.
Each \texttt{@add\_var} call creates a block of variables rather than a single variable: \texttt{h}, for example, comprises one variable for each index in \texttt{0:nh}.
The same holds for the constraints: each \texttt{@add\_con} call creates a block of constraints, one per element of its data iterator (a singleton block when no iterator is given, as for the boundary conditions).
By default, \texttt{@add\_con} declares an equality constraint, since \texttt{lcon} and \texttt{ucon} are both zero.
An inequality constraint is declared by setting the keyword argument \texttt{lcon} or \texttt{ucon}.
The final line instantiates the \texttt{ExaCore} as an \texttt{ExaModel}, the form in which the model is passed to the solver.
The \texttt{SIMDFunction}s and their sparsity structure were already built as the patterns were added.

The example contains no \texttt{for} loop outside the ExaModels.jl syntax.
The objective and constraints are always declared as Julia \texttt{generator} expressions, \texttt{body for i = iterator}, with the loop written inside the macro call.
We follow the Julia convention and call this construct the \texttt{generator}; its iterator component is the \emph{data iterator}.
Because the data iterator is passed to the modeling backend rather than unrolled by the user, each declaration maps directly onto one pattern of \cref{eq:simd}, and the model is stored in the form of the \gls{simd} abstraction.

\paragraph{Solving.} The model is passed to any NLPModels.jl-compatible solver:
\begin{lstlisting}
result = ipopt(model)     # solve via Ipopt (NLPModelsIpopt.jl)
result = madnlp(model)    # solve via MadNLP.jl
result = uno(model)       # solve via Uno (UnoSolver.jl)
\end{lstlisting}
Solving the problem on a GPU requires a GPU-compatible model.
Such a model is obtained by replacing the first line of the example with \texttt{core = ExaCore(minimize = false, backend = CUDABackend())}, with no other changes (\Cref{sec:gpu}).
Passing that model to a GPU-compatible solver such as MadNLP.jl then solves the problem entirely on the GPU.
Solvers in the \gls{jso} ecosystem return their results in a common template defined by SolverCore.jl.
Here \texttt{result} is an instance of a subtype of \texttt{SolverCore.AbstractExecutionStats}.
Its fields carry the solver status and the run statistics, together with the full primal solution (\texttt{result.solution}), the constraint multipliers (\texttt{result.multipliers}), and the bound multipliers (\texttt{result.multipliers\_L} and \texttt{result.multipliers\_U}).
Each of these is a flat vector over all variables or over all constraints.
The accessors below read from this object in the same way regardless of which solver produced it.

\paragraph{Accessing the solution.} From these flat vectors, the ExaModels.jl accessors slice out the part that belongs to a given block, using the block's offset, and return it in the block's shape:
\begin{lstlisting}
# primal: altitude trajectory
h_sol = solution(result, h)

# dual: velocity-dynamics multipliers
lam = multipliers(result, vel)
\end{lstlisting}
The variable-bound multipliers are read back the same way, via \texttt{multipliers\_L(result, h)} and \texttt{multipliers\_U(result, h)}.

\paragraph{Updating and resolving.} Suppose that, after the solve, we wish to re-solve the problem with an updated value of the required final mass, held in the parameter \texttt{mf}.
Such a change leaves the problem structure intact, so the numerical data can be modified in place and the model re-solved without rebuilding it.
Here we raise the required final mass, tighten the thrust upper bound, and warm-start from the previous solution.
\emph{Warm starting} initializes the solver at a previously computed solution, primal or dual or both, so that a re-solve after a small change in the data converges in fewer iterations.
The compiled \texttt{SIMDFunction}s are reused unchanged, and only the parameter, bound, and initial-point arrays are overwritten:
\begin{lstlisting}
# tighten max thrust
set_uvar!(model, tau, fill(0.8 * T_max, nh + 1))

# raise the required final mass
set_value!(model, mf, [1.05 * m_f])

# warm-start from previous solve
set_start!(model, tau, solution(result, tau))
result = ipopt(model)   # resolve, no rebuild
\end{lstlisting}
Because an \texttt{ExaModel} is a standard \texttt{AbstractNLPModel}, such warm-started resolves also work with any repeated-solve interface the chosen solver exposes.
Some solvers reuse their internals across calls through more specialized interfaces, such as a persistent solver object with pre-allocated workspaces; for optimizing re-solve performance, we refer the reader to the documentation of the chosen solver.
The example here demonstrates that the model object itself can be reused, which is most effective when successive solves leave the problem dimensions and the sparsity pattern unchanged.

\subsubsection{Model-Building Macros}
An \texttt{ExaCore} is a caching object: it accumulates the variables, parameters, objective patterns, and constraint patterns as the model is built.
For the full reference to the macros that build an \texttt{ExaCore}, including every signature, keyword argument, and usage example, we refer the reader to the online documentation~\citep{shinExaModelsjl2025}.
Each macro takes the core as its first argument and \emph{rebinds} it.
For example, \texttt{@add\_var(core, x, n)} creates a variable block, binds the handle \texttt{x} to it, and reassigns \texttt{core} to the updated core; the other macros follow the same convention.
The rebinding is needed because an \texttt{ExaCore} is immutable: each call produces a new core (\Cref{sec:data:core}).
The macros are:
\begin{itemize}
  \item \texttt{@add\_var(core, x, dims; start, lvar, uvar)}: adds a variable block \texttt{x} with the given dimensions and optional initial values and bounds.
  \item \texttt{@add\_par(core, p, data)}: adds a parameter block \texttt{p} whose values are appended to the parameter vector $\theta$.
  \item \texttt{@add\_obj(core, expr for i in data)}: adds an objective pattern.
The generator body is the expression tree, and the generator's iterator supplies the data.
  \item \texttt{@add\_con(core, expr for i in data; lcon, ucon)}: adds a constraint pattern that defines a new block with one constraint per element of the data iterator, with bounds \texttt{lcon} and \texttt{ucon}.
  \item \texttt{@add\_con!(core, con[j] += expr for (j, ...) in data)}: augments an existing constraint block \texttt{con} with additive terms.
This corresponds to the $h^{(n)}$ terms in \cref{eq:simd:con} and keeps the number of expression trees small.
In network models, for instance, each contribution type such as generator injection or line flow becomes its own pattern, rather than being inlined into one large expression per constraint.
  \item \texttt{@add\_expr(core, s, expr for i in data)}: adds a named, reusable subexpression block \texttt{s} that can be referenced as \texttt{s[i]} in later objectives and constraints.
This lets the user name an intermediate quantity that recurs across several objectives or constraints, such as a nodal power injection that feeds both an objective term and several balance constraints.
A reference \texttt{s[i]} inlines the subexpression's tree into the consuming expression and introduces no auxiliary variables or constraints.
Because the subexpression is itself a parameterized expression tree (\Cref{sec:graph}), this keeps the model description compact.
The reuse is a modeling convenience rather than a shared runtime computation, since each consuming pattern inlines its own copy of the tree (\Cref{sec:impl:exacore}).
\end{itemize}

\subsubsection{Model Data Accessors}
Once the model is built, its numerical data can be read and modified one block at a time.
In the accessors below, \texttt{x}, \texttt{con}, and \texttt{param} are the handles bound by the user's variable, constraint, and parameter declarations (\Cref{sec:modeling:syntax}):
\begin{itemize}
  \item \texttt{set\_start!(model, x, values)} and \texttt{get\_start(model, x)}: the initial primal values of variable block \texttt{x}.
  \item \texttt{set\_lvar!(model, x, values)}/\allowbreak\texttt{set\_uvar!(model, x, values)} and \texttt{get\_lvar(model, x)}/\allowbreak\texttt{get\_uvar(model, x)}: the variable bounds.
  \item \texttt{set\_lcon!(model, con, values)}/\allowbreak\texttt{set\_ucon!(model, con, values)} and \texttt{get\_lcon(model, con)}/\allowbreak\texttt{get\_ucon(model, con)}: the constraint bounds.
  \item \texttt{set\_value!(model, param, values)} and \texttt{get\_value(model, param)}: the parameter values.
\end{itemize}
Each accessor acts on a single variable, constraint, or parameter block, using the offset and dimension recorded in the block's handle to address the relevant slice without copying.

\subsubsection{Solution Accessors}
The solution reported by the solver is a flat bundle: one vector each for the primal solution $x$, the constraint multipliers $\lambda$, and the bound multipliers, together with the solution status, the objective and constraint values, and solver statistics.
The user, however, usually wants the portion that belongs to one declared variable or constraint block.
ExaModels.jl therefore provides accessors that slice the flat vectors and return the entries of a single block, addressed by the declaration handles \texttt{x} and \texttt{con} and read from the solver's \texttt{result} object introduced in the example above:
\begin{itemize}
  \item \texttt{solution(result, x)}: the primal solution of variable block \texttt{x};
  \item \texttt{multipliers(result, con)}: the multipliers of constraint block \texttt{con};
  \item \texttt{multipliers\_L(result, x)}/\allowbreak\texttt{multipliers\_U(result, x)}: the lower- and upper-bound multipliers of variable block \texttt{x}.
\end{itemize}
The per-block structure of the \gls{simd} abstraction thus carries all the way through to how the user supplies data and reads back the solution.

\section{Core Data Structures}
\label{sec:data}
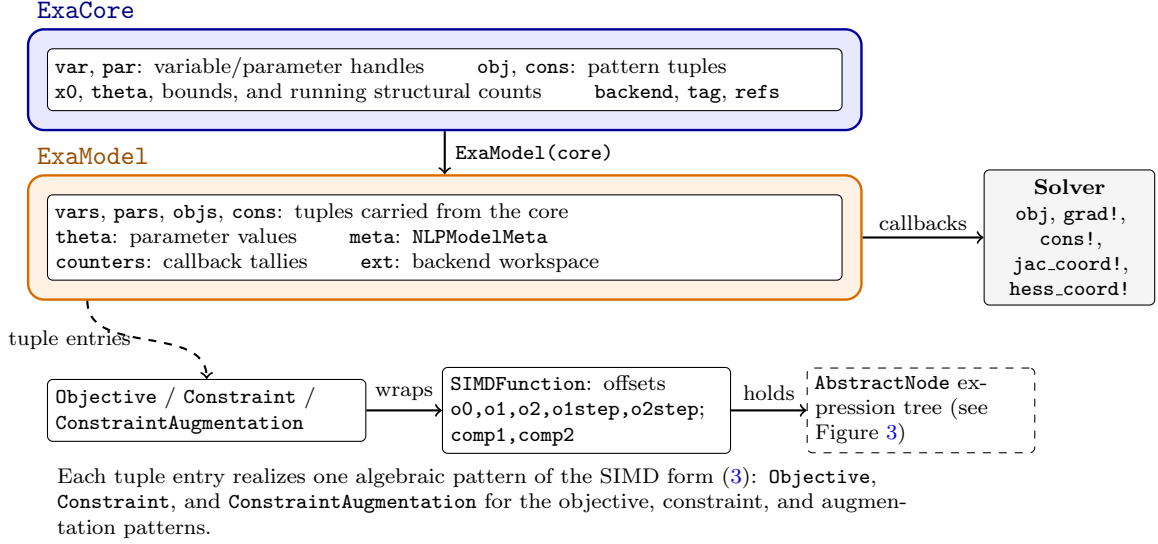
\begin{figure}[t]
  \centering
  \begin{tikzpicture}[
      font=\footnotesize,
      comp/.style={draw, rounded corners=2pt, fill=white, align=left,
        inner sep=3pt, font=\scriptsize},
      lbl/.style={font=\small\bfseries\ttfamily},
      cont/.style={rounded corners=6pt, line width=0.9pt, inner sep=2.5mm,
        fill opacity=0.30, text opacity=1},
      arr/.style={->, thick},
    ]
    \node[comp, text width=103mm] (corefields) {
      \texttt{var}, \texttt{par}: variable/parameter handles \qquad
      \texttt{obj}, \texttt{cons}: pattern tuples\\
      \texttt{x0}, \texttt{theta}, bounds, and running structural counts \qquad
      \texttt{backend}, \texttt{tag}, \texttt{refs}
    };
    \begin{scope}[on background layer]
      \node[cont, draw=blue!60!black, fill=blue!30, fit=(corefields)] (core) {};
    \end{scope}
    \node[lbl, anchor=south west, text=blue!50!black] at (core.north west) {ExaCore};

    \node[comp, text width=103mm, below=11mm of corefields] (modelfields) {
      \texttt{vars}, \texttt{pars}, \texttt{objs}, \texttt{cons}: tuples carried from the core\\
      \texttt{theta}: parameter values \qquad
      \texttt{meta}: \texttt{NLPModelMeta}\\
      \texttt{counters}: callback tallies \qquad
      \texttt{ext}: backend workspace
    };
    \begin{scope}[on background layer]
      \node[cont, draw=orange!85!black, fill=orange!35, fit=(modelfields)] (model) {};
    \end{scope}
    \node[lbl, anchor=south west, text=orange!60!black] at (model.north west) {ExaModel};
    \draw[arr] (core.south) -- node[right, font=\scriptsize] {\texttt{ExaModel(core)}} (model.north);

    \node[draw, rounded corners=2pt, align=center, right=16mm of model.east, font=\scriptsize,
      text width=20mm, fill=black!4] (solver)
      {\textbf{Solver}\\\texttt{obj}, \texttt{grad!},\\\texttt{cons!}, \texttt{jac\_coord!},\\\texttt{hess\_coord!}};
    \draw[arr] (model.east) -- node[midway, above=2pt, fill=white,
      inner xsep=2pt, inner ysep=0.5pt, font=\scriptsize] {callbacks} (solver.west);

    \node[comp, text width=40mm, below=13mm of modelfields.south west, anchor=north west] (entry)
      {\texttt{Objective} / \texttt{Constraint} / \texttt{ConstraintAugmentation}};
    \node[comp, text width=36mm, right=10mm of entry] (simd)
      {\texttt{SIMDFunction}: offsets \texttt{o0,o1,o2,o1step,o2step}; \texttt{comp1,comp2}};
    \node[comp, dashed, text width=28mm, right=10mm of simd] (tree)
      {\texttt{AbstractNode} expression tree (see \Cref{fig:tree})};
    \draw[arr] (entry) -- node[above, font=\scriptsize] {wraps} (simd);
    \draw[arr] (simd) -- node[above, font=\scriptsize] {holds} (tree);
    \draw[arr, dashed] ([xshift=8mm]model.south west) to[out=-90, in=90] node[left, font=\scriptsize, xshift=-2pt] {tuple entries} (entry.north);
    \node[font=\scriptsize, text width=112mm, below=2mm of entry.south west, anchor=north west]
      {Each tuple entry realizes one algebraic pattern of the SIMD form~\cref{eq:simd}: \texttt{Objective}, \texttt{Constraint}, and \texttt{ConstraintAugmentation} for the objective, constraint, and augmentation patterns.};
  \end{tikzpicture}
  \caption{The construction-time and solver-facing containers. Arrows show what \texttt{ExaModel(core)} carries from one to the other, and which wrapper types sit inside both.}
  \label{fig:datastruct}
\end{figure}

The two central data structures in ExaModels.jl both store a model expressed in the \gls{simd} abstraction \cref{eq:simd} of \Cref{sec:modeling:simd} in a standardized form.
An \texttt{ExaCore} is the immutable structure that holds the model as it is being built.
An \texttt{ExaModel} is the immutable completed model, instantiated from a finished \texttt{ExaCore}, and it implements the NLPModels.jl callback interface of \Cref{sec:modeling:nlp} (\Cref{fig:datastruct}).

Inside these two containers sit the data types that represent the individual patterns of \cref{eq:simd}.
Each objective pattern is stored as an \texttt{Objective}, each constraint block as a \texttt{Constraint}, and each additive contribution to an existing constraint block as a \texttt{ConstraintAugmentation}.
Each of the three wrappers holds a \texttt{SIMDFunction} and the data iterator over which the pattern is evaluated.
The \texttt{SIMDFunction} in turn holds the parameterized expression tree for one algebraic pattern, together with the offsets and compressor maps that place its output.
The expression tree itself, and how it is built, is the subject of \Cref{sec:graph}.

Both \texttt{ExaCore} and \texttt{ExaModel} are \emph{parametric structs}.
Their type parameters encode the floating-point type \texttt{T}, the array type \texttt{VT}, the backend \texttt{B}, and the concrete types of all stored variables, parameters, objectives, and constraints.
For example, \texttt{ExaCore\{Float64, Vector\{Float64\}, Nothing, ...\}} is a CPU model in double precision, and \texttt{ExaCore\{Float32, CuArray\{Float32\}, CUDABackend, ...\}} is a single-precision model on an NVIDIA GPU.
Because these choices are captured as type parameters, Julia's compiler specializes the AD passes and solver callbacks for each concrete combination.
A single codebase therefore handles CPU, CUDA, ROCm, oneAPI, and Metal backends with no runtime dispatch and no conditional branching in the evaluation path.

The two structs follow a general principle: ExaModels.jl uses explicit, concrete typing throughout its performance-critical evaluation path.
Every data structure on this path, from the individual expression-tree nodes up to the \texttt{ExaCore} and \texttt{ExaModel} that contain them, has a fully concrete Julia type.
There are no abstractly typed fields, no type-unstable containers, and no runtime polymorphism during callback evaluation.
Concrete types let the compiler specialize and inline the evaluation code completely, and they are also what make \gls{aot} compilation and GPU kernel generation possible (\Cref{sec:aot}).
The price is that even models of moderate complexity carry deeply nested type signatures.
The typing of the expression trees is treated in \Cref{sec:graph:type}.

\subsection{ExaCore: Model-Construction Cache}
\label{sec:data:core}

An \texttt{ExaCore} is the construction-time cache of an optimization model.
It accumulates the variables, parameters, algebraic patterns, numerical arrays, and structural counts needed to instantiate the eventual \texttt{ExaModel}.
It is created by a call to \texttt{ExaCore()}, which optionally takes a GPU backend.
The user then adds variables, parameters, objectives, and constraints with the macros \texttt{@add\_var}, \texttt{@add\_par}, \texttt{@add\_obj}, and \texttt{@add\_con}.
The cache is immutable, so each addition returns an updated \texttt{ExaCore}.
The reconstruction mechanics, and their role in preserving concrete types, are detailed in \Cref{sec:impl:exacore}.

Reflecting the explicit-typing principle, the struct is \texttt{ExaCore\{T, VT, B, S, V, P, O, C, \ldots\}}.
Its key type parameters and fields are:
\begin{itemize}
  \item \texttt{T} and \texttt{VT}: the floating-point type and the array type, respectively, with \texttt{VT} constrained to \texttt{AbstractVector\{T\}}.
  \item \texttt{backend::B}: the KernelAbstractions.jl backend (e.g., \texttt{CUDABackend()}) or \texttt{nothing} for CPU execution.
  \item \texttt{S}: the model tag type, a typed extension point for carrying structural metadata without adding fields to the base core.
For example, a two-stage core uses it to record scenario information and variable/constraint scenario assignments.
  \item \texttt{var::V}: a tuple of \texttt{Variable} handles, each recording the block's dimensions, its offset into the global variable vector $x$, its name, and a tag.
  \item \texttt{par::P}: a tuple of \texttt{Parameter} handles, whose values are concatenated into the global parameter vector $\theta$.
  \item \texttt{obj::O}: a tuple of \texttt{Objective} structs, each bundling a \texttt{SIMDFunction} and the collected data iterator.
  \item \texttt{cons::C}: a tuple of \texttt{Constraint} and \texttt{ConstraintAugmentation} structs.
A \texttt{Constraint} defines a new block of constraints with bounds.
A \texttt{ConstraintAugmentation}, created via \texttt{@add\_con!}, adds terms to an existing constraint block without introducing new rows.
  \item Primal, dual, parameter, and bound arrays, all of type \texttt{VT}: \texttt{x0}, \texttt{\codetheta}, \texttt{lvar}, \texttt{uvar}, \texttt{y0}, \texttt{lcon}, and \texttt{ucon}.
  \item Running dimension counters (all \texttt{Int}): \texttt{nvar} and \texttt{npar} count scalar variables and parameters; \texttt{ncon} counts constraint rows; \texttt{nobj} counts objective-term instances; and \texttt{nconaug} counts augmentation-term instances.
  \item Running sparse-output counts (all \texttt{Int}): \texttt{nnzg} counts objective-gradient contributions, \texttt{nnzj} counts constraint-Jacobian contributions, and \texttt{nnzh} counts Lagrangian-Hessian contributions.
These counts determine output and workspace sizes and are updated as patterns are added.
\end{itemize}
Each \texttt{@add\_obj} or \texttt{@add\_con} call builds its \texttt{SIMDFunction} at once, probing sparsity and building the compressor maps at the moment the pattern is registered.
The updated \texttt{ExaCore} therefore already carries all the information needed for the final model.

\subsection{ExaModel: The Solver-Facing Model}
\label{sec:data:model}

Once all variables, objectives, and constraints have been added, the user calls \texttt{ExaModel(core)} to instantiate the \texttt{ExaCore} as an immutable \texttt{ExaModel}.
The \texttt{ExaModel} carries forward the core's \texttt{var}, \texttt{par}, \texttt{obj}, and \texttt{cons} tuples as the fields \texttt{vars}, \texttt{pars}, \texttt{objs}, and \texttt{cons}, preserving their concrete types while placing them in a form ready for a solver.

The instantiation performs four steps.
\begin{enumerate}
  \item Constructs the \texttt{NLPModelMeta} object.
It stores the problem dimensions $n$ and $m$, the initial primal and dual iterates $x_0$ and $y_0$, the variable and constraint bounds $x^\flat$, $x^\sharp$, $g^\flat$, $g^\sharp$, the bound classification indices (free, lower-bounded, upper-bounded, fixed), and the nonzero counts \texttt{nnzj} and \texttt{nnzh}.
This metadata follows the NLPModels.jl specification, and solvers query it to allocate sparse structures and to select algorithmic options.
  \item Inherits the four tuples and the parameter values \texttt{\codetheta} from the \texttt{ExaCore}.
  \item Builds the backend-specific extension \texttt{ext} when a GPU backend is active.
This is the only step that analyzes the model.
The gradient contributions of all patterns are collected, sorted by target variable, and partitioned by pointer arrays, which yields the sparse-to-dense compression map used by the GPU gradient path (\Cref{sec:impl:gpu}).
The same structure is built for constraint augmentations, and the device-side work buffers are allocated.
The per-iteration callbacks can then launch with no further host-side preparation.
On a CPU, no extension is built, and the instantiation reduces to reorganizing the data above.
  \item Initializes an \texttt{NLPModels.Counters} object to track callback invocations during the solve.
\end{enumerate}

The resulting \texttt{ExaModel} is a subtype of \texttt{NLPModels.AbstractNLPModel\{T,VT\}}, where \texttt{T} is the floating-point type and \texttt{VT} is the array type (\texttt{Vector\{T\}} for CPU, \texttt{CuArray\{T\}} for CUDA, etc.).
Its type has the form \texttt{ExaModel\{T, VT, E, V, P, O, C, \ldots\}}, where \texttt{E} is the type of the backend extension \texttt{ext} (\texttt{Nothing} on CPU) and \texttt{V}, \texttt{P}, \texttt{O}, \texttt{C} are the tuple types carried unchanged from the \texttt{ExaCore}.
As with the core, a few further parameters for tags and references are not shown.
An \texttt{ExaModel} can therefore be passed directly to any solver in the \gls{jso} ecosystem, including Ipopt (via NLPModelsIpopt.jl) and MadNLP.jl~\citep{shinMadNLPjl2025}, without any adapter code.
For each solver callback introduced in \Cref{sec:modeling:nlp}, ExaModels.jl dispatches to loops over the stored \texttt{SIMDFunction}s.

\subsection{Handles and Pattern Wrappers}
\label{sec:data:internal}

The internal types described below connect the two containers above to the expression trees of \Cref{sec:graph} and the AD machinery of \Cref{sec:ad}.
Most of them appear inside both the \texttt{ExaCore} and the \texttt{ExaModel}.
They are the variable and parameter handles, the wrapper types stored in the \texttt{obj} and \texttt{cons} tuples, and the \texttt{SIMDFunction} that each wrapper holds.
The node types underneath the \texttt{SIMDFunction} are the subject of \Cref{sec:graph}.

\subsubsection{Variable and Parameter}
Each call to \texttt{@add\_var} or \texttt{@add\_par} declares a block of variables or parameters, possibly with multiple dimensions.
Indices such as \texttt{x[i]} are local to the declared block, whereas solver callbacks operate on the global vectors $x$ and $\theta$ formed by concatenating all variable and parameter blocks.
Translating a block-local index into a global index therefore requires the block's offset: if a block begins after \texttt{offset} scalar entries, its local entry \texttt{i} is stored at global index \texttt{offset + i}.
The \texttt{var} and \texttt{par} tuples store lightweight \emph{handle} types that record this mapping, rather than the numerical data itself.

A \texttt{Variable\{S,O,T\}}, the handle returned by \texttt{@add\_var}, contains:
\begin{itemize}
  \item \texttt{size::S}: the block's shape (a tuple of dimensions); general tensor blocks of arbitrary dimension are supported, with vectors and matrices as special cases.
  \item \texttt{length::O}: the total number of scalar entries in the block.
  \item \texttt{offset::O}: the offset that translates a block-local variable index into its index in the global primal vector $x$.
  \item \texttt{name::Symbol}: the variable's name, captured automatically from the identifier in \texttt{@add\_var(core, x, ...)} and making the handle accessible as \texttt{core.x}.
  \item \texttt{tag::T}: optional user metadata, e.g.\ a scenario or batch-instance identifier for the two-stage and batch models of \Cref{sec:extensions:twostage,sec:extensions:batch}.
\end{itemize}
Indexing the handle as \texttt{x[i]} inside a later expression produces the symbolic variable reference of \Cref{sec:graph}, and \texttt{solution(result, x)} uses the offset and size to slice the block's values out of the solution vector.

A \texttt{Parameter\{S,O,T\}}, the handle returned by \texttt{@add\_par}, carries the same bookkeeping except for the name:
\begin{itemize}
  \item \texttt{size::S}, \texttt{length::O}: the block's shape and number of entries; as with variables, general tensor blocks of arbitrary dimension are supported.
  \item \texttt{offset::O}: the offset that translates a block-local parameter index into its index in the global parameter vector $\theta$.
  \item \texttt{tag::T}: optional user metadata, as for \texttt{Variable}.
\end{itemize}
Parameter values can be updated with \texttt{set\_value!} without rebuilding the model.
Both handles are pure bookkeeping: they carry offsets and sizes, not arrays, which is consistent with the array-free design of \Cref{sec:graph}.
The numerical data, such as \texttt{x0}, \texttt{lvar}, \texttt{uvar}, and $\theta$, are stored separately in the flat arrays of the \texttt{ExaCore}.

Declaring a block also extends the corresponding numerical arrays by one entry per scalar in the block.
A call to \texttt{@add\_var} appends the block's initial values to \texttt{x0} and its lower and upper bounds to \texttt{lvar} and \texttt{uvar}.
Their defaults are $0$, $-\infty$, and $+\infty$, respectively, and the user may replace each default with a scalar, array, or generator.
Similarly, \texttt{@add\_par} appends the block's values to the global parameter array $\theta$, using $0$ by default.

\subsubsection{Objective, Constraint, and ConstraintAugmentation}
\label{sec:data:wrappers}%
Each entry in the objective and constraint tuples of an \texttt{ExaCore} or \texttt{ExaModel} is one of three wrapper types, each pairing a \texttt{SIMDFunction} with the data over which it is evaluated.
The contents and role of the \texttt{SIMDFunction} are described in \Cref{sec:data:simdfunction}.
These collections are stored as Julia \emph{tuples} rather than arrays, so that the concrete type of every wrapper appears explicitly in the top-level \texttt{ExaCore} or \texttt{ExaModel} type, down to the type of its \texttt{SIMDFunction} and expression tree.
For example, two objective patterns are represented by a heterogeneous tuple such as \texttt{Tuple\{Objective\{F1,I1\}, Objective\{F2,I2\}\}}, which preserves the distinct types \texttt{F1} and \texttt{F2}.
A \texttt{Vector\{Objective\}}, by contrast, would erase those distinctions behind the abstract \texttt{Objective} type and require dynamic dispatch during evaluation.
The \texttt{var} tuple follows the same principle.
Because extending any of these tuples changes its type, each addition reconstructs the core with a new concrete top-level type, as detailed in \Cref{sec:impl:exacore}.

An \texttt{Objective\{F,I\}} corresponds to one objective pattern $\tilde{f}^{(\ell)}$ in \cref{eq:simd:obj} and contains:
\begin{itemize}
  \item \texttt{f::F}: the \texttt{SIMDFunction} constructed from the pattern.
  \item \texttt{itr::I}: the collected data iterator over which the pattern is evaluated.
\end{itemize}
A \texttt{Constraint\{F,I,O,S,T\}} corresponds to one base constraint pattern $\tilde{g}^{(m)}$ in \cref{eq:simd:con} and adds bookkeeping for the block's location:
\begin{itemize}
  \item \texttt{f::F}, \texttt{itr::I}: the \texttt{SIMDFunction} and the data iterator, as above.
  \item \texttt{offset::O}: the block's offset into the global constraint vector.
  \item \texttt{size::S}: the block's dimensions.
  \item \texttt{tag::T}: optional user metadata, as for \texttt{Variable}.
\end{itemize}
A \texttt{ConstraintAugmentation\{F,I,D,T\}} corresponds to an augmentation pattern $\tilde{h}^{(n)}$ in \cref{eq:simd:con}.
Its \texttt{SIMDFunction} wraps an expression tree of the form \texttt{idx => value}, where the index expression selects \emph{which} constraint row to augment and the value expression supplies the additive term:
\begin{itemize}
  \item \texttt{f::F}, \texttt{itr::I}: the \texttt{SIMDFunction} and the data iterator.
  \item \texttt{oa::Int}: the base offset of this augmentation's terms in the GPU augmentation buffer.
Each term is first written to its own position in this buffer and is subsequently accumulated into its target constraint row.
The target constraint block's base row is instead encoded by \texttt{f.o0}.
  \item \texttt{dims::D}: the dimensions of that original constraint block, used to map index expressions to rows.
  \item \texttt{tag::T}: optional user metadata, as for \texttt{Variable}.
\end{itemize}

To see why \texttt{ConstraintAugmentation} is needed, consider the active power balance constraint in \gls{acopf}:
\begin{equation*}
  \sum_{g \in \mathcal{G}_b} p_g - p^d_b = \sum_{(b,j) \in \mathcal{E}} p^{\mathrm{f}}_{bj}  + \sum_{(j,b) \in \mathcal{E}} p^{\mathrm{t}}_{jb}, \quad \forall b \in \mathcal{B},
\end{equation*}
where $\mathcal{G}_b$ is the set of generators at bus $b$, and $p^{\mathrm{f}}_{bj}$ and $p^{\mathrm{t}}_{jb}$ are the flows on the lines incident to $b$, each measured at the $b$ end of its line.
With line losses, the power entering a line at one end differs from the power leaving it at the other, so each end carries its own flow variable.
Each flow variable is the power the bus sends into the line at its end, so every incident line contributes to the balance with the same sign, whether the bus is the from end or the to end of the line.
Each bus has a different number of generators and incident lines, so the full expression cannot be written as a single parameterized tree evaluated over buses.
The tree \emph{structure} varies from bus to bus, which violates the \gls{simd} requirement that all data points share the same tree.

With \texttt{@add\_con!}, the constraint is instead built in layers:
\begin{lstlisting}
# base pattern: one row per bus
@add_con(core, c_p, b.pd for b in bus_data)

# generator injection
@add_con!(core, c_p, g.bus => -pg[g.i] for g in gen_data)

# line flow (from and to ends)
@add_con!(core, c_p, l.from => pf[l.i] for l in line_data)
@add_con!(core, c_p, l.to   => pt[l.i] for l in line_data)
\end{lstlisting}
The accumulated residual is demand plus line flows minus generation, which is the zero-valued form of the balance equation above.
Each \texttt{@add\_con!} call creates a separate \texttt{ConstraintAugmentation} with its own fixed-structure expression tree and its own data iterator.
The augmentation mechanism accumulates all contributions into the correct constraint rows at evaluation time.
Without it, the user would have to write one large expression per bus that inlines every generator and line contribution, which destroys the pattern--data separation and blocks \gls{simd} and \gls{gpu} parallelism.

More broadly, \texttt{@add\_con!} applies whenever a constraint sums over heterogeneous sets.
Reactive power balance in \gls{acopf} has the same structure, as do nodal flow balance constraints in optimal control on networks.
The augmentation mechanism makes the number of distinct patterns scale with the number of \emph{types of contributions} (generator, line-from, line-to), not with the number of constraints or the network topology.

\subsection{SIMDFunction: The Unit of Callback Work}
\label{sec:data:simdfunction}
The parameterized expression tree, described in \Cref{sec:graph}, captures only the algebraic structure of a pattern.
Its \texttt{Var} leaves record which variables the pattern reads (\Cref{sec:graph:indexed}), but the tree knows nothing else about the surrounding \gls{nlp} model.
It does not know whether it contributes to the objective or to the constraints.
It does not know which data iterator it will be evaluated over, where its values and derivative nonzeros belong in the global output arrays, or how repeated derivative contributions should be merged.
Turning the tree into a unit of callback work therefore takes additional bookkeeping.

A \texttt{SIMDFunction\{F,C1,C2\}} attaches the output offsets, the per-data-point strides, and the sparse compressor mappings to the expression tree.
The data iterator is held by the enclosing \texttt{Objective} or \texttt{Constraint}, which pairs it with the \texttt{SIMDFunction}.
Constraints also carry their row and shape metadata.

The \texttt{SIMDFunction} is built at model-construction time by probing the expression tree for sparsity and building the compressor maps.
It is the unit of work for both CPU loops and GPU kernels (\Cref{sec:gpu}), and it contains:
\begin{itemize}
  \item \texttt{f::F}: the parameterized expression tree representing the algebraic pattern.
Its node structure and construction are discussed in \Cref{sec:graph}.
  \item \texttt{comp1::C1}, \texttt{comp2::C2}: \texttt{Compressor} mappings from reverse-pass visits to unique sparse output entries.
The first-order mapping \texttt{comp1} is used for gradient/Jacobian entries, whereas the second-order mapping \texttt{comp2} is used for Lagrangian Hessian entries.
Their construction and role in merging repeated contributions within a pattern are described in \Cref{sec:ad:compressor}.
  \item \texttt{o0}, \texttt{o1}, \texttt{o2}, \texttt{o1step}, \texttt{o2step}: offset and stride integers, one offset per derivative order.
The offsets locate this pattern's entries in the shared output arrays, namely the function values, the gradient and Jacobian nonzeros, and the Lagrangian Hessian nonzeros.
The strides record how far the write position advances from one data point to the next.
Together they let each pattern write directly into its own slice without collisions; the formula for the write position is given in \Cref{sec:ad:callbacks}.
\end{itemize}

\section{Parameterized Expression Tree}
\label{sec:graph}

During model creation, ExaModels.jl captures each user-supplied generator body as a single \emph{parameterized expression tree}: a Julia data structure that represents one algebraic pattern of \cref{eq:simd} and is evaluated at every data point of that pattern.
The word \emph{parameterized} refers to the data-point index \texttt{i}, which is left as an undetermined parameter.
Each leaf that depends on the data point records only \emph{how} to fetch its value once \texttt{i} is known, such as which variable index or which data field, and not the value itself.
Supplying a concrete \texttt{i} at evaluation time instantiates the tree for that data point by a forward pass, without copying it or allocating a new one.
Keeping the decision variables symbolic is standard in symbolic expression-tree generation.
Extending the same treatment to the data-point index is not standard, and it is what lets one tree serve an entire family of expressions rather than one tree per scalar expression.

For a constraint pattern $g^{(m)}$ in \cref{eq:simd:con}, the algebraic template becomes the expression tree stored by a \texttt{SIMDFunction}.
The data points $q_j^{(m)}$ become the data iterator \texttt{itr} stored by the enclosing \texttt{Constraint}, and objective patterns are represented in the same way by \texttt{Objective}.
Evaluating a constraint row instantiates this shared tree at the corresponding element \texttt{i} of \texttt{itr}.
The remaining ingredients of the \texttt{SIMDFunction}, namely the compressor maps and the output offsets, are deferred to \Cref{sec:ad}.

The defining design decision is to keep the array data entirely separate from the tree.
The tree refers to decision variables, parameters, and data-point fields only symbolically, by the index or access path needed to fetch each value.
It therefore carries no array data of its own and is fully immutable.
This is what lets the explicit typing of \Cref{sec:data} reach every node.
Each node type is parameterized by the types of its children, so nesting produces a tree whose type encodes the node itself and its entire subtree, and the type of the root alone determines the algebraic structure of the pattern.
The compiler can then specialize a dispatch-free evaluator for each pattern.
The tree is also a bits type, because it is array-free and concretely typed, so it can be passed directly as a \gls{gpu} kernel argument, with the numerical data supplied alongside as device arrays (\Cref{sec:gpu}).

\subsection{Tree Construction from the Modeling Syntax}
\label{sec:graph:construction:syntax}

A parameterized expression tree is not a distinct data structure but a composition of nested nodes, each holding its children as fields.
The user never writes those nodes.
They are built incrementally by Julia's multiple dispatch as the ordinary expression in an \texttt{@add\_obj} or \texttt{@add\_con} declaration (\Cref{sec:modeling:syntax}) is evaluated.
Each operation in the expression wraps the nodes its operands have already become, so the tree grows outward from the leaves as evaluation proceeds.
\Cref{sec:graph:nodetypes} catalogues the node types in full.

\subsubsection{Building Nodes from the Bottom Up}
\label{sec:graph:bottomup}

Construction starts from the \emph{source nodes}, which stand in for the vectors and the data point that a pattern is written against.
They are \texttt{VarSource} for the decision variables $x$, \texttt{ParameterSource} for the parameters $\theta$, and \texttt{DataSource} for the current data point.
Indexing or field access on a source node does not read a value.
It returns an \emph{indexed node} that records \emph{how} the value will be fetched once a data point is supplied:

\begin{lstlisting}[style=jlrepl]
julia> x = VarSource(); θ = ParameterSource(); i = DataSource();

julia> x[2]
Var
  x[2]

julia> θ[2]
ParameterNode
  θ[2]

julia> i[1]
DataIndexed
  i[1]

julia> i.a
DataIndexed
  i.a
\end{lstlisting}
A node prints as the mathematical expression it stands for, with the decision-variable vector, the parameter vector, and the data point shown as \texttt{x}, \texttt{\codetheta}, and \texttt{i}.
Those names are fixed by the printer rather than taken from the variable a source is bound to, so we bind the sources to the same names here.
The expression \texttt{x[2]} invokes \texttt{getindex(x, 2)}.
Because \texttt{x} is a \texttt{VarSource}, multiple dispatch selects \texttt{getindex(::VarSource, \ldots)} and returns a \texttt{Var} node instead of reading an array element.

The indices used so far have been literals.
The parameterization enters when an indexed node itself serves as the index of a variable or a parameter.
\begin{lstlisting}[style=jlrepl]
julia> x[i[1]]
Var
  x[i[1]]

julia> θ[i.a]
ParameterNode
  θ[i.a]
\end{lstlisting}
Neither expression names a particular entry of $x$ or $\theta$.
Each records the access path by which the entry will be determined once a data point is supplied.
A single node therefore stands for a different variable at every data point, and one tree serves an entire pattern.
This is what the word \emph{parameterized} refers to, and it is the property that the rest of the construction preserves.

These indexed nodes are the ingredients of larger expressions.
Applying an ordinary Julia function or operator to one builds a further node rather than computing a number, because each registered function (\Cref{sec:impl:register}) carries an overload for \texttt{AbstractNode} operands.
The results are the \emph{general nodes} \texttt{Node1} and \texttt{Node2}, for unary and binary operations respectively:
\begin{lstlisting}[style=jlrepl]
julia> x[2]^2
Node1{abs2}
  x[2]^2

julia> sin(x[2])
Node1{sin}
  sin(x[2])

julia> exp(x[1]) + x[2]
Node2{+}
  exp(x[1]) + x[2]

julia> x[i[1]] * θ[i[2]]
Node2{*}
  x[i[1]] * θ[i[2]]
\end{lstlisting}
In each case the function or operator itself is carried in the node's first type parameter, so the algebraic form of the pattern is fixed at compile time.
Applied repeatedly, this yields the complete parameterized expression tree.

\subsubsection{From a Model Declaration to a Tree}
\label{sec:graph:declaration}

The construction above proceeded upward, one expression at a time, with the nodes written out by hand.
From a single \texttt{@add\_obj} or \texttt{@add\_con} statement in the modeling syntax of \Cref{sec:modeling:syntax}, the same machinery produces the tree for a whole pattern without the user naming a single node.

Consider the extended Rosenbrock function:
\begin{equation}\label{eq:rosenbrock}
  \min_{x \in \mathbb{R}^N}\; \sum_{i=2}^{N} \bigl[100(x_{i-1}^2 - x_i)^2 + (x_{i-1} - 1)^2\bigr].
\end{equation}
In ExaModels.jl, this is specified as a single objective pattern:
\begin{lstlisting}
@add_obj(core, 100(x[i-1]^2 - x[i])^2 + (x[i-1] - 1)^2 for i in 2:N)
\end{lstlisting}

Julia lowers the generator into a one-argument function \texttt{f = i -> 100*(x[i-1]\^{}2 - x[i])\^{}2 + (x[i-1] - 1)\^{}2}, paired with the range \texttt{2:N} that supplies its data points.
The name \texttt{x} inside the body is the variable handle declared by \texttt{@add\_var}, which the function captures from the enclosing scope.
To construct a single expression tree for this pattern, ExaModels.jl calls \texttt{f} once, with \texttt{i = DataSource()}.
This one source-node argument turns \texttt{f} into a probe of its own algebraic structure: each operation builds a symbolic node.

The tree is built step by step through multiple dispatch.
When the generator body is evaluated with this source node:
\begin{enumerate}
  \item \texttt{i - 1} dispatches to \texttt{-(::AbstractNode, ::Real)}, returning a \texttt{Node2} with the \texttt{DataSource} and the constant as children;
  \item {\sloppy \texttt{x[i-1]} dispatches to \texttt{getindex} on the variable handle \texttt{x}.
This appends the offset of the variable block within the global decision vector, which is zero in this example, and returns a \texttt{Var} node owning the resulting index subtree;\par}
  \item \texttt{x[i-1]\^{}2} lowers to \texttt{Base.literal\_pow} with the exponent carried as \texttt{Val\{2\}}, applying the simplification rule $(\cdot)^2 \to \operatorname{abs2}(\cdot)$ and returning a \texttt{Node1\{abs2,...\}};
  \item the subtraction, outer square, multiplication by 100, and addition each dispatch similarly, constructing \texttt{Node2} nodes.
\end{enumerate}
No loop is unrolled during model construction: the final result is a single expression tree whose \texttt{Var} leaves resolve variable indices from the data point at evaluation time.
The same tree is evaluated for each $i \in \{2,\ldots,N\}$, with different data producing different variable indices and hence different function values and derivatives.
This is the parameterized aspect: one tree serves all $N-1$ data points, and its evaluation functions are compiled only once.
It is the key to \gls{simd} evaluation and \gls{gpu} acceleration.

The result is a single value whose type spells out the entire algebraic structure of the pattern.
That type is long, so ExaModels.jl abbreviates it when a node is displayed.
The \texttt{fulltype} function prints it in full:
\begin{lstlisting}[style=jlrepl]
julia> N = 10; core = ExaCore(concrete = Val(true));

julia> @add_var(core, x, N); i = DataSource();

julia> fulltype(100*(x[i-1]^2 - x[i])^2 + (x[i-1] - 1)^2)
Node2{typeof(+),Node2{typeof(*),Int64,Node1{typeof(abs2),
Node2{typeof(-),Node1{typeof(abs2),
Var{Node2{typeof(+),Node2{typeof(-),DataSource,Int64},Int64}}},
Var{Node2{typeof(+),DataSource,Int64}}}}},Node1{typeof(abs2),
Node2{typeof(-),Var{Node2{typeof(+),Node2{typeof(-),DataSource,Int64},Int64}},
Int64}}}
\end{lstlisting}
Every operation, every operand, and the index arithmetic inside each \texttt{Var} appear in the type, which is what allows the compiler to specialize the evaluation of this pattern (\Cref{sec:graph:type}).

\begin{figure}[t]
\centering
\begin{tikzpicture}[
  scale=0.8, transform shape,
  x=1.34cm,
  y=1.28cm,
  op/.style={circle, draw=black!48, thick, fill=black!3, align=center,
    minimum size=11mm, inner sep=0.5pt, font=\footnotesize\ttfamily},
  var/.style={circle, draw=j1!85, thick, fill=j1!10, align=center,
    minimum size=11mm, inner sep=0.5pt, font=\footnotesize\ttfamily},
  data/.style={circle, draw=j3!85, thick, fill=j3!10, align=center,
    minimum size=11mm, inner sep=0.5pt, font=\footnotesize\ttfamily},
  const/.style={circle, draw=j4!75, thick, fill=j4!8, align=center,
    minimum size=11mm, inner sep=0.5pt, font=\footnotesize\ttfamily},
  source/.style={circle, draw=j2!85, thick, fill=j2!10, align=center,
    minimum size=11mm, inner sep=0.5pt, font=\footnotesize\ttfamily},
  typ/.style={font=\tiny, text=black!55},
  ann/.style={font=\scriptsize, text=black!60, inner sep=1pt},
  edge/.style={draw=black!38, semithick, line cap=round},
]
\path[use as bounding box] (-5.9,0.7) rectangle (6.1,-9.25);

\node[op]     (root)  at ( 0.0,  0) {\textcolor{black!55}{\tiny Node2}\\[-2pt]$+$};
\node[op]     (mul)   at (-2.6, -1) {\textcolor{black!55}{\tiny Node2}\\[-2pt]$\times$};
\node[op]     (rabs)  at ( 3.3, -1) {\textcolor{black!55}{\tiny Node1}\\[-2pt]\texttt{abs2}};
\node[const]  (c100)  at (-4.1, -2) {\textcolor{black!55}{\tiny Int}\\[-2pt]$100$};
\node[op]     (labs)  at (-1.4, -2) {\textcolor{black!55}{\tiny Node1}\\[-2pt]\texttt{abs2}};
\node[op]     (rsub)  at ( 3.3, -2) {\textcolor{black!55}{\tiny Node2}\\[-2pt]$-$};
\node[op]     (lsub)  at (-1.4, -3) {\textcolor{black!55}{\tiny Node2}\\[-2pt]$-$};
\node[var]    (rvar)  at ( 2.6, -3) {\textcolor{black!55}{\tiny Var}\\[-2pt]$x_{i-1}$};
\node[const]  (c1)    at ( 4.2, -3) {\textcolor{black!55}{\tiny Int}\\[-2pt]$1$};
\node[op]     (lsq)   at (-2.7, -4) {\textcolor{black!55}{\tiny Node1}\\[-2pt]\texttt{abs2}};
\node[var]    (lvar2) at (-0.1, -4) {\textcolor{black!55}{\tiny Var}\\[-2pt]$x_i$};
\node[op]     (ridx)  at ( 2.6, -4) {\textcolor{black!55}{\tiny Node2}\\[-2pt]$+$};
\node[var]    (lvar1) at (-2.7, -5) {\textcolor{black!55}{\tiny Var}\\[-2pt]$x_{i-1}$};
\node[op]     (lidx2) at (-0.1, -5) {\textcolor{black!55}{\tiny Node2}\\[-2pt]$+$};
\node[op]     (rsubi) at ( 1.8, -5) {\textcolor{black!55}{\tiny Node2}\\[-2pt]$-$};
\node[const]  (roff)  at ( 3.5, -5) {\textcolor{black!55}{\tiny Int}\\[-2pt]$0$};
\node[op]     (lidx1) at (-2.7, -6) {\textcolor{black!55}{\tiny Node2}\\[-2pt]$+$};
\node[source] (lds2)  at (-1.1, -6) {\textcolor{black!55}{\fontsize{4pt}{4.6pt}\selectfont DataSource}\\[-2pt]$i$};
\node[const]  (loff2) at ( 0.4, -6) {\textcolor{black!55}{\tiny Int}\\[-2pt]$0$};
\node[source] (rds)   at ( 1.2, -6.4) {\textcolor{black!55}{\fontsize{4pt}{4.6pt}\selectfont DataSource}\\[-2pt]$i$};
\node[const]  (rc1)   at ( 2.5, -6.4) {\textcolor{black!55}{\tiny Int}\\[-2pt]$1$};
\node[op]     (lsubi) at (-3.6, -7) {\textcolor{black!55}{\tiny Node2}\\[-2pt]$-$};
\node[const]  (loff1) at (-2.2, -7) {\textcolor{black!55}{\tiny Int}\\[-2pt]$0$};
\node[source] (lds1)  at (-4.3, -8.2) {\textcolor{black!55}{\fontsize{4pt}{4.6pt}\selectfont DataSource}\\[-2pt]$i$};
\node[const]  (lc1)   at (-2.9, -8.2) {\textcolor{black!55}{\tiny Int}\\[-2pt]$1$};

\node[ann, anchor=west] at ([xshift=3pt]root.east)  {$100(x_{i-1}^2 - x_i)^2 + (x_{i-1} - 1)^2$};
\node[ann, anchor=east] at ([xshift=-3pt]mul.west)  {$100(x_{i-1}^2 - x_i)^2$};
\node[ann, anchor=west] at ([xshift=3pt]labs.east)  {$(x_{i-1}^2 - x_i)^2$};
\node[ann, anchor=east] at ([xshift=-3pt]lsub.west) {$x_{i-1}^2 - x_i$};
\node[ann, anchor=east] at ([xshift=-3pt]lsq.west)  {$x_{i-1}^2$};
\node[ann, anchor=west] at ([xshift=3pt]rabs.east)  {$(x_{i-1} - 1)^2$};
\node[ann, anchor=west] at ([xshift=3pt]rsub.east)  {$x_{i-1} - 1$};
\node[ann, anchor=west] at ([xshift=3pt]ridx.east)  {$(i - 1) + 0$};
\node[ann, anchor=west] at ([xshift=3pt]lidx2.east) {$i + 0$};
\node[ann, anchor=east] at ([xshift=-3pt]lidx1.west) {$(i - 1) + 0$};

\draw[edge] (root) -- (mul);
\draw[edge] (root) -- (rabs);
\draw[edge] (mul) -- (c100);
\draw[edge] (mul) -- (labs);
\draw[edge] (labs) -- (lsub);
\draw[edge] (lsub) -- (lsq);
\draw[edge] (lsub) -- (lvar2);
\draw[edge] (lsq) -- (lvar1);
\draw[edge] (lvar1) -- (lidx1);
\draw[edge] (lidx1) -- (lsubi);
\draw[edge] (lidx1) -- (loff1);
\draw[edge] (lsubi) -- (lds1);
\draw[edge] (lsubi) -- (lc1);
\draw[edge] (lvar2) -- (lidx2);
\draw[edge] (lidx2) -- (lds2);
\draw[edge] (lidx2) -- (loff2);
\draw[edge] (rabs) -- (rsub);
\draw[edge] (rsub) -- (rvar);
\draw[edge] (rsub) -- (c1);
\draw[edge] (rvar) -- (ridx);
\draw[edge] (ridx) -- (rsubi);
\draw[edge] (ridx) -- (roff);
\draw[edge] (rsubi) -- (rds);
\draw[edge] (rsubi) -- (rc1);

\begin{scope}[on background layer]
  \node[draw=j1!48, dashed, rounded corners=5pt, fill=j1!3,
    inner sep=5pt, fit=(lvar1)(lidx1)(lsubi)(loff1)(lds1)(lc1)] {};
\end{scope}
\node[font=\scriptsize\itshape, text=j1!85, anchor=east]
  at ([xshift=-3pt]lvar1.west) {$x[i-1]$};

\end{tikzpicture}
\caption{Parameterized expression tree for the Rosenbrock pattern $100(x_{i-1}^2 - x_i)^2 + (x_{i-1} - 1)^2$. Gray circles are general nodes; colored circles are variable, constant, and source nodes. The dashed box encloses the index subtree of $x[i-1]$.}
\label{fig:tree}
\end{figure}
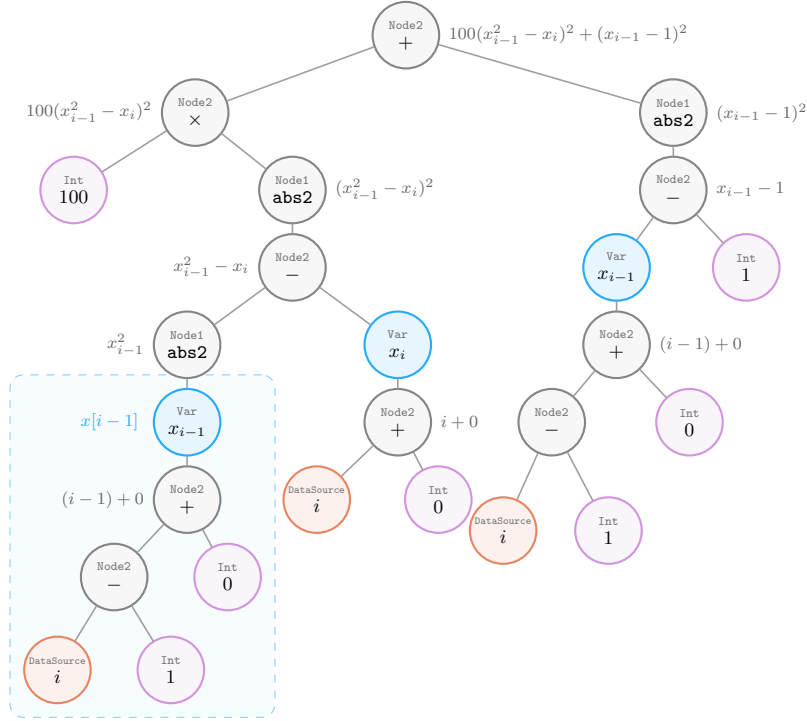

\Cref{fig:tree} shows the expression tree for this pattern.
General nodes carry their operation as a type parameter, and source and constant nodes supply the inputs.
Each \texttt{Var} node contains a subtree that computes the variable index from the data point $i$.
The reference \texttt{x[i-1]} becomes a \texttt{Var} node whose index subtree combines the \texttt{DataSource} leaf with two integer constants and evaluates $(i - 1) + 0$.
The trailing constant is the offset of the variable block \texttt{x} within the global decision vector, which is zero in this example.
At evaluation time the subtree resolves to the index $i - 1$, so the same tree evaluates $100(x_1^2 - x_2)^2 + (x_1 - 1)^2$ when $i = 2$, $100(x_2^2 - x_3)^2 + (x_2 - 1)^2$ when $i = 3$, and so on for the remaining data points.

\subsubsection{Forms of the Data Iterator}
\label{sec:graph:iterator}

The model declarations so far iterate over \texttt{2:N}, the obvious choice when a pattern is indexed by an integer.
The data iterator is not restricted to a range.
It may be any collection whose elements the pattern accesses, because those accesses are captured symbolically.
The \texttt{getindex}, \texttt{getproperty}, and \texttt{indexed\_iterate} methods on \texttt{DataSource} all return \texttt{DataIndexed} nodes, so an element may be addressed by position or by field name.

The elements must be concretely typed, so that the access paths resolve at compile time.
For \gls{gpu} execution they must additionally be bits types, since the data iterator is copied to the device and the kernel receives it as an argument.

A tuple carries the data by position, and the element may be indexed or destructured:
\begin{lstlisting}
data = [(100.0, k-1, k) for k in 2:N]

@add_obj(core, d[1]*(x[d[2]]^2 - x[d[3]])^2 + (x[d[2]] - 1)^2 for d in data)
@add_obj(core, c*(x[i]^2 - x[j])^2 + (x[i] - 1)^2 for (c, i, j) in data)
\end{lstlisting}
The two declarations produce the same tree: destructuring dispatches to \texttt{indexed\_iterate} on the \texttt{DataSource}, which returns the same \texttt{DataIndexed} nodes that explicit indexing does.
A named tuple carries the data by name instead, which is clearer when a pattern has several fields:
\begin{lstlisting}
data = [(coef = 100.0, i = k-1, j = k) for k in 2:N]

@add_obj(core, d.coef*(x[d.i]^2 - x[d.j])^2 + (x[d.i] - 1)^2 for d in data)
\end{lstlisting}
A user-defined \texttt{struct} behaves exactly as the named tuple does, since the field accesses resolve through the same \texttt{getproperty} overload:
\begin{lstlisting}
struct Arc
    coef::Float64
    i::Int
    j::Int
end

data = [Arc(100.0, k-1, k) for k in 2:N]

@add_obj(core, a.coef*(x[a.i]^2 - x[a.j])^2 + (x[a.i] - 1)^2 for a in data)
\end{lstlisting}
Nested combinations of these forms work in the same way.
In each case the field name or position becomes part of the \texttt{DataIndexed} type parameter, so the access path is fixed at compile time and involves no runtime lookup.

\subsection{Node Types}
\label{sec:graph:nodetypes}

The node types fall into three families: the leaves, which either stand for the vectors and the data point or carry a constant value; the indexed nodes obtained from the source nodes; and the general nodes that apply registered functions and operators.

\subsubsection{Leaf Nodes}
\label{sec:graph:sources}

The leaves are of two kinds: the \emph{source nodes}, which record where a value will be fetched once a data point is supplied, and constants, which carry a value directly.
\begin{itemize}
  \item \texttt{VarSource}: the decision-variable vector $x$.
Indexing it produces a \texttt{Var} node.
  \item \texttt{ParameterSource}: the parameter vector $\theta$.
Indexing it produces a \texttt{ParameterNode}, and it appears only in models that carry parameters (\Cref{sec:modeling:nlp}).
  \item \texttt{DataSource}: the current data point.
It admits field access as well as indexing, so a structured data point can be addressed by name, and either form produces a \texttt{DataIndexed} node.
  \item \texttt{Real}: a literal coefficient is not wrapped at all.
It is carried as a plain Julia number in the child field of the general node above it, which is why a general node's children are typed as \texttt{AbstractNode} or \texttt{Real}.
  \item \texttt{Constant\{V\}}: a constant that carries its value \texttt{V} in its type rather than in a field, so the value is available at compile time (\Cref{sec:impl:constants}).
\end{itemize}
\texttt{VarSource} and \texttt{ParameterSource} exist only to be indexed; a pattern never refers to the whole vector $x$ or $\theta$ at once.
A \texttt{DataSource} is different: when the data point is itself the quantity a pattern needs, it is used directly rather than indexed.
The Rosenbrock declaration of \Cref{sec:graph:declaration} is such a case, since its data iterator \texttt{2:N} yields plain integers and the pattern uses \texttt{i} as an index into \texttt{x} rather than reading a field of it.
The \texttt{Var}, \texttt{ParameterNode}, and \texttt{DataIndexed} nodes that indexing produces are the subject of \Cref{sec:graph:indexed}.

\subsubsection{Indexed Nodes}
\label{sec:graph:indexed}

The indexed nodes produced from the source nodes are:
\begin{itemize}
  \item \texttt{Var\{I\}}: a decision-variable access \texttt{x[i]}, where the index type \texttt{I} may be a literal or a data-dependent subtree, so the index itself is computed from the data point.
  \item \texttt{DataIndexed\{I,J\}}: the lookup \texttt{inner.J} or \texttt{inner[J]}, where \texttt{I} is the type of the inner node being indexed into and \texttt{J} is the field name or index itself.
Indexing the data point \texttt{i} at position $1$ gives \texttt{DataIndexed\{DataSource, 1\}}, and accessing its field \texttt{a} gives \texttt{DataIndexed\{DataSource, :a\}}.
The inner node may itself be a \texttt{DataIndexed}, so a nested access path such as \texttt{i.a[2]} is expressed by nesting the type, as \texttt{DataIndexed\{DataIndexed\{DataSource, :a\}, 2\}}.
Since \texttt{J} is carried \emph{in the type}, the access compiles to a direct \texttt{getfield} or \texttt{getindex} with no runtime lookup.
  \item \texttt{ParameterNode\{I\}}: a model-parameter access \texttt{\codetheta[i]}, produced by indexing the \texttt{ParameterSource} leaf.
As with \texttt{Var\{I\}}, the index type \texttt{I} may be a literal or a data-dependent subtree, so a parameter entry can also be selected by the data point.
\end{itemize}

\subsubsection{General Nodes}
\label{sec:graph:nodes}

The general nodes come in two forms, distinguished by the number of arguments:
\begin{itemize}
  \item \texttt{Node1\{F,I\}}: a \emph{univariate} operation \texttt{f} with \texttt{F = typeof(f)} and a single child subtree of type \texttt{I}, evaluating as \texttt{f(inner(i, x, \codetheta))};
  \item \texttt{Node2\{F,I1,I2\}}: a \emph{bivariate} operation \texttt{f} with \texttt{F = typeof(f)} and two child subtrees of types \texttt{I1}, \texttt{I2}, evaluating as \texttt{f(inner1(i, x, \codetheta), inner2(i, x, \codetheta))}.
\end{itemize}
ExaModels.jl supports only univariate and bivariate functions as building blocks, so an $n$-ary operation such as a sum of several terms is represented as a chain of nested nodes rather than a single wide node.

\label{sec:graph:type}The type parameters are what make this representation more than a bookkeeping device.
The operation is carried in \texttt{F}, so two nodes with different operations are different concrete types, and the child types \texttt{I} (or \texttt{I1}, \texttt{I2}) record the subtrees in full.
Together with the parameters carried by \texttt{Var\{I\}}, \texttt{DataIndexed\{I,J\}}, and \texttt{Constant\{V\}}, the entire nesting structure of an expression is therefore part of its type.
The expression $x_{i-1}^2$, for instance, has a type of the form
\begin{center}
\texttt{Node1\{typeof(abs2),\, Var\{Node2\{typeof(+),\, \ldots\}\}\}}.
\end{center}
This is the explicit-typing principle of \Cref{sec:data} carried furthest: a complete pattern is one deeply nested concrete type, with no abstractly typed field anywhere on the evaluation path.
The compiler can therefore inline and specialize the evaluation code for each pattern.
Different patterns are different types, so it produces a separate specialized routine for each, which is what makes \gls{simd} evaluation and \gls{gpu} kernel compilation possible.

\subsubsection{The Root}
\label{sec:graph:root}

The outermost node of the tree is called the \emph{root}, but it is not a distinct node type.
It is simply the general, indexed, or constant node that represents the complete algebraic pattern and is stored in the \texttt{f} field of the \texttt{SIMDFunction} (\Cref{sec:data:internal}).
Its role is to be the entry point: evaluating the pattern at a data point means calling the root, which calls its children in turn.
\Cref{sec:graph:eval} describes that traversal.

\subsection{Node Evaluation}
\label{sec:graph:eval}

Every \texttt{node <: AbstractNode} is callable with the standard evaluation arguments, \texttt{node(i, x, \codetheta)}, and the call returns a scalar.
Here \texttt{i} is the data point, that is, the element of the data iterator, \texttt{x} is the variable vector, and \texttt{\codetheta} is the parameter vector, which is unused when the model carries no parameters.
Every node, and hence every expression tree, is evaluated through calls of this form.
For example, a tree can be constructed and evaluated directly, and the same tree returns a different value at each data point:
\begin{lstlisting}[style=jlrepl]
julia> x = VarSource(); i = DataSource();

julia> node = 100 * (x[i-1]^2 - x[i])^2 + (x[i-1] - 1)^2;

julia> node(2, [2.0, 3.0, 5.0], nothing)
101.0

julia> node(3, [2.0, 3.0, 5.0], nothing)
1604.0
\end{lstlisting}
Arithmetic on the source nodes assembles the tree by operator overloading, and each call evaluates that one tree at the given data point and variable values.
The data point selects which entries of \texttt{x} the \texttt{Var} leaves read, which is why the two calls return different values.
This tree reads no parameters, so the third argument is unused and may be passed as \texttt{nothing}.
Evaluation descends from the root to the leaves: each general node evaluates its children and applies its operation, while the leaves read from the arguments.
A \texttt{Var} node, for example, reads its entry of \texttt{x}, and a \texttt{DataIndexed} node extracts its field or entry from the data point \texttt{i}.

At the individual node level, evaluation applies the function associated with the node type to the recursively evaluated children.
For example, for every registered bivariate function \texttt{f}, the registration macro (\Cref{sec:impl:register}) generates the \texttt{Node2} evaluation method
\begin{lstlisting}
@inline (n::Node2{typeof(f),I1,I2})(i, x, θ) where {I1,I2} =
    f(n.inner1(i, x, θ), n.inner2(i, x, θ))
\end{lstlisting}
so calling a \texttt{Node2} first evaluates its two child subtrees \texttt{inner1} and \texttt{inner2} through the same interface and then applies \texttt{f} to their values.
A \texttt{Node1} evaluates in the same way with a single child.
When a child is a plain \texttt{Real} coefficient, a specialized method uses the stored value directly instead of recursing.
Values therefore flow back up from the leaves to the root, and the call on the root returns the instantiation of the tree at the given data entry \texttt{i}, variables \texttt{x}, and parameters \codetheta.

Every use of the tree goes through this same call; only the type of the arguments varies.
Called with numerical arguments, it computes the function value, which is the regular evaluation path used by the objective and constraint callbacks.
Called with a derivative-carrying \emph{adjoint source} in place of \texttt{x}, the identical traversal builds the computation graph from which the derivatives are evaluated, for both first- and second-order derivatives.
This is the basis of the forward pass of the \gls{ad} algorithm of \Cref{sec:ad}.

\section{Sparse Automatic Differentiation}
\label{sec:ad}

The same evaluation call \texttt{node(i, x, \codetheta)} that instantiates the tree also produces derivatives, once \texttt{x} is replaced by an input type that carries derivative information.
The operator overloads then propagate derivatives through the same expression structure.

The derivative evaluation itself (\Cref{sec:ad:first}) is standard reverse-mode AD.
The one difference is that it runs on the \emph{parameterized} expression tree, so the same compiled code serves every data point of a pattern and the data points can be processed in parallel (\Cref{sec:gpu}).
Second derivatives are obtained by reverse over reverse, rather than by the forward-over-reverse composition that machine learning frameworks generally use.

The larger departure from common practice is in how the sparse Jacobian and Hessian are obtained.
The usual route in a machine learning framework is to expose Jacobian--vector and Hessian--vector products.
The sparse matrix is recovered from a sequence of such products, and the sequence is chosen by a graph coloring of the sparsity pattern.
ExaModels.jl instead evaluates the nonzero entries directly, in a single pass per pattern.
The sparsity pattern is determined \emph{a priori} from the expression tree (\Cref{sec:ad:sparsity}), so no matrix--vector products are formed and no coloring is computed.

\subsection{Derivative Evaluation}
\label{sec:ad:first}

We focus on the first-order case, the gradient.
The second-order case follows the same structure and is described at the end of this subsection.
Consider one term of the extended Rosenbrock function \cref{eq:rosenbrock}, written in the two variables it involves,
\begin{align*}
  f(x_1, x_2) = 100(x_1^2 - x_2)^2 + (x_1 - 1)^2,
\end{align*}
evaluated through the intermediates
\begin{align*}
  s = x_1^2 - x_2, \qquad t = x_1 - 1, \qquad f = 100 s^2 + t^2.
\end{align*}
Reverse-mode AD associates with every intermediate quantity $v$ an \emph{adjoint} $\bar{v} = \partial f / \partial v$, the sensitivity of the output to that intermediate.
The adjoint of the output itself is $\bar{f} = 1$.
The chain rule then propagates adjoints backward through each elementary operation.
If $v = F(u_1, u_2)$, each argument $u_j$ receives the contribution $\bar{v} \cdot \partial F / \partial u_j$, the adjoint of the result times the local partial derivative of the operation.
In the example, the backward traversal proceeds as
\begin{align*}
  \bar{f} &= 1, \qquad \bar{s} = 200 s, \qquad \bar{t} = 2 t,\\
  \bar{x}_1 &= \bar{s} \cdot 2 x_1 + \bar{t} \cdot 1 = 400 x_1 (x_1^2 - x_2) + 2 (x_1 - 1), \qquad \bar{x}_2 = -\bar{s} = -200 (x_1^2 - x_2),
\end{align*}
which are the entries of $\nabla f$.
A single backward traversal of these adjoints therefore yields the full gradient at a cost proportional to one function evaluation.

ExaModels.jl follows this standard scheme, in the same two phases: a \emph{forward pass} that evaluates the parameterized expression tree of \Cref{sec:graph} and records the local derivatives along the way, and a \emph{reverse pass} that propagates adjoints from the root back to the leaves.

\subsubsection{Forward Pass}
\label{sec:ad:forward}
The forward pass uses the node-evaluation machinery of \Cref{sec:graph:eval}.
It is the ordinary evaluation call \texttt{node(i, x, \codetheta)}, with an \texttt{AdjointNodeSource} passed in place of \texttt{x}.
No separate symbolic construction takes place.
The evaluation produces the adjoint tree: a structure with the same shape as the original expression tree, in which each node carries its primal value and its local partial derivatives.
The construction reuses the multiple dispatch mechanism of the primal tree, with the local derivative $F'$ supplied by the function's registered first derivative (\Cref{sec:impl:register}).

The first-order source \texttt{AdjointNodeSource} is a thin wrapper around the primal vector $x$.
The pass evaluates the stored tree as \texttt{node(i, AdjointNodeSource(x), \codetheta)}.
Indexing that source returns an \texttt{AdjointNodeVar(i, x[i])} leaf, a variable node that also carries its primal value.
As the operator overloads run on these leaves, each interior operation builds an \texttt{AdjointNode1} or \texttt{AdjointNode2} node (both \texttt{<: AbstractAdjointNode}) storing, alongside the primal value, the local partial derivatives of its operation:
\begin{itemize}
  \item \texttt{AdjointNodeVar\{I,T\}}: stores the variable index $i$ and the primal value $x_i$.
  \item \texttt{AdjointNode1\{F,T,I\}}: for a unary $F$, stores $v = F(u)$ and $\partial F/\partial u = F'(u)$, where $u$ is the child's primal value.
  \item \texttt{AdjointNode2\{F,T,I1,I2\}}: for a binary $F$, stores $v = F(u_1, u_2)$ and both partial derivatives $\partial F/\partial u_1$, $\partial F/\partial u_2$.
\end{itemize}
The second-order case is identical in structure: a \texttt{SecondAdjointNodeSource} seeds \texttt{SecondAdjointNodeVar} leaves, and the overloads build \texttt{<: AbstractSecondAdjointNode} nodes that also store second derivatives (\Cref{sec:ad:second}).
In both cases the adjoint tree is built during the forward pass and traversed back during the reverse pass.
\Cref{fig:ad} illustrates both passes on the Rosenbrock pattern.

Evaluating $x_1^2 x_2$ on an \texttt{AdjointNodeSource} returns the root of an adjoint tree rather than a number.
Its fields hold the value of the subexpression and the local partial derivatives with respect to its children:
\begin{lstlisting}[style=jlrepl]
julia> node = x[1]^2 * x[2];

julia> g = node(nothing, AdjointNodeSource([3.0, 5.0]), nothing);

julia> g.x, g.y1, g.y2      # value, and the two local partials
(45.0, 5.0, 9.0)

julia> g.inner1.x, g.inner1.y     # the x[1]^2 child: value and d/dx[1]
(9.0, 6.0)
\end{lstlisting}
No global derivative has been computed at this point: each node knows only its own value and how it responds to its immediate children.
Combining these local factors into the gradient is the work of the reverse pass.

\subsubsection{Reverse Pass}
\label{sec:ad:reverse}
The reverse pass traverses the adjoint tree from root to leaves, propagating the adjoint (sensitivity) $\bar{v}$.
The traversal is implemented by a family of reverse-pass functions that share this structure and differ only in where the results are written: \texttt{drpass} (dense gradient), \texttt{grpass} (sparse gradient), \texttt{jrpass} (sparse Jacobian), and the second-order \texttt{hrpass} and \texttt{hdrpass} (sparse Hessian).
Each is discussed where it is used.
Taking \texttt{drpass} as the representative:
\begin{itemize}
  \item At the root, $\bar{v}$ is set to the seed adjoint: $1$ for the gradient and the Jacobian, and the multiplier $\lambda_j$ of constraint $j$ for the Hessian of the Lagrangian.
  \item At \texttt{AdjointNode1\{F\}} with stored derivative $F'$: propagate $\bar{v} \cdot F'$ to the child.
  \item At \texttt{AdjointNode2\{F\}} with stored partials $\partial_1 F$, $\partial_2 F$: propagate $\bar{v} \cdot \partial_1 F$ to child~1 and $\bar{v} \cdot \partial_2 F$ to child~2.
  \item At \texttt{AdjointNodeVar} with index $i$: accumulate $\bar{v}$ into the gradient entry $(\nabla f)_i$.
\end{itemize}

Returning to the pattern above, the reverse pass multiplies these local factors along every path from the root and accumulates the products at the leaves:
\begin{lstlisting}[style=jlrepl]
julia> y = zeros(2);

julia> drpass(g, y, 1.0)    # 1.0 is the seed adjoint at the root

julia> y
2-element Vector{Float64}:
 30.0
  9.0
\end{lstlisting}
matching $\nabla(x_1^2 x_2) = (2 x_1 x_2,\; x_1^2) = (30, 9)$ at $x = (3, 5)$.
Here the gradient is dense and each leaf writes to its own variable index.
For the sparse Jacobian and Hessian, the reverse pass writes instead into the nonzero positions of a \gls{coo} array, which presupposes that those positions are already known.
How they are determined, without graph coloring, is the subject of \Cref{sec:ad:sparsity}.

\subsubsection{Second-Order Derivatives}
\label{sec:ad:second}

The Hessian of the Lagrangian requires second-order derivatives, which are propagated in the same fashion as the first-order case.
The expression tree is evaluated with \texttt{SecondAdjointNodeSource(x)}, producing a second-adjoint tree whose nodes also store the second derivatives of their operations ($h = F''$ for unary nodes; $h_{11}$, $h_{12}$, $h_{22}$ for binary nodes).
The reverse pass then propagates the pair of first and second adjoints $(\bar{v}, \bar{v}_2)$ from root to leaves.
At a unary node with $y = F'$ and $h = F''$, the child receives $\bar{v} \cdot y$ and $\bar{v}_2 \cdot y^2 + \bar{v} \cdot h$, which is the work of \texttt{hrpass}.
The companion \texttt{hdrpass} traverses pairs of subtrees and accumulates the mixed contributions $\bar{v}_2 \cdot y_1 y_2 + \bar{v} \cdot h_{12}$.
At the leaves, \texttt{hrpass} accumulates $\bar{v}_2$ into the diagonal Hessian entry $H_{i,i}$.
The function \texttt{hdrpass} accumulates into the entry $H_{i,j}$ when its two subtrees reach variable leaves with indices $i$ and $j$.
The second-derivative work is skipped on affine nodes, whose $h$ terms are zero, as described in \Cref{sec:impl:affine}.

\subsubsection{Worked Example: Gradient of the Rosenbrock Pattern}
\label{sec:ad:example}

\Cref{fig:ad} recalls the parameterized tree of \Cref{fig:tree} for the pattern of \cref{eq:rosenbrock} and instantiates it at a concrete data point.
For $i = 2$ and $x = [2, 3, \ldots]$ the variable references resolve to $x_1$ and $x_2$, and the tree evaluates to $f_2 = 101$.
We let $s = x_{i-1}^2 - x_i$ and $t = x_{i-1} - 1$, so that one term of the objective is $f_i = 100 s^2 + t^2$.

The forward pass evaluates this tree bottom-up, with $x$ replaced by an \texttt{AdjointNodeSource} so that every node also records its local derivatives (\Cref{fig:ad}).
The inner \texttt{abs2} node over $x_{i-1}$ stores $v = x_{i-1}^2$ and $F' = 2 x_{i-1}$.
The subtraction node $s = x_{i-1}^2 - x_i$ stores $\partial_1 = 1$ and $\partial_2 = -1$.
The outer \texttt{abs2} stores $v = s^2$ and $F' = 2s$.
The multiplication by $100$ stores $\partial_1 = s^2$ and $\partial_2 = 100$.
The $t = x_{i-1} - 1$ and $t^2$ nodes store their partials in the same way.
The root $+$ node holds the primal value $v = 100 s^2 + t^2$.

The reverse pass then traverses the same tree from the root down (\Cref{fig:ad}), seeding the root adjoint $\bar{v} = 1$ and multiplying by each node's stored local derivative as it descends:
\begin{enumerate}
  \item the root $+$ sends $\bar{v} \cdot 1 = 1$ to each of its two children (the $100 s^2$ branch and the $t^2$ branch);
  \item at the multiplication node the $100 s^2$ branch receives $1 \cdot 100 = 100$, which the outer \texttt{abs2} scales by $F' = 2s$ to send $200s$ into the subtraction node;
  \item the subtraction propagates $200s \cdot 1 = 200s$ to the $x_{i-1}^2$ subtree and $200s \cdot (-1) = -200s$ to the leaf $x_i$, accumulating $-200s$ into $(\nabla f)_{x_i}$;
  \item the inner \texttt{abs2} scales by its stored $F' = 2 x_{i-1}$, so the leaf $x_{i-1}$ receives $200s \cdot 2 x_{i-1} = 400 s\, x_{i-1}$;
  \item on the other branch, the \texttt{abs2} node applied to $t$ sends $1 \cdot 2t = 2t$ through the $t = x_{i-1} - 1$ node to the leaf $x_{i-1}$.
\end{enumerate}
Accumulating the two contributions at the shared leaf $x_{i-1}$ yields $\partial f_i / \partial x_{i-1} = 400 s\, x_{i-1} + 2t$ and $\partial f_i / \partial x_i = -200s$, which is the gradient of $f_i$.
At the data point of \Cref{fig:ad}, where $s = t = 1$ and $x_1 = 2$, these are $\partial f_2 / \partial x_1 = 802$ and $\partial f_2 / \partial x_2 = -200$.

\begin{sidewaysfigure}
\centering
\begin{tikzpicture}[
  scale=0.95, transform shape,
  nd/.style={circle, draw=black!48, thick, fill=black!3,
    minimum size=11mm, inner sep=0.5pt, font=\footnotesize\ttfamily},
  ndv/.style={nd, align=center, minimum size=11mm},
  ndvar/.style={nd, draw=j1!85, fill=j1!10},
  ndvarv/.style={ndvar, align=center, minimum size=11mm},
  ndconst/.style={nd, draw=j4!75, fill=j4!8},
  fwd/.style={font=\footnotesize, text=ForestGreen},
  rev/.style={font=\footnotesize, text=Maroon},
  edge from parent/.style={draw, black!38, semithick},
  level 1/.style={sibling distance=40mm, level distance=15mm},
  level 2/.style={sibling distance=24mm, level distance=15mm},
  level 3/.style={sibling distance=20mm, level distance=15mm},
  level 4/.style={sibling distance=15mm, level distance=15mm},
  level 5/.style={sibling distance=13mm, level distance=14mm},
]

\node[font=\bfseries, align=center] at (0,1.2) {Parameterized tree\\{\normalfont\small (same as \Cref{fig:tree}, leaf nodes simplified)}};

\node[nd] (Lplus) at (0,0) {$+$}
  child { node[nd] {$\times$}
    child { node[ndconst] {100} }
    child { node[nd] {abs2}
      child { node[nd] {$-$}
        child { node[nd] {abs2}
          child { node[ndvar] {$x_{i-1}$} }
        }
        child { node[ndvar] {$x_i$} }
      }
    }
  }
  child { node[nd] {abs2}
    child { node[nd] {$-$}
      child { node[ndvar] {$x_{i-1}$} }
      child { node[ndconst] {1} }
    }
  };

\node[font=\bfseries, align=center] at (10,1.7) {Instantiation ($i=2$,\; $x=[2,3,\ldots]$)\\ with forward and reverse passes};

\begin{scope}[
  level 1/.style={sibling distance=68mm, level distance=16mm},
  level 2/.style={sibling distance=34mm, level distance=16mm},
  level 3/.style={sibling distance=28mm, level distance=16mm},
  level 4/.style={sibling distance=20mm, level distance=16mm},
]
\node[ndv] (plus) at (10,0) {$+$\\\textcolor{ForestGreen}{\scriptsize 101}}
  child { node[ndv] (times) {$\times$\\\textcolor{ForestGreen}{\scriptsize 100}}
    child { node[ndconst] (c100) {100} }
    child { node[ndv] (sq1) {abs2\\\textcolor{ForestGreen}{\scriptsize 1}}
      child { node[ndv] (sub1) {$-$\\\textcolor{ForestGreen}{\scriptsize 1}}
        child { node[ndv] (sq2) {abs2\\\textcolor{ForestGreen}{\scriptsize 4}}
          child { node[ndvarv] (va1) {$x_1$\\\textcolor{ForestGreen}{\scriptsize 2}} }
        }
        child { node[ndvarv] (vb) {$x_2$\\\textcolor{ForestGreen}{\scriptsize 3}} }
      }
    }
  }
  child { node[ndv] (sq3) {abs2\\\textcolor{ForestGreen}{\scriptsize 1}}
    child { node[ndv] (sub2) {$-$\\\textcolor{ForestGreen}{\scriptsize 1}}
      child { node[ndvarv] (va2) {$x_1$\\\textcolor{ForestGreen}{\scriptsize 2}} }
      child { node[ndconst] (c1) {1} }
    }
  };

\node[fwd, left=4pt of times, anchor=east] {$\partial_1 = 1,\; \partial_2 = 100$};
\node[fwd, right=4pt of sq1, anchor=west] {$F' = 2$};
\node[fwd, left=4pt of sub1, anchor=east] {$\partial_1 = 1,\; \partial_2 = {-}1$};
\node[fwd, left=4pt of sq2, anchor=east] {$F' = 4$};
\node[fwd, right=4pt of sq3, anchor=west] {$F' = 2$};
\node[fwd, right=4pt of sub2, anchor=west] {$\partial_1 = 1,\; \partial_2 = {-}1$};

\node[rev, above=2pt of plus, anchor=south] {$\bar{v} = 1$};

\def\revL#1#2#3#4{\draw[->, Maroon, semithick, dashed]
  ($($(#1)!6.5mm!(#2)$)!3.5pt!90:(#2)$) --
  node[auto, rev, pos=0.4, inner sep=2pt, #4] {#3}
  ($($(#2)!6.5mm!(#1)$)!3.5pt!-90:(#1)$);}
\def\revR#1#2#3#4{\draw[->, Maroon, semithick, dashed]
  ($($(#1)!6.5mm!(#2)$)!3.5pt!-90:(#2)$) --
  node[auto, swap, rev, pos=0.4, inner sep=2pt, #4] {#3}
  ($($(#2)!6.5mm!(#1)$)!3.5pt!90:(#1)$);}
\def\fwdL#1#2{\draw[->, ForestGreen, semithick]
  ($($(#2)!6.5mm!(#1)$)!3.5pt!90:(#1)$) --
  ($($(#1)!6.5mm!(#2)$)!3.5pt!-90:(#2)$);}
\def\fwdR#1#2{\draw[->, ForestGreen, semithick]
  ($($(#2)!6.5mm!(#1)$)!3.5pt!-90:(#1)$) --
  ($($(#1)!6.5mm!(#2)$)!3.5pt!90:(#2)$);}
\fwdL{plus}{times}
\fwdR{plus}{sq3}
\fwdR{times}{sq1}
\fwdR{sq1}{sub1}
\fwdL{sub1}{sq2}
\fwdR{sub1}{vb}
\fwdL{sq2}{va1}
\fwdR{sq3}{sub2}
\fwdL{sub2}{va2}

\revL{plus}{times}{$1{\cdot}1$}{}
\revR{plus}{sq3}{$1{\cdot}1$}{}
\revR{times}{sq1}{$1{\cdot}100$}{}
\revR{sq1}{sub1}{$100{\cdot}2$}{}
\revL{sub1}{sq2}{$200{\cdot}1$}{swap, xshift=-3pt}
\revR{sub1}{vb}{$200{\cdot}({-}1)$}{swap, xshift=3pt}
\revL{sq2}{va1}{$200{\cdot}4$}{swap, xshift=-3pt}
\revR{sq3}{sub2}{$1{\cdot}2$}{}
\revL{sub2}{va2}{$2{\cdot}1$}{}

\node[rev, font=\footnotesize\bfseries, below=4pt of va1, anchor=north] {$(\nabla f)_{x_1} \mathrel{+}= 800$};
\node[rev, font=\footnotesize\bfseries, below=4pt of vb, anchor=north] {$(\nabla f)_{x_2} \mathrel{+}= {-}200$};
\node[rev, font=\footnotesize\bfseries, below=4pt of va2, anchor=north] {$(\nabla f)_{x_1} \mathrel{+}= 2$};

\node[anchor=north, font=\scriptsize] at ([yshift=-5mm]current bounding box.south) {
  \textcolor{ForestGreen}{\rule{8pt}{2pt}}~Forward (leaves to root: primal $v$, local derivatives)\qquad
  \textcolor{Maroon}{\rule{8pt}{2pt}}~\textcolor{Maroon}{- -}~Reverse (adjoint $\bar{v}$ propagation)
};
\end{scope}

\end{tikzpicture}
\caption{The parameterized tree of one Rosenbrock pattern (left), and the same tree instantiated at the data point $i = 2$ with $x = [2, 3, \ldots]$, carrying the first-order \gls{ad} passes (right). \textcolor{ForestGreen}{Green} marks the forward pass and its stored local derivatives; \textcolor{Maroon}{maroon} dashed arrows mark the reverse pass.}
\label{fig:ad}
\end{sidewaysfigure}

\subsection{Sparsity Detection Without Coloring}
\label{sec:ad:sparsity}

Sparse Jacobian and Hessian evaluation is commonly implemented with graph coloring~\citep{gebremedhinWhatColorYour2005}.
Structurally independent columns are grouped into colors, so that the whole derivative is recovered from one compressed evaluation per color.
This is the mechanism behind the sparse \gls{ad} of general modeling and machine learning frameworks, among them ADNLPModels.jl~\citep{Migot_ADNLPModels_jl_Automatic_Differentiation} in the \gls{jso} ecosystem and the coloring-based sparse Jacobian and Hessian routines available for JAX and PyTorch.
All of them first detect a sparsity pattern, color it, and then evaluate one matrix--vector product per color.
The number of evaluations is therefore the chromatic number of the pattern, which the framework must compute.
The entries must afterwards be scattered from the compressed result back to their positions in the sparse matrix.

ExaModels.jl takes a different route.
The \gls{simd} abstraction fixes the structure of each pattern \emph{a priori} at model-construction time.
The exact structural sparsity layout can therefore be determined without coloring, before any derivative is computed.
The reverse pass then writes each nonzero straight to its position in the \gls{coo} array.
No coloring is computed, no compressed evaluation is performed, and no scatter is needed.
The entries so produced are not fully deduplicated.
A variable that appears at more than one data point, or in more than one pattern, contributes a separate entry for each, and the solver sums them when it assembles the matrix.
This is not peculiar to our route, since a coloring-based evaluation likewise produces a compressed array from which the entries must afterwards be recovered.

The bookkeeping for this layout takes place at the \texttt{SIMDFunction} level.
The stored expression tree serves as the template of the pattern's algebraic structure, the \texttt{Compressor} merges repeated visits to the same variable into unique nonzero entries, and the offsets place each pattern's entries in the global arrays.

\subsubsection{Probing the Tree with Sentinel Inputs}
\label{sec:ad:probe}
For each pattern, ExaModels.jl performs a single probe evaluation of the expression tree using \emph{sentinel} arguments: placeholder sources that stand in for the ordinary inputs.
For the Rosenbrock pattern of \cref{eq:rosenbrock}, with \texttt{obj} the created \texttt{Objective} and \texttt{f = obj.f.f} its lifted function, the first-order probe proceeds as follows in the Julia REPL:
\begin{lstlisting}[style=jlrepl]
julia> d = f(Identity(),
             AdjointNodeSource(NaNSource{Float64}()),
             NaNSource{Float64}());

julia> raw = Any[];

julia> grpass(d, nothing, nothing, nothing, raw, NaN);

julia> raw
3-element Vector{Any}:
 x[i - 1]
 x[i]
 x[i - 1]
\end{lstlisting}
The primal vector $x$ and the parameter vector $\theta$ are replaced by \texttt{NaNSource\{T\}} sentinels, and the data point $i$ by an \texttt{Identity} sentinel.
The numerical values these return are immaterial, since the probe records only which entries are touched.
\texttt{NaN} is used so that any result that did depend on one would be visibly poisoned rather than plausible.
Under these sentinels, the data-access nodes (\texttt{DataSource}, \texttt{DataIndexed}) evaluate to \texttt{NaN}, since no real data is available, while the \texttt{Var} nodes record which symbolic variable references are visited.
The probe runs the \emph{same} reverse-pass infrastructure used for ordinary derivative evaluation, \texttt{grpass} for first order and \texttt{hrpass0} for second order.
The only difference is at the leaves.
With the \texttt{Compressor} argument set to \texttt{nothing}, a separate method makes each variable leaf record the symbolic variable reference it resolved, instead of accumulating a numeric value.
Because the probe reuses that same infrastructure, every general node is visited in the identical order.
The recorded sequence of nonzero positions therefore matches, entry for entry, the order in which the numeric passes will later write them.
The three entries displayed above are the symbolic variable references visited in order: the first two from the term $100(x_{i-1}^2 - x_i)^2$ and the third from $(x_{i-1} - 1)^2$.

\subsubsection{Compressor: Merging Repeated Leaf Visits}
\label{sec:ad:compressor}
The sequence \texttt{raw} collected by the probe may contain duplicates, for instance when $x_i$ appears more than once in the expression.
Without merging them, the number of stored nonzeros would grow with the number of leaf visits rather than with the number of unique entries.
The \texttt{Compressor} is a mapping from the position in \texttt{raw} to the unique nonzero index:
\begin{enumerate}
  \item Compute the list of unique variable references (for the gradient and Jacobian) or unique $(i,j)$ pairs (for the Hessian) from \texttt{raw}.
  \item Build a mapping tuple: for each entry of \texttt{raw}, record which unique index it corresponds to.
  \item Store this mapping in a \texttt{Compressor\{I\}}, whose type parameter \texttt{I} is an \texttt{NTuple\{N,Int\}} and so fixes the length $N$ at compile time.
\end{enumerate}
Because the length is fixed in the type, looking up \texttt{comp(cnt)} is an inlined, constant-time operation with no runtime allocation.
Continuing the Rosenbrock example in the REPL, the recorded sequence contains two unique references, which determine \texttt{o1step} $= 2$:
\begin{lstlisting}[style=jlrepl]
julia> unique(raw)
2-element Vector{Any}:
 x[i - 1]
 x[i]
\end{lstlisting}
and mapping each position of \texttt{raw} to the unique entry it corresponds to gives the stored compressor:
\begin{lstlisting}[style=jlrepl]
julia> obj.f.comp1
Compressor{Tuple{Int64, Int64, Int64}}((1, 2, 1))
\end{lstlisting}
The first and third leaf visits accumulate into the first unique nonzero, and the second visit into the second.
The second-order probe records four entries over three unique nonzeros, the two diagonal entries and one off-diagonal, yielding \texttt{Compressor((1, 2, 3, 1))} with \texttt{o2step} $= 3$.
Throughout, only the lower triangle ($i \geq j$) of the symmetric Hessian is stored, following the NLPModels.jl convention.

The resulting Jacobian and Hessian are stored in a \emph{partially compressed} coordinate format, and the qualifier \emph{partially} reflects where uniqueness has been resolved.
At the level of an individual expression tree, the \texttt{Compressor} has already resolved it: duplicate entries arising from the same tree, such as those from a variable that appears more than once in one expression, are merged.
Duplicates can still occur \emph{across} different expression trees, and we deliberately do \emph{not} compress across different patterns or across different data points of the same pattern.
If two patterns contribute to the same $(i,j)$ entry of the Hessian, they write to separate positions in the COO array.
The final matrix is therefore only partially compressed, and a solver that requires a fully compressed format, such as \gls{csc}, must perform the remaining compression during its sparse matrix assembly.
This design choice has two benefits.
First, different patterns, and different data points within a pattern, never write to the same memory location, which removes write conflicts and allows lock-free parallel evaluation.
Second, the same compiled kernel runs identically for every data point in a pattern, with only the offset advancing by a fixed stride, both on the CPU and the GPU.
The trade-off is a larger COO array, proportional to the total number of pattern$\times$data-point contributions rather than to the number of unique structural nonzeros, which remains moderate for most applications.
When a fully compressed format is needed, the conversion from COO to CSC carries a one-time cost for analyzing the sparsity structure.
The repeated conversions that follow reuse a fixed structure, the typical case in nonlinear optimization, and are inexpensive.

\subsection{Callback Evaluation}
\label{sec:ad:callbacks}

Each callback is implemented internally as a loop over patterns and data points, built from the expression trees, compressors, and offsets described above.
The callback interface itself was introduced in \Cref{sec:modeling:nlp}.
The zero-order callbacks \texttt{NLPModels.obj} and \texttt{NLPModels.cons!} are the simplest: the pattern's tree is instantiated and evaluated at each data point.
Objective terms are accumulated into the scalar objective value, and constraint values are written to their rows, located by the offset \texttt{o0}.

\subsubsection{Structure Callbacks}
\label{sec:ad:structure}
The sparsity patterns themselves are exposed to the solver through the structure callbacks introduced in \Cref{sec:modeling:nlp}:
\begin{itemize}
  \item \texttt{NLPModels.jac\_structure!}: fills the row and column index vectors \texttt{i} and \texttt{j} of the constraint Jacobian's COO representation;
  \item \texttt{NLPModels.hess\_structure!}: fills the corresponding index vectors of the Lagrangian Hessian.
\end{itemize}
Both reuse the reverse-pass infrastructure of \Cref{sec:ad:first} once more.
The passes run with NaN-backed sentinel sources (\Cref{sec:ad:sparsity}), and instead of accumulating numeric values, each variable leaf evaluates its variable index and stores it at the position determined by the \texttt{Compressor} (\Cref{sec:ad:compressor}) and the offsets.
The difference from the probe lies in what is instantiated: the probe is applied once per parameterized expression tree, whereas the structure evaluation is performed for each instantiation of the pattern.
For data points $k = 1, \ldots, N$, the symbolic indices are resolved at the concrete data entry \texttt{itr[k]}.
The resulting global row and column indices are written to the block \texttt{[o1 + (k-1)*o1step + 1, ..., o1 + k*o1step]}, and analogously with \texttt{o2} and \texttt{o2step} for the Hessian.
This loop runs over all data points of every pattern, but it is evaluated only once, at the beginning of the NLP solve.
The numeric callbacks then refill only the value vector \texttt{v} at every iterate.

\subsubsection{Derivative Callbacks}
\label{sec:ad:derivative}
Everything described so far concerns a single algebraic pattern.
A model holds several, stored as the \texttt{obj} and \texttt{cons} tuples of the \texttt{ExaModel} (\Cref{sec:data:model}).
Each derivative callback applies the machinery of this section to each of them in turn, accumulating into the same global output array.
What distinguishes one pattern's contributions from another's is where they are written, which is the role of the offsets carried by each \texttt{SIMDFunction}.

\begin{figure}[t]
\centering
\begin{tikzpicture}[
  cell/.style={draw=black!55, thick, minimum width=10mm, minimum height=8mm,
    inner sep=0pt, font=\scriptsize\ttfamily, text=black!45},
  lb/.style={font=\footnotesize\ttfamily},
  visit/.style={font=\footnotesize\ttfamily, text=Maroon},
]
\node[cell, fill=black!8]  (c1) at (0,0) {};
\node[cell, fill=black!8,  right=0pt of c1] (c2) {};
\node[cell, fill=j1!12, right=0pt of c2] (c3) {1};
\node[cell, fill=j1!12, right=0pt of c3] (c4) {2};
\node[cell, fill=j1!28, right=0pt of c4] (c5) {1};
\node[cell, fill=j1!28, right=0pt of c5] (c6) {2};
\node[cell, fill=j1!12, right=0pt of c6] (c7) {1};
\node[cell, fill=j1!12, right=0pt of c7] (c8) {2};
\node[cell, fill=black!8,  right=0pt of c8] (c9) {$\cdots$};
\node[lb, left=4pt of c1] {v:};

\draw[->, thick, black!70] ($(c3.south west)+(-7mm,-9mm)$) node[left, lb] {o1} -- ($(c3.south west)+(0,-0.5mm)$);

\draw[decorate, decoration={brace, amplitude=3pt}] ($(c3.north west)+(0,1mm)$) -- ($(c4.north east)+(0,1mm)$)
  node[midway, above=4pt, lb] {o1step};

\draw[decorate, decoration={brace, mirror, amplitude=3pt}] ($(c3.south west)+(0,-1mm)$) -- ($(c4.south east)+(0,-1mm)$) node[midway, below=4pt, lb] {k = 1};
\draw[decorate, decoration={brace, mirror, amplitude=3pt}] ($(c5.south west)+(0,-1mm)$) -- ($(c6.south east)+(0,-1mm)$) node[midway, below=4pt, lb] {k = 2};
\draw[decorate, decoration={brace, mirror, amplitude=3pt}] ($(c7.south west)+(0,-1mm)$) -- ($(c8.south east)+(0,-1mm)$) node[midway, below=4pt, lb] {k = 3};

\node[visit] (v1) at ($(c5.north)+(-14mm,16mm)$) {cnt = 1};
\node[visit] (v2) at ($(c5.north)+(1mm,20mm)$) {cnt = 2};
\node[visit] (v3) at ($(c5.north)+(16mm,16mm)$) {cnt = 3};
\draw[->, Maroon, semithick, dashed] (v1.south) to[bend right=12] ($(c5.north)+(-1.5mm,0)$);
\draw[->, Maroon, semithick, dashed] (v2.south) -- (c6.north);
\draw[->, Maroon, semithick, dashed] (v3.south) to[bend left=18] ($(c5.north)+(1.5mm,0)$);
\node[visit, right=4mm of v3] {comp1 = (1, 2, 1)};

\node[lb, anchor=north] at ($(c6.south)+(6mm,-13mm)$) {write position: o1 + (k-1)*o1step + comp1(cnt)};
\end{tikzpicture}
\caption{Offset arithmetic for the sparse gradient of the Rosenbrock pattern, showing where each data point's contributions are written in the global \gls{coo} array.}
\label{fig:offset}
\end{figure}
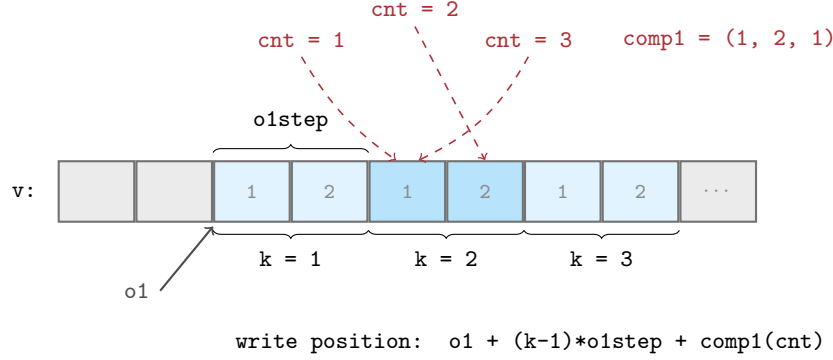

The derivative callbacks first introduced in \Cref{sec:modeling:nlp} are the dense gradient \texttt{NLPModels.grad!}, the sparse Jacobian \texttt{NLPModels.jac\_coord!}, and the sparse Hessian \texttt{NLPModels.hess\_coord!}.
All three share the same mechanics.
The reverse pass traverses the tree and \emph{counts} the entry updates, incrementing a running counter \texttt{cnt} each time a leaf node is reached.
The entry to update is \texttt{offset + comp(cnt)}, where \texttt{offset} locates the pattern's block within the global COO array and \texttt{comp(cnt)} maps the counter to the correct slot among the pattern's unique nonzeros.
The offset advances by \texttt{o1step} (first order) or \texttt{o2step} (second order) for each data point, so data point \texttt{k} writes to positions \texttt{[o1 + (k-1)*o1step + 1, ..., o1 + k*o1step]}.
Each data point has the same number of unique nonzeros, since that number is set by the pattern and not by the data.
The offsets are therefore perfectly regular, and no dynamic bookkeeping is needed.
\Cref{fig:offset} illustrates this offset calculation for the gradient of the Rosenbrock pattern.

\begin{itemize}
  \item \texttt{NLPModels.grad!} (dense gradient): on the CPU, it runs the dense reverse pass \texttt{drpass} (\Cref{sec:ad:first}), which accumulates each contribution $\bar{v}$ directly in place into the dense gradient vector.
On the GPU, this direct accumulation would create race conditions, because threads evaluating different patterns may contribute to the same gradient entry at the same time.
The GPU path therefore first evaluates a \emph{sparse but uncompressed} gradient with \texttt{grpass}, each thread writing to its own conflict-free slot through the compressor and offsets.
It then reduces that buffer into the dense gradient vector (\Cref{sec:impl:gpu}).

  \item \texttt{NLPModels.jac\_coord!} (sparse Jacobian): it loops over the collection of \texttt{SIMDFunction}s stored in the model (\Cref{sec:data:simdfunction}).
For each pattern, the expression tree is instantiated at every data point and run through the sparse reverse pass \texttt{jrpass} (\Cref{sec:ad:first}), which adds each contribution $\bar{v}$ to the COO value-vector entry \texttt{v[offset + comp1(cnt)]}.
The \texttt{Compressor} (\Cref{sec:ad:compressor}) merges the writes into the partially compressed form, and the offsets and strides stored in the \texttt{SIMDFunction} place each pattern's block and advance the write position across data points.

  \item \texttt{NLPModels.hess\_coord!} (sparse Hessian): it proceeds in the same way with the second-order passes \texttt{hrpass} and \texttt{hdrpass} (\Cref{sec:ad:second}).
Each diagonal and off-diagonal contribution is written into the COO value vector at the position determined by \texttt{comp2} and the second-order offsets \texttt{o2} and \texttt{o2step} (\Cref{sec:data:simdfunction}).
\end{itemize}

\section{GPU Kernel Generation}
\label{sec:gpu}

ExaModels.jl evaluates models and their derivatives natively on the GPU.

\subsection{Selecting a GPU Backend}
\label{sec:gpu:frontend}

Moving a model to the GPU takes one keyword argument at \texttt{ExaCore} construction.
Nothing else in the model code changes:
\begin{lstlisting}
using ExaModels, CUDA, MadNLPGPU

N = 10_000

# CPU (default)
core = ExaCore()

# NVIDIA GPU: the only change on the user side
core = ExaCore(backend = CUDABackend())

# everything that follows remains the same
@add_var(core, x, N; start = (mod(i,2)==1 ? -1.2 : 1.0 for i=1:N))
@add_obj(core, 100(x[i]^2 - x[i+1])^2 + (x[i] - 1)^2 for i = 1:N-1)
model  = ExaModel(core)
result = madnlp(model)   # GPU-capable solver (MadNLPGPU.jl)
\end{lstlisting}
The supported backends, provided through KernelAbstractions.jl~\citep{churavyKernelAbstractionsjl2025}, are CUDA.jl (NVIDIA, \texttt{CUDABackend()}), AMDGPU.jl (AMD, \texttt{ROCBackend()}), oneAPI.jl (Intel), Metal.jl (Apple), and the \texttt{CPU()} backend for multithreaded execution on the host.
The \texttt{ExaCore} then keeps all of its numerical arrays in the array type of the selected device.
When the user supplies host arrays, for example as initial values or as parameter data, they are converted to that type as they are stored.
Supplying the data in the device array type from the start avoids the conversion, but it is not required.
The modeling syntax, the callbacks, and the solution access are the same as on the CPU.

\subsection{SIMDFunction as the Unit of Kernel Generation}
\label{sec:gpu:simdfunction}

\begin{figure}[t]
\centering
\begin{tikzpicture}[
  font=\footnotesize,
  box/.style={rectangle, rounded corners=3pt, draw=black!50, thick, fill=black!3,
    align=center, inner sep=4pt},
  thread/.style={rectangle, rounded corners=3pt, draw=j1!85, thick, fill=j1!10,
    align=center, minimum width=24mm, minimum height=7mm, font=\footnotesize\ttfamily},
  slot/.style={rectangle, draw=j3!85, thick, fill=j3!10, align=center,
    text width=48mm, minimum height=7mm, inner sep=2pt, font=\scriptsize\ttfamily},
  minicirc/.style={circle, draw=black!48, fill=black!3, thick, minimum size=4.5mm, inner sep=0pt},
  arr/.style={->, thick, draw=black!55},
  lbl/.style={font=\scriptsize, text=black!60},
]
\node[box, minimum width=30mm, minimum height=26mm] (kernel) at (0,0) {};
\node[minicirc] (t0) at (0, 0.75) {};
\node[minicirc] (t1) at (-0.5, 0.15) {};
\node[minicirc] (t2) at (0.5, 0.15) {};
\draw[black!48] (t0) -- (t1) (t0) -- (t2);
\node[font=\scriptsize, align=center] at (0,-0.55) {expression tree\\(one pattern)};
\node[font=\scriptsize\bfseries, anchor=south] at (kernel.north) {compiled kernel};

\node[thread] (th1) at (4.6,  1.65) {thread 1};
\node[thread] (th2) at (4.6,  0.55) {thread 2};
\node[thread] (th3) at (4.6, -0.55) {thread 3};
\node[font=\footnotesize] at (4.6, -1.35) {$\vdots$};
\node[thread] (thN) at (4.6, -2.2) {thread N};

\node[slot] (s1) at (9.2,  1.65) {v[o1+1 : o1+o1step]};
\node[slot] (s2) at (9.2,  0.55) {v[o1+o1step+1 : o1+2o1step]};
\node[slot] (s3) at (9.2, -0.55) {v[o1+2o1step+1 : o1+3o1step]};
\node[font=\footnotesize] at (9.2, -1.35) {$\vdots$};
\node[slot] (sN) at (9.2, -2.2) {v[o1+(N-1)o1step+1 : o1+No1step]};

\draw[arr] (kernel.east) -- node[lbl, above, sloped] {itr[1]} (th1.west);
\draw[arr] (kernel.east) -- node[lbl, above, sloped] {itr[2]} (th2.west);
\draw[arr] (kernel.east) -- node[lbl, above, sloped] {itr[3]} (th3.west);
\draw[arr] (kernel.east) -- node[lbl, below, sloped] {itr[N]} (thN.west);
\draw[arr] (th1.east) -- (s1.west);
\draw[arr] (th2.east) -- (s2.west);
\draw[arr] (th3.east) -- (s3.west);
\draw[arr] (thN.east) -- (sN.west);

\node[lbl, align=center] at (4.6, 2.45) {one thread per\\data point};
\node[lbl, align=center] at (9.2, 2.45) {disjoint output blocks\\(no race conditions)};
\end{tikzpicture}
\caption{Mapping of the \gls{simd} abstraction onto a \gls{gpu}. One kernel is compiled per algebraic pattern (left), one thread runs per data point (middle), and each thread writes to its own block of the output array (right).}
\label{fig:gpumap}
\end{figure}
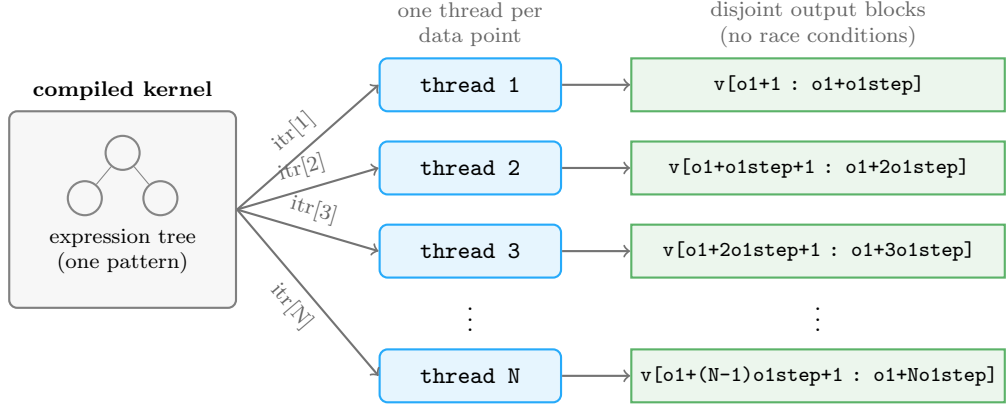

The \gls{simd} abstraction of \Cref{sec:modeling} expresses each pattern as a single function evaluated over a data iterator, and the parameterized expression tree of \Cref{sec:graph} is type-stable end to end, so the Julia compiler specializes a dedicated evaluator per pattern.
The AD passes of \Cref{sec:ad} therefore map onto the device directly: each algebraic pattern becomes a GPU kernel, and each data point becomes a thread (\Cref{fig:gpumap}).
All threads within a pattern execute the same compiled kernel and instantiate the same expression tree on different data indices.
The expression-tree structure therefore introduces no thread divergence, and every thread follows the same control flow, up to any branching inside the elementary functions themselves.
This is well suited to GPU hardware.
The implementation uses KernelAbstractions.jl for backend-agnostic kernels, written once with the \texttt{@kernel} macro and executed on any of the supported backends.

The \texttt{SIMDFunction} of \Cref{sec:data} is the unit of kernel generation.
One kernel is compiled per pattern from its expression tree, and the \texttt{Compressor} mappings and offsets that the \texttt{SIMDFunction} carries tell each thread where to write.
Within a pattern, the entries produced by data point $k$ begin at $\texttt{o1} + (k-1)\cdot\texttt{o1step}$ in the first-order output and at $\texttt{o2} + (k-1)\cdot\texttt{o2step}$ in the second-order output.
The strides \texttt{o1step} and \texttt{o2step} are the numbers of unique entries per data point, determined by the sparsity analysis of \Cref{sec:ad:sparsity}.
A reverse-pass visit writes at that base plus the index returned by the corresponding compressor.
Distinct data points therefore write to disjoint slices of the output arrays, and no two threads target the same position.
The derivative kernels need no atomics and no locks.

\subsection{CPU vs.\ GPU Evaluation}
\label{sec:gpu:cpugpu}

A GPU model exposes the same callback names as a CPU model, namely \texttt{NLPModels.obj}, \texttt{NLPModels.cons!}, \texttt{NLPModels.grad!}, \texttt{NLPModels.jac\_coord!}, and \texttt{NLPModels.hess\_coord!}.
Julia's multiple dispatch, keyed on the model's backend and array types, selects the implementation.
The two implementations of \texttt{NLPModels.cons!} follow, simplified to show the structure.

\begin{center}
\noindent
\begin{minipage}[t]{0.485\linewidth}
\begin{lstlisting}[basicstyle=\small\ttfamily, numbers=left, numbersep=4pt, xleftmargin=1.4em, framexleftmargin=1.4em]
# CPU (simplified)
function NLPModels.cons!(
    model, x, g)
  for con in model.cons
    @simd for i in
        eachindex(con.itr)
      g[offset0(con, i)] +=
        con.f(con.itr[i],
          x, model.θ)
    end
  end
end
\end{lstlisting}
\end{minipage}\hfill
\begin{minipage}[t]{0.485\linewidth}
\begin{lstlisting}[basicstyle=\small\ttfamily, numbers=left, numbersep=4pt, xleftmargin=1.4em, framexleftmargin=1.4em, showlines=true]
# GPU (simplified)
function NLPModels.cons!(
    model, x, g)
  for con in model.cons
    kerf(model.ext.backend)(
      g, con.f, con.itr,
      x, model.θ;
      ndrange =
        length(con.itr))
  end
end
 
\end{lstlisting}
\end{minipage}
\end{center}

Both implementations loop over the constraint patterns in the same way.
They differ only in the inner iteration over data points.
On the CPU, a \texttt{@simd for} loop evaluates the parameterized expression tree \texttt{con.f(con.itr[i], x, model.\codetheta)} one data point at a time, and the annotation lets the compiler vectorize across data points using CPU \gls{simd} instructions.
On the GPU, the loop body becomes the body of the KernelAbstractions.jl kernel \texttt{kerf}, launched over the data iterator.

\noindent\begin{minipage}{\linewidth}
\begin{lstlisting}[basicstyle=\small\ttfamily, numbers=left, numbersep=4pt, xleftmargin=1.4em, framexleftmargin=1.4em]
# one thread per data point
@kernel function kerf(y, @Const(f), @Const(itr),
    @Const(x), @Const(θ))
  I = @index(Global)
  @inbounds y[offset0(f, itr, I)] = f(itr[I], x, θ)
end
\end{lstlisting}
\end{minipage}

\texttt{@index(Global)} provides the thread's global index $I$, and the launch with \texttt{ndrange = length(con.itr)} assigns one thread per data point.
Each thread instantiates the expression tree at its own data point \texttt{itr[I]} and writes the result at \texttt{offset0(f, itr, I)}, the per-data-point base position described above, so distinct threads write to disjoint positions.
Because \texttt{f} is a compile-time-specialized expression tree and a bits type (\Cref{sec:graph:type}), it is a valid kernel argument.
The kernel also contains no dynamic dispatch, which GPU compilation requires.

Multithreaded execution on the host is obtained the same way, by passing \texttt{CPU()} as the backend.
No separate implementation and no thread-specific specialization is written for it.
KernelAbstractions.jl compiles the same kernel body for the host and, at launch, divides the index range evenly among the available Julia threads, each of which runs the kernel over its own contiguous block of data points~\citep{churavyKernelAbstractionsjl2025}.
The kernel is therefore data-parallel in the same sense as on a \gls{gpu}, and only the number of threads differs.
This is the backend used for the multithreaded \gls{cpu} results of \Cref{sec:numerics:ad}.

The same comparison applies to the sparse gradient:

\noindent
\begin{minipage}[t]{0.485\linewidth}
\begin{lstlisting}[basicstyle=\small\ttfamily, numbers=left, numbersep=4pt, xleftmargin=1.4em, framexleftmargin=1.4em]
# CPU (simplified):
# dense reverse pass
function NLPModels.grad!(
    model, x, g)
  for obj in model.objs
    @simd for i in
        eachindex(obj.itr)
      drpass(
        obj.f(obj.itr[i],
          AdjointNodeSource(x),
          model.θ),
        g, 1.0)
    end
  end
end
\end{lstlisting}
\end{minipage}\hfill
\begin{minipage}[t]{0.485\linewidth}
\begin{lstlisting}[basicstyle=\small\ttfamily, numbers=left, numbersep=4pt, xleftmargin=1.4em, framexleftmargin=1.4em, showlines=true]
# GPU (simplified):
# sparse reverse pass
function NLPModels.grad!(
    model, x, g)
  buf = model.ext.gradbuffer
  for obj in model.objs
    kerg(model.ext.backend)(
      buf, obj.f, obj.itr,
      x, model.θ, 1.0;
      ndrange =
        length(obj.itr))
  end
  compress_to_dense(g, buf)
end
 
\end{lstlisting}
\end{minipage}

On the CPU, each data point's forward evaluation builds the adjoint tree, and the dense reverse pass \texttt{drpass} (\Cref{sec:ad:first}) accumulates contributions directly into the dense gradient vector.
On the GPU, the wrapper launches the kernel \texttt{kerg} into a sparse buffer, which is then compressed into the dense gradient (\Cref{sec:impl:gpu}).

\noindent\begin{minipage}{\linewidth}
\begin{lstlisting}[basicstyle=\small\ttfamily, numbers=left, numbersep=4pt, xleftmargin=1.4em, framexleftmargin=1.4em]
# the AD pass in each thread
@kernel function kerg(y, @Const(f), @Const(itr),
    @Const(x), @Const(θ), @Const(adj))
  I = @index(Global)
  @inbounds grpass(f(itr[I], AdjointNodeSource(x), θ),
    f.comp1, y, offset1(f, I), 0, adj)
end
\end{lstlisting}
\end{minipage}

Each thread runs the full forward and reverse pass: the forward evaluation on an \texttt{AdjointNodeSource} builds the adjoint tree, and \texttt{grpass} writes the contributions to the thread's own slot of the buffer, at positions determined by the \texttt{Compressor} \texttt{f.comp1} and the base offset \texttt{offset1(f, I)}.
The Jacobian and Hessian kernels follow the same structure, calling \texttt{jrpass} and \texttt{hrpass0} respectively.

The only thing that changes is the top-level wrapper that distributes the work.
The pattern is instantiated over its data points, and each instantiation, that is, each forward and reverse pass, is mapped to a GPU thread.
The expression trees are bits types, meaning concretely typed, immutable, and free of references to host memory (\Cref{sec:graph:type}), so they are passed to the kernels directly.
The partially compressed COO output format (\Cref{sec:ad:sparsity}) keeps the evaluation free of race conditions.
The GPU is well utilized as long as the model carries many data points.

Several further implementation choices, namely the native COO output format, the objective's buffer-then-reduce evaluation, the gradient's sparse-then-compress path, and the restriction to data-level parallelism, are described in \Cref{sec:impl:gpu}.

\section{Implementation Notes}
\label{sec:impl}

This section collects implementation details deferred from the preceding sections.
Each item records a design or engineering choice that matters for using or extending ExaModels.jl but is not needed to follow the main development.

\subsection{Model Construction}
\label{sec:impl:construction}

\subsubsection{Immutable ExaCore and Reconstruction}
\label{sec:impl:exacore}

The central design choice that enables type stability is that \texttt{ExaCore} is an \emph{immutable} struct rather than a \texttt{mutable struct}.
As described in \Cref{sec:data:core}, each call to \texttt{@add\_var}, \texttt{@add\_obj}, or \texttt{@add\_con} extends the component tuples and returns a new \texttt{ExaCore}.
The type parameters of the new core encode the exact set of patterns added so far, so its concrete type differs from that of the previous core.
The Julia compiler sees these as distinct concrete types and generates specialized code for each construction step.
By the time \texttt{ExaModel(core)} is called, the final \texttt{ExaCore} already holds all expression trees, \texttt{Compressor} tuples, and offset computations.
The instantiation builds no new expression trees or sparsity structures, only the solver-facing metadata and backend workspaces described in \Cref{sec:data:model}.

This design contrasts with JuMP~\citep{lubinJuMP10Recent2023} and similar \glspl{ams}, where the model object is a mutable container that accumulates variables, constraints, and objectives behind a single stable type.
In JuMP, each \texttt{@variable}, \texttt{@constraint}, and \texttt{@objective} mutates internal dictionaries within the same \texttt{Model} object:
\begin{lstlisting}
model = Model(Ipopt.Optimizer)

@variable(model, x[1:N])                    # mutates model

@constraint(model, [i = 1:N-2], g(x, i) == 0)  # mutates model

@objective(model, Min, sum(f(x, i) for i = 1:N-1))  # mutates model
\end{lstlisting}
This API is natural and convenient, but the model's concrete type does not change across these calls.
All structure is hidden in runtime data structures.
The derivative evaluation code must therefore dispatch dynamically on the stored expression types, which rules out the compile-time specialization described in \Cref{sec:graph:type,sec:gpu:cpugpu}.

The ExaModels.jl macros \texttt{@add\_var}, \texttt{@add\_obj}, and \texttt{@add\_con} are designed to \emph{look} like the JuMP convention while achieving the opposite effect.
Each macro call rebinds the name \texttt{core} to a new \texttt{ExaCore} of a different concrete type, so the compiler can track the growing model structure statically.
The user writes:
\begin{lstlisting}
core = ExaCore(concrete = Val(true))
# core :: ExaCore{..., Tuple{}, Tuple{}}

@add_var(core, x, N; start = (x0(i) for i = 1:N))
# core :: ExaCore{..., Tuple{Variable}, Tuple{}}

@add_con(core, g(x, i) for i = 1:N-2)
@add_obj(core, f(x, i) for i = 1:N-1)
# core :: ExaCore{..., Tuple{Variable}, Tuple{Objective{...}, ...}}

model = ExaModel(core)
\end{lstlisting}
This reads like a conventional modeling script, but each line produces a new concrete type that the compiler can fully specialize.

The mutating appearance is achieved purely by macro expansion: \texttt{@add\_var(core, x, n)} unrolls to
\begin{lstlisting}
local _var
core, _var = add_var(core, n; name = Val(:x))
x = _var
\end{lstlisting}
The underlying function \texttt{add\_var} returns a new core together with the variable handle.
The macro rebinds the name \texttt{core} in the calling scope to the new object and binds the handle to \texttt{x}.
No object is ever mutated.

Creating a new \texttt{ExaCore} at every declaration may appear expensive, but the cost is negligible.
The struct holds only counters, handles, and references to the numerical arrays, so each rebind copies the references while the arrays themselves are never recreated.
For the example above with $N = 10^6$, the struct occupies 144 bytes when empty and 176, 304, and 392 bytes after the variable, constraint, and objective declarations.
Each increment is the declared handle or wrapper stored inline in the core's tuples.
These figures were measured on a 64-bit machine with ExaModels.jl v0.11.2 and may change slightly across versions.
The struct therefore grows with the number of declared blocks and patterns, not with the number of data points.
It stays small for the typical large-scale model, which has few patterns.

\subsubsection{Reusable Expression Objects}
\label{sec:impl:expressions}

The \texttt{@add\_expr} macro defines a subexpression that can be reused in objectives and constraints.
For example,
\begin{lstlisting}
@add_expr(core, s, x[i]^2 for i = 1:N)
@add_obj(core, (s[i] - 1)^2 for i = 1:N)
@add_con(core, s[i] + s[i+1] for i = 1:N-1; ucon = 1.0)
\end{lstlisting}
The first line returns an \texttt{Expression} named \texttt{s}.
It adds no objective, constraint, callback, or independent algebraic pattern to the model.
Each later reference \texttt{s[i]} substitutes the corresponding expression tree directly into the objective or constraint that uses it.
The object can therefore be reused at several sites without introducing auxiliary variables or constraint rows.

An \texttt{Expression\{S,F,I,T\}} stores the information needed to reconstruct a reusable subexpression at a selected index:
\begin{itemize}
  \item \texttt{size::S} and \texttt{length::Int}: the shape and number of scalar entries;
  \item \texttt{f::F}: the generator function that constructs the subexpression tree;
  \item \texttt{iter::I}: the collected data iterator used when the expression is indexed; and
  \item \texttt{tag::T}: optional user metadata.
\end{itemize}
The \texttt{add\_expr} function returns an updated \texttt{ExaCore} together with the \texttt{Expression}.
An anonymous expression is returned to the caller without being stored in the core.
The named form \texttt{@add\_expr(core, s, ...)} also records the object in \texttt{refs}, making it accessible as \texttt{core.s} or \texttt{model.s}.
This entry is only a name-to-object association.
The expression is never added to the \texttt{obj} or \texttt{cons} tuples, is never compiled into its own \texttt{SIMDFunction}, and is never evaluated as a callback.
Indexing it as \texttt{s[i]} instead invokes its generator for that index and substitutes the resulting expression tree into the objective or constraint being built.
One consequence is that a subexpression referenced at several sites is re-evaluated at each site rather than computed once and shared.
We have found this sufficient for the use cases of ExaModels.jl, since the subexpressions that arise in practice are inexpensive to re-evaluate.

\subsubsection{Compile-Time Constants and Algebraic Simplification}
\label{sec:impl:constants}

The \texttt{Constant\{V\}} node embeds its value \texttt{V} as a type parameter rather than a struct field.
This means the constant occupies zero runtime storage and is fully resolved at compile time.
The Julia compiler (and the AOT compiler \texttt{juliac --trim=safe}) can propagate the value concretely through all downstream operations.

This design also enables algebraic simplification at model-construction time.
Trivial branches are collapsed before the expression tree is finalized.
The identities applied automatically include the following:
\begin{center}
\begin{tabular}{ll}
  $x + 0 \to x$, \quad $x \cdot 1 \to x$, \quad $x / 1 \to x$, & $x^0 \to 1$, \quad $x^1 \to x$, \\
  $x^2 \to \operatorname{abs2}(x)$, \quad $x^{-1} \to \operatorname{inv}(x)$, & $0 \cdot x \to 0$, \quad $0 / x \to 0$.
\end{tabular}
\end{center}
These rules reduce the depth of the expression tree and remove unnecessary operations.
This matters most for GPU kernels, where every instruction counts.
The simplification is exact, with no floating-point approximation, because it operates on the symbolic tree structure rather than on numeric values.

Integer literal exponents reach the power rules through a lowering step in Julia rather than through the registered operators of \Cref{sec:impl:register}.
Such an exponent is lowered by the compiler to \texttt{Base.literal\_pow}, with the exponent carried as a \texttt{Val} type parameter.
The expression \texttt{x[i]\^{}2} therefore arrives at ExaModels.jl as \texttt{literal\_pow(\^{}, x[i], Val\{2\})} rather than as a call to \texttt{\^{}} with the runtime value $2$.
The exponent is thus known to the method that builds the node, which is what allows $x^2$ to be recorded as an \texttt{abs2} node and $x^1$ to be replaced by $x$.
The derivative rules of \texttt{abs2} are then used in place of the general power rule, and no exponent is stored or dispatched on at evaluation time.
The other rules in the table fire when an operand is a \texttt{Constant} node, as in $x^{-1} \to \operatorname{inv}(x)$.

\subsubsection{Function Registration}
\label{sec:impl:register}

For a function to appear in an ExaModels.jl expression, it must be \emph{registered}.
The system needs the function itself together with its first and second derivatives, which are required for the gradient, Jacobian, and Lagrangian Hessian computations of \Cref{sec:ad}.
ExaModels.jl provides two macros for this purpose:
\begin{itemize}
  \item \texttt{@register\_univariate(f, df, ddf)}: registers a unary function with its first and second derivatives.
  \item \texttt{@register\_bivariate(f, df1, df2, ddf11, ddf12, ddf22)}: registers a binary function with its partial derivatives up to second order.
\end{itemize}
Each macro generates the following method groups via multiple dispatch:
\begin{itemize}
  \item a primal method that constructs a \texttt{Node1} or \texttt{Node2} when applied to \texttt{AbstractNode} arguments;
  \item a constant-folding method for \texttt{Constant} arguments;
  \item a first-order adjoint method that constructs an \texttt{AdjointNode1} or \texttt{AdjointNode2} with the derivative values;
  \item a second-order adjoint method that constructs a \texttt{SecondAdjointNode1} or \texttt{SecondAdjointNode2}; and
  \item an evaluation method that applies the function when the constructed node is called on a data point.
\end{itemize}

ExaModels.jl v0.11.2 ships with 61 pre-registered function signatures (52 univariate and 9 bivariate), which cover the operations needed by most algebraic models:
\begin{itemize}
  \item arithmetic: $+$, $-$ (unary and binary), $*$, $/$, and \texttt{\^{}} (power);
  \item exponential and logarithmic: \texttt{exp}, \texttt{exp2}, \texttt{exp10}, \texttt{expm1}, \texttt{log}, \texttt{log2}, \texttt{log10}, \texttt{log1p}, \texttt{sqrt}, \texttt{cbrt};
  \item trigonometric: \texttt{sin}, \texttt{cos}, \texttt{tan}, \texttt{csc}, \texttt{sec}, \texttt{cot}; the inverses \texttt{asin}, \texttt{acos}, \texttt{atan} (including the two-argument form), \texttt{acot}; the degree variants \texttt{sind}, \texttt{cosd}, \texttt{tand}, \texttt{cscd}, \texttt{secd}, \texttt{cotd}, \texttt{atand}, \texttt{acotd}; and \texttt{sinpi}, \texttt{cospi}, \texttt{sinc};
  \item hyperbolic: \texttt{sinh}, \texttt{cosh}, \texttt{tanh}, \texttt{csch}, \texttt{sech}, \texttt{coth}, and the inverses \texttt{asinh}, \texttt{acosh}, \texttt{atanh}, \texttt{acoth};
  \item utility: \texttt{abs}, \texttt{abs2}, \texttt{inv}, \texttt{sign}, \texttt{signbit}, \texttt{floor}, \texttt{ceil}, \texttt{max}, \texttt{min}, \texttt{hypot}, \texttt{deg2rad}, \texttt{rad2deg}.
\end{itemize}
Each signature is registered together with its first and second derivatives.

Less common functions are provided through extension packages, so that the dependency is incurred only when needed.
The \texttt{ExaModelsSpecialFunctions} extension is loaded when SpecialFunctions.jl is present.
It registers \texttt{erf}, \texttt{erfc}, \texttt{erfinv}, \texttt{gamma}, \texttt{digamma}, \texttt{trigamma}, \texttt{beta}, the Bessel functions \texttt{besselj0}, \texttt{besselj1}, \texttt{bessely0}, and \texttt{bessely1}, and the Airy functions \texttt{airyai} and \texttt{airybi}, among others, all with the derivative rules required for AD.
Beyond these, users can register custom functions with the same two macros.
The registered derivatives enter the type system through multiple dispatch, so no modification of ExaModels.jl internals is required.

\subsection{Callback Implementation}
\label{sec:impl:callbacks}

\subsubsection{Compiler Specialization of the Callbacks}
\label{sec:impl:specialization}

The expression tree type is fully concrete and known at compile time (\Cref{sec:graph:type}), and the AD traversal is annotated with \texttt{@inline} to force aggressive inlining.
The Julia compiler therefore produces specialized machine code for each pattern with no interpretive overhead.
Forcing inlining in this way is also what makes \gls{aot} compilation possible (\Cref{sec:aot}).
The generated code for a complete gradient kernel is reproduced in \Cref{app:llvm}.
As detailed there, the kernel uses 13 floating-point operations against the best hand-written sequence of 12 that we could construct, and no unnecessary overhead is observable in the LLVM code.
Each user-facing callback, namely the objective, constraint, gradient, Jacobian, and Hessian evaluations, thus compiles to a single piece of LLVM code, compiled once per model.
This lets the compiler optimize the entire evaluation as a whole.
All callbacks are also non-allocating on the CPU, so no garbage-collection pauses occur during the solve.
The user writes natural mathematical expressions, dispatch builds the tree, and the compiler produces code whose performance is competitive with hand-written source transformations.

\subsubsection{Affine-Node Specialization}
\label{sec:impl:affine}

Many expression trees contain long chains of affine operations ($+$, $-$, multiplication by a constant) whose second derivatives are identically zero.
When a binary operation combines one subtree with a plain scalar, the adjoint construction reduces it to a unary node whose function is recorded as \texttt{FirstFixed\{F\}} or \texttt{SecondFixed\{F\}}.
These wrapper types indicate that the first or the second argument of the binary function \texttt{F} is the fixed scalar.
ExaModels.jl exploits this in the second-order reverse pass of \Cref{sec:ad} with \texttt{hrpass0}, a specialized entry point that dispatches on the function type.
For affine \texttt{SecondAdjointNode1} nodes, such as \texttt{FirstFixed\{typeof(*)\}} and \texttt{SecondFixed\{typeof(+)\}}, the $h$ term is known to be zero at compile time.
The recursion then propagates only $(\bar{v} \cdot y,\; \bar{v}_2 \cdot y^2)$ and never computes the $\bar{v} \cdot h$ contribution.
This removes unnecessary multiplications and additions from the inner loop, in particular for the common case of weighted sums of nonlinear terms.
The specialization is selected by dispatch on the type of the expression tree, so the choice is resolved at compile time and no branch on the affine case remains in the generated code.

\subsubsection{GPU Evaluation Details}
\label{sec:impl:gpu}

\paragraph{Jacobian and Hessian evaluation with partial compression only.}
For Jacobian and Hessian evaluation, each data point in a pattern writes to its own contiguous block of the COO value array (\Cref{sec:ad:sparsity}).
Different data points write to disjoint positions, so all writes are conflict-free and no atomics or synchronization are needed.
This follows from the \gls{simd} structure.
Within a pattern, data point $k$ has first-order block base offset $\texttt{o1} + (k{-}1) \cdot \texttt{o1step}$, and a reverse-pass visit writes at that base plus the one-based index returned by \texttt{comp1}.
The second-order path uses \texttt{o2}, \texttt{o2step}, and \texttt{comp2} in the same way.
ExaModels.jl therefore stores Jacobian and Hessian entries in the partially compressed COO format of \Cref{sec:ad:sparsity} and does \emph{not} further compress them to \gls{csc} format.
Full compression would require sorting and merging entries.
Avoiding it also suits solvers such as Ipopt~\citep{wachterImplementationInteriorpointFilter2006} and MadNLP.jl~\citep{shinMadNLPjl2025}, which accept COO input directly.

\paragraph{Objective evaluation by buffering and reduction.}
\label{sec:gpu:objective}
The objective callback evaluates $f(x) = \sum_\ell \sum_{i=1}^{I_\ell} \tilde{f}^{(\ell)}(p_i, x, \theta)$.
On the CPU, this sum is computed by direct scalar accumulation while iterating over patterns and data points (\Cref{sec:gpu:cpugpu}).
No array of objective terms is formed, so the CPU objective path does not use \texttt{o0}.
On the GPU, direct accumulation by many threads would require atomics, so ExaModels.jl first evaluates the objective terms into a buffer of length $\sum_\ell I_\ell$.
Each pattern's \texttt{o0} identifies the beginning of its disjoint slice of this buffer, and thread $i$ writes its value at \texttt{o0 + i}.
Once all pattern kernels have completed, a parallel \texttt{sum} reduction over the buffer produces the scalar objective value.
Thus, unlike the CPU objective path, the GPU implementation requires the zero-order offset to place contributions before the subsequent reduction.

\paragraph{Gradient evaluation by sparse pass and compression.}
On the CPU, the gradient is computed by the dense reverse pass \texttt{drpass} (\Cref{sec:ad:first}), which accumulates each contribution $\bar{v}$ directly into the dense gradient vector.
On the GPU, the sparse reverse pass \texttt{grpass} is used instead.
Each thread writes its gradient contributions into a sparse buffer indexed by variable index and position, and a subsequent KernelAbstractions.jl kernel \texttt{compress\_to\_dense} reduces that buffer into the dense gradient vector.
Two phases are needed because threads from different patterns may contribute to the same gradient entry, and atomic additions to a dense vector would cause contention:
\begin{lstlisting}
@kernel function compress_to_dense(
        y, @Const(y0), @Const(ptr), @Const(sparsity))
    I = @index(Global)
    @inbounds for j = ptr[I]:(ptr[I+1]-1)
        (k, l) = sparsity[j]
        y[k] += y0[l]
    end
end
\end{lstlisting}
The \texttt{ptr} array partitions the sparse buffer entries by target variable, so each thread of the compression kernel handles a distinct set of target indices without conflicts.
This partition, together with the sorted sparsity array, is computed once at \texttt{ExaModel} instantiation (\Cref{sec:data:model}).
The evaluation of \texttt{ConstraintAugmentation}s uses the same strategy.
Several augmentations may contribute to the same constraint entry, so their contributions are written into a separate buffer and reduced into the constraint vector by the same \texttt{compress\_to\_dense} kernel.

\paragraph{Parallelism over data points only.}
ExaModels.jl parallelizes at the data-point level, so each thread executes the full forward--reverse AD pass for one data point of one pattern.
The expression tree itself is never split across threads.
This design keeps the kernel simple and avoids the synchronization overhead that splitting the tree traversal would incur.
For typical \gls{nlp} patterns, the expression trees are small, on the order of tens of nodes, so the parallelism from thousands of data points is more than enough to saturate GPU hardware.
The result is high occupancy with minimal launch overhead, at one kernel launch per pattern per callback.

\subsection{Compile-Time Design}
\label{sec:impl:compiletime}

\subsubsection{Type Stability and Ahead-of-Time Compilation}
\label{sec:aot}

A Julia function is \emph{type-stable} if the compiler can infer the concrete type of every variable and return value from the types of the arguments alone, with no runtime dispatch.
GPU kernels \emph{require} it, since dynamic dispatch is not supported on GPU hardware, and Julia's \gls{aot} compiler \texttt{juliac --trim=safe} relies on it to produce standalone binaries with unused code eliminated.

Type stability originates in the construction design of \Cref{sec:data}.
The \texttt{ExaCore} is immutable, and each macro call produces a core of a new concrete type whose parameters encode the model built so far.
The final type therefore determines the model structure completely, and every function invoked on the model is compiled specifically for that type.
The construction-level details, including the contrast with mutable model containers such as JuMP's, are given in \Cref{sec:impl:exacore}.

Because all ExaModels.jl code paths are type-stable, the entire model construction and solve pipeline is compatible with \texttt{juliac --trim=safe}.
The compiler produces a standalone binary or shared library that runs without a Julia installation, since the required runtime libraries are bundled with the artifact.
Combined with NLPModelsIpoptLite.jl, a lightweight AOT-compatible wrapper around the Ipopt C library, a complete optimization application can be packaged as a single executable:
\begin{lstlisting}
# RosenbrockApp.jl: compiles to a binary
module RosenbrockApp
using ExaModels, NLPModelsIpoptLite
function (@main)(ARGS)
    N = parse(Int, ARGS[1])

    core = ExaCore(concrete = Val(true))
    @add_var(core, x, N; start = (mod(i,2)==1 ? -1.2 : 1.0 for i=1:N))
    @add_obj(core, 100(x[i]^2 - x[i+1])^2 + (x[i] - 1)^2 for i = 1:N-1)

    model = ExaModel(core)
    result = ipopt(model; print_level = 5)
    return result.status == 0 ? 0 : 1
end
end
\end{lstlisting}
Running \texttt{juliac --trim=safe --output-exe rosenbrock RosenbrockApp.jl} produces a standalone binary that solves the Rosenbrock problem for any $N$ specified at the command line, with no JIT compilation latency.
As of this writing, AOT compilation is supported for CPU backends only.
GPU backends still require Julia's JIT runtime for kernel compilation.

A model, or an entire model library, developed with ExaModels.jl can be compiled together with a solver into a standalone binary, or into a shared library whose entry points are exposed to other programs.
Developers can therefore build model libraries on ExaModels.jl without imposing a Julia-runtime dependence on their users.
The library ships as a conventional compiled artifact that can be called from Python, C++, or any language with a C foreign-function interface, which eases deployment in production environments.

\subsubsection{Generic Numeric Precision}
\label{sec:aot:precision}

Julia's multiple dispatch and JIT compilation make the entire ExaModels.jl pipeline generic over the floating-point type.
This covers expression tree construction, AD evaluation, and the solver callbacks.
All internal code operates on the abstract type \texttt{T <: AbstractFloat}, and the compiler specializes automatically:
\begin{lstlisting}
core32 = ExaCore(Float32)  # all arrays use Float32
core64 = ExaCore(Float64)  # all arrays use Float64
\end{lstlisting}
This extends to user-defined numeric types, including simulated higher precisions.
One example is the double-double type \texttt{Double64} of DoubleFloats.jl~\citep{juliamathDoubleFloatsjl}, which simulates roughly quadruple precision with a pair of \texttt{Float64}s.
Another is Julia's built-in arbitrary-precision \texttt{BigFloat}.
Both run through the entire pipeline with no code changes beyond specifying the type at \texttt{ExaCore} construction.
The \texttt{Constant\{V\}} type parameter (\Cref{sec:impl:constants}) ensures that compile-time constants are converted to the appropriate precision by \texttt{replace\_T}, so no mixed-precision artifacts arise.

\texttt{Float32} arithmetic is much faster than \texttt{Float64} on consumer GPUs, and some hardware does not support double precision at all, not even through software emulation.
Apple GPUs, targeted by the Metal backend, execute only single-precision arithmetic, and ExaModels.jl models run there entirely in \texttt{Float32}.
The same genericity also serves applications that require extended precision, such as sensitivity analysis near degeneracy.

\paragraph{Implications for model libraries.}
Portability also propagates to model libraries built on top of ExaModels.jl, such as ExaModelsPower.jl and COPSBenchmark.jl.
Such a library need only accept the numeric type and the backend as parameters, pass them to \texttt{ExaCore(T; backend = backend, concrete = Val(true))}, convert its numerical data to \texttt{T}, and avoid type-unstable constructs.
It then inherits arbitrary precision, every GPU backend, and AOT compatibility without further effort:
\begin{lstlisting}
function my_model(data; T = Float64, backend = nothing)
    core = ExaCore(T; backend = backend, concrete = Val(true))

    @add_var(core, x, length(data))
    @add_obj(core, f(x, d) for d in convert.(T, data))

    return ExaModel(core)
end
\end{lstlisting}

\subsection{Two-Stage and Batch Structures}
\label{sec:extensions}

The extensively studied classes of structured \glspl{nlp} include two-stage stochastic programs~\citep{IntroductionStochasticProgramming1997,lubinScalableStochasticOptimization2011,chiangStructuredNonconvexOptimization2014}, graph-structured problems~\citep{jalvingGraphbasedModelingAbstraction2022,colePlasmoDatajlJuliaFramework2024}, and batch optimization~\citep{amosOptNetDifferentiableOptimization2017,kimLeveragingGPUBatching2021}.
Of these, the two-stage and batch structures are the ones aligned with the \gls{simd} abstraction.
Each adds an axis of repetition over which the same patterns are evaluated, so capturing it at the modeling level extends \gls{simd} parallelism to the scenario and instance axes.
Graph-induced structure is less amenable, since an irregular graph yields heterogeneous algebraic patterns.

ExaModels.jl accordingly targets both structures within a single formulation.
No separate data structures are needed, since every problem can be regarded as a batched two-stage problem.
A standard \gls{nlp} is the special case of one batch instance and no recourse variables.
A further mechanism, the external oracle, accommodates components that cannot be expressed as parameterized expression trees.
It is described in \Cref{sec:extensions:oracles}.

\subsubsection{Unified Formulation}
\label{sec:extensions:formulation}

The most general formulation targeted by ExaModels.jl is a \emph{batched two-stage parametric} \gls{nlp}.
Given $B \geq 1$ independent batch instances and $S \geq 0$ scenarios per instance, the problem is
\begin{equation}\label{eq:extensions}
\begin{aligned}
  \min_{\{x_{b,s}\}}\quad & \sum_{b=1}^{B} \left( f_0(x_{b,0};\, \theta_{b,0}) + \sum_{s=1}^{S} f_s(x_{b,0}, x_{b,s};\, \theta_{b,s}) \right) \\
  \st\quad & g_0(x_{b,0};\, \theta_{b,0}) \in \mathcal{C}_0, \quad b = 1, \ldots, B, \\
           & g_s(x_{b,0}, x_{b,s};\, \theta_{b,s}) \in \mathcal{C}_s, \quad s = 1, \ldots, S,\; b = 1, \ldots, B,
\end{aligned}
\end{equation}
where the scenario index $s$ runs over $0, 1, \ldots, S$ and $s = 0$ denotes the root.
The vectors $x_{b,0}$ are the root, or first-stage, variables of instance~$b$, and $x_{b,s}$ for $s \geq 1$ are the scenario, or second-stage, variables of that instance.
The vectors $\theta_{b,s}$ are parameters, which are not optimized over.
Each $f_s(\cdot)$ and $g_s(\cdot)$ has the \gls{simd} structure of \cref{eq:simd}, with scenario probabilities absorbed into $f_s(\cdot)$.
Formulation \cref{eq:extensions} has two structural dimensions.
The first is the scenarios, over which the root variables are shared while the scenario variables and constraints are replicated.
The second is the batch instances, which share the same expression trees, and hence the same compiled functions, but differ in their numerical data.
The instances couple through neither the objective nor the constraints, so minimizing the aggregated objective is equivalent to solving the $B$ instances independently.
The batched form exists so that all instances are evaluated by one set of kernels.
Both dimensions are optional and can be active at the same time: $B = 1$ gives a two-stage stochastic program, $S = 0$ gives a batch of standard \glspl{nlp}, and $S = 0$ with $B = 1$ recovers \cref{eq:simd}.

All scenarios and batch instances share the same algebraic patterns.
The \texttt{SIMDFunction}s, \texttt{Compressor} maps, and sparsity structures are therefore built once and reused across every scenario and instance.
Their concrete types coincide, so the functions specialized for them are compiled only once.
The only quantities that vary are numerical: parameter values, variable and constraint bounds, and the primal and dual initial points.
This is what makes the extensions efficient.
It is also why ExaModels.jl provides accessors that address a single scenario, modifying those arrays without touching the model structure and hence without triggering recompilation.

The goal of ExaModels.jl is not to provide decomposition algorithms, which are the domain of algorithm and software developers.
It is to provide what such algorithms need, namely efficient callbacks for structured models together with the structural information required to exploit them.
ExaModels.jl therefore tags each variable and constraint with the stage it belongs to, and reports the scenario that each variable and constraint belongs to as well as the number of scenarios.
With these, a developer can implement methods such as Benders decomposition or progressive hedging on top of the model.

\subsubsection{Two-Stage Stochastic Programs}
\label{sec:extensions:twostage}
\label{sec:impl:accessors}

The two-stage structure is declared with \texttt{TwoStageExaCore(nscen)} together with a marker \texttt{EachScenario()} that distinguishes root from scenario declarations, as in the following example:
\begin{lstlisting}
S = 4                                  # number of scenarios
core = TwoStageExaCore(S)

@add_var(core, x, 2)                   # root (shared)
@add_var(core, z, EachScenario(), 3)   # per scenario
@add_par(core, theta, EachScenario(), theta_vals)

@add_obj(core, x[1]^2 + x[2]^2)        # root objective
@add_obj(core, (z[j, s] - theta[j, s])^2 for (j, s) in idx)
@add_con(core, EachScenario(), (z[j, s] + x[2] for (j, s) in idx))

model = ExaModel(core)
\end{lstlisting}
where \texttt{theta\_vals} holds the per-scenario parameter values and \texttt{idx} collects the index pairs $(j, s)$ over which the scenario patterns are evaluated.
The marker also applies to \texttt{@add\_par} and \texttt{@add\_con}, and appends the scenario index as the last dimension of the declared block.
Objectives are declared with the plain \texttt{@add\_obj}, with the scenario structure carried by the data iterator.
Internally, the scenario blocks are replicated $S$ times and tagged by stage.
The scenario patterns are therefore still constructed once and evaluated across all scenarios in parallel, with the scenario index acting as an additional data dimension.

Two-stage models also provide counterparts of the model data accessors of \Cref{sec:modeling:simd} that address a single scenario.
Parameter values, bounds, and initial points can thus be set and queried one scenario at a time, again without touching the model structure or the compiled functions.
Their signatures are given in the online documentation~\citep{shinExaModelsjl2025}.

\subsubsection{Batch Models}
\label{sec:extensions:batch}

The batch structure solves $B$ independent instances that share the same expression trees, and hence the same compiled functions, but differ in parameter values.
All instances are evaluated together, with each pair of a data point and an instance mapped to a separate thread.
This adds a second level of parallelism on top of the parallelism across data points, and it is useful for parametric sweeps, ensemble optimization, and warm-starting.
At the time of writing, this extension is under active development in the ExaModels.jl repository and is not part of a released version.
It follows the scenario mechanism of \Cref{sec:extensions:twostage}: declarations are replicated across instances, the internal storage generalizes from vectors to matrices, and a flattening wrapper exposes the batch model as a standard vector-valued \texttt{NLPModels.AbstractNLPModel}.

\subsection{External Oracles}
\label{sec:extensions:oracles}

ExaModels.jl does not aim to provide efficient callbacks for every class of function.
Components such as neural-network layers, PDE solvers, and custom GPU kernels are handled more efficiently by external libraries that already supply derivatives.
For these cases, ExaModels.jl provides an \emph{external oracle} interface.
Let $x$ be a variable block of length $n_x$ and let $f(\cdot): \mathbb{R}^{n_x} \to \mathbb{R}^{m}$ be an opaque map.
The call \texttt{embed\_oracle(core, x, m; \ldots)} adds a block $z$ of $m$ auxiliary output variables, registers $z - f(x) = 0$ as an oracle constraint, and returns the updated core, the block $z$, and the oracle object.
The user supplies matrix-free callbacks that operate on local vectors, without global offsets: \texttt{f!} (function evaluation $y = f(x)$), \texttt{jvp!} (Jacobian--vector product $J_f(x)\, v$), \texttt{vjp!} (transposed-Jacobian--vector product $J_f(x)^\top w$), and an optional \texttt{hvp!} (Hessian--vector product $(\sum_{i=1}^{m} w_i \nabla^2 f_i(x))\,v$).
\begin{lstlisting}
core, z, _ = embed_oracle(
    core, x, N;
    f!   = (y, xv)     -> (y   .= xv .^ 2;   nothing),
    jvp! = (Jv, xv, v) -> (Jv  .= 2 .* xv .* v; nothing),
    vjp! = (Jtv, xv, w)-> (Jtv .= 2 .* xv .* w; nothing),
    hvp! = (Hv, xv, w, v) -> (Hv .= 2 .* w .* v; nothing),
)
@add_obj(core, z[i] + sin(x[i]) for i in 1:N)
\end{lstlisting}
The oracle output $z$ is a regular \texttt{Variable} and can be combined with \gls{simd} expressions in any subsequent objective or constraint.

\section{Numerical Experiments}
\label{sec:numerics}

This section validates the preceding claims numerically.
It measures the speedup those callbacks obtain on a \gls{gpu}, and shows their evaluation time to be $O(1)$ in the number of data points, as claimed earlier, for as long as the device is not saturated.
It then exercises the callbacks inside a \gls{gpu}-resident optimization solver, to quantify what that speedup is worth in a complete solve.

\subsection{Benchmark Problems}
\label{sec:numerics:problems}

We evaluate ExaModels.jl on three benchmark suites that cover a range of problem structures and scales: the Luk\v{s}an--Vl\v{c}ek collection, the COPS test set, and PGLIB-OPF instances.

\paragraph{Luk\v{s}an--Vl\v{c}ek problems.}
The Luk\v{s}an--Vl\v{c}ek collection~\citep{lukVsanSparsePartiallySeparable} consists of 18 scalable, sparse, equality-constrained \glspl{nlp}.
They were designed for testing unconstrained and equality-constrained optimization algorithms.
Most are synthetic instances with no direct application context, but their partially separable structure makes them useful for benchmarking \gls{ad} systems.
Each problem is parameterized by the number of variables~$N$.
The objectives and constraints are sums of partially separable functions, which fits the SIMD abstraction directly.
The collection includes generalized Rosenbrock, Wood, Cragg--Levy, Powell, Broyden (tridiagonal and banded), and chained variants of the Hock--Schittkowski problems (HS46--HS53), among others.
Implementations for ExaModels.jl are provided in LuksanVlcekBenchmark.jl~\citep{LuksanVlcekBenchmarkjl}.

\paragraph{COPS problems.}
The COPS test set~\citep{dolanBenchmarkingOptimizationSoftware2001} is a collection of 18 scalable \glspl{nlp} from applications in optimal control, fluid dynamics, parameter estimation, PDE-constrained optimization, and mesh generation.
The problems include, among others, the Goddard rocket, robot arm, glider trajectory, journal bearing, minimum surface, torsion, catenary chain, and several chemical engineering models (gas--oil, methanol, pinene, catalytic mixing).
Each problem is parameterized by a sizing parameter, which for the PDE-based problems is a discretization resolution.
The suite spans a wide range of problem structures: some problems are unconstrained, others have equality or inequality constraints, and several involve nonlinear PDEs discretized on 2D grids.
Implementations for ExaModels.jl are provided in COPSBenchmark.jl~\citep{COPSBenchmarkjl}.

\paragraph{PGLIB-OPF problems.}
The PGLib-OPF benchmark library~\citep{babaeinejadsarookolaeePowerGridLibrary2021} provides standardized \gls{acopf} instances ranging from 3 to 78{,}484 buses.
These formulations, and the case data format the library distributes them in, have a lineage: MATPOWER~\citep{zimmermanMATPOWERSteadyState2011} established the reference implementations and the case format that the field still uses, and PowerModels.jl~\citep{coffrinPowerModelsJLOpenSource2018} later expressed the polar, rectangular, and relaxed power flow formulations over that format in a common framework.
We test all 66 instances in the polar (ACP) voltage formulation using ExaModelsPower.jl~\citep{johnsonExaModelsPowerjl2025}, which provides ExaModels implementations in that tradition.
Every instance has the same algebraic structure of 15 algebraic patterns.
The instances differ only in the network data, that is, the bus parameters, the line admittances, the power demands, and the number of data points over which each pattern is repeated.
ExaModels.jl therefore compiles the expression trees and the \gls{ad} kernels once for the whole set of 66 instances.
The first instance triggers compilation, and every later one reuses the compiled functions and changes only the data arrays.
The pattern--data separation of \Cref{sec:modeling:simd} thus applies at the scale of a benchmark library.
The OPF problems also carry inequality constraints and variable bounds, which the LV problems do not.
The OPF instances are parameterized by the network topology rather than a scalar sizing parameter, making them complementary to the scalable LV and COPS problems.

\paragraph{Exclusion of CUTEst.}
CUTEst~\citep{gouldCUTEstConstrainedUnconstrained2015} is the standard test environment for mathematical optimization.
We do not use it for three reasons.
First, most instances in the collection are small enough that data-parallel evaluation has nothing to work with, so \gls{gpu} acceleration is neither necessary nor useful for them.
Second, those instances that are scalable are each parameterized by their own set of admissible dimensions rather than by a common scalable size, so the collection does not offer a size sweep of the kind our experiments require.
Third, the large CUTEst instances are drawn largely from the same families as the LV and COPS suites, which we benchmark directly.

\paragraph{Instance sizes.}
For each LV and COPS problem, we select three values of the sizing parameter targeting approximately $10^2$, $10^4$, and $10^6$ Hessian nonzeros (nnzh), producing \emph{small}, \emph{medium}, and \emph{large} instances.
The relation between the sizing parameter and nnzh varies across problems, being linear for chain, quadratic for bearing, and cubic for polygon, so the sizes are chosen per problem to reach the target ranges.
For the LV suite, we use $N \in \{20,\; 2{,}000,\; 200{,}000\}$ across all 18 problems.
For COPS, the sizes are problem-specific: for example, bearing uses grid sizes $(5\times 5)$, $(35 \times 35)$, and $(320 \times 320)$, while catmix uses $N \in \{3,\; 200,\; 20{,}000\}$.
The PGLIB-OPF instances span a natural range from 3 to 78{,}484 buses.
In total, the benchmark comprises $18 \times 3 = 54$ LV instances, $18 \times 3 = 54$ COPS instances, and $66$ OPF instances, for 174 benchmark instances.

\subsection{Experimental Setup}
\label{sec:numerics:setup}

\paragraph{Reproducibility.}
The benchmark implementation, the generators that produce every table and figure in this section, and the exact package environment are available at \url{https://github.com/exanauts/exa-models-benchmark}.
The batch scripts that ran each leg of the campaign are there too, in \texttt{slurm/} and named by site as \texttt{orcd-*} and \texttt{jlse-*}, alongside a \texttt{Makefile} whose targets drive the runs, fetch the archived results, and rebuild the paper from them.
Each run is archived on the \texttt{results} branch of that repository as a self-contained bundle, one directory per run, holding the benchmarking script, the hardware description of the machine, the raw per-instance timings in machine-readable form, and the resolved \texttt{Manifest.toml}.
The hardware used for benchmarking is listed in \Cref{tab:hardware}, where each platform is named by a label whose letter gives the vendor: C is \gls{cpu} only, N NVIDIA, A AMD, and M Apple.
The AMD and B200 systems are provided by the Joint Laboratory for System Evaluation at Argonne National Laboratory, and all other clusters are part of MIT's Office of Research Computing and Data.
The manifest fixes the software environment of that run, recording the Julia version and the resolved version of every package in the dependency graph, including the commit of ExaModels.jl used to generate the results, KernelAbstractions.jl, CUDA.jl, and the CUDA runtime and driver artifacts.
The results reported here were produced with Julia~1.12 and ExaModels.jl~v0.11.2~\citep{shinExaModelsjl2025}, in Float64 arithmetic unless noted.
Every number reported in this section can therefore be traced to the run that produced it and to the environment in which it ran.

\begin{table}[t]
\centering
\caption{Hardware platforms used in the benchmark experiments.}
\label{tab:hardware}
\scriptsize
\begin{tabular*}{\textwidth}{@{\extracolsep{\fill}}ll rr l}
  \toprule
  \textbf{Platform} & \textbf{CPU} & \textbf{Cores} & \textbf{RAM} & \textbf{GPU} \\
  \midrule
  \textbf{A1} & EPYC 9654 & 192 & 2266.6 GiB & Instinct MI300X \\
  \textbf{C1} & Xeon Gold 5420+ & 56 & 251.1 GiB & -- \\
  \textbf{C3} & EPYC 9474F & 96 & 376.9 GiB & -- \\
  \textbf{I1} & Xeon CPU Max 9470C & 104 & 1088.0 GiB & Intel Data Center GPU Max 1550 \\
  \textbf{M1} & M2 Pro & 12 & 16.0 GiB & M2 Pro \\
  \textbf{N1} & Xeon 6960P & 144 & 2266.6 GiB & B200 \\
  \textbf{N2} & EPYC 7763 & 64 & 503.2 GiB & A100-SXM4-80GB \\
  \textbf{N3} & Xeon Platinum 8580 & 120 & 2014.9 GiB & H100 80GB \\
  \textbf{N4} & Xeon Platinum 8562Y+ & 64 & 1006.9 GiB & L40S \\
  \textbf{N6} & Xeon 6767P & 128 & 2014.9 GiB & RTX PRO 6000 Blackwell \\
  \bottomrule
\end{tabular*}

\end{table}

\paragraph{Compilation.}
Julia compiles the model-construction code and the evaluation callbacks the first time a given model structure is encountered.
That first-time cost is substantial, and it is specific to the problem, since a new set of algebraic patterns requires a new set of compiled functions.
The concern of this paper is runtime performance, so every timing reported in this section is a steady-state measurement and excludes that compilation.
Model creation and compilation are measured jointly in \Cref{sec:numerics:compilation}.
This holds for the model creation times as well: they measure steady-state construction, not the first call that triggers compilation.

\paragraph{Aggregate metric.}
We report the \gls{sgm} with a shift of $10\,\mu$s as the aggregate performance metric.
Denoting the wall times of a modeling system on problems $i = 1, \ldots, n$ by $t_i$, the SGM is
\begin{equation*}
  \operatorname{SGM}(t; \sigma) = \exp\!\Bigl(\frac{1}{n}\sum_{i=1}^n \ln(t_i + \sigma)\Bigr) - \sigma,
\end{equation*}
where the shift is $\sigma = 10^{-5}$\,s.
The SGM is less sensitive to outliers than the arithmetic mean and is standard in optimization benchmarking.
Instances are classified into the \emph{Small}, \emph{Medium}, and \emph{Large} size classes by $\text{nnz} = \max(\text{nnzj}, \text{nnzh})$.
Full per-instance timings are provided in \Cref{app:tables}.

\paragraph{\gls{gpu} timing protocol.}
A \gls{gpu} callback carries two costs that a \gls{cpu} callback does not.
The first is the kernel launch cost, a few microseconds per call, which a solver pays as well and which every timing we report includes.
The second is host--device synchronization, paid only when the host reads a device value.
A device-resident solver does not do that after each callback, since the results stay in device arrays that the following linear algebra consumes, so this cost is purely an artifact of measuring from the host rather than a property of the callback.
We exclude it not by measuring and subtracting it, but by making it negligible.
The timing routine of the benchmarking script is, in simplified form:
\begin{lstlisting}[basicstyle=\small\ttfamily]
f(); sync()      # warm up
sync()           # quiesce the device
t0 = time_ns()
for i = 1:N      # N back-to-back launches
    f()
end
sync()           # wait for the last kernel
(time_ns() - t0) / 1e9 / N    # apparent per-call time
\end{lstlisting}
The batch size $N$ is calibrated from three synchronized single calls so that the batch spans $0.5$\,s, and the reported value is the best of three batches.
The apparent per-call time satisfies
\begin{equation*}
  \hat{t} = \frac{N t + t_{\mathrm{sync}}}{N} = t + \frac{t_{\mathrm{sync}}}{N},
\end{equation*}
where $\hat{t}$ is the reported estimate of the average callback time, $t$ is the true average per-call cost, $t_{\mathrm{sync}}$ is the one synchronization inside the timed interval, and $N$ is the batch size.
The second term is the bias.
It decreases as $1/N$ and is determined by the batch duration $N t$ rather than by the speed of the individual kernel.

A $0.5$\,s batch against a synchronization of order $10\,\mu$s therefore yields a relative bias of approximately $2 \times 10^{-5}$ for every callback, irrespective of its speed.
The bias is positive, so the reported time overestimates the \gls{gpu} callback time.

\subsection{NLP Function Evaluation}
\label{sec:numerics:ad}

We measure the \gls{nlp} function evaluation performance of ExaModels.jl, on the \gls{cpu} and on a range of \gls{gpu} platforms.
The callbacks are \texttt{NLPModels.obj}, \texttt{NLPModels.cons!}, \texttt{NLPModels.grad!}, \texttt{NLPModels.jac\_coord!}, and \texttt{NLPModels.hess\_coord!}.
Every timing is the minimum wall time over repeated trials, which is robust to system noise.
Speedups are reported against single-threaded execution, the standard sequential baseline.
The multithreaded row of each summary table shows the best thread count among 4, 8, and 16; the full ladder appears in the appendix tables, and the summary row's choice reflects where adding threads stops helping.
Each summary table reports one NVIDIA \gls{gpu}, the one with the highest aggregate speedup on the large size class, together with the AMD and Intel \glspl{gpu}.
The composite row of each summary table aggregates the five callbacks into one evaluation cost per solve.
The composite approximates the time a modeling system spends on \gls{nlp} function evaluation when it is deployed inside an optimization solver.
Each instance's callback times are weighted by the call counts of a \gls{gpu}-resident solve of the largest \gls{opf} instance, \cntObj{} objective, \cntCons{} constraint, \cntGrad{} gradient, \cntJac{} Jacobian, and \cntHess{} Hessian evaluations, and the row reports the \gls{sgm} of the weighted sum over the instances of the class.
The weights are derived from that one \gls{opf} solve and reused across all three suites, so on the LV and COPS suites the row is a weighting of measured callback times, not a measurement of a solve.

The full per-instance results appear in \Cref{app:tables}, and \Cref{tab:gpu_lv,tab:gpu_cops,tab:gpu_opf_polar} summarize them as \glspl{sgm}, one table per suite.

\begin{table}[t]
\centering
\caption{GPU speedup over single-threaded CPU: Luk\v{s}an--Vl\v{c}ek (SGM, $\sigma = 10^{-5}$\,s).}
\label{tab:gpu_lv}
\scriptsize
\setlength{\tabcolsep}{3pt}
\begin{tabular*}{\textwidth}{@{\extracolsep{\fill}}ll rrrr}
  \toprule
  & & \multicolumn{4}{c@{}}{SGM time, and speedup over single-thread ExaModels CPU (C3)} \\
  \cmidrule(l){3-6}
  \textbf{callback} & \textbf{backend} & \shortstack{Small (18)\\[-2pt]{\tiny nnz$<10^3$}} & \shortstack{Medium (18)\\[-2pt]{\tiny $10^3{\leq}$nnz${<}10^5$}} & \shortstack{Large (18)\\[-2pt]{\tiny $10^5{\leq}$nnz}} & Total (54) \\
  \midrule
  \texttt{obj} & CPU (C3) & \textbf{0.1\,\textmu s ($1.0\times$)} & \textbf{5.6\,\textmu s ($1.0\times$)} & 371.4\,\textmu s ($1.0\times$) & 29.1\,\textmu s ($1.0\times$) \\
   & CPU-8T (C3) & 0.2\,\textmu s ($0.4\times$) & 6.0\,\textmu s ($0.9\times$) & 270.9\,\textmu s ($1.4\times$) & \textbf{25.8\,\textmu s ($1.1\times$)} \\
   & CUDA (N1) & 41.9\,\textmu s ($0.0\times$) & 48.0\,\textmu s ($0.1\times$) & \textbf{52.9\,\textmu s ($7.0\times$)} & 47.4\,\textmu s ($0.6\times$) \\
   & AMDGPU (A1) & 51.3\,\textmu s ($0.0\times$) & 68.3\,\textmu s ($0.1\times$) & 89.1\,\textmu s ($4.2\times$) & 68.1\,\textmu s ($0.4\times$) \\
   & oneAPI (I1) & 114.0\,\textmu s ($0.0\times$) & 143.4\,\textmu s ($0.0\times$) & 273.9\,\textmu s ($1.4\times$) & 165.5\,\textmu s ($0.2\times$) \\
  \midrule
  \texttt{cons!} & CPU (C3) & \textbf{0.1\,\textmu s ($1.0\times$)} & \textbf{3.1\,\textmu s ($1.0\times$)} & 149.3\,\textmu s ($1.0\times$) & 17.6\,\textmu s ($1.0\times$) \\
   & CPU-8T (C3) & 0.4\,\textmu s ($0.2\times$) & 4.7\,\textmu s ($0.7\times$) & 49.0\,\textmu s ($3.0\times$) & \textbf{10.8\,\textmu s ($1.6\times$)} \\
   & CUDA (N1) & 13.9\,\textmu s ($0.0\times$) & 14.3\,\textmu s ($0.2\times$) & \textbf{14.4\,\textmu s ($10\times$)} & 14.2\,\textmu s ($1.2\times$) \\
   & AMDGPU (A1) & 28.9\,\textmu s ($0.0\times$) & 28.6\,\textmu s ($0.1\times$) & 28.7\,\textmu s ($5.2\times$) & 28.7\,\textmu s ($0.6\times$) \\
   & oneAPI (I1) & 24.0\,\textmu s ($0.0\times$) & 24.4\,\textmu s ($0.1\times$) & 25.6\,\textmu s ($5.8\times$) & 24.7\,\textmu s ($0.7\times$) \\
  \midrule
  \texttt{grad!} & CPU (C3) & \textbf{0.1\,\textmu s ($1.0\times$)} & \textbf{11.5\,\textmu s ($1.0\times$)} & 1.05\,ms ($1.0\times$) & 51.4\,\textmu s ($1.0\times$) \\
   & CPU-8T (C3) & 0.4\,\textmu s ($0.4\times$) & 27.8\,\textmu s ($0.4\times$) & 1.15\,ms ($0.9\times$) & 66.9\,\textmu s ($0.8\times$) \\
   & CUDA (N1) & 20.7\,\textmu s ($0.0\times$) & 21.8\,\textmu s ($0.5\times$) & \textbf{22.3\,\textmu s ($47\times$)} & \textbf{21.6\,\textmu s ($2.4\times$)} \\
   & AMDGPU (A1) & 47.3\,\textmu s ($0.0\times$) & 46.6\,\textmu s ($0.2\times$) & 46.7\,\textmu s ($23\times$) & 46.9\,\textmu s ($1.1\times$) \\
   & oneAPI (I1) & 42.1\,\textmu s ($0.0\times$) & 41.5\,\textmu s ($0.3\times$) & 42.1\,\textmu s ($25\times$) & 41.9\,\textmu s ($1.2\times$) \\
  \midrule
  \texttt{jac\_coord!} & CPU (C3) & \textbf{0.1\,\textmu s ($1.0\times$)} & 7.5\,\textmu s ($1.0\times$) & 349.4\,\textmu s ($1.0\times$) & 29.9\,\textmu s ($1.0\times$) \\
   & CPU-8T (C3) & 0.6\,\textmu s ($0.2\times$) & \textbf{7.5\,\textmu s ($1.0\times$)} & 198.7\,\textmu s ($1.8\times$) & 23.8\,\textmu s ($1.3\times$) \\
   & CUDA (N1) & 19.2\,\textmu s ($0.0\times$) & 19.9\,\textmu s ($0.4\times$) & \textbf{20.4\,\textmu s ($17\times$)} & \textbf{19.9\,\textmu s ($1.5\times$)} \\
   & AMDGPU (A1) & 41.6\,\textmu s ($0.0\times$) & 40.8\,\textmu s ($0.2\times$) & 40.9\,\textmu s ($8.5\times$) & 41.1\,\textmu s ($0.7\times$) \\
   & oneAPI (I1) & 44.8\,\textmu s ($0.0\times$) & 46.5\,\textmu s ($0.2\times$) & 50.5\,\textmu s ($6.9\times$) & 47.2\,\textmu s ($0.6\times$) \\
  \midrule
  \texttt{hess\_coord!} & CPU (C3) & \textbf{0.4\,\textmu s ($1.0\times$)} & 29.9\,\textmu s ($1.0\times$) & 2.73\,ms ($1.0\times$) & 94.3\,\textmu s ($1.0\times$) \\
   & CPU-8T (C3) & 1.0\,\textmu s ($0.4\times$) & \textbf{25.3\,\textmu s ($1.2\times$)} & 689.7\,\textmu s ($4.0\times$) & 54.8\,\textmu s ($1.7\times$) \\
   & CUDA (N1) & 26.9\,\textmu s ($0.0\times$) & 28.3\,\textmu s ($1.1\times$) & \textbf{35.7\,\textmu s ($76\times$)} & \textbf{30.1\,\textmu s ($3.1\times$)} \\
   & AMDGPU (A1) & 57.1\,\textmu s ($0.0\times$) & 56.4\,\textmu s ($0.5\times$) & 61.2\,\textmu s ($45\times$) & 58.2\,\textmu s ($1.6\times$) \\
   & oneAPI (I1) & 71.8\,\textmu s ($0.0\times$) & 73.8\,\textmu s ($0.4\times$) & 99.5\,\textmu s ($27\times$) & 80.9\,\textmu s ($1.2\times$) \\
  \midrule
  composite & CPU (C3) & \textbf{79.8\,\textmu s ($1.0\times$)} & \textbf{5.96\,ms ($1.0\times$)} & 615.17\,ms ($1.0\times$) & \textbf{6.9\,ms ($1.0\times$)} \\
   & CPU-8T (C3) & 277.0\,\textmu s ($0.3\times$) & 8.06\,ms ($0.7\times$) & 293.5\,ms ($2.1\times$) & 8.78\,ms ($0.8\times$) \\
   & CUDA (N1) & 14.18\,ms ($0.0\times$) & 15.35\,ms ($0.4\times$) & \textbf{16.86\,ms ($36\times$)} & 15.42\,ms ($0.4\times$) \\
   & AMDGPU (A1) & 26.26\,ms ($0.0\times$) & 27.98\,ms ($0.2\times$) & 31.16\,ms ($20\times$) & 28.4\,ms ($0.2\times$) \\
   & oneAPI (I1) & 34.7\,ms ($0.0\times$) & 38.66\,ms ($0.2\times$) & 57.71\,ms ($11\times$) & 42.62\,ms ($0.2\times$) \\
  \bottomrule
\end{tabular*}

\end{table}

\begin{table}[t]
\centering
\caption{GPU speedup over single-threaded CPU: COPS (SGM, $\sigma = 10^{-5}$\,s).}
\label{tab:gpu_cops}
\scriptsize
\setlength{\tabcolsep}{3pt}
\begin{tabular*}{\textwidth}{@{\extracolsep{\fill}}ll rrrr}
  \toprule
  & & \multicolumn{4}{c@{}}{SGM time, and speedup over single-thread ExaModels CPU (C3)} \\
  \cmidrule(l){3-6}
  \textbf{callback} & \textbf{backend} & \shortstack{Small (16)\\[-2pt]{\tiny nnz$<10^3$}} & \shortstack{Medium (18)\\[-2pt]{\tiny $10^3{\leq}$nnz${<}10^5$}} & \shortstack{Large (20)\\[-2pt]{\tiny $10^5{\leq}$nnz}} & Total (54) \\
  \midrule
  \texttt{obj} & CPU (C3) & \textbf{0.5\,\textmu s ($1.0\times$)} & \textbf{3.5\,\textmu s ($1.0\times$)} & 19.9\,\textmu s ($1.0\times$) & \textbf{6.8\,\textmu s ($1.0\times$)} \\
   & CPU-8T (C3) & 0.7\,\textmu s ($0.7\times$) & 5.6\,\textmu s ($0.6\times$) & \textbf{18.9\,\textmu s ($1.1\times$)} & 7.5\,\textmu s ($0.9\times$) \\
   & CUDA (N1) & 44.5\,\textmu s ($0.0\times$) & 47.2\,\textmu s ($0.1\times$) & 48.1\,\textmu s ($0.4\times$) & 46.7\,\textmu s ($0.1\times$) \\
   & AMDGPU (A1) & 56.4\,\textmu s ($0.0\times$) & 60.6\,\textmu s ($0.1\times$) & 65.5\,\textmu s ($0.3\times$) & 61.0\,\textmu s ($0.1\times$) \\
   & oneAPI (I1) & 129.4\,\textmu s ($0.0\times$) & 157.8\,\textmu s ($0.0\times$) & 1.03\,ms ($0.0\times$) & 303.4\,\textmu s ($0.0\times$) \\
  \midrule
  \texttt{cons!} & CPU (C3) & \textbf{0.4\,\textmu s ($1.0\times$)} & \textbf{13.8\,\textmu s ($1.0\times$)} & 614.5\,\textmu s ($1.0\times$) & 52.5\,\textmu s ($1.0\times$) \\
   & CPU-8T (C3) & 3.3\,\textmu s ($0.1\times$) & 20.1\,\textmu s ($0.7\times$) & 179.9\,\textmu s ($3.4\times$) & \textbf{36.7\,\textmu s ($1.4\times$)} \\
   & CUDA (N1) & 43.8\,\textmu s ($0.0\times$) & 52.2\,\textmu s ($0.3\times$) & \textbf{65.2\,\textmu s ($9.4\times$)} & 53.9\,\textmu s ($1.0\times$) \\
   & AMDGPU (A1) & 84.6\,\textmu s ($0.0\times$) & 97.8\,\textmu s ($0.1\times$) & 109.3\,\textmu s ($5.6\times$) & 97.7\,\textmu s ($0.5\times$) \\
   & oneAPI (I1) & 86.2\,\textmu s ($0.0\times$) & 103.2\,\textmu s ($0.1\times$) & 113.4\,\textmu s ($5.4\times$) & 101.4\,\textmu s ($0.5\times$) \\
  \midrule
  \texttt{grad!} & CPU (C3) & \textbf{1.0\,\textmu s ($1.0\times$)} & \textbf{5.7\,\textmu s ($1.0\times$)} & 60.0\,\textmu s ($1.0\times$) & \textbf{14.6\,\textmu s ($1.0\times$)} \\
   & CPU-8T (C3) & 1.1\,\textmu s ($0.9\times$) & 8.4\,\textmu s ($0.7\times$) & 70.5\,\textmu s ($0.9\times$) & 17.4\,\textmu s ($0.8\times$) \\
   & CUDA (N1) & 23.0\,\textmu s ($0.0\times$) & 22.8\,\textmu s ($0.2\times$) & \textbf{24.8\,\textmu s ($2.4\times$)} & 23.6\,\textmu s ($0.6\times$) \\
   & AMDGPU (A1) & 52.3\,\textmu s ($0.0\times$) & 47.3\,\textmu s ($0.1\times$) & 49.0\,\textmu s ($1.2\times$) & 49.4\,\textmu s ($0.3\times$) \\
   & oneAPI (I1) & 51.6\,\textmu s ($0.0\times$) & 47.2\,\textmu s ($0.1\times$) & 47.8\,\textmu s ($1.3\times$) & 48.7\,\textmu s ($0.3\times$) \\
  \midrule
  \texttt{jac\_coord!} & CPU (C3) & \textbf{0.5\,\textmu s ($1.0\times$)} & \textbf{21.2\,\textmu s ($1.0\times$)} & 937.4\,\textmu s ($1.0\times$) & 69.9\,\textmu s ($1.0\times$) \\
   & CPU-8T (C3) & 3.2\,\textmu s ($0.1\times$) & 23.5\,\textmu s ($0.9\times$) & 359.8\,\textmu s ($2.6\times$) & \textbf{51.9\,\textmu s ($1.3\times$)} \\
   & CUDA (N1) & 52.7\,\textmu s ($0.0\times$) & 57.0\,\textmu s ($0.4\times$) & \textbf{70.6\,\textmu s ($13\times$)} & 60.4\,\textmu s ($1.2\times$) \\
   & AMDGPU (A1) & 101.6\,\textmu s ($0.0\times$) & 103.4\,\textmu s ($0.2\times$) & 122.0\,\textmu s ($7.7\times$) & 109.4\,\textmu s ($0.6\times$) \\
   & oneAPI (I1) & 155.9\,\textmu s ($0.0\times$) & 164.2\,\textmu s ($0.1\times$) & 182.0\,\textmu s ($5.2\times$) & 168.1\,\textmu s ($0.4\times$) \\
  \midrule
  \texttt{hess\_coord!} & CPU (C3) & \textbf{1.9\,\textmu s ($1.0\times$)} & \textbf{50.2\,\textmu s ($1.0\times$)} & 2.83\,ms ($1.0\times$) & 145.2\,\textmu s ($1.0\times$) \\
   & CPU-8T (C3) & 4.9\,\textmu s ($0.4\times$) & 52.8\,\textmu s ($0.9\times$) & 754.3\,\textmu s ($3.7\times$) & 93.5\,\textmu s ($1.6\times$) \\
   & CUDA (N1) & 71.9\,\textmu s ($0.0\times$) & 75.6\,\textmu s ($0.7\times$) & \textbf{92.9\,\textmu s ($30\times$)} & \textbf{80.4\,\textmu s ($1.8\times$)} \\
   & AMDGPU (A1) & 140.5\,\textmu s ($0.0\times$) & 140.6\,\textmu s ($0.4\times$) & 172.0\,\textmu s ($16\times$) & 151.5\,\textmu s ($1.0\times$) \\
   & oneAPI (I1) & 211.5\,\textmu s ($0.0\times$) & 224.8\,\textmu s ($0.2\times$) & 229.3\,\textmu s ($12\times$) & 222.4\,\textmu s ($0.7\times$) \\
  \midrule
  composite & CPU (C3) & \textbf{287.9\,\textmu s ($1.0\times$)} & \textbf{10.72\,ms ($1.0\times$)} & 816.87\,ms ($1.0\times$) & \textbf{18.45\,ms ($1.0\times$)} \\
   & CPU-8T (C3) & 1.1\,ms ($0.3\times$) & 14.44\,ms ($0.7\times$) & 252.04\,ms ($3.2\times$) & 19.48\,ms ($0.9\times$) \\
   & CUDA (N1) & 28.06\,ms ($0.0\times$) & 30.78\,ms ($0.3\times$) & \textbf{38.1\,ms ($21\times$)} & 32.41\,ms ($0.6\times$) \\
   & AMDGPU (A1) & 52.03\,ms ($0.0\times$) & 54.72\,ms ($0.2\times$) & 67.64\,ms ($12\times$) & 58.31\,ms ($0.3\times$) \\
   & oneAPI (I1) & 77.92\,ms ($0.0\times$) & 88.97\,ms ($0.1\times$) & 249.29\,ms ($3.3\times$) & 125.61\,ms ($0.1\times$) \\
  \bottomrule
\end{tabular*}

\end{table}

\begin{table}[t]
\centering
\caption{GPU speedup over single-threaded CPU: PGLIB-OPF (polar) (SGM, $\sigma = 10^{-5}$\,s).}
\label{tab:gpu_opf_polar}
\scriptsize
\setlength{\tabcolsep}{3pt}
\begin{tabular*}{\textwidth}{@{\extracolsep{\fill}}ll rrrr}
  \toprule
  & & \multicolumn{4}{c@{}}{SGM time, and speedup over single-thread ExaModels CPU (C3)} \\
  \cmidrule(l){3-6}
  \textbf{callback} & \textbf{backend} & \shortstack{Small (3)\\[-2pt]{\tiny nnz$<10^3$}} & \shortstack{Medium (19)\\[-2pt]{\tiny $10^3{\leq}$nnz${<}10^5$}} & \shortstack{Large (44)\\[-2pt]{\tiny $10^5{\leq}$nnz}} & Total (66) \\
  \midrule
  \texttt{obj} & CPU (C3) & \textbf{0.0\,\textmu s ($1.0\times$)} & \textbf{0.0\,\textmu s ($1.0\times$)} & \textbf{0.0\,\textmu s ($1.0\times$)} & \textbf{0.0\,\textmu s ($1.0\times$)} \\
   & CPU-4T (C3) & 0.2\,\textmu s ($0.1\times$) & 0.3\,\textmu s ($0.1\times$) & 2.3\,\textmu s ($0.0\times$) & 1.6\,\textmu s ($0.0\times$) \\
   & CUDA (N1) & 40.6\,\textmu s ($0.0\times$) & 43.3\,\textmu s ($0.0\times$) & 46.4\,\textmu s ($0.0\times$) & 45.2\,\textmu s ($0.0\times$) \\
   & AMDGPU (A1) & 50.2\,\textmu s ($0.0\times$) & 52.0\,\textmu s ($0.0\times$) & 66.2\,\textmu s ($0.0\times$) & 61.1\,\textmu s ($0.0\times$) \\
   & oneAPI (I1) & 360.6\,\textmu s ($0.0\times$) & 365.3\,\textmu s ($0.0\times$) & 464.5\,\textmu s ($0.0\times$) & 428.6\,\textmu s ($0.0\times$) \\
  \midrule
  \texttt{cons!} & CPU (C3) & \textbf{0.5\,\textmu s ($1.0\times$)} & \textbf{12.3\,\textmu s ($1.0\times$)} & 378.2\,\textmu s ($1.0\times$) & \textbf{134.8\,\textmu s ($1.0\times$)} \\
   & CPU-4T (C3) & 41.7\,\textmu s ($0.0\times$) & 62.7\,\textmu s ($0.2\times$) & 360.4\,\textmu s ($1.0\times$) & 202.0\,\textmu s ($0.7\times$) \\
   & CUDA (N1) & 185.9\,\textmu s ($0.0\times$) & 199.2\,\textmu s ($0.1\times$) & \textbf{203.2\,\textmu s ($1.9\times$)} & 201.2\,\textmu s ($0.7\times$) \\
   & AMDGPU (A1) & 293.8\,\textmu s ($0.0\times$) & 284.2\,\textmu s ($0.0\times$) & 280.1\,\textmu s ($1.4\times$) & 281.9\,\textmu s ($0.5\times$) \\
   & oneAPI (I1) & 321.3\,\textmu s ($0.0\times$) & 454.0\,\textmu s ($0.0\times$) & 434.9\,\textmu s ($0.9\times$) & 434.3\,\textmu s ($0.3\times$) \\
  \midrule
  \texttt{grad!} & CPU (C3) & \textbf{0.0\,\textmu s ($1.0\times$)} & \textbf{0.2\,\textmu s ($1.0\times$)} & \textbf{4.7\,\textmu s ($1.0\times$)} & \textbf{3.0\,\textmu s ($1.0\times$)} \\
   & CPU-4T (C3) & 0.2\,\textmu s ($0.1\times$) & 0.6\,\textmu s ($0.4\times$) & 7.5\,\textmu s ($0.6\times$) & 4.8\,\textmu s ($0.6\times$) \\
   & CUDA (N1) & 19.9\,\textmu s ($0.0\times$) & 20.9\,\textmu s ($0.0\times$) & 21.9\,\textmu s ($0.2\times$) & 21.5\,\textmu s ($0.1\times$) \\
   & AMDGPU (A1) & 45.8\,\textmu s ($0.0\times$) & 45.5\,\textmu s ($0.0\times$) & 45.4\,\textmu s ($0.1\times$) & 45.4\,\textmu s ($0.1\times$) \\
   & oneAPI (I1) & 32.2\,\textmu s ($0.0\times$) & 40.1\,\textmu s ($0.0\times$) & 39.6\,\textmu s ($0.1\times$) & 39.4\,\textmu s ($0.1\times$) \\
  \midrule
  \texttt{jac\_coord!} & CPU (C3) & \textbf{1.0\,\textmu s ($1.0\times$)} & \textbf{24.6\,\textmu s ($1.0\times$)} & 792.5\,\textmu s ($1.0\times$) & 257.1\,\textmu s ($1.0\times$) \\
   & CPU-4T (C3) & 22.8\,\textmu s ($0.0\times$) & 46.7\,\textmu s ($0.5\times$) & 475.2\,\textmu s ($1.7\times$) & 221.4\,\textmu s ($1.2\times$) \\
   & CUDA (N1) & 138.0\,\textmu s ($0.0\times$) & 155.5\,\textmu s ($0.2\times$) & \textbf{158.9\,\textmu s ($5.0\times$)} & \textbf{156.9\,\textmu s ($1.6\times$)} \\
   & AMDGPU (A1) & 253.4\,\textmu s ($0.0\times$) & 243.5\,\textmu s ($0.1\times$) & 241.5\,\textmu s ($3.3\times$) & 242.6\,\textmu s ($1.1\times$) \\
   & oneAPI (I1) & 360.2\,\textmu s ($0.0\times$) & 559.1\,\textmu s ($0.0\times$) & 529.9\,\textmu s ($1.5\times$) & 528.8\,\textmu s ($0.5\times$) \\
  \midrule
  \texttt{hess\_coord!} & CPU (C3) & \textbf{1.5\,\textmu s ($1.0\times$)} & \textbf{37.3\,\textmu s ($1.0\times$)} & 1.25\,ms ($1.0\times$) & 385.6\,\textmu s ($1.0\times$) \\
   & CPU-4T (C3) & 23.9\,\textmu s ($0.1\times$) & 63.1\,\textmu s ($0.6\times$) & 631.6\,\textmu s ($2.0\times$) & 290.4\,\textmu s ($1.3\times$) \\
   & CUDA (N1) & 153.4\,\textmu s ($0.0\times$) & 168.0\,\textmu s ($0.2\times$) & \textbf{172.3\,\textmu s ($7.3\times$)} & \textbf{170.1\,\textmu s ($2.3\times$)} \\
   & AMDGPU (A1) & 267.3\,\textmu s ($0.0\times$) & 262.1\,\textmu s ($0.1\times$) & 260.7\,\textmu s ($4.8\times$) & 261.4\,\textmu s ($1.5\times$) \\
   & oneAPI (I1) & 399.2\,\textmu s ($0.0\times$) & 613.3\,\textmu s ($0.1\times$) & 585.1\,\textmu s ($2.1\times$) & 582.9\,\textmu s ($0.7\times$) \\
  \midrule
  composite & CPU (C3) & \textbf{272.7\,\textmu s ($1.0\times$)} & \textbf{7.08\,ms ($1.0\times$)} & 274.29\,ms ($1.0\times$) & 70.05\,ms ($1.0\times$) \\
   & CPU-4T (C3) & 10.15\,ms ($0.0\times$) & 19.48\,ms ($0.4\times$) & 167.95\,ms ($1.6\times$) & 79.52\,ms ($0.9\times$) \\
   & CUDA (N1) & 61.48\,ms ($0.0\times$) & 67.18\,ms ($0.1\times$) & \textbf{68.81\,ms ($4.0\times$)} & \textbf{67.99\,ms ($1.0\times$)} \\
   & AMDGPU (A1) & 103.89\,ms ($0.0\times$) & 101.24\,ms ($0.1\times$) & 102.03\,ms ($2.7\times$) & 101.88\,ms ($0.7\times$) \\
   & oneAPI (I1) & 219.23\,ms ($0.0\times$) & 285.82\,ms ($0.0\times$) & 293.25\,ms ($0.9\times$) & 287.26\,ms ($0.2\times$) \\
  \bottomrule
\end{tabular*}

\end{table}

ExaModels.jl compiles a derivative evaluation function specialized to the type of each expression tree, and pays that compilation cost upfront, once per model structure (\Cref{sec:numerics:compilation}).
\Cref{app:llvm} shows the compiled gradient kernel that ExaModels.jl produces, in which almost nothing remains beyond the arithmetic the pattern requires.

\label{sec:numerics:ad:gpu}

The \gls{gpu} rows of the same summary tables run the same callbacks on the same instances under the timing protocol of \Cref{sec:numerics:setup}.
The \gls{gpu} speedups grow with problem size across the three suites, reaching $\gsLvHessLarge\times$ for the Hessian on the large size class of the LV suite and $\gsCopsHessLarge\times$ on the large size class of the COPS suite.
On small instances the \gls{gpu} is slower than the \gls{cpu}, and on medium instances the two are within a factor of two.

The growth of the speedup reflects \gls{gpu} times that are nearly constant in the problem size rather than \gls{gpu} performance that improves with size.
On the LV suite, the Hessian takes $\gsLvHessAbsSmall$, $\gsLvHessAbsMedium$, and $\gsLvHessAbsLarge\,\mu$s on the small, medium, and large size classes, and the Jacobian $\gsLvJacAbsSmall$, $\gsLvJacAbsMedium$, and $\gsLvJacAbsLarge\,\mu$s, while the \gls{cpu} times grow by more than three orders of magnitude across the same classes.

The absolute level of the \gls{gpu} times, a few tens of microseconds, is set by the kernel launch overhead.
The flatness across four orders of magnitude in problem size is the constant-time parallel evaluation argued in \Cref{sec:intro}, realized in practice for as long as the device is not saturated.

The spread across devices is modest.
On the large-instance LV Hessian, the five NVIDIA cards measured span $\gsLvSpanLo$--$\gsLvSpanHi\times$, the AMD \platAOneName{} reaches $\gsLvAmdHessLarge\times$, the Intel \platIOneName{} reaches $\gsLvIntelHessLarge\times$, and the Apple \platMOneName{} reaches $\gsLvAppleHessLarge\times$.
The COPS suite shows the same modest spread.
The same kernels run across four vendors without change.
The Apple timings, in the per-instance figures and the appendix tables, are in \texttt{Float32}, while all other series are in \texttt{Float64}.

The LV problems show larger \gls{gpu} speedups than the COPS problems.
The absolute Hessian times show the same thing, $\gsCopsHessAbsLarge\,\mu$s on large COPS instances against $\gsLvHessAbsLarge\,\mu$s on large LV instances.
On the PGLIB-OPF instances the acceleration is more moderate still, reaching $\gsOpfPolarHessLarge\times$ for the Hessian on the large size class, and the objective and gradient callbacks show no speedup there.
The LV expression trees are simpler and more uniform, which gives a higher arithmetic intensity per data point.
The COPS problems have more expression structures, as in the rocket dynamics with its exponentials and divisions, and this raises register pressure and lowers \gls{gpu} occupancy.
The OPF objective is a simple sum over generators and offers little parallelism, which is why the objective and gradient callbacks do not accelerate there.

The multi-threaded \gls{cpu} rows use the KernelAbstractions.jl \texttt{CPU()} backend, so they run the same data-parallel kernels as the \gls{gpu} rows, one thread per data point, on host threads (\Cref{sec:gpu:cpugpu}).
On the large instances of the LV and COPS suites they stay an order of magnitude below the \gls{gpu} on the same callbacks.

\begin{figure}[t]
\centering
\vspace{0.1cm}
\pgfplotslegendfromname{cops-legend}\\[2pt]
\begin{tikzpicture}
\begin{axis}[
  width=.99\textwidth,
  height=3.3cm,
  ymode=log,
  xmin=-0.5, xmax=53.5,
  grid=both,
  grid style={line width=.1pt, draw=gray!20},
  xtick={0,1,2,3,4,5,6,7,8,9,10,11,12,13,14,15,16,17,18,19,20,21,22,23,24,25,26,27,28,29,30,31,32,33,34,35,36,37,38,39,40,41,42,43,44,45,46,47,48,49,50,51,52,53},
  ylabel style={font=\scriptsize},
  ylabel={\texttt{obj}},
  yticklabel style={font=\scriptsize},
  xticklabels={},
  legend to name=cops-legend,
  legend columns=4,
  legend style={font=\scriptsize, draw=black, fill=white, /tikz/every even column/.append style={column sep=2pt}},
]
\addplot[gray, dashed, thin, forget plot] coordinates {(-0.5,1) (53.5,1)};
\addplot[only marks, mark=o,mark size=1.8pt,black,thick] coordinates {
  (0, 1.0)
  (1, 1.0)
  (2, 1.0)
  (3, 1.0)
  (4, 1.0)
  (5, 1.0)
  (6, 1.0)
  (7, 1.0)
  (8, 1.0)
  (9, 1.0)
  (10, 1.0)
  (11, 1.0)
  (12, 1.0)
  (13, 1.0)
  (14, 1.0)
  (15, 1.0)
  (16, 1.0)
  (17, 1.0)
  (18, 1.0)
  (19, 1.0)
  (20, 1.0)
  (21, 1.0)
  (22, 1.0)
  (23, 1.0)
  (24, 1.0)
  (25, 1.0)
  (26, 1.0)
  (27, 1.0)
  (28, 1.0)
  (29, 1.0)
  (30, 1.0)
  (31, 1.0)
  (32, 1.0)
  (33, 1.0)
  (34, 1.0)
  (35, 1.0)
  (36, 1.0)
  (37, 1.0)
  (38, 1.0)
  (39, 1.0)
  (40, 1.0)
  (41, 1.0)
  (42, 1.0)
  (43, 1.0)
  (44, 1.0)
  (45, 1.0)
  (46, 1.0)
  (47, 1.0)
  (48, 1.0)
  (49, 1.0)
  (50, 1.0)
  (51, 1.0)
  (52, 1.0)
  (53, 1.0)
};
\addlegendentry{ExaModels (CPU, C3)}
\addplot[only marks, mark=square,mark size=1.8pt,black!55,thick] coordinates {
  (0, 0.15384615384615385)
  (1, 0.3333333333333333)
  (2, 0.14285714285714285)
  (3, 0.15384615384615385)
  (4, 0.16666666666666669)
  (5, 0.35294117647058826)
  (6, 0.125)
  (7, 0.017391304347826087)
  (8, 0.8401826484018264)
  (9, 0.5595238095238095)
  (10, 0.0975609756097561)
  (11, 0.9032602886157135)
  (12, 0.942857142857143)
  (13, 0.8962264150943396)
  (14, 0.8658536585365854)
  (15, 0.891089108910891)
  (16, 0.05405405405405406)
  (17, 0.01694915254237288)
  (18, 0.6128356543615331)
  (19, 0.30340136054421774)
  (20, 1.0151705799957265)
  (21, 0.15384615384615385)
  (22, 0.15384615384615385)
  (23, 0.16666666666666669)
  (24, 0.875)
  (25, 0.891089108910891)
  (26, 0.001606425702811245)
  (27, 0.2963235294117647)
  (28, 0.11764705882352942)
  (29, 0.9739130434782609)
  (30, 0.942857142857143)
  (31, 0.8554216867469879)
  (32, 0.8878504672897197)
  (33, 0.05555555555555556)
  (34, 0.15384615384615385)
  (35, 1.7619047619047619)
  (36, 0.15384615384615385)
  (37, 2.291513593327837)
  (38, 1.2682926829268293)
  (39, 1.7653205513301522)
  (40, 3.14088837413244)
  (41, 0.01680672268907563)
  (42, 0.16666666666666669)
  (43, 0.15384615384615385)
  (44, 0.892156862745098)
  (45, 0.16666666666666669)
  (46, 0.00010487621984153204)
  (47, 0.15384615384615385)
  (48, 0.93006993006993)
  (49, 0.11764705882352942)
  (50, 0.962979455343968)
  (51, 0.9024390243902439)
  (52, 0.9074074074074073)
  (53, 0.05555555555555556)
};
\addlegendentry{ExaModels (CPU-4T, C3)}
\addplot[only marks, mark=triangle,mark size=1.8pt,red!80!black] coordinates {
  (0, 0.0003979920061554399)
  (1, 0.0010012384278064)
  (2, 0.0003839183188898655)
  (3, 0.0003913441491424763)
  (4, 0.0003939115649825815)
  (5, 0.0011458195731542508)
  (6, 0.00041472150874225064)
  (7, 0.00039303955946506854)
  (8, 0.02376692273842427)
  (9, 0.005153807817107528)
  (10, 0.0006027333517506546)
  (11, 0.027334836564955563)
  (12, 0.024718741359442235)
  (13, 0.01817448595068713)
  (14, 0.013811721380377593)
  (15, 0.01631959121121253)
  (16, 0.00039228411817466484)
  (17, 0.00041709407426085963)
  (18, 0.5990502088777571)
  (19, 0.14409672614130575)
  (20, 0.45850360466773515)
  (21, 0.00038495323363882503)
  (22, 0.00038852372030072477)
  (23, 0.0004149994229121775)
  (24, 0.005211925907268748)
  (25, 0.017880178522616037)
  (26, 0.0003051748089523299)
  (27, 0.05598628061055971)
  (28, 0.0003962297173326121)
  (29, 0.014103559333402118)
  (30, 0.024932552346063647)
  (31, 0.013804190834904643)
  (32, 0.017958503078994804)
  (33, 0.00040404801483573806)
  (34, 0.0004096120321338919)
  (35, 1.0041894664941007)
  (36, 0.0003846206309887223)
  (37, 34.31538954986155)
  (38, 0.03146888624467537)
  (39, 8.555023371009494)
  (40, 28.338979444040916)
  (41, 0.00039093434054073825)
  (42, 0.0004096742350164787)
  (43, 0.0003990247993513773)
  (44, 0.017496886563627436)
  (45, 0.00039404108609268737)
  (46, 0.00020216581715731095)
  (47, 0.00041990181154606405)
  (48, 0.022597570549802153)
  (49, 0.0003902915004865627)
  (50, 1.0644386269815505)
  (51, 0.013989738069623687)
  (52, 0.019396221589117232)
  (53, 0.0003861475042959538)
};
\addlegendentry{ExaModels (AMDGPU, A1)}
\addplot[only marks, mark=diamond,mark size=1.8pt,green!55!black] coordinates {
  (0, 0.0004415032264768814)
  (1, 0.001342564607081981)
  (2, 0.0005341771822872637)
  (3, 0.00048213293101590193)
  (4, 0.0005261146240091456)
  (5, 0.001264809773579741)
  (6, 0.0004638506797496308)
  (7, 0.0005062182764712433)
  (8, 0.034186324046075266)
  (9, 0.00870070215480923)
  (10, 0.0007948731808554738)
  (11, 0.033276237938025824)
  (12, 0.033330449744423625)
  (13, 0.019748134455061184)
  (14, 0.015293822406596215)
  (15, 0.02023946908040596)
  (16, 0.00047764037515680546)
  (17, 0.0004719957009309971)
  (18, 0.8458545264603947)
  (19, 0.22284866188193986)
  (20, 0.6984816612774079)
  (21, 0.0005297364609048148)
  (22, 0.0004381000988366966)
  (23, 0.00047683061220487336)
  (24, 0.007112420537084835)
  (25, 0.021618992959827787)
  (26, 0.00039087287142067255)
  (27, 0.083274968876342)
  (28, 0.00040713579663749353)
  (29, 0.024537422859408853)
  (30, 0.029101067367626136)
  (31, 0.015433087116324977)
  (32, 0.02101593596170884)
  (33, 0.00044778044558153017)
  (34, 0.0004988562510717771)
  (35, 1.1862815037712267)
  (36, 0.0004846272011313138)
  (37, 60.77954985682688)
  (38, 0.0488860766859903)
  (39, 16.932826003930355)
  (40, 65.37445218235098)
  (41, 0.0004474718910966896)
  (42, 0.00040287112797008604)
  (43, 0.0004660871081662195)
  (44, 0.022988652848952042)
  (45, 0.00045989512606732573)
  (46, 0.00037933819649200404)
  (47, 0.0004586604907970426)
  (48, 0.030693936147284627)
  (49, 0.0005154473195906935)
  (50, 1.9479443724057341)
  (51, 0.01477644168260236)
  (52, 0.02337874850052464)
  (53, 0.00045633842689310776)
};
\addlegendentry{ExaModels (CUDA, N1)}
\addplot[only marks, mark=pentagon,mark size=1.8pt,violet] coordinates {
  (0, 6.542741689696298e-5)
  (1, 0.00015199585419202932)
  (2, 6.809762844173731e-5)
  (3, 6.589200431184949e-5)
  (4, 6.491385167692472e-5)
  (5, 0.00019220540927340573)
  (6, 7.162266561561348e-5)
  (7, 6.973793100441404e-5)
  (8, 0.003597845514039508)
  (9, 0.0010897130145088192)
  (10, 0.000129563809797709)
  (11, 0.005646493476878732)
  (12, 0.004460408710131398)
  (13, 0.0023066319927773017)
  (14, 0.002113547670424522)
  (15, 0.0030761940932898376)
  (16, 6.494437118324374e-5)
  (17, 6.984143905642247e-5)
  (18, 0.16362356189126243)
  (19, 0.031918847610380205)
  (20, 0.13639292107433887)
  (21, 6.0313528021716435e-5)
  (22, 6.486372475723276e-5)
  (23, 6.076005160172448e-5)
  (24, 0.000860044374204975)
  (25, 0.0031567020978921175)
  (26, 6.103851550873135e-5)
  (27, 0.013703268059027695)
  (28, 7.045874814937872e-5)
  (29, 0.00303419022207236)
  (30, 0.004428555202877641)
  (31, 0.0021331951169864093)
  (32, 0.0024106262980867774)
  (33, 7.160922033337844e-5)
  (34, 6.174374737658191e-5)
  (35, 0.14520074179164674)
  (36, 7.045377380665113e-5)
  (37, 9.514116503525656)
  (38, 0.0060299004494735065)
  (39, 2.013955257373103)
  (40, 9.213515296593943)
  (41, 6.705482617855089e-5)
  (42, 6.558485890469159e-5)
  (43, 6.517175854534561e-5)
  (44, 0.003241824251639319)
  (45, 6.003185436066419e-5)
  (46, 4.1338352436222035e-5)
  (47, 6.30454971648815e-5)
  (48, 0.00466793294942959)
  (49, 7.305819642061572e-5)
  (50, 0.2722031186704258)
  (51, 0.0022424817097999082)
  (52, 0.003198205951511659)
  (53, 7.012669377880068e-5)
};
\addlegendentry{ExaModels (Metal, M1\textsuperscript{*})}
\addplot[only marks, mark=star,mark size=1.8pt,blue!70!black] coordinates {
  (0, 0.00020763436490304254)
  (1, 0.0004134181994023762)
  (2, 0.00018300954936142168)
  (3, 0.00016045128348286957)
  (4, 0.00020842771287551861)
  (5, 0.0005471491579197548)
  (6, 0.0002096815234986456)
  (8, 0.012533804769875234)
  (9, 0.001133137983931824)
  (10, 0.00033084535749705195)
  (11, 0.01039158729067062)
  (12, 0.01058188979034093)
  (13, 0.00803465638576619)
  (14, 0.005044903556784764)
  (15, 0.007874603029449381)
  (16, 0.00017846689509322744)
  (18, 0.23145329947120408)
  (19, 0.031364781577843955)
  (20, 0.172097295429746)
  (21, 0.00020611823774422452)
  (22, 0.00017376060813182462)
  (23, 0.0001245420137134183)
  (24, 0.0019739021847971373)
  (25, 0.007816203218025179)
  (26, 9.236561474206678e-5)
  (27, 0.018542316109754722)
  (28, 0.00017662583257575866)
  (29, 0.008298850282839082)
  (30, 0.008597788744092314)
  (31, 0.0041460975434114535)
  (32, 0.00820646867032005)
  (33, 0.00020957808503832207)
  (34, 0.00020688910129089514)
  (35, 0.29381899521050586)
  (36, 0.00021228465393599843)
  (37, 11.641373800646788)
  (38, 0.0005603908765248933)
  (39, 0.1593939112311677)
  (40, 0.6143511987715272)
  (42, 6.7277810152915365e-6)
  (43, 7.243828313446071e-6)
  (44, 0.00033614311622058053)
  (45, 5.65589440360978e-6)
  (46, 7.021213069665226e-5)
  (47, 0.00016332016169839247)
  (48, 0.0003515225797082885)
  (49, 0.00020517687139859835)
  (50, 0.43092057791733507)
  (51, 0.00021161895273274887)
  (52, 0.00033317281672523097)
  (53, 7.157478399550557e-6)
};
\addlegendentry{ExaModels (oneAPI, I1)}
\end{axis}
\end{tikzpicture}%
\vspace{-0.2cm}
\begin{tikzpicture}
\begin{axis}[
  width=.99\textwidth,
  height=3.3cm,
  ymode=log,
  xmin=-0.5, xmax=53.5,
  grid=both,
  grid style={line width=.1pt, draw=gray!20},
  xtick={0,1,2,3,4,5,6,7,8,9,10,11,12,13,14,15,16,17,18,19,20,21,22,23,24,25,26,27,28,29,30,31,32,33,34,35,36,37,38,39,40,41,42,43,44,45,46,47,48,49,50,51,52,53},
  ylabel style={font=\scriptsize},
  ylabel={\texttt{cons}},
  yticklabel style={font=\scriptsize},
  xticklabels={},
]
\addplot[gray, dashed, thin, forget plot] coordinates {(-0.5,1) (53.5,1)};
\addplot[only marks, mark=o,mark size=1.8pt,black,thick] coordinates {
  (0, 1.0)
  (1, 1.0)
  (2, 1.0)
  (3, 1.0)
  (4, 1.0)
  (5, 1.0)
  (6, 1.0)
  (7, 1.0)
  (8, 1.0)
  (9, 1.0)
  (10, 1.0)
  (11, 1.0)
  (12, 1.0)
  (13, 1.0)
  (14, 1.0)
  (15, 1.0)
  (16, 1.0)
  (17, 1.0)
  (18, 1.0)
  (19, 1.0)
  (20, 1.0)
  (21, 1.0)
  (22, 1.0)
  (23, 1.0)
  (24, 1.0)
  (25, 1.0)
  (26, 1.0)
  (27, 1.0)
  (28, 1.0)
  (29, 1.0)
  (30, 1.0)
  (31, 1.0)
  (32, 1.0)
  (33, 1.0)
  (34, 1.0)
  (35, 1.0)
  (36, 1.0)
  (37, 1.0)
  (38, 1.0)
  (39, 1.0)
  (40, 1.0)
  (41, 1.0)
  (42, 1.0)
  (43, 1.0)
  (44, 1.0)
  (45, 1.0)
  (46, 1.0)
  (47, 1.0)
  (48, 1.0)
  (49, 1.0)
  (50, 1.0)
  (51, 1.0)
  (52, 1.0)
  (53, 1.0)
};
\addplot[only marks, mark=square,mark size=1.8pt,black!55,thick] coordinates {
  (0, 0.1134020618556701)
  (1, 0.24561403508771928)
  (2, 0.03225806451612903)
  (3, 0.002163372825492167)
  (4, 0.11764705882352942)
  (5, 0.3846153846153846)
  (6, 0.10526315789473685)
  (7, 0.2057335581787521)
  (8, 0.10344827586206896)
  (9, 0.47368421052631576)
  (10, 0.3003533568904594)
  (11, 0.21794871794871792)
  (12, 0.2682926829268293)
  (13, 0.39622641509433965)
  (14, 0.3993710691823899)
  (15, 0.3424124513618677)
  (16, 0.176759410801964)
  (17, 0.7448504559658184)
  (18, 0.1935483870967742)
  (19, 0.15196599362380447)
  (20, 1.0186335403726707)
  (21, 0.05293721691709164)
  (22, 1.010989010989011)
  (23, 0.8944444444444445)
  (24, 0.5103123158515027)
  (25, 1.0485232067510548)
  (26, 0.1817116060961313)
  (27, 0.5714285714285714)
  (28, 0.2883509833585477)
  (29, 0.8316766070245196)
  (30, 0.6746482350037029)
  (31, 0.973055824459503)
  (32, 1.225801382778127)
  (33, 0.9808458826522577)
  (34, 1.5621705012730467)
  (35, 1.25)
  (36, 1.856995759237733)
  (37, 0.608)
  (38, 3.886602602629642)
  (39, 2.3330339097237816)
  (40, 1.9257743190661478)
  (41, 3.5417222109127207)
  (42, 3.0271214642262896)
  (43, 1.8018393893321705)
  (44, 3.57119608889129)
  (45, 2.998128016901333)
  (46, 3.7052845239954015)
  (47, 2.2444554182137417)
  (48, 3.8950192239949177)
  (49, 3.4222931945144297)
  (50, 4.06251537014503)
  (51, 3.900409671048847)
  (52, 4.012062480710771)
  (53, 1.8283273675917473)
};
\addplot[only marks, mark=triangle,mark size=1.8pt,red!80!black] coordinates {
  (0, 0.0012516374931529314)
  (1, 0.00290309671564063)
  (2, 0.0008232221597509414)
  (3, 0.00036407399993730307)
  (4, 0.0017877460810832124)
  (5, 0.004110749004366314)
  (6, 0.001391990531028161)
  (7, 0.007302621479914898)
  (8, 0.0012250537181461455)
  (9, 0.007469282265463189)
  (10, 0.00852876248388829)
  (11, 0.0024406995925981365)
  (12, 0.004420610610705247)
  (13, 0.008047039567572327)
  (14, 0.00936986388282905)
  (15, 0.006683094644886688)
  (16, 0.009159677062693758)
  (17, 0.048861203177083375)
  (18, 0.00252735302047403)
  (19, 0.12821753268240252)
  (20, 0.012028441746315665)
  (21, 0.009194707007585607)
  (22, 0.04363220476292064)
  (23, 0.07487039478348848)
  (24, 0.17788281125584096)
  (25, 0.19425382129243124)
  (26, 0.030620435209233385)
  (27, 0.006494538661258016)
  (28, 0.09766093954086898)
  (29, 0.5256009138185437)
  (30, 0.25357744258575954)
  (31, 0.4079926389820599)
  (32, 0.7082467952697281)
  (33, 0.678554340848483)
  (34, 0.7812874037873229)
  (35, 0.016931987668128604)
  (36, 2.1585443091112198)
  (37, 0.01663220619439733)
  (38, 11.832023407576527)
  (39, 8.759969376786234)
  (40, 1.07776078233306)
  (41, 5.506804794398026)
  (42, 3.613372681737951)
  (43, 1.0169032196978378)
  (44, 16.163628250780384)
  (45, 8.128040428310463)
  (46, 2.4619977521661824)
  (47, 7.123647939542176)
  (48, 20.18666534608434)
  (49, 9.396320327293648)
  (50, 41.16932079984097)
  (51, 45.256433581000536)
  (52, 54.109532905770955)
  (53, 81.52752863241375)
};
\addplot[only marks, mark=diamond,mark size=1.8pt,green!55!black] coordinates {
  (0, 0.0028387050036289235)
  (1, 0.006820530992855781)
  (2, 0.001976535970494974)
  (3, 0.0005986453663647431)
  (4, 0.004512512575491199)
  (5, 0.008676590383877173)
  (6, 0.0024936484414137066)
  (7, 0.017171155303309343)
  (8, 0.0027586435451754866)
  (9, 0.020480437848094907)
  (10, 0.01607391595927017)
  (11, 0.005210042543337923)
  (12, 0.006875974142892986)
  (13, 0.010698985823391625)
  (14, 0.014140436774398378)
  (15, 0.015764996556432098)
  (16, 0.016660399528289786)
  (17, 0.1356373325472066)
  (18, 0.004839056040660394)
  (19, 0.2733962261914426)
  (20, 0.040522906297996715)
  (21, 0.01569141632714648)
  (22, 0.08696370282031843)
  (23, 0.14840167193840179)
  (24, 0.4636176801747995)
  (25, 0.2988752096084158)
  (26, 0.05922243759193973)
  (27, 0.01719843104733006)
  (28, 0.13391007216410977)
  (29, 0.9635567519321876)
  (30, 0.6073566922717182)
  (31, 0.6923318924434589)
  (32, 1.0505939021065382)
  (33, 1.1690003589579054)
  (34, 1.0782091612553903)
  (35, 0.03373826903516512)
  (36, 2.257536196124951)
  (37, 0.03240290426037094)
  (38, 26.091770512448562)
  (39, 18.807455268751365)
  (40, 3.904789214456584)
  (41, 6.574004437082183)
  (42, 5.47437512549975)
  (43, 1.57247877715667)
  (44, 34.331322555893315)
  (45, 15.495866412101213)
  (46, 5.044331687221267)
  (47, 10.298261977106742)
  (48, 55.46955112448738)
  (49, 21.262986146094992)
  (50, 69.66241764603878)
  (51, 69.26037075758008)
  (52, 50.12519892908915)
  (53, 111.28861251928167)
};
\addplot[only marks, mark=pentagon,mark size=1.8pt,violet] coordinates {
  (0, 0.0023171380099637033)
  (1, 0.0040161308592685015)
  (2, 0.0015656251060898296)
  (3, 0.0006805477679857927)
  (4, 0.0033090063326497074)
  (5, 0.008796832994082937)
  (6, 0.002582967239043911)
  (7, 0.010149639736825268)
  (8, 0.001994579288726479)
  (9, 0.014906959590184664)
  (10, 0.015099908340160778)
  (11, 0.004480285840668276)
  (12, 0.009932632711164258)
  (13, 0.010921906818359441)
  (14, 0.017273642648313058)
  (15, 0.01273738649029852)
  (16, 0.018158320776868436)
  (17, 0.10372419258921611)
  (18, 0.004892139443763745)
  (19, 0.2492154460009858)
  (20, 0.043179112083534504)
  (21, 0.018563423744464597)
  (22, 0.07582545552718224)
  (23, 0.12065332812453101)
  (24, 0.3201278767939846)
  (25, 0.31130983601355233)
  (26, 0.05105102901694865)
  (27, 0.013965415994480696)
  (28, 0.17543082404718408)
  (29, 0.6881745842314219)
  (30, 0.5896647846612825)
  (31, 0.8276321993905786)
  (32, 1.3942608091439304)
  (33, 1.1097209410298396)
  (34, 0.7824189237403653)
  (35, 0.03241793299217417)
  (36, 1.5276819397940309)
  (37, 0.030306250902523556)
  (38, 21.779873241099388)
  (39, 12.446082677202513)
  (40, 3.14879357896402)
  (41, 6.852646957935981)
  (42, 6.409921817823322)
  (43, 1.7678950066530545)
  (44, 20.64741802351464)
  (45, 13.48582447468615)
  (46, 4.660247376321669)
  (47, 5.265600382729653)
  (48, 30.28267705503051)
  (49, 13.182098336378342)
  (50, 35.108695452206256)
  (51, 44.1856374859386)
  (52, 47.60842805983491)
  (53, 17.487858490376095)
};
\addplot[only marks, mark=star,mark size=1.8pt,blue!70!black] coordinates {
  (0, 0.0011050887580215769)
  (1, 0.002339208696525629)
  (2, 0.0009188202775194353)
  (3, 0.000264878740711063)
  (4, 0.0014827291115363079)
  (5, 0.005700428916232858)
  (6, 0.001702425897684099)
  (8, 0.0015518255095946678)
  (9, 0.006069756310209706)
  (10, 0.00960685717120948)
  (11, 0.002628450651661565)
  (12, 0.0032361705559174024)
  (13, 0.009740588733613907)
  (14, 0.004107457757600077)
  (15, 0.007415720142546523)
  (16, 0.009342046636010619)
  (18, 0.003087823362587869)
  (19, 0.06600557513551707)
  (20, 0.021544215739793367)
  (21, 0.007385317017288276)
  (22, 0.03817912217178309)
  (23, 0.059433103177902985)
  (24, 0.1233442095627122)
  (25, 0.19621785524031873)
  (26, 0.032903713264178326)
  (27, 0.009350040359449212)
  (28, 0.10288604322743179)
  (29, 0.5931773308227027)
  (30, 0.1849696095236113)
  (31, 0.19021463979172842)
  (32, 0.8928089520686684)
  (33, 0.6832399916894522)
  (34, 0.376382542564386)
  (35, 0.02330585487588665)
  (36, 1.1227865712302103)
  (37, 0.01870980328010621)
  (38, 9.862525985497056)
  (39, 7.124855183152121)
  (40, 1.8045010601383569)
  (42, 3.0025307841244655)
  (43, 0.8840820309796071)
  (44, 16.16133721106715)
  (45, 6.182505945375832)
  (46, 2.9342613962812814)
  (47, 5.36488963042434)
  (48, 17.6861459841255)
  (49, 11.042139481765703)
  (50, 59.80326044112891)
  (51, 18.698724057098026)
  (52, 74.36992519834577)
  (53, 69.81214520810668)
};
\end{axis}
\end{tikzpicture}%
\vspace{-0.2cm}
\begin{tikzpicture}
\begin{axis}[
  width=.99\textwidth,
  height=3.3cm,
  ymode=log,
  xmin=-0.5, xmax=53.5,
  grid=both,
  grid style={line width=.1pt, draw=gray!20},
  xtick={0,1,2,3,4,5,6,7,8,9,10,11,12,13,14,15,16,17,18,19,20,21,22,23,24,25,26,27,28,29,30,31,32,33,34,35,36,37,38,39,40,41,42,43,44,45,46,47,48,49,50,51,52,53},
  ylabel style={font=\scriptsize},
  ylabel={\texttt{grad}},
  yticklabel style={font=\scriptsize},
  xticklabels={},
]
\addplot[gray, dashed, thin, forget plot] coordinates {(-0.5,1) (53.5,1)};
\addplot[only marks, mark=o,mark size=1.8pt,black,thick] coordinates {
  (0, 1.0)
  (1, 1.0)
  (2, 1.0)
  (3, 1.0)
  (4, 1.0)
  (5, 1.0)
  (6, 1.0)
  (7, 1.0)
  (8, 1.0)
  (9, 1.0)
  (10, 1.0)
  (11, 1.0)
  (12, 1.0)
  (13, 1.0)
  (14, 1.0)
  (15, 1.0)
  (16, 1.0)
  (17, 1.0)
  (18, 1.0)
  (19, 1.0)
  (20, 1.0)
  (21, 1.0)
  (22, 1.0)
  (23, 1.0)
  (24, 1.0)
  (25, 1.0)
  (26, 1.0)
  (27, 1.0)
  (28, 1.0)
  (29, 1.0)
  (30, 1.0)
  (31, 1.0)
  (32, 1.0)
  (33, 1.0)
  (34, 1.0)
  (35, 1.0)
  (36, 1.0)
  (37, 1.0)
  (38, 1.0)
  (39, 1.0)
  (40, 1.0)
  (41, 1.0)
  (42, 1.0)
  (43, 1.0)
  (44, 1.0)
  (45, 1.0)
  (46, 1.0)
  (47, 1.0)
  (48, 1.0)
  (49, 1.0)
  (50, 1.0)
  (51, 1.0)
  (52, 1.0)
  (53, 1.0)
};
\addplot[only marks, mark=square,mark size=1.8pt,black!55,thick] coordinates {
  (0, 0.09090909090909091)
  (1, 0.23333333333333336)
  (2, 0.0625)
  (3, 0.08333333333333334)
  (4, 0.09090909090909091)
  (5, 0.3055555555555556)
  (6, 0.08333333333333334)
  (7, 0.5984848484848484)
  (8, 0.848360655737705)
  (9, 0.6511627906976745)
  (10, 0.07692307692307691)
  (11, 1.028720626631854)
  (12, 1.9324324324324322)
  (13, 0.8587570621468927)
  (14, 1.2222222222222223)
  (15, 0.9745222929936306)
  (16, 0.5)
  (17, 0.5985401459854014)
  (18, 0.692018135807676)
  (19, 0.4263693669393912)
  (20, 1.1512314884362282)
  (21, 0.3333333333333333)
  (22, 0.25)
  (23, 0.42857142857142855)
  (24, 0.5862068965517241)
  (25, 0.9050279329608938)
  (26, 0.013399503722084368)
  (27, 0.2756482525366404)
  (28, 0.576271186440678)
  (29, 0.4105263157894737)
  (30, 1.6334519572953736)
  (31, 1.0625)
  (32, 0.8392857142857142)
  (33, 0.9803921568627451)
  (34, 0.71875)
  (35, 0.4147348269570159)
  (36, 0.8813559322033899)
  (37, 1.3349518417490138)
  (38, 0.7369402985074627)
  (39, 0.6735376513632257)
  (40, 2.4620618596118886)
  (41, 0.8534031413612566)
  (42, 0.8921775898520086)
  (43, 0.9397363465160075)
  (44, 0.8177777777777777)
  (45, 1.0124300807955253)
  (46, 0.05397025784084434)
  (47, 0.9467831612390786)
  (48, 0.9445154691894656)
  (49, 1.0032649757150567)
  (50, 0.5720216940015433)
  (51, 0.9231473230413526)
  (52, 0.9101821811829298)
  (53, 0.9694912374199627)
};
\addplot[only marks, mark=triangle,mark size=1.8pt,red!80!black] coordinates {
  (0, 0.0004310501733565259)
  (1, 0.0015098649876139735)
  (2, 0.00042259444561680364)
  (3, 0.000433436435219376)
  (4, 0.000424621903175667)
  (5, 0.0023732890767564387)
  (6, 0.0004564512568049865)
  (7, 0.016768852658544515)
  (8, 0.02707672304846167)
  (9, 0.012484083657191086)
  (10, 0.0008524868923235161)
  (11, 0.06799300048990166)
  (12, 0.08911590497485052)
  (13, 0.030507680641777557)
  (14, 0.0332766981936948)
  (15, 0.028077825021644427)
  (16, 0.0018844910980786408)
  (17, 0.017148088803678223)
  (18, 0.8593522305776152)
  (19, 0.5206431470150593)
  (20, 1.7331279226291703)
  (21, 0.0023609000305513344)
  (22, 0.0015819042233556925)
  (23, 0.0033027115823552103)
  (24, 0.010626178686732954)
  (25, 0.03496268321673476)
  (26, 0.006165187893451588)
  (27, 0.20512098421561725)
  (28, 0.007661001236080016)
  (29, 0.027582295234459396)
  (30, 0.10226233857664278)
  (31, 0.03960531964675147)
  (32, 0.03679397341019318)
  (33, 0.08944941310093278)
  (34, 0.014744455822843098)
  (35, 0.9907458293285327)
  (36, 0.04000652306358552)
  (37, 65.80293020796698)
  (38, 0.08598699110865995)
  (39, 37.51229509119868)
  (40, 117.54621729544829)
  (41, 0.0629526567576623)
  (42, 0.09456673162818491)
  (43, 0.21893467650005988)
  (44, 0.19499543247573856)
  (45, 0.37237330153285153)
  (46, 0.5995353603088909)
  (47, 0.2679761805325313)
  (48, 0.7468640162566781)
  (49, 0.8395054516532556)
  (50, 2.4743402190376416)
  (51, 0.8464294477879488)
  (52, 0.7147208116415864)
  (53, 11.569962982036126)
};
\addplot[only marks, mark=diamond,mark size=1.8pt,green!55!black] coordinates {
  (0, 0.00098506412769934)
  (1, 0.003747746647251126)
  (2, 0.0012500177424393318)
  (3, 0.0009867110398514975)
  (4, 0.001110155131856955)
  (5, 0.004870505587658755)
  (6, 0.0009112249772780357)
  (7, 0.042444907933287486)
  (8, 0.062129895476715796)
  (9, 0.032215885459824935)
  (10, 0.0017598129192160342)
  (11, 0.13423959597561277)
  (12, 0.21484245995538118)
  (13, 0.05102802965304603)
  (14, 0.07449177370654868)
  (15, 0.06971301484667597)
  (16, 0.0036025172012541015)
  (17, 0.03423077277426433)
  (18, 1.7177586289882936)
  (19, 1.050458464696327)
  (20, 3.2379376629697214)
  (21, 0.006046326436470157)
  (22, 0.003296066002083076)
  (23, 0.006771262417417514)
  (24, 0.03001052257187613)
  (25, 0.08598167963937373)
  (26, 0.013682552825029842)
  (27, 0.32866034444471115)
  (28, 0.011050994782922872)
  (29, 0.06667305730770237)
  (30, 0.21555866974433727)
  (31, 0.09405657224171458)
  (32, 0.06578680132562365)
  (33, 0.19228060662608584)
  (34, 0.036841188055273116)
  (35, 1.3450304982836052)
  (36, 0.12059509318510404)
  (37, 125.18415680348882)
  (38, 0.1747264645033699)
  (39, 77.29997075824764)
  (40, 331.37749892918185)
  (41, 0.12135417485815103)
  (42, 0.1706988250524812)
  (43, 0.47568138045138847)
  (44, 0.5367495935842928)
  (45, 0.7559831191057855)
  (46, 1.3060410042182247)
  (47, 0.4672900988249171)
  (48, 1.6889184955092098)
  (49, 2.262332084392891)
  (50, 5.731196994957983)
  (51, 1.3351859802253419)
  (52, 0.8439170420962901)
  (53, 18.188934193896465)
};
\addplot[only marks, mark=pentagon,mark size=1.8pt,violet] coordinates {
  (0, 0.000936992997006385)
  (1, 0.0030528256536813915)
  (2, 0.0006882020026121028)
  (3, 0.0009235213759108714)
  (4, 0.0008773947503164386)
  (5, 0.0046747795870520845)
  (6, 0.0008909958879560653)
  (7, 0.03145590851439136)
  (8, 0.043998404421758154)
  (9, 0.023003233276065276)
  (10, 0.001679226852852768)
  (11, 0.1293942254523681)
  (12, 0.19647500516993802)
  (13, 0.019691126769833353)
  (14, 0.06252356710713645)
  (15, 0.06603301480625202)
  (16, 0.0044661938660293)
  (17, 0.034575942149206)
  (18, 1.4659790324944126)
  (19, 0.8416177381488692)
  (20, 2.8537053549311366)
  (21, 0.004764113317317615)
  (22, 0.0032787586040633424)
  (23, 0.006842666964246043)
  (24, 0.023189642702155358)
  (25, 0.06954615606978898)
  (26, 0.01228805108163578)
  (27, 0.29020611866722223)
  (28, 0.015743363578057924)
  (29, 0.05178212973629858)
  (30, 0.21330892722257708)
  (31, 0.07945704581295839)
  (32, 0.0635551469527253)
  (33, 0.2208456178708274)
  (34, 0.03302315321870985)
  (35, 0.6661532282506074)
  (36, 0.0951652663519549)
  (37, 25.135628978855188)
  (38, 0.17934094365811296)
  (39, 13.312547650807698)
  (40, 70.95981876974318)
  (41, 0.140676810353655)
  (42, 0.18482085307566148)
  (43, 0.45091413488154636)
  (44, 0.40254189940242235)
  (45, 0.7424241197802173)
  (46, 0.8894109882849929)
  (47, 0.5586941943584403)
  (48, 1.5804112459473472)
  (49, 1.6986849787094167)
  (50, 3.771703781898837)
  (51, 1.37004261082076)
  (52, 1.1446097831156123)
  (53, 2.338791322318252)
};
\addplot[only marks, mark=star,mark size=1.8pt,blue!70!black] coordinates {
  (0, 0.0005458404640908112)
  (1, 0.00083250703640299)
  (2, 0.0005171277721342044)
  (3, 0.0005350161202497096)
  (4, 0.0004771905403121597)
  (5, 0.0032275613931617357)
  (6, 0.0005775830238414865)
  (8, 0.031153407177332948)
  (9, 0.003950925864839198)
  (10, 0.0010780204005227407)
  (11, 0.07253202081155752)
  (12, 0.10942897575460749)
  (13, 0.03754809979516178)
  (14, 0.02909886102976531)
  (15, 0.03499377099440427)
  (16, 0.002595967061331539)
  (18, 0.9615261649896122)
  (19, 0.1418653241039662)
  (20, 1.8265066387098081)
  (21, 0.002899016721523179)
  (22, 0.0018586639585189326)
  (23, 0.0035519929203382066)
  (24, 0.008047004782852476)
  (25, 0.04414050154879211)
  (26, 0.007224272708653077)
  (27, 0.1855992213000589)
  (28, 0.008274607450493052)
  (29, 0.03338298181638396)
  (30, 0.11392449994050734)
  (31, 0.034376716560898274)
  (32, 0.04759934308247887)
  (33, 0.13019212373584424)
  (34, 0.01874038643638295)
  (35, 0.8399818071888077)
  (36, 0.05771370338301858)
  (37, 68.17131022199482)
  (38, 0.06652644083117132)
  (39, 11.33128156796756)
  (40, 159.18727792310654)
  (42, 0.11236946958044891)
  (43, 0.26316511782963686)
  (44, 0.24150989296974573)
  (45, 0.3871283148135251)
  (46, 0.7961079902008135)
  (47, 0.3169034289759225)
  (48, 0.9782114365675524)
  (49, 1.0042655287673685)
  (50, 3.2725073313833137)
  (51, 0.5191312937921206)
  (52, 0.8835339826305566)
  (53, 9.339689996342356)
};
\end{axis}
\end{tikzpicture}%
\vspace{-0.2cm}
\begin{tikzpicture}
\begin{axis}[
  width=.99\textwidth,
  height=3.3cm,
  ymode=log,
  xmin=-0.5, xmax=53.5,
  grid=both,
  grid style={line width=.1pt, draw=gray!20},
  xtick={0,1,2,3,4,5,6,7,8,9,10,11,12,13,14,15,16,17,18,19,20,21,22,23,24,25,26,27,28,29,30,31,32,33,34,35,36,37,38,39,40,41,42,43,44,45,46,47,48,49,50,51,52,53},
  ylabel style={font=\scriptsize},
  ylabel={\texttt{jac}},
  yticklabel style={font=\scriptsize},
  xticklabels={},
]
\addplot[gray, dashed, thin, forget plot] coordinates {(-0.5,1) (53.5,1)};
\addplot[only marks, mark=o,mark size=1.8pt,black,thick] coordinates {
  (0, 1.0)
  (1, 1.0)
  (2, 1.0)
  (3, 1.0)
  (4, 1.0)
  (5, 1.0)
  (6, 1.0)
  (7, 1.0)
  (8, 1.0)
  (9, 1.0)
  (10, 1.0)
  (11, 1.0)
  (12, 1.0)
  (13, 1.0)
  (14, 1.0)
  (15, 1.0)
  (16, 1.0)
  (17, 1.0)
  (18, 1.0)
  (19, 1.0)
  (20, 1.0)
  (21, 1.0)
  (22, 1.0)
  (23, 1.0)
  (24, 1.0)
  (25, 1.0)
  (26, 1.0)
  (27, 1.0)
  (28, 1.0)
  (29, 1.0)
  (30, 1.0)
  (31, 1.0)
  (32, 1.0)
  (33, 1.0)
  (34, 1.0)
  (35, 1.0)
  (36, 1.0)
  (37, 1.0)
  (38, 1.0)
  (39, 1.0)
  (40, 1.0)
  (41, 1.0)
  (42, 1.0)
  (43, 1.0)
  (44, 1.0)
  (45, 1.0)
  (46, 1.0)
  (47, 1.0)
  (48, 1.0)
  (49, 1.0)
  (50, 1.0)
  (51, 1.0)
  (52, 1.0)
  (53, 1.0)
};
\addplot[only marks, mark=square,mark size=1.8pt,black!55,thick] coordinates {
  (0, 0.09836065573770492)
  (1, 0.25)
  (2, 0.056105610561056105)
  (3, 0.006172839506172839)
  (4, 0.1369047619047619)
  (5, 0.125)
  (6, 0.08898305084745764)
  (7, 0.2332814930015552)
  (8, 0.04687499999999999)
  (9, 0.18749999999999997)
  (10, 0.32298136645962733)
  (11, 0.05319148936170213)
  (12, 0.22631578947368422)
  (13, 0.401673640167364)
  (14, 0.41477272727272724)
  (15, 0.3767123287671233)
  (16, 0.29192902117916425)
  (17, 0.8971062052505967)
  (18, 0.0909090909090909)
  (19, 0.06043956043956045)
  (20, 0.34057971014492755)
  (21, 0.2643823264201984)
  (22, 1.0934911242603549)
  (23, 1.3258785942492013)
  (24, 0.7100271002710028)
  (25, 1.2268602540834845)
  (26, 0.4506235254465791)
  (27, 0.38095238095238093)
  (28, 0.37241252976735667)
  (29, 0.99880161559971)
  (30, 0.8000971581248482)
  (31, 1.0424377652360326)
  (32, 1.3303421200224341)
  (33, 1.0504334568692169)
  (34, 3.385525566354934)
  (35, 0.8461538461538461)
  (36, 4.523923650192486)
  (37, 0.4158415841584158)
  (38, 3.33585113456903)
  (39, 0.4556561736380624)
  (40, 0.3539091816710848)
  (41, 3.5222044454193866)
  (42, 2.345704597374976)
  (43, 2.2129939226809143)
  (44, 3.16220451762261)
  (45, 2.586932967668623)
  (46, 2.2177565137060125)
  (47, 5.289969830164633)
  (48, 3.0082721421232437)
  (49, 1.8534738568165063)
  (50, 3.7229477743602213)
  (51, 3.1686227591131004)
  (52, 3.553888991657936)
  (53, 4.2951648684441395)
};
\addplot[only marks, mark=triangle,mark size=1.8pt,red!80!black] coordinates {
  (0, 0.0011762813038156048)
  (1, 0.0030687843923512302)
  (2, 0.0013614346100866573)
  (3, 0.0006531231564334875)
  (4, 0.002256696842317434)
  (5, 0.0008562350339481137)
  (6, 0.0017006941195162595)
  (7, 0.008146155236928132)
  (8, 0.0005042702648183363)
  (9, 0.0012236963112727073)
  (10, 0.00899702962461374)
  (11, 0.0005851969460659368)
  (12, 0.0033782351563896337)
  (13, 0.008016979983647816)
  (14, 0.009007981428780232)
  (15, 0.007052276218797857)
  (16, 0.011403997842727344)
  (17, 0.16386759996647818)
  (18, 0.0009837958155002509)
  (19, 0.01985788945124677)
  (20, 0.005206659839930615)
  (21, 0.04231945798366513)
  (22, 0.0917196051663569)
  (23, 0.16727672555656556)
  (24, 0.286395150363794)
  (25, 0.23374299207346913)
  (26, 0.11293020794390748)
  (27, 0.003321662326741833)
  (28, 0.17369860937072762)
  (29, 0.6165162374961767)
  (30, 0.2407361863422148)
  (31, 0.4433595571591191)
  (32, 0.8014150096949594)
  (33, 0.9438188571843871)
  (34, 3.3762471963122715)
  (35, 0.01456389463468858)
  (36, 7.55313376444855)
  (37, 0.006900333469458292)
  (38, 20.998788129981918)
  (39, 1.3337629929500192)
  (40, 0.21430868903489012)
  (41, 14.042015815249131)
  (42, 8.160351639137206)
  (43, 4.23191326395185)
  (44, 20.05417193808566)
  (45, 17.389225798537634)
  (46, 9.340720356026575)
  (47, 44.687266898050595)
  (48, 19.137948478457645)
  (49, 17.91487436012725)
  (50, 50.2452738258275)
  (51, 46.09922181869244)
  (52, 64.27522384393775)
  (53, 42.404467146221435)
};
\addplot[only marks, mark=diamond,mark size=1.8pt,green!55!black] coordinates {
  (0, 0.0025627202995789273)
  (1, 0.006905498543194221)
  (2, 0.0033790080342316247)
  (3, 0.0010432007620342991)
  (4, 0.00582457979816266)
  (5, 0.0016343235400036232)
  (6, 0.002788014241277644)
  (7, 0.01815898681066301)
  (8, 0.0011950923246864592)
  (9, 0.0033772784048870922)
  (10, 0.016009568697539566)
  (11, 0.0011071106455213878)
  (12, 0.005215296992946595)
  (13, 0.012266529812400917)
  (14, 0.012614080709159887)
  (15, 0.016609687987675537)
  (16, 0.017799840645356047)
  (17, 0.2757975190017896)
  (18, 0.0020523161800769804)
  (19, 0.043094211389643367)
  (20, 0.010285459581470207)
  (21, 0.07095534360786779)
  (22, 0.1933496213370079)
  (23, 0.3299691022100882)
  (24, 0.754070076793033)
  (25, 0.4111377121062726)
  (26, 0.21962743748509506)
  (27, 0.008113555044608934)
  (28, 0.26387465322126535)
  (29, 1.1521977737323705)
  (30, 0.6229873056066838)
  (31, 0.6396724373121974)
  (32, 1.0190776797933252)
  (33, 1.523521495224574)
  (34, 5.407242809591829)
  (35, 0.0259201149910556)
  (36, 24.946807242440865)
  (37, 0.01353897362311048)
  (38, 40.650316951155204)
  (39, 3.0636435907426005)
  (40, 0.8281336284478639)
  (41, 13.387420059859974)
  (42, 13.826515951915765)
  (43, 6.228633587830157)
  (44, 33.977945326801105)
  (45, 33.26804854901772)
  (46, 18.06121393966933)
  (47, 36.58076257932695)
  (48, 55.533083361002824)
  (49, 42.19866261018716)
  (50, 82.550498111328)
  (51, 61.23488184571421)
  (52, 47.858407114288646)
  (53, 96.92442994793264)
};
\addplot[only marks, mark=pentagon,mark size=1.8pt,violet] coordinates {
  (0, 0.002017568546143949)
  (1, 0.005560858774790422)
  (2, 0.0026310768981289375)
  (3, 0.0009964453234767225)
  (4, 0.0038420330906714342)
  (5, 0.0017493490957708818)
  (6, 0.0032578556129457662)
  (7, 0.010248190058496622)
  (8, 0.0009673278840707934)
  (9, 0.002620977080060313)
  (10, 0.016178955874103125)
  (11, 0.001029314357849887)
  (12, 0.00790336591129274)
  (13, 0.014266447604442695)
  (14, 0.016594999475755788)
  (15, 0.014465421100284924)
  (16, 0.021469268150720155)
  (17, 0.12999108500098117)
  (18, 0.001783125015457739)
  (19, 0.03869811621528026)
  (20, 0.009889015492106112)
  (21, 0.048130942023403796)
  (22, 0.16493893779499058)
  (23, 0.2842856325750407)
  (24, 0.532615094934632)
  (25, 0.36344457955891574)
  (26, 0.19170019190806104)
  (27, 0.007503831002339714)
  (28, 0.32055588468906876)
  (29, 0.7578118317646613)
  (30, 0.6255100624048023)
  (31, 0.9032486420974533)
  (32, 1.6384849740103633)
  (33, 1.6405855582592188)
  (34, 6.864712983985454)
  (35, 0.020119503309071928)
  (36, 18.27515159068589)
  (37, 0.013198538811202448)
  (38, 30.782352274077837)
  (39, 2.7372737760983052)
  (40, 0.6454276824065073)
  (41, 12.836730136470667)
  (42, 12.616866494744466)
  (43, 5.659353634150901)
  (44, 18.00985150579148)
  (45, 18.782674895298715)
  (46, 13.700138130514993)
  (47, 35.59096126949642)
  (48, 21.878098283694438)
  (49, 13.206774128875391)
  (50, 36.39525765845516)
  (51, 32.40262395321114)
  (52, 38.47210653935235)
  (53, 9.348851890416327)
};
\addplot[only marks, mark=star,mark size=1.8pt,blue!70!black] coordinates {
  (0, 0.0007634521286179708)
  (1, 0.0017970776757060044)
  (2, 0.0010138062127851679)
  (3, 0.0004206055339314874)
  (4, 0.001081955616910025)
  (5, 0.0009434238393628455)
  (6, 0.0013989532677249389)
  (8, 0.00030930889843624874)
  (9, 0.0005091354383849172)
  (10, 0.006164231256740305)
  (11, 0.0004110478825748794)
  (12, 0.0017024317458552804)
  (13, 0.006025202264988763)
  (14, 0.0028105535567315627)
  (15, 0.005256846546363878)
  (16, 0.007150475634251302)
  (18, 0.0008608924125534183)
  (19, 0.007091036161398049)
  (20, 0.003881307149925192)
  (21, 0.02651361185111262)
  (22, 0.05796287532536306)
  (23, 0.07861329371451262)
  (24, 0.15991989579284754)
  (25, 0.1310830114852637)
  (26, 0.07252497860788372)
  (27, 0.0038173004867859757)
  (28, 0.13203407088216948)
  (29, 0.4064676029181584)
  (30, 0.11872543221193357)
  (31, 0.12909696545432875)
  (32, 0.6030971232787811)
  (33, 0.5494095048485471)
  (34, 2.058932974230634)
  (35, 0.015857916375411608)
  (36, 6.508618257210791)
  (37, 0.005805903591790325)
  (38, 11.381503100688938)
  (39, 0.5303190458266855)
  (40, 0.24884877244413497)
  (42, 4.946828483431141)
  (43, 2.8749711605206114)
  (44, 11.830899275554833)
  (45, 8.046873101574073)
  (46, 6.408213424333051)
  (47, 31.18316164805931)
  (48, 11.41447732791781)
  (49, 13.230388305238518)
  (50, 43.10913500391207)
  (51, 13.940884240011814)
  (52, 58.03715899471945)
  (53, 20.20207527863888)
};
\end{axis}
\end{tikzpicture}%
\vspace{-0.2cm}
\begin{tikzpicture}
\begin{axis}[
  width=.99\textwidth,
  height=3.3cm,
  ymode=log,
  xmin=-0.5, xmax=53.5,
  grid=both,
  grid style={line width=.1pt, draw=gray!20},
  xtick={0,1,2,3,4,5,6,7,8,9,10,11,12,13,14,15,16,17,18,19,20,21,22,23,24,25,26,27,28,29,30,31,32,33,34,35,36,37,38,39,40,41,42,43,44,45,46,47,48,49,50,51,52,53},
  ylabel style={font=\scriptsize},
  ylabel={\texttt{hess}},
  yticklabel style={font=\scriptsize},
  xticklabels={},
]
\addplot[gray, dashed, thin, forget plot] coordinates {(-0.5,1) (53.5,1)};
\addplot[only marks, mark=o,mark size=1.8pt,black,thick] coordinates {
  (0, 1.0)
  (1, 1.0)
  (2, 1.0)
  (3, 1.0)
  (4, 1.0)
  (5, 1.0)
  (6, 1.0)
  (7, 1.0)
  (8, 1.0)
  (9, 1.0)
  (10, 1.0)
  (11, 1.0)
  (12, 1.0)
  (13, 1.0)
  (14, 1.0)
  (15, 1.0)
  (16, 1.0)
  (17, 1.0)
  (18, 1.0)
  (19, 1.0)
  (20, 1.0)
  (21, 1.0)
  (22, 1.0)
  (23, 1.0)
  (24, 1.0)
  (25, 1.0)
  (26, 1.0)
  (27, 1.0)
  (28, 1.0)
  (29, 1.0)
  (30, 1.0)
  (31, 1.0)
  (32, 1.0)
  (33, 1.0)
  (34, 1.0)
  (35, 1.0)
  (36, 1.0)
  (37, 1.0)
  (38, 1.0)
  (39, 1.0)
  (40, 1.0)
  (41, 1.0)
  (42, 1.0)
  (43, 1.0)
  (44, 1.0)
  (45, 1.0)
  (46, 1.0)
  (47, 1.0)
  (48, 1.0)
  (49, 1.0)
  (50, 1.0)
  (51, 1.0)
  (52, 1.0)
  (53, 1.0)
};
\addplot[only marks, mark=square,mark size=1.8pt,black!55,thick] coordinates {
  (0, 0.17261904761904762)
  (1, 0.3596491228070175)
  (2, 0.06206896551724137)
  (3, 0.010783055198973043)
  (4, 0.14778325123152708)
  (5, 0.6101694915254238)
  (6, 0.08333333333333334)
  (7, 0.5158132530120482)
  (8, 0.6787878787878788)
  (9, 0.6028368794326241)
  (10, 0.2154255319148936)
  (11, 0.8209549071618037)
  (12, 1.2767441860465114)
  (13, 0.6532769556025371)
  (14, 0.543859649122807)
  (15, 0.641304347826087)
  (16, 0.03928571428571429)
  (17, 0.9267576496976495)
  (18, 0.892242019302153)
  (19, 0.4555503041575826)
  (20, 1.0728623659993066)
  (21, 0.41515946431036865)
  (22, 1.0376737381126555)
  (23, 1.1072319201995013)
  (24, 0.9778869778869779)
  (25, 0.9981916817359856)
  (26, 0.506751676641304)
  (27, 1.1617021276595743)
  (28, 0.5405511041625952)
  (29, 0.8249834541564699)
  (30, 0.7814768644129941)
  (31, 0.966644758462044)
  (32, 1.3493364850914857)
  (33, 0.35485984052765923)
  (34, 3.2434078727809093)
  (35, 2.6993609469349864)
  (36, 4.09077918674428)
  (37, 3.314326746104024)
  (38, 3.120084944523968)
  (39, 2.8800757661225678)
  (40, 3.759709236638446)
  (41, 3.5280398287492396)
  (42, 2.4246081498541865)
  (43, 2.494363709783354)
  (44, 2.719653335558598)
  (45, 2.6572200586502426)
  (46, 2.5826921782822274)
  (47, 4.599704696566523)
  (48, 3.1021141860178227)
  (49, 2.859846117417025)
  (50, 3.263298959486198)
  (51, 3.111216783934973)
  (52, 4.014827325928366)
  (53, 2.6107664128433967)
};
\addplot[only marks, mark=triangle,mark size=1.8pt,red!80!black] coordinates {
  (0, 0.0024640351887879695)
  (1, 0.005451725768790268)
  (2, 0.0012124747346451384)
  (3, 0.001075136347889272)
  (4, 0.00246866055991361)
  (5, 0.00945426351121798)
  (6, 0.0019723731490069686)
  (7, 0.015917975193563486)
  (8, 0.021548131884225316)
  (9, 0.010942028648197991)
  (10, 0.005261177042778109)
  (11, 0.05203729870481728)
  (12, 0.03501053405430231)
  (13, 0.02092168302503384)
  (14, 0.01478537568166798)
  (15, 0.016664731965237656)
  (16, 0.0014343183650530938)
  (17, 0.42689079135347235)
  (18, 0.6901106696978887)
  (19, 0.41323548448126496)
  (20, 1.2879734611400173)
  (21, 0.08150352067450653)
  (22, 0.2401686121476981)
  (23, 0.20084680904971391)
  (24, 0.4557106231741503)
  (25, 0.2905993866775101)
  (26, 0.11999644263602391)
  (27, 1.1558971661830306)
  (28, 0.21126220579996416)
  (29, 0.37377716924964044)
  (30, 0.18255477837607226)
  (31, 0.29066150820568265)
  (32, 0.6187852215946301)
  (33, 0.1014759798178267)
  (34, 3.7287888613459885)
  (35, 17.862223762799918)
  (36, 13.504550478555059)
  (37, 48.55329311428722)
  (38, 31.113321907964078)
  (39, 30.722400048341825)
  (40, 78.74001648844305)
  (41, 7.67628766473354)
  (42, 14.970658789400794)
  (43, 8.685865379257914)
  (44, 12.891861976836218)
  (45, 20.400062692261436)
  (46, 10.977920524492163)
  (47, 49.07966450856178)
  (48, 12.519308369620676)
  (49, 21.1769819944029)
  (50, 31.505089542333526)
  (51, 31.583879316243863)
  (52, 47.201155146035234)
  (53, 10.56983463031394)
};
\addplot[only marks, mark=diamond,mark size=1.8pt,green!55!black] coordinates {
  (0, 0.005106169815801526)
  (1, 0.007962194676269622)
  (2, 0.0028360345893177308)
  (3, 0.0018322011745495821)
  (4, 0.006140839729333482)
  (5, 0.01792919581281561)
  (6, 0.0030178518443915066)
  (7, 0.04213371209593346)
  (8, 0.04709728915194389)
  (9, 0.026315373967779457)
  (10, 0.010079239935715604)
  (11, 0.10637191850089489)
  (12, 0.041323339808234)
  (13, 0.0308902561629375)
  (14, 0.03425078378959606)
  (15, 0.045963223243636275)
  (16, 0.0024631199954048034)
  (17, 1.2814018885726925)
  (18, 1.267421605106935)
  (19, 0.7974394656442474)
  (20, 1.904581165195566)
  (21, 0.13415158447656)
  (22, 0.5026256503135695)
  (23, 0.39437958637255804)
  (24, 1.1878815232039621)
  (25, 0.5988843844055478)
  (26, 0.22558948436722373)
  (27, 2.487745922493307)
  (28, 0.32042254851285135)
  (29, 0.7096554329345846)
  (30, 0.5304641385170472)
  (31, 0.48193798070625554)
  (32, 0.7334628064488234)
  (33, 0.18531377460430212)
  (34, 5.184742355281165)
  (35, 36.37922368515351)
  (36, 26.589463834873623)
  (37, 86.45404675208901)
  (38, 57.88107646549756)
  (39, 61.2983810734226)
  (40, 240.90937016576112)
  (41, 58.701300276608876)
  (42, 33.33831654008691)
  (43, 13.817975483248393)
  (44, 18.86862880927326)
  (45, 41.739600207048255)
  (46, 21.15988925351569)
  (47, 41.122945474277415)
  (48, 40.30341133655382)
  (49, 44.148812139771174)
  (50, 43.06190536945117)
  (51, 43.18467119296831)
  (52, 37.967065705950205)
  (53, 15.297998033786602)
};
\addplot[only marks, mark=pentagon,mark size=1.8pt,violet] coordinates {
  (0, 0.004305994896567615)
  (1, 0.009248793138866064)
  (2, 0.002348899994640656)
  (3, 0.0016319028083910737)
  (4, 0.004245453116870163)
  (5, 0.018241870843248888)
  (6, 0.0035525813096398515)
  (7, 0.00812100559665011)
  (8, 0.03601907667527615)
  (9, 0.01930014878459126)
  (10, 0.009824110839284122)
  (11, 0.09315687770026097)
  (12, 0.08499743312110519)
  (13, 0.03649405828719457)
  (14, 0.02863615003302143)
  (15, 0.042249072788639076)
  (16, 0.0017822773587373164)
  (17, 0.26446619646555336)
  (18, 1.1802551102992025)
  (19, 0.6260914905525806)
  (20, 2.417823132017696)
  (21, 0.09156297322859198)
  (22, 0.40205750172652877)
  (23, 0.322870454771434)
  (24, 0.7942921961463124)
  (25, 0.3384116809146227)
  (26, 0.23156296988147795)
  (27, 0.8198212727119673)
  (28, 0.3774350769142594)
  (29, 0.6307557202716082)
  (30, 0.5039577509079619)
  (31, 0.6232088477805968)
  (32, 1.2110818442092735)
  (33, 0.1927357733789613)
  (34, 6.7057414638161505)
  (35, 7.643139967426188)
  (36, 20.55566566020239)
  (37, 43.43279542907261)
  (38, 17.784200394878294)
  (39, 23.423195055205746)
  (40, 85.30639733519368)
  (41, 5.033750469276806)
  (42, 8.100107368869052)
  (43, 9.967505616600377)
  (44, 7.138832537568704)
  (45, 20.33162069030172)
  (46, 10.507189309341328)
  (47, 37.32395165970945)
  (48, 19.62736030804226)
  (49, 22.034734045714636)
  (50, 39.376101411831435)
  (51, 37.45724637961827)
  (52, 68.09780814836341)
  (53, 15.050544807285778)
};
\addplot[only marks, mark=star,mark size=1.8pt,blue!70!black] coordinates {
  (0, 0.0015260096229503797)
  (1, 0.0028496991632279297)
  (2, 0.0008672043805466007)
  (3, 0.0006667635947453226)
  (4, 0.0011730528899026517)
  (5, 0.009852476731200443)
  (6, 0.0014984967422744504)
  (8, 0.018113497187360418)
  (9, 0.0030781942203854692)
  (10, 0.0035351085566628438)
  (11, 0.03585944980296979)
  (12, 0.019596200022500412)
  (13, 0.014935184896039523)
  (14, 0.00503238347354021)
  (15, 0.0145825325117846)
  (16, 0.0008774085441212731)
  (18, 0.4725532394321097)
  (19, 0.10118394164356467)
  (20, 0.9311007022183814)
  (21, 0.049393157262165074)
  (22, 0.1470353683656803)
  (23, 0.08697418053557737)
  (24, 0.22678022490797753)
  (25, 0.14910735666662464)
  (26, 0.0741846609201733)
  (27, 1.1618453559998438)
  (28, 0.14897022794986958)
  (29, 0.2485032918670942)
  (30, 0.10550327117125342)
  (31, 0.09955879664201181)
  (32, 0.45370934031219057)
  (33, 0.06266708433863112)
  (34, 2.3696772721478383)
  (35, 17.812404384211124)
  (36, 6.924925208721237)
  (37, 36.41515794730757)
  (38, 15.159116989651913)
  (39, 8.637231545148468)
  (40, 79.78040386320556)
  (42, 12.453651272162219)
  (43, 5.262402700081551)
  (44, 13.071884021803566)
  (45, 8.967672522233583)
  (46, 7.428690157631338)
  (47, 32.05514269793962)
  (48, 8.907700615677651)
  (49, 15.919883220984566)
  (50, 25.69119277216599)
  (51, 10.118010837034177)
  (52, 43.695934863820746)
  (53, 6.124447871847154)
};
\end{axis}
\end{tikzpicture}%
\vspace{-0.2cm}
\begin{tikzpicture}
\begin{axis}[
  width=.99\textwidth,
  height=3.3cm,
  ymode=log,
  xmin=-0.5, xmax=53.5,
  grid=both,
  grid style={line width=.1pt, draw=gray!20},
  xtick={0,1,2,3,4,5,6,7,8,9,10,11,12,13,14,15,16,17,18,19,20,21,22,23,24,25,26,27,28,29,30,31,32,33,34,35,36,37,38,39,40,41,42,43,44,45,46,47,48,49,50,51,52,53},
  ylabel style={font=\scriptsize},
  ylabel={composite},
  yticklabel style={font=\scriptsize},
  xticklabels={roc-2,poly-10,cam-20,rob-2,steer-5,elec-5,chain-50,glid-2,bear-5x5,tor-5x5,mix-3,surf-5x5,mar-2,gas-2,pin-2,meth-2,channel-5,glid-80,bear-35x35,tor-35x35,surf-30x30,rob-150,roc-200,steer-400,poly-100,meth-50,cam-2000,elec-50,chain-5000,mix-200,mar-150,pin-100,gas-200,channel-400,dir-5,elec-200,dir-20,bear-320x320,poly-800,tor-320x320,surf-280x280,glid-8000,roc-17000,rob-15000,meth-4000,steer-42000,cam-170000,dir-100,mar-14000,chain-500000,mix-20000,pin-10000,gas-18000,channel-42000},
  x tick label style={rotate=45, anchor=east, font=\scriptsize},
  xlabel={COPS instance (sorted by $\text{nnz}_J + \text{nnz}_H$)},
  xlabel style={font=\small},
]
\addplot[gray, dashed, thin, forget plot] coordinates {(-0.5,1) (53.5,1)};
\addplot[only marks, mark=o,mark size=1.8pt,black,thick] coordinates {
  (0, 1.0)
  (1, 1.0)
  (2, 1.0)
  (3, 1.0)
  (4, 1.0)
  (5, 1.0)
  (6, 1.0)
  (7, 1.0)
  (8, 1.0)
  (9, 1.0)
  (10, 1.0)
  (11, 1.0)
  (12, 1.0)
  (13, 1.0)
  (14, 1.0)
  (15, 1.0)
  (16, 1.0)
  (17, 1.0)
  (18, 1.0)
  (19, 1.0)
  (20, 1.0)
  (21, 1.0)
  (22, 1.0)
  (23, 1.0)
  (24, 1.0)
  (25, 1.0)
  (26, 1.0)
  (27, 1.0)
  (28, 1.0)
  (29, 1.0)
  (30, 1.0)
  (31, 1.0)
  (32, 1.0)
  (33, 1.0)
  (34, 1.0)
  (35, 1.0)
  (36, 1.0)
  (37, 1.0)
  (38, 1.0)
  (39, 1.0)
  (40, 1.0)
  (41, 1.0)
  (42, 1.0)
  (43, 1.0)
  (44, 1.0)
  (45, 1.0)
  (46, 1.0)
  (47, 1.0)
  (48, 1.0)
  (49, 1.0)
  (50, 1.0)
  (51, 1.0)
  (52, 1.0)
  (53, 1.0)
};
\addplot[only marks, mark=square,mark size=1.8pt,black!55,thick] coordinates {
  (0, 0.13252861602497398)
  (1, 0.2942308849184864)
  (2, 0.05221919640724284)
  (3, 0.005523984392162138)
  (4, 0.1348610006653792)
  (5, 0.4247429105640387)
  (6, 0.09009797690287456)
  (7, 0.36822539192205955)
  (8, 0.6819877801970131)
  (9, 0.592460881934566)
  (10, 0.25684207127056663)
  (11, 0.8045264646476525)
  (12, 1.0426531206423635)
  (13, 0.6095181275146788)
  (14, 0.558478111333388)
  (15, 0.5845849598014632)
  (16, 0.16568402360767542)
  (17, 0.9095206370782338)
  (18, 0.7204785582959803)
  (19, 0.3846493695041185)
  (20, 1.0848822176967368)
  (21, 0.22800210647436459)
  (22, 1.0386128073409926)
  (23, 1.1281786723735037)
  (24, 0.7877431215067809)
  (25, 1.0677354752114534)
  (26, 0.2975575451387566)
  (27, 0.6461828839463415)
  (28, 0.4140311822226473)
  (29, 0.8793396561586317)
  (30, 0.7737721921211862)
  (31, 0.995887524333052)
  (32, 1.296062667203809)
  (33, 0.9045062809535306)
  (34, 2.744807768221753)
  (35, 1.312556005899183)
  (36, 3.441291762836905)
  (37, 2.0349788220136795)
  (38, 3.2621890686288406)
  (39, 1.0859226568197633)
  (40, 3.0910617665022238)
  (41, 3.5261627086520764)
  (42, 2.4430645854224244)
  (43, 2.317106627370165)
  (44, 3.017788354181596)
  (45, 2.658607212130165)
  (46, 1.5839195515612303)
  (47, 4.019412059595979)
  (48, 3.234846092941451)
  (49, 2.4274229622778116)
  (50, 3.4282904438279855)
  (51, 3.33286579854779)
  (52, 3.8124989764691777)
  (53, 2.727654630857701)
};
\addplot[only marks, mark=triangle,mark size=1.8pt,red!80!black] coordinates {
  (0, 0.0013812138996377665)
  (1, 0.0030444757475727655)
  (2, 0.0009941923786305374)
  (3, 0.0006558643711426934)
  (4, 0.0017736876639548823)
  (5, 0.0034646476951749037)
  (6, 0.001454530708986466)
  (7, 0.011770779989828437)
  (8, 0.01701836396123136)
  (9, 0.008658174442899465)
  (10, 0.005645347851676546)
  (11, 0.03046742537861468)
  (12, 0.024606078377481865)
  (13, 0.015503402665987454)
  (14, 0.013292895536513298)
  (15, 0.01301828788958244)
  (16, 0.006159673751696775)
  (17, 0.2618939112434467)
  (18, 0.5066810779869383)
  (19, 0.30346300934665843)
  (20, 0.6120549048043054)
  (21, 0.03732591309953611)
  (22, 0.10314335541886906)
  (23, 0.11620050966323868)
  (24, 0.21610978507015574)
  (25, 0.19633081563069274)
  (26, 0.07131298491364511)
  (27, 0.2844782642538448)
  (28, 0.13101815509106848)
  (29, 0.35895254754117806)
  (30, 0.18962946552493543)
  (31, 0.31671839948227476)
  (32, 0.5555845776891543)
  (33, 0.4824384755322497)
  (34, 1.9008110379368617)
  (35, 2.9490985644359364)
  (36, 5.2799126515023715)
  (37, 34.58574011700605)
  (38, 14.056447485626213)
  (39, 20.285844013510747)
  (40, 41.21639340418433)
  (41, 7.708919385936471)
  (42, 8.06640876864584)
  (43, 3.9952787168839063)
  (44, 12.595386088081137)
  (45, 12.186827716192838)
  (46, 5.83064822576578)
  (47, 21.245964428971757)
  (48, 13.75743980633748)
  (49, 13.164080918251871)
  (50, 29.38434821518614)
  (51, 33.50600066129489)
  (52, 44.38102912641402)
  (53, 38.16892368242656)
};
\addplot[only marks, mark=diamond,mark size=1.8pt,green!55!black] coordinates {
  (0, 0.002682832607164297)
  (1, 0.0055178169013805674)
  (2, 0.002251427156345461)
  (3, 0.0010770832018089846)
  (4, 0.004048717356176652)
  (5, 0.005515205369503222)
  (6, 0.0023036388708108373)
  (7, 0.02763015474878688)
  (8, 0.034307157962632784)
  (9, 0.019001689976265734)
  (10, 0.010069645782696215)
  (11, 0.0554768757566334)
  (12, 0.03525374632002557)
  (13, 0.021940564201779787)
  (14, 0.021892412761012336)
  (15, 0.029864748560725192)
  (16, 0.010003388698417267)
  (17, 0.6311241294131433)
  (18, 0.8979907086146349)
  (19, 0.5486062724767141)
  (20, 1.170264958180731)
  (21, 0.06290418742084836)
  (22, 0.19260993170268187)
  (23, 0.2112102824149325)
  (24, 0.4883898821739271)
  (25, 0.34123951758535875)
  (26, 0.129269618756505)
  (27, 0.4988550903181354)
  (28, 0.18340304524484027)
  (29, 0.6810126587148403)
  (30, 0.43598287644163497)
  (31, 0.49939199163701015)
  (32, 0.7301072009470844)
  (33, 0.7877402320115646)
  (34, 2.7575572446456182)
  (35, 4.164932774800798)
  (36, 7.948655021438408)
  (37, 63.65723619033777)
  (38, 26.21761907446075)
  (39, 41.130324996175226)
  (40, 126.09450607626088)
  (41, 32.53445295218648)
  (42, 13.544846303658296)
  (43, 6.123854731727162)
  (44, 20.837569499705058)
  (45, 22.03635226459284)
  (46, 11.44414982266168)
  (47, 23.939045213410196)
  (48, 35.92321074793864)
  (49, 27.53486409532494)
  (50, 48.15670710752883)
  (51, 45.82168800870473)
  (52, 38.52308244505304)
  (53, 65.30232682149857)
};
\addplot[only marks, mark=pentagon,mark size=1.8pt,violet] coordinates {
  (0, 0.0011086958097458743)
  (1, 0.0018333056785064597)
  (2, 0.0009150429171049304)
  (3, 0.0008068232135018539)
  (4, 0.0014229382844979117)
  (5, 0.0015934589835072822)
  (6, 0.001321749539039206)
  (7, 0.007284777977387592)
  (8, 0.009092486907916192)
  (9, 0.004692840791435778)
  (10, 0.005129915056902889)
  (11, 0.024741689592120742)
  (12, 0.024690044106322487)
  (13, 0.0101856670051348)
  (14, 0.012476008240208334)
  (15, 0.01450302019609268)
  (16, 0.0050596963663504255)
  (17, 0.18112877796055765)
  (18, 0.3750437052991781)
  (19, 0.147780717567128)
  (20, 0.614533661361472)
  (21, 0.03974415363334726)
  (22, 0.08114010142199558)
  (23, 0.08472119871251212)
  (24, 0.13872389151238032)
  (25, 0.1561970766242429)
  (26, 0.06027233314990938)
  (27, 0.1377587912832116)
  (28, 0.11571938821166715)
  (29, 0.26789566333885273)
  (30, 0.20252291412309317)
  (31, 0.31910520772991685)
  (32, 0.4320096681215646)
  (33, 0.4383863648973546)
  (34, 1.7101477011149129)
  (35, 1.1064876111158726)
  (36, 4.388875169993408)
  (37, 17.241419950252126)
  (38, 7.348022116024856)
  (39, 7.115184591526456)
  (40, 34.315138614986076)
  (41, 5.116009699993334)
  (42, 4.799933317827042)
  (43, 4.308142396320949)
  (44, 7.36743345563381)
  (45, 7.797312500750649)
  (46, 3.756995066196182)
  (47, 13.826553375189622)
  (48, 12.482356209524832)
  (49, 8.955616350583304)
  (50, 19.12258450029382)
  (51, 21.93491135480155)
  (52, 28.65664653769804)
  (53, 10.22111961452416)
};
\addplot[only marks, mark=star,mark size=1.8pt,blue!70!black] coordinates {
  (0, 0.0009634737865597409)
  (1, 0.0016667982564851057)
  (2, 0.0007717637994001451)
  (3, 0.000432791173011375)
  (4, 0.0010214446259065358)
  (5, 0.0028336275853810392)
  (6, 0.001214776007723855)
  (8, 0.01318302485704403)
  (9, 0.0024368598291324662)
  (10, 0.004270821678549296)
  (11, 0.020823542235392276)
  (12, 0.014339528831189646)
  (13, 0.012024017980861456)
  (14, 0.00492626914478674)
  (15, 0.010890124417859188)
  (16, 0.004111708669743238)
  (18, 0.3367118425033019)
  (19, 0.07624033431667504)
  (20, 0.4571015761353146)
  (21, 0.025585422875759056)
  (22, 0.06841535372501838)
  (23, 0.05821230028761571)
  (24, 0.11553126093749422)
  (25, 0.12241981552643787)
  (26, 0.04404558392048691)
  (27, 0.1622654326692075)
  (28, 0.0983697745264783)
  (29, 0.268417202705591)
  (30, 0.10641864263236024)
  (31, 0.11457468792417295)
  (32, 0.4439235267822741)
  (33, 0.3218615612046407)
  (34, 1.0632882613936063)
  (35, 1.840333245786439)
  (36, 3.1923417980174627)
  (37, 20.302149265150145)
  (38, 1.010856292408514)
  (39, 0.8273055534210331)
  (40, 4.149834195669346)
  (42, 1.0223941219044805)
  (43, 1.1259129071158482)
  (44, 2.06265542004775)
  (45, 1.1477810797957682)
  (46, 3.774808425768581)
  (47, 14.933791064512787)
  (48, 1.8413429431089767)
  (49, 10.57670585287043)
  (50, 23.184354822990507)
  (51, 3.846562878836933)
  (52, 6.978953631252895)
  (53, 6.592112990787474)
};
\end{axis}
\end{tikzpicture}%
\caption{Per-callback and composite speedup over single-threaded CPU (C3) on COPS instances, sorted by $\text{nnz}_J + \text{nnz}_H$. \textsuperscript{*}Float32 precision.}
\label{fig:cops_speedup}
\end{figure}
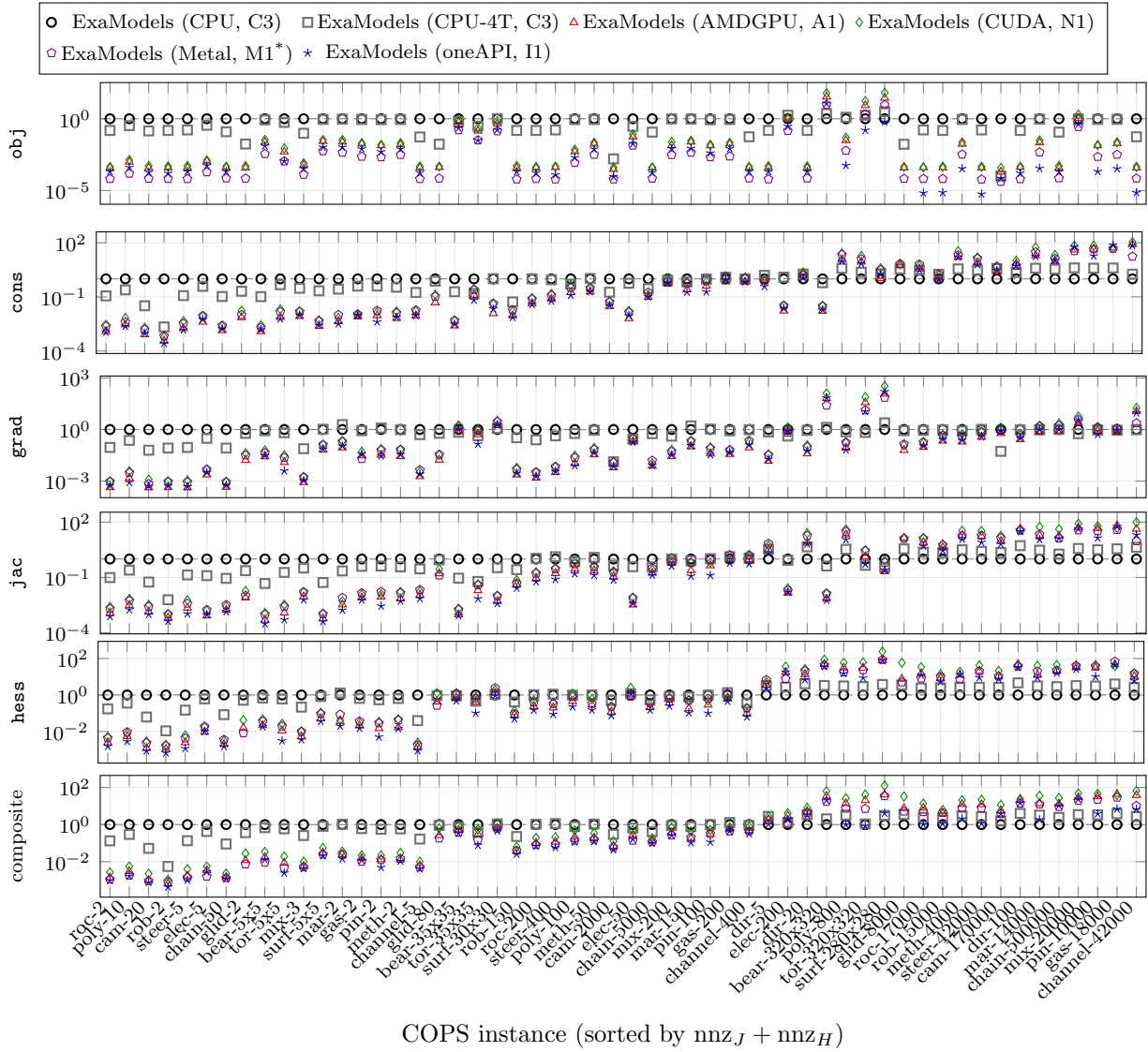

\Cref{fig:opf_speedup,fig:lv_speedup,fig:cops_speedup} show the same comparison per instance.
The crossover at which the \gls{gpu} overtakes the \gls{cpu} is visible in each suite, and the curves flatten to the right, where the \gls{gpu} time has stopped growing with the problem.

\providecommand{\figdir}{results/figures}
\begin{figure}[t]
\centering
\vspace{0.1cm}
\pgfplotslegendfromname{lv-legend}\\[2pt]
\begin{tikzpicture}
\begin{axis}[
  width=.99\textwidth,
  height=3.3cm,
  ymode=log,
  xmin=-0.5, xmax=53.5,
  grid=both,
  grid style={line width=.1pt, draw=gray!20},
  xtick={0,1,2,3,4,5,6,7,8,9,10,11,12,13,14,15,16,17,18,19,20,21,22,23,24,25,26,27,28,29,30,31,32,33,34,35,36,37,38,39,40,41,42,43,44,45,46,47,48,49,50,51,52,53},
  ylabel style={font=\scriptsize},
  ylabel={\texttt{obj}},
  yticklabel style={font=\scriptsize},
  xticklabels={},
  legend to name=lv-legend,
  legend columns=4,
  legend style={font=\scriptsize, draw=black, fill=white, /tikz/every even column/.append style={column sep=2pt}},
]
\addplot[gray, dashed, thin, forget plot] coordinates {(-0.5,1) (53.5,1)};
\addplot[only marks, mark=o,mark size=1.8pt,black,thick] coordinates {
  (0, 1.0)
  (1, 1.0)
  (2, 1.0)
  (3, 1.0)
  (4, 1.0)
  (5, 1.0)
  (6, 1.0)
  (7, 1.0)
  (8, 1.0)
  (9, 1.0)
  (10, 1.0)
  (11, 1.0)
  (12, 1.0)
  (13, 1.0)
  (14, 1.0)
  (15, 1.0)
  (16, 1.0)
  (17, 1.0)
  (18, 1.0)
  (19, 1.0)
  (20, 1.0)
  (21, 1.0)
  (22, 1.0)
  (23, 1.0)
  (24, 1.0)
  (25, 1.0)
  (26, 1.0)
  (27, 1.0)
  (28, 1.0)
  (29, 1.0)
  (30, 1.0)
  (31, 1.0)
  (32, 1.0)
  (33, 1.0)
  (34, 1.0)
  (35, 1.0)
  (36, 1.0)
  (37, 1.0)
  (38, 1.0)
  (39, 1.0)
  (40, 1.0)
  (41, 1.0)
  (42, 1.0)
  (43, 1.0)
  (44, 1.0)
  (45, 1.0)
  (46, 1.0)
  (47, 1.0)
  (48, 1.0)
  (49, 1.0)
  (50, 1.0)
  (51, 1.0)
  (52, 1.0)
  (53, 1.0)
};
\addlegendentry{ExaModels (CPU, C3)}
\addplot[only marks, mark=square,mark size=1.8pt,black!55,thick] coordinates {
  (0, 0.3333333333333333)
  (1, 0.3125)
  (2, 0.3333333333333333)
  (3, 0.47619047619047616)
  (4, 0.6666666666666666)
  (5, 0.3888888888888889)
  (6, 0.2777777777777778)
  (7, 0.2631578947368421)
  (8, 0.5)
  (9, 0.18749999999999997)
  (10, 0.3888888888888889)
  (11, 0.411764705882353)
  (12, 0.368421052631579)
  (13, 0.391304347826087)
  (14, 0.7222222222222222)
  (15, 0.13333333333333333)
  (16, 0.3888888888888889)
  (17, 0.6)
  (18, 1.3353658536585364)
  (19, 0.7178098676293622)
  (20, 1.0068587105624143)
  (21, 1.140969162995595)
  (22, 1.0836820083682008)
  (23, 1.1460176991150444)
  (24, 0.844402277039848)
  (25, 0.9393939393939393)
  (26, 0.4705882352941177)
  (27, 0.9277108433734941)
  (28, 1.3376623376623376)
  (29, 0.9307535641547863)
  (30, 1.0672097759674135)
  (31, 0.8455497382198953)
  (32, 1.0535062829347386)
  (33, 0.0017271157167530226)
  (34, 0.4416013206768469)
  (35, 0.4638888888888889)
  (36, 1.411198619006313)
  (37, 1.9307058275539166)
  (38, 0.9147216310626601)
  (39, 1.0751323591477198)
  (40, 1.065726129429154)
  (41, 1.078139460867453)
  (42, 1.3361396704524626)
  (43, 1.8333733493397357)
  (44, 0.6232890499194848)
  (45, 1.8277625461914413)
  (46, 1.0763545304997983)
  (47, 1.293633987456076)
  (48, 1.7933210497348997)
  (49, 1.2937301839892674)
  (50, 2.1800582562540294)
  (51, 5.732745153680566e-5)
  (52, 0.8341155419064631)
  (53, 1.3381340819323706)
};
\addlegendentry{ExaModels (CPU-4T, C3)}
\addplot[only marks, mark=triangle,mark size=1.8pt,red!80!black] coordinates {
  (0, 0.0009957964702508382)
  (1, 0.0009707809514736326)
  (2, 0.0009591044131589785)
  (3, 0.0019388893750438565)
  (4, 0.004676310793881618)
  (5, 0.0014145129509703968)
  (6, 0.0009704155467650744)
  (7, 0.0010298410703123413)
  (8, 0.0017716448001702477)
  (9, 0.0005909195135653777)
  (10, 0.0013453020086833059)
  (11, 0.0013397554472418143)
  (12, 0.001363488265790767)
  (13, 0.0017465623088442133)
  (14, 0.005401966408021485)
  (15, 0.0003956676082441279)
  (16, 0.0022200846155121484)
  (17, 0.0030785839208492625)
  (18, 0.09933415078938908)
  (19, 0.3628140240721752)
  (20, 0.10952975937589214)
  (21, 0.037037341523918806)
  (22, 0.040017131596697265)
  (23, 0.038769728672493835)
  (24, 0.06526390454159309)
  (25, 0.022620655078934977)
  (26, 0.009192721880778038)
  (27, 0.021636147564562034)
  (28, 0.07823016925970562)
  (29, 0.06725297585367697)
  (30, 0.07288701678926282)
  (31, 0.050292306123412914)
  (32, 0.3676157251014021)
  (33, 0.0003097853840780607)
  (34, 0.12928099173048457)
  (35, 0.18003798570482996)
  (36, 6.588920633171603)
  (37, 28.619473084366476)
  (38, 7.243381817075898)
  (39, 3.592607409778649)
  (40, 3.1236381558565816)
  (41, 3.1093203508566503)
  (42, 5.416346184358578)
  (43, 2.1727313977895983)
  (44, 0.6282312125964674)
  (45, 1.7825657354000861)
  (46, 5.3752932406589915)
  (47, 4.613287873843468)
  (48, 7.767326859097117)
  (49, 3.294513623526481)
  (50, 25.464515193931607)
  (51, 0.00019957614298160166)
  (52, 9.477437711565118)
  (53, 14.069568663769065)
};
\addlegendentry{ExaModels (AMDGPU, A1)}
\addplot[only marks, mark=diamond,mark size=1.8pt,green!55!black] coordinates {
  (0, 0.0011294384664383304)
  (1, 0.0013222617291107678)
  (2, 0.0013366384374153455)
  (3, 0.002093787115124711)
  (4, 0.004990820965261041)
  (5, 0.0018348061964854692)
  (6, 0.001157890051897072)
  (7, 0.0010915465723659499)
  (8, 0.0023960075868675076)
  (9, 0.0007129987413718417)
  (10, 0.0018662898801585686)
  (11, 0.0014841710223923476)
  (12, 0.0018675797787127133)
  (13, 0.0020224396556084586)
  (14, 0.005886481868718911)
  (15, 0.0005145528195364262)
  (16, 0.0029040717257425696)
  (17, 0.003957588894372348)
  (18, 0.12453279824287823)
  (19, 0.5211064822106095)
  (20, 0.12284392551064355)
  (21, 0.05720195177564941)
  (22, 0.06017368951938401)
  (23, 0.06301588446509108)
  (24, 0.11520306859139187)
  (25, 0.03901746855112112)
  (26, 0.014250622541229038)
  (27, 0.031737636941514664)
  (28, 0.09891045220185911)
  (29, 0.10920152895808655)
  (30, 0.11249181673019559)
  (31, 0.06061114047811565)
  (32, 0.4274114834230216)
  (33, 0.000393318437747897)
  (34, 0.1811338498232455)
  (35, 0.25169944656502247)
  (36, 11.57722401913556)
  (37, 43.220638254580315)
  (38, 15.296136472863651)
  (39, 5.221318574885627)
  (40, 4.527457885569945)
  (41, 4.822299382259047)
  (42, 7.99280907484662)
  (43, 2.881471115833161)
  (44, 1.0957229922395546)
  (45, 3.2509916355624195)
  (46, 9.457547563819114)
  (47, 9.532741716815188)
  (48, 11.169696815789846)
  (49, 5.8995934971731545)
  (50, 53.67229091678998)
  (51, 0.0003758063747847452)
  (52, 16.525677733932064)
  (53, 22.325324399244145)
};
\addlegendentry{ExaModels (CUDA, N1)}
\addplot[only marks, mark=pentagon,mark size=1.8pt,violet] coordinates {
  (0, 0.00016249783903656895)
  (1, 0.00017042924008612518)
  (2, 0.00015656083823496352)
  (3, 0.0003414731196678499)
  (4, 0.0008005116588949476)
  (5, 0.00023652077050446404)
  (6, 0.00016837197042193682)
  (7, 0.00016133595166038572)
  (8, 0.0002915263826025486)
  (9, 9.996657711295192e-5)
  (10, 0.00023833096247247905)
  (11, 0.0002448391868364903)
  (12, 0.0002470633032288562)
  (13, 0.00028530230266720925)
  (14, 0.0008968297109153484)
  (15, 6.905300183568895e-5)
  (16, 0.0004581464924411213)
  (17, 0.0005064803529515107)
  (18, 0.01884251680479835)
  (19, 0.07549395164565982)
  (20, 0.02443789304235)
  (21, 0.008697222675451017)
  (22, 0.006725134172954432)
  (23, 0.007984885554716744)
  (24, 0.015276428822805664)
  (25, 0.005199247516511465)
  (26, 0.0021652434169487523)
  (27, 0.005367647248904106)
  (28, 0.01593051935959396)
  (29, 0.016200454057853364)
  (30, 0.015744330494899025)
  (31, 0.010579890797835539)
  (32, 0.08777647758587072)
  (33, 6.8761151336811e-5)
  (34, 0.03286848160045331)
  (35, 0.039971545301122)
  (36, 1.439811039203441)
  (37, 6.436198301895802)
  (38, 1.9829084902152014)
  (39, 0.7140066830514306)
  (40, 0.7296681656139447)
  (41, 0.7390758619875502)
  (42, 1.342365341980783)
  (43, 0.3958288418855994)
  (44, 0.17142434097672302)
  (45, 0.3365339140862938)
  (46, 1.4574841704332635)
  (47, 1.2444716718939142)
  (48, 1.4246491353823394)
  (49, 0.876679657529715)
  (50, 7.068505222959938)
  (51, 5.594395462391685e-5)
  (52, 2.9308449069845706)
  (53, 3.2201898743423145)
};
\addlegendentry{ExaModels (Metal, M1\textsuperscript{*})}
\addplot[only marks, mark=star,mark size=1.8pt,blue!70!black] coordinates {
  (0, 0.0004557343519504091)
  (1, 0.000461970375552802)
  (2, 0.00046621561201338346)
  (3, 0.0008654667671545128)
  (4, 0.0021787132262960454)
  (5, 0.0006453969500295622)
  (6, 0.0004477536896780117)
  (7, 0.00041732266241956394)
  (8, 0.0008231234751346101)
  (9, 0.00026795510868035015)
  (10, 0.0006159669795378251)
  (11, 0.0006424776558248924)
  (12, 0.0006024791870509403)
  (13, 0.0007573213078159451)
  (14, 0.002219868461540628)
  (15, 0.00018608908656053238)
  (16, 0.0009661678281439411)
  (17, 0.001283019707519714)
  (18, 0.0475679919971784)
  (19, 0.09029472164599797)
  (20, 0.06369105472065559)
  (21, 0.016865612663074986)
  (22, 0.022899083388505465)
  (23, 0.021020851946560514)
  (24, 0.03881169732645683)
  (25, 0.013786322393343065)
  (26, 0.005275378108563328)
  (27, 0.010409551804194775)
  (28, 0.04494154878167541)
  (29, 0.040400384000964545)
  (30, 0.0355742148618707)
  (31, 0.02441641096435398)
  (32, 0.21786890146677323)
  (33, 9.557207051871812e-5)
  (34, 0.04819531386872089)
  (35, 0.053608392993108066)
  (36, 2.5524377218569674)
  (37, 7.552437210993065)
  (38, 2.8503281780804492)
  (39, 0.8568279093428972)
  (40, 1.0339721758452947)
  (41, 1.0629239487297713)
  (42, 1.7820186971875187)
  (43, 0.6442027520277449)
  (44, 0.25504637479676673)
  (45, 0.5137327313297009)
  (46, 1.928696809155392)
  (47, 1.4091334652315823)
  (48, 1.9286001610322068)
  (49, 1.2944916858337188)
  (50, 9.185380167949543)
  (51, 7.012009150627016e-5)
  (52, 3.387483466293959)
  (53, 3.8498651343580907)
};
\addlegendentry{ExaModels (oneAPI, I1)}
\end{axis}
\end{tikzpicture}%
\vspace{-0.2cm}
\begin{tikzpicture}
\begin{axis}[
  width=.99\textwidth,
  height=3.3cm,
  ymode=log,
  xmin=-0.5, xmax=53.5,
  grid=both,
  grid style={line width=.1pt, draw=gray!20},
  xtick={0,1,2,3,4,5,6,7,8,9,10,11,12,13,14,15,16,17,18,19,20,21,22,23,24,25,26,27,28,29,30,31,32,33,34,35,36,37,38,39,40,41,42,43,44,45,46,47,48,49,50,51,52,53},
  ylabel style={font=\scriptsize},
  ylabel={\texttt{cons}},
  yticklabel style={font=\scriptsize},
  xticklabels={},
]
\addplot[gray, dashed, thin, forget plot] coordinates {(-0.5,1) (53.5,1)};
\addplot[only marks, mark=o,mark size=1.8pt,black,thick] coordinates {
  (0, 1.0)
  (1, 1.0)
  (2, 1.0)
  (3, 1.0)
  (4, 1.0)
  (5, 1.0)
  (6, 1.0)
  (7, 1.0)
  (8, 1.0)
  (9, 1.0)
  (10, 1.0)
  (11, 1.0)
  (12, 1.0)
  (13, 1.0)
  (14, 1.0)
  (15, 1.0)
  (16, 1.0)
  (17, 1.0)
  (18, 1.0)
  (19, 1.0)
  (20, 1.0)
  (21, 1.0)
  (22, 1.0)
  (23, 1.0)
  (24, 1.0)
  (25, 1.0)
  (26, 1.0)
  (27, 1.0)
  (28, 1.0)
  (29, 1.0)
  (30, 1.0)
  (31, 1.0)
  (32, 1.0)
  (33, 1.0)
  (34, 1.0)
  (35, 1.0)
  (36, 1.0)
  (37, 1.0)
  (38, 1.0)
  (39, 1.0)
  (40, 1.0)
  (41, 1.0)
  (42, 1.0)
  (43, 1.0)
  (44, 1.0)
  (45, 1.0)
  (46, 1.0)
  (47, 1.0)
  (48, 1.0)
  (49, 1.0)
  (50, 1.0)
  (51, 1.0)
  (52, 1.0)
  (53, 1.0)
};
\addplot[only marks, mark=square,mark size=1.8pt,black!55,thick] coordinates {
  (0, 0.1842105263157895)
  (1, 0.19444444444444445)
  (2, 0.16216216216216214)
  (3, 0.22222222222222218)
  (4, 0.10204081632653061)
  (5, 0.23076923076923075)
  (6, 0.17142857142857143)
  (7, 0.15)
  (8, 0.04895104895104895)
  (9, 0.17647058823529413)
  (10, 0.2857142857142857)
  (11, 0.3)
  (12, 0.15)
  (13, 0.6111111111111112)
  (14, 0.3333333333333333)
  (15, 0.7837837837837837)
  (16, 0.4705882352941177)
  (17, 0.0891089108910891)
  (18, 0.04861111111111111)
  (19, 0.10204081632653061)
  (20, 0.25925925925925924)
  (21, 0.727810650887574)
  (22, 0.6449704142011834)
  (23, 0.6257309941520468)
  (24, 0.7241379310344828)
  (25, 0.6855670103092784)
  (26, 0.6636363636363636)
  (27, 0.6137566137566137)
  (28, 0.318359375)
  (29, 1.1915708812260535)
  (30, 0.6234782608695653)
  (31, 0.11054852320675106)
  (32, 0.4201183431952663)
  (33, 1.2433668801463862)
  (34, 0.6624548736462094)
  (35, 0.7399868680236376)
  (36, 0.04895104895104895)
  (37, 0.10416666666666667)
  (38, 0.25925925925925924)
  (39, 1.639519931418774)
  (40, 1.632288047382387)
  (41, 1.5819994045171626)
  (42, 1.9982268653855368)
  (43, 1.559506085312537)
  (44, 1.920066932681169)
  (45, 1.6622517482517483)
  (46, 6.489773135042462)
  (47, 4.005635393867818)
  (48, 4.476783433398535)
  (49, 1.8299900570848486)
  (50, 6.5961355035229685)
  (51, 4.154641020041733)
  (52, 7.793678646934461)
  (53, 5.493618110844864)
};
\addplot[only marks, mark=triangle,mark size=1.8pt,red!80!black] coordinates {
  (0, 0.001867384534098803)
  (1, 0.0019211708893252643)
  (2, 0.001619397038402972)
  (3, 0.0024211436342976814)
  (4, 0.0010138540581861833)
  (5, 0.002635682902167312)
  (6, 0.0015848306183485548)
  (7, 0.0016764053383485121)
  (8, 0.0009173366733194675)
  (9, 0.002509619077133761)
  (10, 0.005001817160174291)
  (11, 0.003668174295812327)
  (12, 0.0016500712220241564)
  (13, 0.008865631775559085)
  (14, 0.005071521769697379)
  (15, 0.024574623020622157)
  (16, 0.006877423475080929)
  (17, 0.0017290444258273323)
  (18, 0.0009161180508987495)
  (19, 0.001063541031318333)
  (20, 0.002942569851595399)
  (21, 0.03434361220707565)
  (22, 0.031514601946753526)
  (23, 0.030723614660944823)
  (24, 0.05182936473282394)
  (25, 0.03743306036498376)
  (26, 0.05946693747638455)
  (27, 0.0307155731885541)
  (28, 0.1371429006054977)
  (29, 0.2613916168615478)
  (30, 0.5993685063448423)
  (31, 0.07810986877497327)
  (32, 0.17585504023001905)
  (33, 2.2746848688991905)
  (34, 0.2966726577028788)
  (35, 0.10505479723475192)
  (36, 0.0008994144119629024)
  (37, 0.0010165887211835804)
  (38, 0.0029012447616575208)
  (39, 3.091640991632249)
  (40, 3.1196600268809678)
  (41, 2.986559737085957)
  (42, 5.119932424799603)
  (43, 3.762214723029849)
  (44, 6.102419287496955)
  (45, 3.1866784245397133)
  (46, 13.130716146128599)
  (47, 26.84117053174757)
  (48, 65.51559787254371)
  (49, 7.270637299864135)
  (50, 17.797188246651576)
  (51, 238.2581233303888)
  (52, 32.349019561214085)
  (53, 10.836647964426211)
};
\addplot[only marks, mark=diamond,mark size=1.8pt,green!55!black] coordinates {
  (0, 0.003936966112946631)
  (1, 0.004415654851423528)
  (2, 0.004332437113357517)
  (3, 0.004302398402823509)
  (4, 0.0016295022875262718)
  (5, 0.006183796531986739)
  (6, 0.0033509534283487244)
  (7, 0.0029307347654411093)
  (8, 0.0019855235816040716)
  (9, 0.005441668833584543)
  (10, 0.0136679104161414)
  (11, 0.006302441521728427)
  (12, 0.004648432035066035)
  (13, 0.021140025049392226)
  (14, 0.010845420795819623)
  (15, 0.06452902622406222)
  (16, 0.013771616876648247)
  (17, 0.003970135800369838)
  (18, 0.0017191524904369936)
  (19, 0.0025956419296591835)
  (20, 0.0045367775753669265)
  (21, 0.06241388659654581)
  (22, 0.062028842216579616)
  (23, 0.061913984498370024)
  (24, 0.1308111938502294)
  (25, 0.08602962411348007)
  (26, 0.1298869994669918)
  (27, 0.06465099899866192)
  (28, 0.29393439683606853)
  (29, 0.5621230262455241)
  (30, 1.1933644542366535)
  (31, 0.15092152858162852)
  (32, 0.2768771651095619)
  (33, 4.177162273548587)
  (34, 0.6842903310745786)
  (35, 0.25774614648356153)
  (36, 0.0013553249967886542)
  (37, 0.002020239584414414)
  (38, 0.0077920053735450094)
  (39, 6.8405237491270725)
  (40, 5.81951790202146)
  (41, 6.578464887567473)
  (42, 10.79968055509966)
  (43, 7.741631556490627)
  (44, 10.512606148991559)
  (45, 8.69748383068301)
  (46, 27.07459027164412)
  (47, 67.6187251824347)
  (48, 156.1872765339872)
  (49, 16.38779613634159)
  (50, 44.79111044977081)
  (51, 339.7910546771181)
  (52, 56.08518925223618)
  (53, 21.47351855160762)
};
\addplot[only marks, mark=pentagon,mark size=1.8pt,violet] coordinates {
  (0, 0.003864954959460412)
  (1, 0.003547536619670867)
  (2, 0.0032522555956951667)
  (3, 0.004438400667789937)
  (4, 0.0013530581518082338)
  (5, 0.005125053330176873)
  (6, 0.003363093126346991)
  (7, 0.002439720779971467)
  (8, 0.0004710332293192035)
  (9, 0.0033864319928706694)
  (10, 0.006662024667076308)
  (11, 0.006848163209042663)
  (12, 0.003237819278660886)
  (13, 0.013321561269911342)
  (14, 0.009299375352196783)
  (15, 0.042895494343672645)
  (16, 0.011870920909068779)
  (17, 0.0022452127836999586)
  (18, 0.0014195603230160556)
  (19, 0.0014510314693349452)
  (20, 0.005017044893744639)
  (21, 0.0618887413747964)
  (22, 0.0591243080086713)
  (23, 0.055550123156708216)
  (24, 0.10451103975139957)
  (25, 0.0743184729620849)
  (26, 0.11824309548390548)
  (27, 0.04835185409021121)
  (28, 0.23455340395106358)
  (29, 0.4026402479179586)
  (30, 0.8908565297164199)
  (31, 0.13805491726609181)
  (32, 0.3353017223567933)
  (33, 4.026081482642965)
  (34, 0.524166962837197)
  (35, 0.14159672822472305)
  (36, 0.0014505543389514367)
  (37, 0.0014840218494340217)
  (38, 0.005232894255766555)
  (39, 5.833832147926679)
  (40, 5.687778673401881)
  (41, 5.957531507206959)
  (42, 9.110292061771878)
  (43, 7.190516487902399)
  (44, 10.575740194862265)
  (45, 4.723090839455677)
  (46, 14.261958077182335)
  (47, 37.867345447117046)
  (48, 42.137236763863555)
  (49, 6.90837156354297)
  (50, 15.551785932561199)
  (51, 96.2937026299681)
  (52, 19.640937539819685)
  (53, 14.914041628276731)
};
\addplot[only marks, mark=star,mark size=1.8pt,blue!70!black] coordinates {
  (0, 0.0025335196227899096)
  (1, 0.0025080863660400065)
  (2, 0.0021473773030572367)
  (3, 0.0030877136067407667)
  (4, 0.0011951244048973636)
  (5, 0.0032080795526307393)
  (6, 0.002094294309818067)
  (7, 0.0019498864746846142)
  (8, 0.0009704264642239503)
  (9, 0.003168677278658412)
  (10, 0.0067272041036752295)
  (11, 0.004866557269343138)
  (12, 0.0018979924165523198)
  (13, 0.012843741880346931)
  (14, 0.006360842382894957)
  (15, 0.018958549106682317)
  (16, 0.007507697384050769)
  (17, 0.0019977679693191142)
  (18, 0.001021494179842097)
  (19, 0.0011192294719318883)
  (20, 0.0035583996364006924)
  (21, 0.043476381657140684)
  (22, 0.0321185357605081)
  (23, 0.036170008138218024)
  (24, 0.06452485062471476)
  (25, 0.04584743888967293)
  (26, 0.07520289044472796)
  (27, 0.03915485739009402)
  (28, 0.1755260270046685)
  (29, 0.34042820931976076)
  (30, 0.8199875030244902)
  (31, 0.0808587929690724)
  (32, 0.2230907586517555)
  (33, 2.3705425188024227)
  (34, 0.3391334329455396)
  (35, 0.11609250163589735)
  (36, 0.0010272564597995216)
  (37, 0.0009738453744012913)
  (38, 0.0035320206813938978)
  (39, 3.9566278966899526)
  (40, 4.079451072522458)
  (41, 3.022686093499266)
  (42, 6.620593827322892)
  (43, 3.057329034647309)
  (44, 8.157961765437397)
  (45, 4.035256813330353)
  (46, 17.80570838080574)
  (47, 31.943830142534143)
  (48, 80.25227047118597)
  (49, 7.385898371294505)
  (50, 21.863013794744727)
  (51, 251.6963252099492)
  (52, 32.50911434870934)
  (53, 10.689895712705848)
};
\end{axis}
\end{tikzpicture}%
\vspace{-0.2cm}
\begin{tikzpicture}
\begin{axis}[
  width=.99\textwidth,
  height=3.3cm,
  ymode=log,
  xmin=-0.5, xmax=53.5,
  grid=both,
  grid style={line width=.1pt, draw=gray!20},
  xtick={0,1,2,3,4,5,6,7,8,9,10,11,12,13,14,15,16,17,18,19,20,21,22,23,24,25,26,27,28,29,30,31,32,33,34,35,36,37,38,39,40,41,42,43,44,45,46,47,48,49,50,51,52,53},
  ylabel style={font=\scriptsize},
  ylabel={\texttt{grad}},
  yticklabel style={font=\scriptsize},
  xticklabels={},
]
\addplot[gray, dashed, thin, forget plot] coordinates {(-0.5,1) (53.5,1)};
\addplot[only marks, mark=o,mark size=1.8pt,black,thick] coordinates {
  (0, 1.0)
  (1, 1.0)
  (2, 1.0)
  (3, 1.0)
  (4, 1.0)
  (5, 1.0)
  (6, 1.0)
  (7, 1.0)
  (8, 1.0)
  (9, 1.0)
  (10, 1.0)
  (11, 1.0)
  (12, 1.0)
  (13, 1.0)
  (14, 1.0)
  (15, 1.0)
  (16, 1.0)
  (17, 1.0)
  (18, 1.0)
  (19, 1.0)
  (20, 1.0)
  (21, 1.0)
  (22, 1.0)
  (23, 1.0)
  (24, 1.0)
  (25, 1.0)
  (26, 1.0)
  (27, 1.0)
  (28, 1.0)
  (29, 1.0)
  (30, 1.0)
  (31, 1.0)
  (32, 1.0)
  (33, 1.0)
  (34, 1.0)
  (35, 1.0)
  (36, 1.0)
  (37, 1.0)
  (38, 1.0)
  (39, 1.0)
  (40, 1.0)
  (41, 1.0)
  (42, 1.0)
  (43, 1.0)
  (44, 1.0)
  (45, 1.0)
  (46, 1.0)
  (47, 1.0)
  (48, 1.0)
  (49, 1.0)
  (50, 1.0)
  (51, 1.0)
  (52, 1.0)
  (53, 1.0)
};
\addplot[only marks, mark=square,mark size=1.8pt,black!55,thick] coordinates {
  (0, 0.19999999999999998)
  (1, 0.1935483870967742)
  (2, 0.19999999999999998)
  (3, 0.358974358974359)
  (4, 0.6571428571428571)
  (5, 0.2857142857142857)
  (6, 0.2903225806451613)
  (7, 0.2727272727272727)
  (8, 0.3235294117647059)
  (9, 0.1515151515151515)
  (10, 0.5111111111111111)
  (11, 0.2571428571428572)
  (12, 0.2564102564102564)
  (13, 0.3181818181818182)
  (14, 0.5689655172413793)
  (15, 0.23333333333333336)
  (16, 0.3968253968253968)
  (17, 0.5416666666666667)
  (18, 0.3634751773049645)
  (19, 0.9928543800004073)
  (20, 0.3795005807200929)
  (21, 0.20678285445124828)
  (22, 0.20832409910908203)
  (23, 0.19753629564452266)
  (24, 0.28060671722643554)
  (25, 0.38286713286713286)
  (26, 0.14471403812824957)
  (27, 0.37688219663418954)
  (28, 0.5901213171577123)
  (29, 0.28523862375138737)
  (30, 0.414996288047513)
  (31, 0.264277035236938)
  (32, 0.6423889293517844)
  (33, 0.20478210566910912)
  (34, 0.5845500398194851)
  (35, 0.5952474934870136)
  (36, 0.896003592569557)
  (37, 1.509593104229823)
  (38, 0.8536658827019769)
  (39, 0.43974575906900765)
  (40, 0.47143527861936607)
  (41, 0.44294888443622377)
  (42, 0.5683571305342017)
  (43, 0.8887873861620976)
  (44, 0.2549902644015996)
  (45, 0.8809343259137362)
  (46, 1.6783923699939105)
  (47, 0.5681483192337482)
  (48, 1.0833860041296433)
  (49, 0.4720207598297087)
  (50, 1.1890534617021398)
  (51, 0.315500415932449)
  (52, 0.5932211579818207)
  (53, 1.141887432827378)
};
\addplot[only marks, mark=triangle,mark size=1.8pt,red!80!black] coordinates {
  (0, 0.0012544533485891586)
  (1, 0.0013066016770541754)
  (2, 0.0012570996141133693)
  (3, 0.0029514551674307923)
  (4, 0.009569775294373932)
  (5, 0.0022902681606968398)
  (6, 0.0018592291929524234)
  (7, 0.0019440083643681488)
  (8, 0.002404908418081304)
  (9, 0.0010801730910839798)
  (10, 0.004913159006690409)
  (11, 0.0019115218348519684)
  (12, 0.0021086072491817054)
  (13, 0.002815520862832618)
  (14, 0.007115470344362633)
  (15, 0.0015546053841994116)
  (16, 0.004418199186653005)
  (17, 0.005700611838336945)
  (18, 0.17742045778394339)
  (19, 1.0789737383212403)
  (20, 0.27978129117575995)
  (21, 0.097681412411633)
  (22, 0.10425962919197745)
  (23, 0.10253244029359025)
  (24, 0.17070241821791043)
  (25, 0.1939583252595782)
  (26, 0.07047603402220837)
  (27, 0.17905490160352439)
  (28, 0.44145804928741905)
  (29, 0.17602525160087149)
  (30, 0.23411898128079905)
  (31, 0.19931358189345247)
  (32, 0.7430909789285937)
  (33, 0.1168381397287512)
  (34, 0.3614039851557892)
  (35, 0.4759677877542221)
  (36, 17.48727282050872)
  (37, 104.36245050815062)
  (38, 35.66103914906802)
  (39, 10.22068189297332)
  (40, 10.582054212377837)
  (41, 9.468301472709898)
  (42, 16.68276697956204)
  (43, 20.623665572322672)
  (44, 7.01068203548543)
  (45, 19.030729304928148)
  (46, 43.08539496162845)
  (47, 17.328475260421182)
  (48, 25.865150559551115)
  (49, 19.053917965772918)
  (50, 75.27365723883308)
  (51, 12.726399384707888)
  (52, 39.09815019699187)
  (53, 51.197904118341654)
};
\addplot[only marks, mark=diamond,mark size=1.8pt,green!55!black] coordinates {
  (0, 0.0026965581940179015)
  (1, 0.0034930642531238293)
  (2, 0.003581838224808022)
  (3, 0.005160107578365089)
  (4, 0.01878037892149877)
  (5, 0.005111077723279832)
  (6, 0.004513859874452404)
  (7, 0.003582442734752397)
  (8, 0.006319010596894604)
  (9, 0.0024041570182836863)
  (10, 0.014152375526073027)
  (11, 0.0034135302569101764)
  (12, 0.006020093204889026)
  (13, 0.006435118330817455)
  (14, 0.014148974613733958)
  (15, 0.0037497289682509862)
  (16, 0.009419565968253576)
  (17, 0.015881272680480726)
  (18, 0.3790259938846189)
  (19, 2.8986254571872903)
  (20, 0.6114927903760187)
  (21, 0.1744841348730407)
  (22, 0.21836160454058334)
  (23, 0.21641662794022054)
  (24, 0.4168965253649102)
  (25, 0.5177680805411656)
  (26, 0.147668831953024)
  (27, 0.38118715521123486)
  (28, 0.8890035856782558)
  (29, 0.40478494870327303)
  (30, 0.4181618164596226)
  (31, 0.4121012987434008)
  (32, 1.232210713117257)
  (33, 0.22762703102629736)
  (34, 0.8205042659738281)
  (35, 1.2092689017897769)
  (36, 28.564994375937307)
  (37, 245.327361446694)
  (38, 98.65110447080474)
  (39, 22.335669793529107)
  (40, 21.762418886042745)
  (41, 21.42802456861543)
  (42, 35.79706853306856)
  (43, 41.83873408003087)
  (44, 13.045185523839258)
  (45, 47.91947827924037)
  (46, 93.63235947108588)
  (47, 38.988524039246556)
  (48, 59.973351297141534)
  (49, 40.748525537128835)
  (50, 154.5195373917314)
  (51, 19.91011448647458)
  (52, 73.98327653813689)
  (53, 103.83825286785174)
};
\addplot[only marks, mark=pentagon,mark size=1.8pt,violet] coordinates {
  (0, 0.0026451760007108805)
  (1, 0.002334710341264621)
  (2, 0.0027306562701078225)
  (3, 0.006108129236105996)
  (4, 0.018432168100675175)
  (5, 0.0045925296084596235)
  (6, 0.0037479892672701264)
  (7, 0.0034986349112066093)
  (8, 0.0030345344789468338)
  (9, 0.001862441165290839)
  (10, 0.009126390465357747)
  (11, 0.0037656650312895426)
  (12, 0.004340404026384716)
  (13, 0.0053002042343907)
  (14, 0.014326154705832587)
  (15, 0.0031244755875784206)
  (16, 0.008112128513950441)
  (17, 0.01086589173093447)
  (18, 0.37241151361483615)
  (19, 1.887501634945921)
  (20, 0.5295527147373279)
  (21, 0.19765450466947448)
  (22, 0.18434553292260536)
  (23, 0.19175746561973195)
  (24, 0.339895703986602)
  (25, 0.38592900352547493)
  (26, 0.14959476765709442)
  (27, 0.30575979344168125)
  (28, 0.8835355980287485)
  (29, 0.3155399310664504)
  (30, 0.40548837731668363)
  (31, 0.36316184141888364)
  (32, 1.491814286304103)
  (33, 0.23312112297755322)
  (34, 0.6910391861194342)
  (35, 0.9363109888954004)
  (36, 18.493397949244077)
  (37, 45.26606033544858)
  (38, 20.248534810311284)
  (39, 8.097792502962164)
  (40, 8.650895281602212)
  (41, 8.071951716179235)
  (42, 11.890255757005825)
  (43, 15.052017836621575)
  (44, 5.170114957034069)
  (45, 14.549330555195056)
  (46, 44.706417002645814)
  (47, 11.589772032774047)
  (48, 14.770511749925838)
  (49, 9.830562305442784)
  (50, 39.6397923458285)
  (51, 8.920305846582513)
  (52, 21.080089939333256)
  (53, 27.487946108383277)
};
\addplot[only marks, mark=star,mark size=1.8pt,blue!70!black] coordinates {
  (0, 0.001515511843400484)
  (1, 0.0012363763866262542)
  (2, 0.0015050810709159534)
  (3, 0.0034584296582160948)
  (4, 0.011791332853072347)
  (5, 0.0025541478253460538)
  (6, 0.0020900613957857304)
  (7, 0.002165253275136843)
  (8, 0.0025797208556294308)
  (9, 0.0009142921977145467)
  (10, 0.006212059240487276)
  (11, 0.0023737269695668526)
  (12, 0.0024245868331284537)
  (13, 0.0028386621814993156)
  (14, 0.0083333026305845)
  (15, 0.0018798973036810223)
  (16, 0.004466694331346454)
  (17, 0.007204286530535336)
  (18, 0.21500516088424532)
  (19, 1.4470575023295564)
  (20, 0.33749660921778896)
  (21, 0.10050714830992549)
  (22, 0.11728228420453975)
  (23, 0.10859688994325102)
  (24, 0.19286186754738896)
  (25, 0.21498498239093844)
  (26, 0.08561071088830009)
  (27, 0.19804850082113912)
  (28, 0.45667209906777617)
  (29, 0.20773120640216047)
  (30, 0.21812299485486292)
  (31, 0.18657358694842874)
  (32, 0.8598674545515814)
  (33, 0.14708502208151983)
  (34, 0.3601033834246025)
  (35, 0.6156251911411437)
  (36, 22.542551390014207)
  (37, 135.00322846656087)
  (38, 42.032114153313415)
  (39, 10.907906891429715)
  (40, 12.228433714521765)
  (41, 11.283840909742297)
  (42, 19.469744549913695)
  (43, 15.732538933738542)
  (44, 8.775156694095845)
  (45, 21.372513508021097)
  (46, 46.725450229180886)
  (47, 19.781344408083573)
  (48, 22.872546844209054)
  (49, 21.151874249327534)
  (50, 79.14898025516248)
  (51, 16.076515132757102)
  (52, 36.98885432394119)
  (53, 62.60575706107337)
};
\end{axis}
\end{tikzpicture}%
\vspace{-0.2cm}
\begin{tikzpicture}
\begin{axis}[
  width=.99\textwidth,
  height=3.3cm,
  ymode=log,
  xmin=-0.5, xmax=53.5,
  grid=both,
  grid style={line width=.1pt, draw=gray!20},
  xtick={0,1,2,3,4,5,6,7,8,9,10,11,12,13,14,15,16,17,18,19,20,21,22,23,24,25,26,27,28,29,30,31,32,33,34,35,36,37,38,39,40,41,42,43,44,45,46,47,48,49,50,51,52,53},
  ylabel style={font=\scriptsize},
  ylabel={\texttt{jac}},
  yticklabel style={font=\scriptsize},
  xticklabels={},
]
\addplot[gray, dashed, thin, forget plot] coordinates {(-0.5,1) (53.5,1)};
\addplot[only marks, mark=o,mark size=1.8pt,black,thick] coordinates {
  (0, 1.0)
  (1, 1.0)
  (2, 1.0)
  (3, 1.0)
  (4, 1.0)
  (5, 1.0)
  (6, 1.0)
  (7, 1.0)
  (8, 1.0)
  (9, 1.0)
  (10, 1.0)
  (11, 1.0)
  (12, 1.0)
  (13, 1.0)
  (14, 1.0)
  (15, 1.0)
  (16, 1.0)
  (17, 1.0)
  (18, 1.0)
  (19, 1.0)
  (20, 1.0)
  (21, 1.0)
  (22, 1.0)
  (23, 1.0)
  (24, 1.0)
  (25, 1.0)
  (26, 1.0)
  (27, 1.0)
  (28, 1.0)
  (29, 1.0)
  (30, 1.0)
  (31, 1.0)
  (32, 1.0)
  (33, 1.0)
  (34, 1.0)
  (35, 1.0)
  (36, 1.0)
  (37, 1.0)
  (38, 1.0)
  (39, 1.0)
  (40, 1.0)
  (41, 1.0)
  (42, 1.0)
  (43, 1.0)
  (44, 1.0)
  (45, 1.0)
  (46, 1.0)
  (47, 1.0)
  (48, 1.0)
  (49, 1.0)
  (50, 1.0)
  (51, 1.0)
  (52, 1.0)
  (53, 1.0)
};
\addplot[only marks, mark=square,mark size=1.8pt,black!55,thick] coordinates {
  (0, 0.08333333333333334)
  (1, 0.07407407407407407)
  (2, 0.07407407407407407)
  (3, 0.2058823529411765)
  (4, 0.08955223880597013)
  (5, 0.14705882352941177)
  (6, 0.10869565217391304)
  (7, 0.0847457627118644)
  (8, 0.04142011834319527)
  (9, 0.17647058823529413)
  (10, 0.5)
  (11, 0.2888888888888889)
  (12, 0.2777777777777778)
  (13, 0.5714285714285714)
  (14, 0.5652173913043478)
  (15, 0.9076923076923076)
  (16, 0.6)
  (17, 0.1782608695652174)
  (18, 0.04166666666666667)
  (19, 0.08823529411764706)
  (20, 0.21212121212121213)
  (21, 0.7464788732394366)
  (22, 0.7363636363636363)
  (23, 0.7318181818181818)
  (24, 1.54468085106383)
  (25, 1.250965250965251)
  (26, 1.3260869565217392)
  (27, 0.7272727272727272)
  (28, 0.7262479871175523)
  (29, 1.2886699507389163)
  (30, 0.845403212543172)
  (31, 0.50066401062417)
  (32, 0.8604972375690607)
  (33, 1.5377545951316443)
  (34, 0.8614623419461243)
  (35, 1.14002828854314)
  (36, 0.04142011834319527)
  (37, 0.08823529411764706)
  (38, 0.1818181818181818)
  (39, 0.7138259130077965)
  (40, 0.7245126543062772)
  (41, 0.7202953111591223)
  (42, 1.5582273290127528)
  (43, 1.163312016851264)
  (44, 1.281366620275147)
  (45, 0.7652823198481488)
  (46, 3.0379691426102577)
  (47, 3.255601353038217)
  (48, 2.9614815894502056)
  (49, 3.0047721295037304)
  (50, 3.5283015746533475)
  (51, 3.941577465495294)
  (52, 2.706902591971142)
  (53, 3.200090651610572)
};
\addplot[only marks, mark=triangle,mark size=1.8pt,red!80!black] coordinates {
  (0, 0.000789801710684442)
  (1, 0.0008117613125529392)
  (2, 0.0007877588187753452)
  (3, 0.0019214268027944123)
  (4, 0.000933139252188642)
  (5, 0.001458741869797007)
  (6, 0.0009911314450316863)
  (7, 0.0010155202219556743)
  (8, 0.0007632351603318234)
  (9, 0.0016603686990934586)
  (10, 0.004162562102929038)
  (11, 0.0035574115003106895)
  (12, 0.0027335581989288093)
  (13, 0.006576972572910397)
  (14, 0.0055235994319853265)
  (15, 0.025631134429187703)
  (16, 0.006649066705303359)
  (17, 0.0034262333382899506)
  (18, 0.000766151408866249)
  (19, 0.0009964915326819027)
  (20, 0.001971508591502102)
  (21, 0.033130995260661135)
  (22, 0.03417377589492554)
  (23, 0.03440633190710468)
  (24, 0.10163119182786666)
  (25, 0.0696331097854758)
  (26, 0.09839811498247202)
  (27, 0.037775915879835)
  (28, 0.3845274984715778)
  (29, 0.3876086332462048)
  (30, 0.6419255198058603)
  (31, 0.3351694514506516)
  (32, 0.5251198320506628)
  (33, 2.6615281475753774)
  (34, 0.6529215826220826)
  (35, 0.3281491601600374)
  (36, 0.000750926621727978)
  (37, 0.0009667445513672869)
  (38, 0.0016422848810654415)
  (39, 3.5622405065979073)
  (40, 3.477520409715798)
  (41, 3.3222052485373212)
  (42, 10.70502594755096)
  (43, 7.395554140083682)
  (44, 11.173253624891196)
  (45, 4.020825569865678)
  (46, 43.70982034103773)
  (47, 38.1640654701019)
  (48, 69.74999382182786)
  (49, 31.148342222012936)
  (50, 53.01168098391417)
  (51, 277.9852812890159)
  (52, 70.45299447662119)
  (53, 35.98599276142785)
};
\addplot[only marks, mark=diamond,mark size=1.8pt,green!55!black] coordinates {
  (0, 0.001656934798207296)
  (1, 0.0021737302729290626)
  (2, 0.0021689497320739773)
  (3, 0.003327431072716909)
  (4, 0.0015812153263691285)
  (5, 0.0031311389305568676)
  (6, 0.0021809996328287116)
  (7, 0.0020548869594970726)
  (8, 0.0017481404025357498)
  (9, 0.0035902719212713686)
  (10, 0.010925515733507443)
  (11, 0.006672907785988842)
  (12, 0.007332871861041375)
  (13, 0.013883860689654226)
  (14, 0.010874901102185449)
  (15, 0.0677206043451009)
  (16, 0.013846950450561767)
  (17, 0.007945969777674498)
  (18, 0.0015086593588055051)
  (19, 0.0026080049164022277)
  (20, 0.004746713116464194)
  (21, 0.05525206879904208)
  (22, 0.07033635422040985)
  (23, 0.0654585777916447)
  (24, 0.23915693425318932)
  (25, 0.14938607992289907)
  (26, 0.21712374150763641)
  (27, 0.08161244565312475)
  (28, 0.7918391407755239)
  (29, 0.8036064283697545)
  (30, 1.2098359339771796)
  (31, 0.6347002798248187)
  (32, 0.8914079916572224)
  (33, 4.789650295261323)
  (34, 1.4363730200784965)
  (35, 0.7703209947273133)
  (36, 0.0012006765709562204)
  (37, 0.002053792129318144)
  (38, 0.003997988212331554)
  (39, 8.57176148227925)
  (40, 7.124193610420074)
  (41, 7.406603059127937)
  (42, 22.58608081670237)
  (43, 14.270237059648286)
  (44, 20.866134933684066)
  (45, 10.276810484960158)
  (46, 99.85022848762257)
  (47, 84.62115274404516)
  (48, 151.70434910001586)
  (49, 70.97499171702978)
  (50, 131.62794218905677)
  (51, 384.9172441516679)
  (52, 93.83923880682782)
  (53, 65.464417777531)
};
\addplot[only marks, mark=pentagon,mark size=1.8pt,violet] coordinates {
  (0, 0.0015659232749967269)
  (1, 0.0015852762905095962)
  (2, 0.0015743538875116524)
  (3, 0.003571023024947619)
  (4, 0.0010817507374029458)
  (5, 0.002831244057261839)
  (6, 0.0017195871337357268)
  (7, 0.0016841567964614988)
  (8, 0.0004713921940986119)
  (9, 0.002373507546744762)
  (10, 0.008404887479024025)
  (11, 0.0064159011548888525)
  (12, 0.005167517740534978)
  (13, 0.011850340646290394)
  (14, 0.010569910239608868)
  (15, 0.046660151411813115)
  (16, 0.01211701746363288)
  (17, 0.005052451835472519)
  (18, 0.0012964904039772095)
  (19, 0.0014623168591630345)
  (20, 0.003516572896373235)
  (21, 0.06326350765958907)
  (22, 0.05573784136109727)
  (23, 0.05725598824683038)
  (24, 0.1968979581570834)
  (25, 0.1299913197120416)
  (26, 0.20280497559330538)
  (27, 0.05822344934968723)
  (28, 0.7620285019911474)
  (29, 0.6457192166475448)
  (30, 1.1105423168631665)
  (31, 0.580288912591704)
  (32, 0.965649107325845)
  (33, 4.664059924238263)
  (34, 1.153420646332605)
  (35, 0.4838318637406493)
  (36, 0.0012748991739524964)
  (37, 0.001476309109601368)
  (38, 0.0030093544687926303)
  (39, 5.132215928852673)
  (40, 5.1392123873043545)
  (41, 5.1059111867181075)
  (42, 9.070456043353905)
  (43, 8.42735641804534)
  (44, 8.95637844984755)
  (45, 5.521480483347794)
  (46, 21.378798389350266)
  (47, 34.76769301971341)
  (48, 37.93174838582443)
  (49, 21.33051355343145)
  (50, 25.325045564832955)
  (51, 83.31090491469169)
  (52, 14.329395804069716)
  (53, 13.059528072710021)
};
\addplot[only marks, mark=star,mark size=1.8pt,blue!70!black] coordinates {
  (0, 0.0007548582051324304)
  (1, 0.0007273844310385434)
  (2, 0.0007604916409629682)
  (3, 0.001776487774694593)
  (4, 0.0008026264903130212)
  (5, 0.0013429531646990806)
  (6, 0.0009306909792882931)
  (7, 0.0009699780427622426)
  (8, 0.0005595038887970099)
  (9, 0.0015735943441683247)
  (10, 0.004557689871081184)
  (11, 0.003493670869883278)
  (12, 0.0026717999852912066)
  (13, 0.00757979263099098)
  (14, 0.005816318097407543)
  (15, 0.02266647070542298)
  (16, 0.006029003122420717)
  (17, 0.002558896216896377)
  (18, 0.0005822254098083792)
  (19, 0.0008240755778721143)
  (20, 0.0018253076866408143)
  (21, 0.030207152157295355)
  (22, 0.030077444964866758)
  (23, 0.02888889906479923)
  (24, 0.09399187858382876)
  (25, 0.060484335484909245)
  (26, 0.06961404374546129)
  (27, 0.033374612514230226)
  (28, 0.3662201475413277)
  (29, 0.35162327611245175)
  (30, 0.7224130704626089)
  (31, 0.27539416739843964)
  (32, 0.5574575750514058)
  (33, 2.395297358466382)
  (34, 0.6330870712162272)
  (35, 0.23721314871586907)
  (36, 0.0005950553363379253)
  (37, 0.000776924207407562)
  (38, 0.0015279880034605874)
  (39, 3.085742665085325)
  (40, 3.1101575368793317)
  (41, 3.106766694473001)
  (42, 10.358090183797849)
  (43, 5.965318760530479)
  (44, 11.398804119311617)
  (45, 3.684190502402362)
  (46, 48.819198644981924)
  (47, 33.91487426860642)
  (48, 67.99714893323072)
  (49, 30.953682630878284)
  (50, 42.64915711610129)
  (51, 138.60519701356813)
  (52, 24.071976713693942)
  (53, 20.65658583126892)
};
\end{axis}
\end{tikzpicture}%
\vspace{-0.2cm}
\begin{tikzpicture}
\begin{axis}[
  width=.99\textwidth,
  height=3.3cm,
  ymode=log,
  xmin=-0.5, xmax=53.5,
  grid=both,
  grid style={line width=.1pt, draw=gray!20},
  xtick={0,1,2,3,4,5,6,7,8,9,10,11,12,13,14,15,16,17,18,19,20,21,22,23,24,25,26,27,28,29,30,31,32,33,34,35,36,37,38,39,40,41,42,43,44,45,46,47,48,49,50,51,52,53},
  ylabel style={font=\scriptsize},
  ylabel={\texttt{hess}},
  yticklabel style={font=\scriptsize},
  xticklabels={},
]
\addplot[gray, dashed, thin, forget plot] coordinates {(-0.5,1) (53.5,1)};
\addplot[only marks, mark=o,mark size=1.8pt,black,thick] coordinates {
  (0, 1.0)
  (1, 1.0)
  (2, 1.0)
  (3, 1.0)
  (4, 1.0)
  (5, 1.0)
  (6, 1.0)
  (7, 1.0)
  (8, 1.0)
  (9, 1.0)
  (10, 1.0)
  (11, 1.0)
  (12, 1.0)
  (13, 1.0)
  (14, 1.0)
  (15, 1.0)
  (16, 1.0)
  (17, 1.0)
  (18, 1.0)
  (19, 1.0)
  (20, 1.0)
  (21, 1.0)
  (22, 1.0)
  (23, 1.0)
  (24, 1.0)
  (25, 1.0)
  (26, 1.0)
  (27, 1.0)
  (28, 1.0)
  (29, 1.0)
  (30, 1.0)
  (31, 1.0)
  (32, 1.0)
  (33, 1.0)
  (34, 1.0)
  (35, 1.0)
  (36, 1.0)
  (37, 1.0)
  (38, 1.0)
  (39, 1.0)
  (40, 1.0)
  (41, 1.0)
  (42, 1.0)
  (43, 1.0)
  (44, 1.0)
  (45, 1.0)
  (46, 1.0)
  (47, 1.0)
  (48, 1.0)
  (49, 1.0)
  (50, 1.0)
  (51, 1.0)
  (52, 1.0)
  (53, 1.0)
};
\addplot[only marks, mark=square,mark size=1.8pt,black!55,thick] coordinates {
  (0, 0.136986301369863)
  (1, 0.15492957746478875)
  (2, 0.14084507042253522)
  (3, 0.4054054054054054)
  (4, 0.4451612903225806)
  (5, 0.2903225806451613)
  (6, 0.1388888888888889)
  (7, 0.12345679012345678)
  (8, 0.1)
  (9, 0.13793103448275865)
  (10, 0.6987951807228915)
  (11, 0.4383561643835617)
  (12, 0.37037037037037035)
  (13, 0.6103896103896104)
  (14, 0.7425742574257426)
  (15, 0.8188976377952757)
  (16, 0.6160714285714285)
  (17, 0.32142857142857145)
  (18, 1.1076672104404568)
  (19, 1.4929160499048422)
  (20, 1.1889053254437871)
  (21, 1.1365409622886866)
  (22, 1.1558441558441557)
  (23, 1.1485411140583555)
  (24, 1.1635311143270621)
  (25, 1.1139705882352942)
  (26, 1.0340501792114696)
  (27, 1.0885189437428244)
  (28, 1.0733597692862293)
  (29, 1.1971496437054632)
  (30, 1.0407980273481283)
  (31, 0.8036598493003229)
  (32, 1.0334707131260148)
  (33, 1.4494296577946768)
  (34, 1.2059848876582542)
  (35, 1.2791420995513587)
  (36, 3.721462555114068)
  (37, 3.7846261654416766)
  (38, 4.229028379809848)
  (39, 2.255435152486748)
  (40, 2.2670727592042033)
  (41, 2.238986272545342)
  (42, 2.9071072597507643)
  (43, 2.0581033052028372)
  (44, 1.569334510632483)
  (45, 1.953876749904464)
  (46, 3.2070616505827316)
  (47, 3.276437939943037)
  (48, 3.3287805268377872)
  (49, 3.0357803706814526)
  (50, 3.4502204780580463)
  (51, 3.61148935546476)
  (52, 3.3020887396033833)
  (53, 3.3885684179827864)
};
\addplot[only marks, mark=triangle,mark size=1.8pt,red!80!black] coordinates {
  (0, 0.0015260882450577606)
  (1, 0.0016909632204585151)
  (2, 0.0015151581611861768)
  (3, 0.00582126520053135)
  (4, 0.008936318365516392)
  (5, 0.0037774826500536676)
  (6, 0.001546323381968693)
  (7, 0.0015860682404972273)
  (8, 0.0019364756391973196)
  (9, 0.0015859514611585909)
  (10, 0.015852016722948385)
  (11, 0.0062208540352339525)
  (12, 0.0059748621536565094)
  (13, 0.012213918931599552)
  (14, 0.020631935213544696)
  (15, 0.02931116259541664)
  (16, 0.014515917411355493)
  (17, 0.00794570134886388)
  (18, 0.12594732993028532)
  (19, 0.9499701516283553)
  (20, 0.4849793511256275)
  (21, 0.1399752718787832)
  (22, 0.1427104427069229)
  (23, 0.14039431006503264)
  (24, 0.32487482858736294)
  (25, 0.1498814840607645)
  (26, 0.11267897404876293)
  (27, 0.14455207104986317)
  (28, 1.624349746269614)
  (29, 0.7401785295507591)
  (30, 1.2071205082728702)
  (31, 0.7955585337391206)
  (32, 2.2252132255111086)
  (33, 3.1838999123968144)
  (34, 1.417743896291426)
  (35, 0.8235299895327106)
  (36, 12.577700826000655)
  (37, 93.62031299929505)
  (38, 63.87393315571343)
  (39, 14.574371555386369)
  (40, 14.538929077970625)
  (41, 13.349648470794813)
  (42, 33.33136548241595)
  (43, 15.399817895328708)
  (44, 11.596565959909496)
  (45, 15.416389631779312)
  (46, 161.1671454619382)
  (47, 68.72215761817779)
  (48, 66.86644535206413)
  (49, 74.49833951342015)
  (50, 174.5578136394989)
  (51, 240.56404958625322)
  (52, 126.95323067950777)
  (53, 87.500003035042)
};
\addplot[only marks, mark=diamond,mark size=1.8pt,green!55!black] coordinates {
  (0, 0.00324750038271792)
  (1, 0.004152547752638136)
  (2, 0.0038366050084852653)
  (3, 0.010577503147873752)
  (4, 0.0148562600643357)
  (5, 0.008292705526137338)
  (6, 0.003157921499395781)
  (7, 0.0029700080949649635)
  (8, 0.004499260800374318)
  (9, 0.0035466999816471576)
  (10, 0.04181910481432998)
  (11, 0.011721803564169688)
  (12, 0.015447069948024575)
  (13, 0.026944982730989234)
  (14, 0.0442484352824067)
  (15, 0.07223923807441839)
  (16, 0.02868021947251149)
  (17, 0.015494443754160599)
  (18, 0.2493630765727426)
  (19, 2.0930257385156086)
  (20, 1.0379651409539097)
  (21, 0.26625180695185474)
  (22, 0.25539256449915543)
  (23, 0.2592170853201692)
  (24, 0.7913696415988767)
  (25, 0.35984630928502087)
  (26, 0.2495185675458179)
  (27, 0.28317029964862483)
  (28, 3.4518517802321265)
  (29, 1.7001260517210826)
  (30, 1.956713164051839)
  (31, 1.6243756781920755)
  (32, 3.360971703668515)
  (33, 6.0060838541319495)
  (34, 2.8264282329859896)
  (35, 1.6862729668225458)
  (36, 19.72353356174957)
  (37, 192.50449167258216)
  (38, 123.85236190920763)
  (39, 28.8535088103709)
  (40, 26.073220988491766)
  (41, 27.849334333875476)
  (42, 60.32868870310512)
  (43, 29.73511192676641)
  (44, 20.86873180430763)
  (45, 33.04200035182617)
  (46, 172.99175597120635)
  (47, 135.19760633402984)
  (48, 134.17836528498088)
  (49, 115.9326736722142)
  (50, 213.3667996066657)
  (51, 341.52206739784214)
  (52, 179.26671674433058)
  (53, 151.56360958323097)
};
\addplot[only marks, mark=pentagon,mark size=1.8pt,violet] coordinates {
  (0, 0.00292206969642614)
  (1, 0.002191920293807498)
  (2, 0.002913236386952752)
  (3, 0.01040947951021733)
  (4, 0.012142604703038955)
  (5, 0.006806969661151578)
  (6, 0.0029358026561840346)
  (7, 0.0022571913517355244)
  (8, 0.0030810671961455424)
  (9, 0.0030580332058320397)
  (10, 0.025905534724256826)
  (11, 0.011662872385237635)
  (12, 0.010671086023680734)
  (13, 0.01363027573888456)
  (14, 0.03553438453440785)
  (15, 0.05027608248858741)
  (16, 0.022552049488891914)
  (17, 0.011461201328344203)
  (18, 0.18088349088677053)
  (19, 1.2064391084453907)
  (20, 0.7781504821561804)
  (21, 0.23450008621662435)
  (22, 0.22193421751962564)
  (23, 0.2483948939517252)
  (24, 0.5888650158300175)
  (25, 0.250230254939281)
  (26, 0.21640118108169587)
  (27, 0.18710431409928555)
  (28, 2.9908569042929516)
  (29, 1.2484289152489607)
  (30, 0.8945675516861037)
  (31, 0.9844145484332835)
  (32, 2.6599777032643024)
  (33, 4.03976676072982)
  (34, 2.0372434683392004)
  (35, 1.141133197649997)
  (36, 19.925215988485437)
  (37, 94.82723834812833)
  (38, 26.51234359860197)
  (39, 14.24086839472968)
  (40, 14.599969905463771)
  (41, 14.309140091110487)
  (42, 26.325110133029)
  (43, 9.987386649139582)
  (44, 9.317250904206125)
  (45, 9.151891245450324)
  (46, 76.73402240124068)
  (47, 38.74800316227653)
  (48, 16.223919419971452)
  (49, 21.742495610747174)
  (50, 46.585816723001706)
  (51, 46.6192080877936)
  (52, 29.4729942417929)
  (53, 22.64312325716824)
};
\addplot[only marks, mark=star,mark size=1.8pt,blue!70!black] coordinates {
  (0, 0.0012324622683586859)
  (1, 0.001344300828754006)
  (2, 0.001244763414128827)
  (3, 0.004606403673138918)
  (4, 0.006922520711271506)
  (5, 0.003056310451974219)
  (6, 0.001244684162426872)
  (7, 0.0012305612136725977)
  (8, 0.001435126860598152)
  (9, 0.0013244619108594752)
  (10, 0.014691292119329756)
  (11, 0.005213211045738532)
  (12, 0.004969182661299276)
  (13, 0.008449167813687639)
  (14, 0.01688952423926194)
  (15, 0.024080077426895725)
  (16, 0.01142473128071688)
  (17, 0.005465529375169767)
  (18, 0.09472257683085968)
  (19, 0.7753478268282593)
  (20, 0.4267884128038074)
  (21, 0.09750447101508462)
  (22, 0.09810143218082433)
  (23, 0.1038316442546582)
  (24, 0.26488515847122773)
  (25, 0.08898291854528304)
  (26, 0.09539047958079412)
  (27, 0.11331472549609936)
  (28, 1.4671684292356826)
  (29, 0.5927872590097161)
  (30, 0.840590156786666)
  (31, 0.5768389334305912)
  (32, 1.8377995556642635)
  (33, 2.6954189396487287)
  (34, 1.2071290980893372)
  (35, 0.5055178430012907)
  (36, 9.659234223358268)
  (37, 75.90072496282441)
  (38, 40.46645327763583)
  (39, 10.68542617320052)
  (40, 10.705590392600772)
  (41, 10.69319455814933)
  (42, 24.853965799742795)
  (43, 11.613181441230294)
  (44, 9.832210396137588)
  (45, 11.385632062013011)
  (46, 100.20603570737852)
  (47, 58.353569386251635)
  (48, 31.07430743803105)
  (49, 39.40889983347103)
  (50, 74.84740998882782)
  (51, 73.18226829405495)
  (52, 49.67017554898764)
  (53, 36.43939242547901)
};
\end{axis}
\end{tikzpicture}%
\vspace{-0.2cm}
\input{\figdir/LV_composite}%
\caption{Per-callback and composite speedup over single-threaded CPU (C3) on Luk\v{s}an--Vl\v{c}ek instances, sorted by $\text{nnz}_J + \text{nnz}_H$. \textsuperscript{*}Float32 precision.}
\label{fig:lv_speedup}
\end{figure}

\providecommand{\figdir}{results/figures}

\providecommand{\figdir}{results/figures}
\begin{figure}[t]
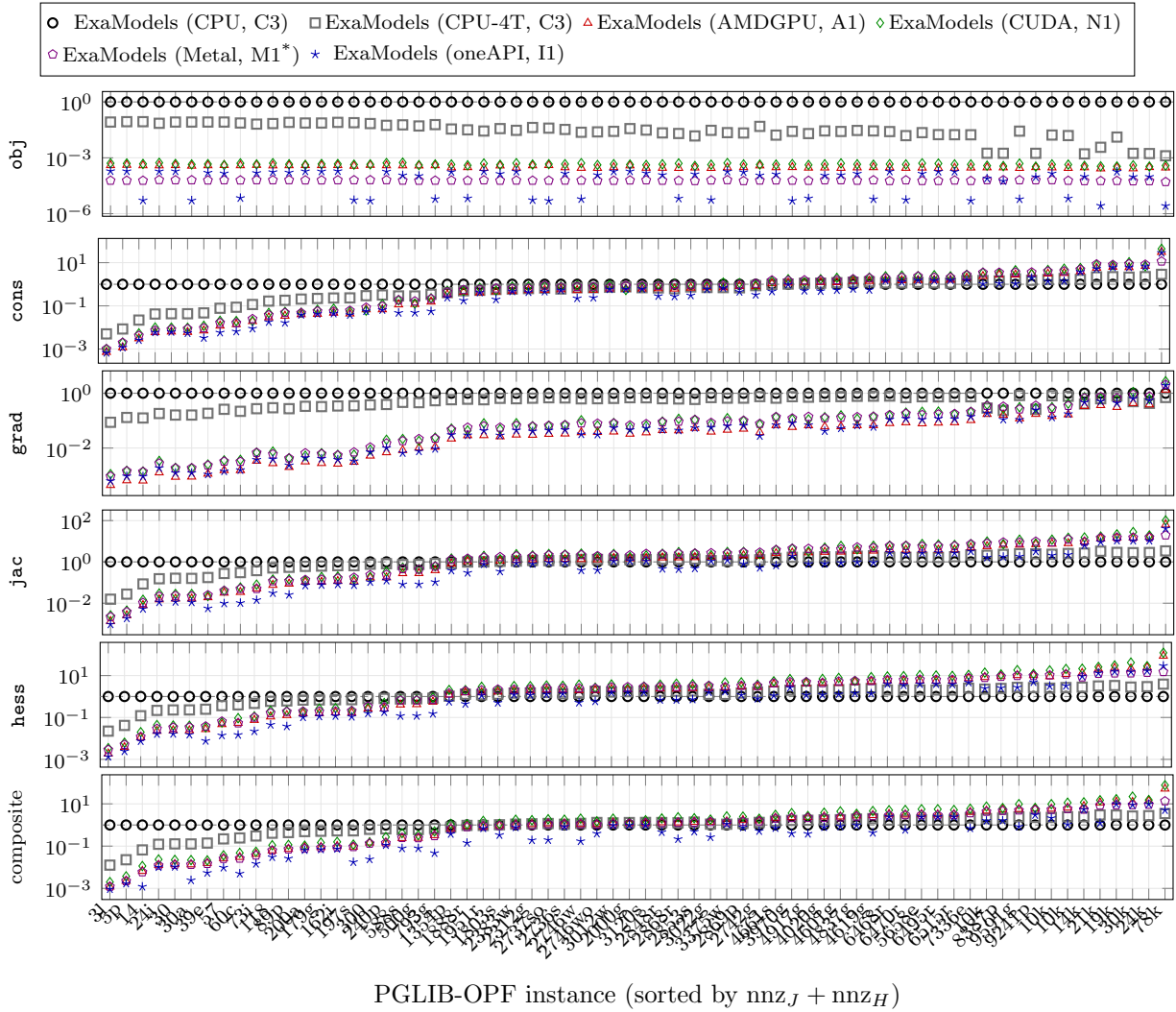

\centering
\vspace{0.1cm}
\pgfplotslegendfromname{opfpolar-legend}\\[2pt]
\input{\figdir/OPF_obj}%
\vspace{-0.2cm}
\input{\figdir/OPF_cons}%
\vspace{-0.2cm}
\input{\figdir/OPF_grad}%
\vspace{-0.2cm}
\input{\figdir/OPF_jac}%
\vspace{-0.2cm}
\input{\figdir/OPF_hess}%
\vspace{-0.2cm}
\input{\figdir/OPF_composite}%
\caption{Per-callback and composite speedup over single-threaded CPU (C3) on PGLIB-OPF (polar formulation) instances, sorted by $\text{nnz}_J + \text{nnz}_H$. \textsuperscript{*}Float32 precision.}
\label{fig:opf_speedup}
\end{figure}

\subsection{Per-Instance GPU Scaling}
\label{sec:numerics:comparison}

The class-level summaries above aggregate over instances; this subsection follows single instances as they grow, on the \gls{cpu} and CUDA backends of platform N6.
The first ladder is the Luk\v{s}an--Vl\v{c}ek Rosenbrock problem at $N = \cmpLvSizeA$, $\cmpLvSizeB$, and $\cmpLvSizeC$; the second is the \gls{acopf} polar formulation at \texttt{case118\_ieee}, \texttt{case1354\_pegase}, and \texttt{case78484\_epigrids}.

The ExaModels.jl \gls{gpu} Hessian takes \cmpLvHessGpuExaA, \cmpLvHessGpuExaB, and \cmpLvHessGpuExaC{} over the Rosenbrock ladder, and the Jacobian takes \cmpLvJacGpuExaA, \cmpLvJacGpuExaB, and \cmpLvJacGpuExaC.
The times grow with $N$, while the problem grows by four orders of magnitude.
The growth shows the saturation boundary of the constant-time argument of \Cref{sec:numerics:ad:gpu}: this ladder runs on the \platNSixName{} of platform N6, whose reported FP64 peak is $1.9$\,TFLOPS~\citep{nvidiaRTXProBlackwellArch2025}, against $40$\,TFLOPS for the \platNOneName{} of the summary tables~\citep{nvidiaHGXB200Datasheet2024}, where the Hessian stays within \gsLvHessAbsSmall{}--\gsLvHessAbsLarge{}\,\textmu s across the size classes.
The launch cost makes the \gls{cpu} the faster device at $N = \cmpLvSizeA$, and the objective callback, a single scalar reduction, gains little from the \gls{gpu} at any size.

A callback issues one kernel per algebraic pattern it touches, and one more to zero its output buffer.
The polar formulation has one objective pattern and fourteen constraint patterns.
The Hessian touches all fifteen patterns, so \texttt{NLPModels.hess\_coord!} issues sixteen kernels.
The objective does not enter the Jacobian, so \texttt{NLPModels.jac\_coord!} issues fifteen.
These counts are read from the model and are exact.
Dividing each measured \gls{gpu} time by its kernel count gives a per-kernel cost.
For the Hessian it is \cmpOpfPerKernelHessA, \cmpOpfPerKernelHessB, and \cmpOpfPerKernelHessC{} at the three \gls{acopf} cases.
For the Jacobian it is \cmpOpfPerKernelJacA, \cmpOpfPerKernelJacB, and \cmpOpfPerKernelJacC.
The Jacobian per-kernel cost varies by at most $10\%$ across those three cases, while the network grows by a factor of $665$ in bus count.
The Hessian per-kernel cost is flat between the two smaller cases and roughly doubles at the largest.
The kernel count therefore accounts for the Jacobian cost across this range, and does not account for the Hessian cost at the largest case.
That the fixed part is the kernel launch is an inference, supported by its magnitude and by the Jacobian, the cheaper kernel, coming out below the Hessian at every size.
The floor is a property of the granularity of the abstraction, one kernel per pattern, and fusing the patterns of a model into a single generated kernel, or issuing them on concurrent streams, would lower it.
The following section shows, however, that lowering it further carries little practical benefit, since \gls{nlp} function evaluation is already a small share of a complete solve.

\subsection{Measuring the Benefits of GPU-Resident Callbacks}
\label{sec:numerics:resident}

\Cref{sec:numerics:ad,sec:numerics:comparison} measure callbacks in isolation.
What a user experiences is a solve, in which the callbacks are one term among several.
Speeding up the callbacks can therefore save only the share of the solve that the callbacks occupy.
We run one interior-point method on two sets of instances and move only the device on which each of its two dominant costs is paid.

\subsubsection{Solver Configurations and Reported Quantities}
\label{sec:numerics:resident:setup}

Two of the four configurations pair a model on one device with a solver on the other.
A user reaches for such a configuration when only one half of the solution procedure has been ported, and its cost is what the experiment is meant to expose.
The bridge is \texttt{WrapperNLPModel}, provided by ExaModels.jl, which presents a model living on one device to a solver living on the other and forwards every callback through intermediate copy buffers.

The solver is MadNLP.jl~\citep{shinMadNLPjl2025} with the LiftedKKT formulation~\citep{pacaudCondensedspaceMethodsNonlinear2024}, which relaxes the equality constraints and condenses the KKT system so that the resulting matrix admits a sparse Cholesky-like factorization on a \gls{gpu}.
MadNLP.jl itself runs on either device, and the linear solver is chosen to match: LDLFactorizations.jl runs only on the \gls{cpu}, and NVIDIA cuDSS only on the \gls{gpu}.
Each configuration is named by the device of the \gls{nlp} solver and the device of the model.
The four configurations are as close to one another as two devices allow.
Each solves the same problem with the same formulation at the same tolerance, and each factorizes the same condensed KKT system with a limited-pivoting LDL factorization.
LDLFactorizations.jl is chosen as the \gls{cpu} linear solver because it has limited pivoting, similar to cuDSS, because it is purely single-threaded, and because it performs an LDL rather than an LBL factorization.
The \gls{cpu} and \gls{gpu} configurations therefore run essentially the same algorithm, and the comparison isolates where the solver and the model run rather than which algorithm runs.
A better \gls{cpu} alternative exists, solving the augmented system with an LBL factorization with pivoting, but the algorithms are deliberately synchronized so that the experiment measures the effect of running the \gls{nlp} solver and the \gls{nlp} functions on the \gls{cpu} against the \gls{gpu}.
The augmented route is preferred on the \gls{cpu} because pivoting handles the indefinite KKT system stably without the conditioning penalty that condensation introduces; condensation is the price of a pivoting-free factorization, which the \gls{gpu} requires.
The Julia code that runs is nearly identical across the configurations.
What differs is the array types it is dispatched on, and hence the low-level kernels those dispatches select, e.g., BLAS operations, map and map-reduce operations, and the \gls{nlp} function evaluation kernels.
The comparison is therefore between the same software on different hardware, together with the bridge between the two.

\Cref{tab:breakdown,fig:breakdown} separate the \gls{nlp} function evaluation out of the solve, wherever that quantity can be measured without disturbing it.
The evaluation column of \Cref{tab:breakdown} is the \gls{nlp} function evaluation time contained in the solution time.
In the fully \gls{gpu}-resident configuration a callback only launches its kernels and nothing forces it to complete, so its evaluation time cannot be measured without altering the solve, and that entry is left blank.
All four configurations were measured on platform~N6, so the rows differ only in where the model and the \gls{nlp} solver run.
Iteration counts are reported alongside the times because they are not identical across configurations.
The \gls{cpu} and the \gls{gpu} evaluate the same expressions in different orders and therefore round differently, so the iterate sequences diverge slightly.

\subsubsection{Solve Time and Evaluation Share}
\label{sec:numerics:resident:times}

\begin{table}[t]
\centering
\caption{Solution time for one interior-point solve, across the four configurations of \Cref{sec:numerics:resident}.}
\label{tab:breakdown}
\scriptsize
\begin{tabular*}{\textwidth}{@{\extracolsep{\fill}}llrrrr}
  \toprule
  Solver & Model & solution (s) & evaluation (s) & total (s) & it. \\
  \midrule
  \multicolumn{6}{@{}l}{\emph{Rosenbrock ($N = 200{,}000$)}} \\
  LDLFactorizations (CPU) & CPU & 0.889 & 0.154 & 0.903 & 7 \\
  LDLFactorizations (CPU) & GPU & 0.732 & 0.051 & 0.763 & 7 \\
  cuDSS (GPU) & CPU & 0.496 & 0.174 & 0.523 & 7 \\
  cuDSS (GPU) & GPU & 0.293 & -- & 0.308 & 7 \\
  \addlinespace
  \multicolumn{6}{@{}l}{\emph{PGLIB-OPF \texttt{case78484\_epigrids}}} \\
  LDLFactorizations (CPU) & CPU & 147 & 9.91 & 147 & 113 \\
  LDLFactorizations (CPU) & GPU & 142 & 3.52 & 142 & 114 \\
  cuDSS (GPU) & CPU & 15.2 & 7.90 & 15.4 & 114 \\
  cuDSS (GPU) & GPU & 7.37 & -- & 7.48 & 113 \\
  \bottomrule
\end{tabular*}

\end{table}

\providecommand{\figdir}{results/figures}
\begin{figure}[t]
\centering
\resizebox{\textwidth}{!}{%
\begin{tikzpicture}
\pgfplotsset{tick scale binop=\times}
\begin{groupplot}[
  group style={group size=3 by 2, horizontal sep=0.9cm, vertical sep=1.15cm},
  width=4.9cm, height=3.7cm,
  ybar stacked, /pgf/bar width=9pt,
  xtick={1,2,3,4},
  xticklabels={C/C,C/G,G/C,G/G},
  xmin=0.4, xmax=4.6,
  x tick label style={font=\tiny},
  yticklabel style={font=\tiny},
  title style={font=\scriptsize, yshift=-2pt},
  ylabel style={font=\scriptsize},
  legend columns=4,
  legend style={draw=black, fill=white, font=\small},
  ymin=0,
  enlarge x limits=0.18,
]
\nextgroupplot[title={$N = 20$},ylabel={wall time (s)},legend to name=breakdownlegend,]
\addplot+[fill=Goldenrod!75, draw=black!60] coordinates {(1,1.7573e-05) (2,0.00185622) (3,0.0016851) (4,0)};
\addplot+[fill=NavyBlue!70, draw=black!65] coordinates {(1,0.000105419) (2,0.000766023) (3,0.0294228) (4,0)};
\addplot+[fill=Plum!55, draw=black!60] coordinates {(1,0) (2,0) (3,0) (4,0.0261143)};
\addplot+[fill=black!20, draw=black!60] coordinates {(1,2.7094e-05) (2,0.00104216) (3,0.000889105) (4,0.00119279)};
\legend{NLP function evaluation, rest of solve, evaluation and solve, model creation}
\nextgroupplot[title={$N = 2{,}000$},]
\addplot+[fill=Goldenrod!75, draw=black!60] coordinates {(1,0.00158606) (2,0.00318143) (3,0.00388878) (4,0)};
\addplot+[fill=NavyBlue!70, draw=black!65] coordinates {(1,0.00825513) (2,0.00679347) (3,0.0411224) (4,0)};
\addplot+[fill=Plum!55, draw=black!60] coordinates {(1,0) (2,0) (3,0) (4,0.0362354)};
\addplot+[fill=black!20, draw=black!60] coordinates {(1,0.000188221) (2,0.0018317) (3,0.0014283) (4,0.0015964)};
\nextgroupplot[title={$N = 200{,}000$},]
\addplot+[fill=Goldenrod!75, draw=black!60] coordinates {(1,0.154072) (2,0.0506939) (3,0.174225) (4,0)};
\addplot+[fill=NavyBlue!70, draw=black!65] coordinates {(1,0.734725) (2,0.681002) (3,0.322219) (4,0)};
\addplot+[fill=Plum!55, draw=black!60] coordinates {(1,0) (2,0) (3,0) (4,0.293043)};
\addplot+[fill=black!20, draw=black!60] coordinates {(1,0.0146549) (2,0.0311681) (3,0.0263698) (4,0.0146903)};
\nextgroupplot[title={1{,}354 buses},ylabel={wall time (s)},]
\addplot+[fill=Goldenrod!75, draw=black!60] coordinates {(1,0.0404767) (2,0.0921019) (3,0.0606242) (4,0)};
\addplot+[fill=NavyBlue!70, draw=black!65] coordinates {(1,0.370225) (2,0.371056) (3,0.223621) (4,0)};
\addplot+[fill=Plum!55, draw=black!60] coordinates {(1,0) (2,0) (3,0) (4,0.280172)};
\addplot+[fill=black!20, draw=black!60] coordinates {(1,0.00162626) (2,0.00664672) (3,0.0036621) (4,0.00588917)};
\nextgroupplot[title={9{,}241 buses},]
\addplot+[fill=Goldenrod!75, draw=black!60] coordinates {(1,0.484059) (2,0.30149) (3,0.558543) (4,0)};
\addplot+[fill=NavyBlue!70, draw=black!65] coordinates {(1,6.52802) (2,5.68002) (3,0.908023) (4,0)};
\addplot+[fill=Plum!55, draw=black!60] coordinates {(1,0) (2,0) (3,0) (4,0.988002)};
\addplot+[fill=black!20, draw=black!60] coordinates {(1,0.00904555) (2,0.0233128) (3,0.0167178) (4,0.0159066)};
\nextgroupplot[title={78{,}478 buses},]
\addplot+[fill=Goldenrod!75, draw=black!60] coordinates {(1,9.90882) (2,3.52166) (3,7.89514) (4,0)};
\addplot+[fill=NavyBlue!70, draw=black!65] coordinates {(1,137.14) (2,138.389) (3,7.34783) (4,0)};
\addplot+[fill=Plum!55, draw=black!60] coordinates {(1,0) (2,0) (3,0) (4,7.36853)};
\addplot+[fill=black!20, draw=black!60] coordinates {(1,0.0819351) (2,0.182482) (3,0.154644) (4,0.109136)};
\end{groupplot}
\end{tikzpicture}%
}
\par\vspace{3pt}
\pgfplotslegendfromname{breakdownlegend}
\caption{Breakdown of the solution time across the four configurations of \Cref{tab:breakdown}: Luk\v{s}an--Vl\v{c}ek Rosenbrock on the top row and PGLIB-OPF on the bottom, at three instances each.}
\label{fig:breakdown}
\end{figure}
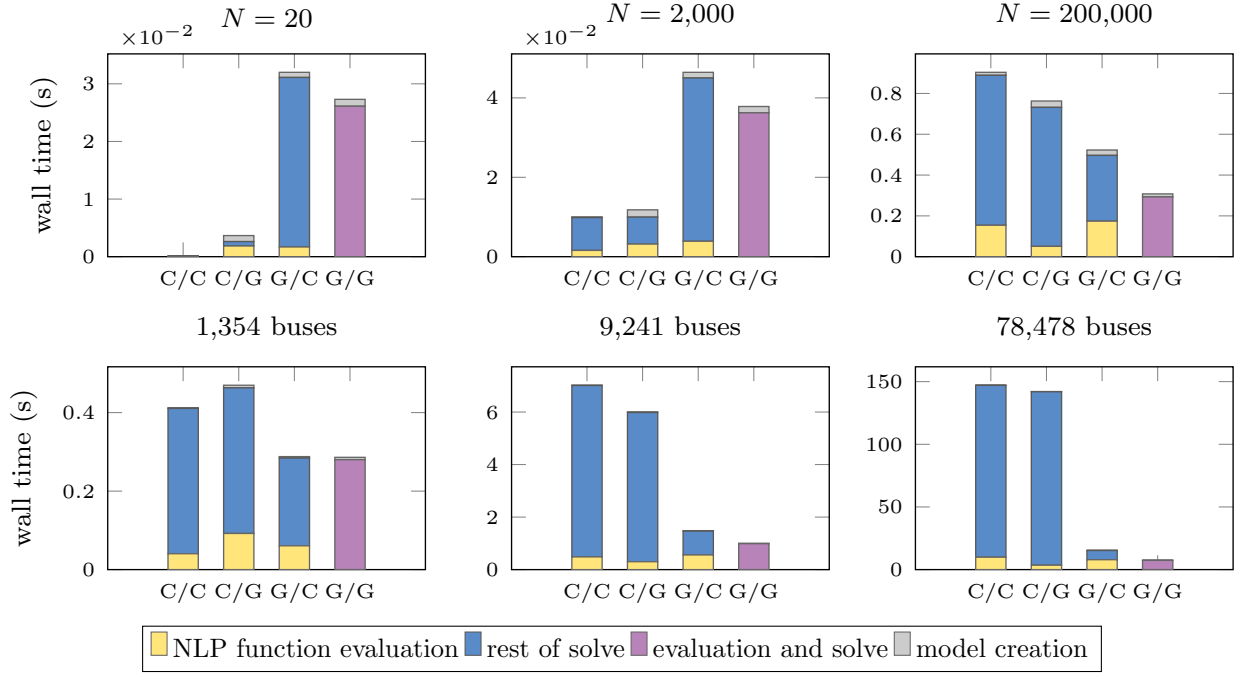

\Cref{tab:breakdown} reports the largest instance of each problem.
When the \gls{nlp} solver runs on the \gls{cpu}, the factorization dominates the solution time, and the device on which the model is evaluated has little effect: making the model \gls{gpu}-resident changes the \bdOpfCase{} solution time by a factor of \bdOpfSpeedupResidentCpuLs{}.
Moving the \gls{nlp} solver to the \gls{gpu} reduces the solution time from \bdOpfTotCC{}\,s to \bdOpfTotGC{}\,s, a factor of \bdOpfSpeedupLinalg{}.
Making the model \gls{gpu}-resident reduces the solution time further, from \bdOpfTotGC{}\,s to \bdOpfTotGG{}\,s, a factor of \bdOpfSpeedupResident{}.
The overall reduction against the all-\gls{cpu} configuration is a factor of \bdOpfSpeedupTotal{}.
On the Rosenbrock instance, the two contributions are comparable: moving the \gls{nlp} solver to the \gls{gpu} gains a factor of \bdRosSpeedupLinalg{}, and making the model \gls{gpu}-resident gains a factor of \bdRosSpeedupResident{}.

The evaluation share of the solve behaves differently on the two devices.
\Cref{fig:breakdown} shows the breakdown across instance sizes, and the largest instances show the contrast most clearly.
Each bar group is labeled by the device of the \gls{nlp} solver over the device of the model evaluation, and each panel carries its own linear axis.
The fully \gls{gpu}-resident configuration (\texttt{G/G}) is drawn as a single segment, since its evaluation share is the quantity that cannot be measured without altering the solve.
With the solver on the \gls{cpu}, evaluation is a small part of the solve, \bdOpfEvalShareCC{}\% of the all-\gls{cpu} configuration at \bdOpfCase{}.
The factorization dominates so completely there that accelerating the callbacks could not repay the transfers and synchronization it would cost.
Moving the solver to the \gls{gpu} removes that dominance and leaves evaluation as the exposed term.
With the model still on the host, evaluation accounts for \bdOpfEvalShareGC{}\% of the same solve, and every callback pays for a transfer across the bus and for the completion barrier the transfer implies.
Making the model \gls{gpu}-resident removes both at once, which is why that configuration is the one where \gls{gpu} acceleration of the callbacks actually pays.

The configuration pairing a \gls{gpu}-resident model with a solver on the \gls{cpu} is the one that does not pay.
The callbacks themselves do become faster, which the evaluation column of \Cref{tab:breakdown} shows.
The solve does not follow, because the function evaluation is a small share of the total with a \gls{cpu} solver, so cutting it changes the total by little.
The achievable improvement is further limited by the bridge: each iteration must synchronize and transfer data between the host and the device, as every callback copies its arguments to the device and its results back.
The practical reading is that the benefit of \gls{gpu}-resident callbacks depends on the rest of the solution procedure already being on the device, and that porting the modeling layer alone is not a shortcut to it.

\subsection{Model Compilation Time}
\label{sec:numerics:compilation}

The preceding subsections measure warm evaluation; this subsection measures the compilation cost paid before it.
ExaModels.jl compiles specialized code separately for model creation and for each of the five \gls{nlp} function callbacks of \Cref{sec:numerics:setup}, one set of kernels per algebraic pattern.
The reported time covers all of it: the model creation and the compilation of every kernel, not of a single one.
\Cref{tab:compile} reports Julia's own accounting of that time on the \gls{cpu} and CUDA backends, during model construction and additionally through the first derivative evaluation and solve, varying the problem size at a fixed set of constraint kernels and the number of distinct constraint kernels at a fixed size.
The constraint kernels are sixteen synthetic algebraic patterns over consecutive variable pairs, such as $x_i^2 + x_{i+1}$, $\sin(x_i)\,x_{i+1}$, $e^{x_i - x_{i+1}}$, $\sqrt{x_i^2 + 1} + x_{i+1}$, and $\tanh(x_i) - x_{i+1}$, each a distinct expression compiled as its own pattern.
The kernel-count sweep splits the same total constraint count among the first $k$ patterns, so the problem size stays fixed while only the number of distinct patterns varies.
In the table, \texttt{model} is the compilation during model construction, \texttt{total} additionally includes the first derivative evaluation and the solve, and times are in milliseconds unless a unit is shown.

On both backends the compilation time is flat in the problem size, within \cmpSizeSpreadGpu\% on the \gls{gpu} across a factor of $1{,}000$, and grows with the number of distinct constraint kernels, by $\cmpKindsRatioCpu\times$ on the \gls{cpu} and $\cmpKindsRatioGpu\times$ on the \gls{gpu} from \cmpKindsLoGpu{} kernel to \cmpKindsHiGpu{}.
The \gls{gpu} increment is proportionally smaller because the measured time on the \gls{gpu} includes compilation beyond the ExaModels.jl kernels, in CUDA.jl and its supporting compilation stack.
That fixed cost, about $\cmpGpuFixedS$\,s, dominates the per-kernel cost, so the growth is masked; the per-kernel growth is present on both backends.

\begin{table}[t]
\caption{Julia compilation time for ExaModels.jl, varying the problem size and the number of distinct constraint kernels. Both sweeps ran on platform N6.}
\label{tab:compile}
\scriptsize
\begin{tabular*}{\textwidth}{@{\extracolsep{\fill}}llrrrrr}
  \toprule
 &  & \multicolumn{2}{c}{\textbf{CPU}} & \multicolumn{2}{c}{\textbf{GPU}} & \\
  \textbf{varying} & \textbf{value} & model & total & model & total & \\
  \midrule
  \multicolumn{7}{@{}l}{\emph{Julia compilation, ExaModels.jl: problem size $n$}} \\
   & 1000 & 318 & 4.7\,s & 16.2\,s & 67.0\,s & \\
   & 10000 & 339 & 5.0\,s & 16.1\,s & 66.7\,s & \\
   & 100000 & 324 & 4.8\,s & 16.3\,s & 66.9\,s & \\
   & 1000000 & 371 & 5.0\,s & 16.1\,s & 66.2\,s & \\
  \multicolumn{7}{@{}l}{\emph{Julia compilation, ExaModels.jl: constraint kernels}} \\
   & 1 & 357 & 5.0\,s & 15.9\,s & 66.7\,s & \\
   & 2 & 504 & 5.2\,s & 16.2\,s & 67.7\,s & \\
   & 4 & 809 & 5.5\,s & 16.9\,s & 71.7\,s & \\
   & 8 & 1.5\,s & 6.4\,s & 17.4\,s & 74.6\,s & \\
   & 16 & 2.8\,s & 8.0\,s & 19.2\,s & 80.9\,s & \\
  \bottomrule
\end{tabular*}
\end{table}

Compilation is paid per algebraic pattern rather than per data point, so a model with few patterns compiles in seconds however large its data, and the cost scales only when the model gains structurally new constraints.
The cost is also paid once per pattern set: recreating the same model, or creating it at a different size, triggers no further compilation.
An \gls{acopf} model compiled on a three-bus network can be recreated for the 78{,}484-bus network at no additional compilation cost, since every instance of the formulation shares the same algebraic patterns (\Cref{sec:numerics:problems}).

\section{Conclusions}
\label{sec:conclusions}

We have presented ExaModels.jl, an \gls{ams} that exploits the partially separable, repetitive structure of large-scale \glspl{nlp} to provide scalable \gls{nlp} function evaluation and \gls{gpu} acceleration.
Two ideas underpin the system.
The first is a \emph{\gls{simd} abstraction} that represents an \gls{nlp} as a small number of algebraic patterns, each evaluated over many data points.
It is realized through a parameterized expression tree that is fully encoded in the Julia type system.
Because the problem structure is recognized at compile time, \gls{nlp} function evaluation reduces to embarrassingly parallel loops over data points, so the same code runs natively on \glspl{gpu}.
The second is a \emph{coloring-free sparse \gls{ad} algorithm}.
It performs the sparsity analysis on the parameterized expression tree itself, and assembles the Jacobian and Hessian directly into a partially compressed \gls{coo} format.

The parameterized expression tree is what makes the \gls{gpu} acceleration possible.
Each algebraic pattern compiles to a single kernel, so every \gls{nlp} function callback becomes a data-parallel loop over the data points, with no device code in the model.
Data-parallel execution on \glspl{gpu} yields up to $\gsLvHessLarge\times$ over single-threaded \gls{cpu} evaluation, on the largest Luk\v{s}an--Vl\v{c}ek instances.
\gls{nlp} function evaluation also takes constant time.
Once the device supplies enough threads, the cost is set by the number of algebraic patterns and not by the number of data points, so it does not grow as the problem grows.
This is the $O(1)$ behavior the \gls{simd} abstraction was designed to deliver, and \Cref{sec:numerics:ad} demonstrates it.
For the class of problems considered here, \gls{nlp} function evaluation on a \gls{gpu} is therefore effectively instantaneous, and what remains is the runtime efficiency of the optimization solver itself.

The implications of this work are as follows.
\begin{itemize}
  \item \emph{For modelers and modeling library developers.}
\gls{gpu} compatibility is the benefit worth pursuing, and either templatization or vectorization delivers it.
A model that runs on the device makes a fully device-resident framework possible, and that is where the large computational gains are.
The remaining differences are largely a matter of preference.
Templatization keeps the constraint and objective syntax that mathematical programming languages have long used, and it can avoid the overhead of graph coloring.
Vectorization integrates seamlessly with the surrounding machine learning workflows.
  \item \emph{For solver developers.}
For classical applications, the major effort now belongs in reducing the number of factorizations the algorithm requires.
Function evaluation is negligible beside them, especially on a \gls{gpu}, since its cost does not grow with the number of data points.
A solver can therefore assume in-place device-resident callbacks, exactly as it already assumes them on the \gls{cpu}, and spend its effort on the linear algebra and on matching the reliability of \gls{cpu} solvers.
Progress there, in either reliability or efficiency, will make the results of this work more useful in practice.
\end{itemize}

\section*{Acknowledgments}

ExaModels.jl is built on the broader Julia ecosystem, and we are indebted to several communities within it.
The cross-platform GPU portability described in \Cref{sec:gpu} rests entirely on the JuliaGPU stack.
We thank Valentin Churavy and the developers of KernelAbstractions.jl for the backend-agnostic kernel layer.
We also thank Tim Besard and the wider JuliaGPU community for CUDA.jl, AMDGPU.jl, oneAPI.jl, Metal.jl, and the GPUCompiler.jl infrastructure, which lets a single ExaModels.jl codebase run unmodified on NVIDIA, AMD, Intel, and Apple hardware.
The \gls{aot} compilation capability of \Cref{sec:aot} relies on Julia's static compilation and trimming effort.
We thank the developers of \texttt{juliac} and the Julia compiler team for making it possible to package complete optimization applications as standalone binaries.
ExaModels.jl also connects to a broad solver ecosystem through the NLPModels.jl callback interface.
We thank Dominique Orban and the JuliaSmoothOptimizers community for NLPModels.jl and the surrounding optimization tooling.

We thank the Argonne Leadership Computing Facility for providing access to the Joint Laboratory for System Evaluation, on which some of the benchmark experiments in \Cref{sec:numerics} were performed.
We thank MIT's Office of Research Computing and Data for the computing resources on which the remaining benchmark experiments were performed.
This work was supported by the U.S. Department of Energy, Office of Science (SC), Advanced Scientific Computing Research (ASCR), Competitive Portfolios Project on Energy Efficient Computing: A Holistic Methodology, under Contract DE-AC02-06CH11357.

\printbibliography

\appendix

\section{Compiled LLVM Code for a Gradient Kernel}
\label{app:llvm}

This appendix reproduces the LLVM code that the Julia compiler generates for an ExaModels.jl gradient kernel.
The listing supports the claim of \Cref{sec:intro:related}: in a language built on multiple dispatch and \gls{jit} compilation, operator overloading can produce derivative evaluation code as free of overhead as source-code transformation.
The model is the extended Rosenbrock objective \cref{eq:rosenbrock}, declared over a data iterator as \texttt{@add\_obj(core, 100*(x[i-1]\^{}2 - x[i])\^{}2 + (x[i-1]-1)\^{}2 for i in 2:N)}.
The function lowered is the per-data-point gradient kernel of \Cref{sec:ad:first}.
It receives one data point \texttt{p}, computes the gradient contribution of that point, and accumulates it into the output array.
The kernel depends only on the argument types, not on the data iterator contents, so the same code serves every data point and every $N$.
The listing was produced with \texttt{@code\_llvm debuginfo=:none} on Julia 1.12.6 with ExaModels.jl v0.11.2.
The annotation comments, marked \texttt{;;}, are ours, and everything else is reproduced verbatim.

In the signature, the second argument is the parameterized expression tree, and its fully concrete type spells out the entire algebraic structure of the pattern (\Cref{sec:graph:type}).
In the body, the global indices of \texttt{x[i-1]} and \texttt{x[i]} are computed from the data point \texttt{p} and the offsets stored in the \texttt{Var} leaves.
The body is a single basic block of straight-line code: the forward evaluation and the reverse pass are fully inlined, and there are no function calls, no branches, no dynamic dispatch, and no memory allocation.

The generated code is also nearly minimal in its operation count.
The body contains 13 floating-point operations, namely 8 multiplications, 3 subtractions, and 2 additions.
The two \texttt{sitofp} instructions convert the integer constants embedded in the tree and involve no arithmetic on the data.
For comparison, write $u = x_{i-1}$, $v = x_i$, $a$ for the seed adjoint, and $y$ for the gradient output.
The same computation, both gradient contributions and their accumulation into $y$, can be hand-written with 12 operations:
\begin{align*}
t &= u \cdot u, & s &= t - v, & h &= a \cdot s, & m &= 200\, h,\\
p &= u \cdot m, & p_2 &= p + p, & r &= u - 1, & z &= a \cdot r,\\
z_2 &= z + z, & q &= p_2 + z_2, & y_{i-1} &\leftarrow y_{i-1} + q, & y_{i} &\leftarrow y_{i} - m,
\end{align*}
that is, 5 multiplications, 3 subtractions, and 4 additions.
We are not aware of a shorter sequence.
The generated kernel is thus one operation off the best hand-written count we could construct.
The difference is a factoring choice rather than interpretive overhead.
The kernel forms the local derivatives $2u$, $2s$, and $2(u-1)$ as multiplications during the forward pass, whereas the hand-written sequence replaces two of them with the doublings $p + p$ and $z + z$.
The listing shows the code shape that source-code transformation aims to produce, obtained here from the operator-overloading construction.
The only code generation involved is the Julia compiler's own lowering to the LLVM intermediate representation, which any compiled language performs.
There is no separate source-level code-generation step.

\begin{lstlisting}[style=llvm]
; Function Signature: gradient!(Array{Float64, 1}, ExaModels.Node2{typeof(Base.:(+)), ExaModels.Node2{typeof(Base.:(*)), Int64, ExaModels.Node1{typeof(Base.abs2), ExaModels.Node2{typeof(Base.:(-)), ExaModels.Node1{typeof(Base.abs2), ExaModels.Var{ExaModels.Node2{typeof(Base.:(+)), ExaModels.Node2{typeof(Base.:(-)), ExaModels.DataSource, Int64}, Int64}}}, ExaModels.Var{ExaModels.Node2{typeof(Base.:(+)), ExaModels.DataSource, Int64}}}}}, ExaModels.Node1{typeof(Base.abs2), ExaModels.Node2{typeof(Base.:(-)), ExaModels.Var{ExaModels.Node2{typeof(Base.:(+)), ExaModels.Node2{typeof(Base.:(-)), ExaModels.DataSource, Int64}, Int64}}, Int64}}}, Array{Float64, 1}, Array{Float64, 1}, Int64, Float64)
define nonnull ptr @"julia_gradient!_3566"(ptr noundef nonnull align 8 dereferenceable(24) %"y::Array", ptr nocapture noundef nonnull readonly align 8 dereferenceable(56) %"f::Node2", ptr noundef nonnull align 8 dereferenceable(24) %"x::Array", ptr noundef nonnull align 8 dereferenceable(24) %"\CE\B8::Array", i64 signext %"p::Int64", double %"adj::Float64") #0 {
top:

;; forward pass: compute the global index of x[i-1] from the data point p
;; and the offsets in the Var leaf; load it; evaluate abs2 and its derivative 2*x[i-1]
  %"f::Node2.inner2_ptr" = getelementptr inbounds i8, ptr %"f::Node2", i64 8
  %"f::Node2.inner2_ptr.unbox" = load i64, ptr %"f::Node2.inner2_ptr", align 8
  %0 = sub i64 %"p::Int64", %"f::Node2.inner2_ptr.unbox"
  %"f::Node2.inner2_ptr.inner2_ptr" = getelementptr inbounds i8, ptr %"f::Node2", i64 16
  %"f::Node2.inner2_ptr.inner2_ptr.unbox" = load i64, ptr %"f::Node2.inner2_ptr.inner2_ptr", align 8
  %1 = add i64 %0, %"f::Node2.inner2_ptr.inner2_ptr.unbox"
  %memoryref_data = load ptr, ptr %"x::Array", align 8
  %memoryref_offset = shl i64 %1, 3
  %2 = getelementptr i8, ptr %memoryref_data, i64 %memoryref_offset
  %memoryref_data1 = getelementptr i8, ptr %2, i64 -8
  %3 = load double, ptr %memoryref_data1, align 8
  %4 = fmul double %3, %3
  %5 = fmul double %3, 2.000000e+00

;; likewise locate and load x[i]; evaluate x[i-1]^2 - x[i] and its local derivative
  %"f::Node2.inner2_ptr.inner2_ptr2" = getelementptr inbounds i8, ptr %"f::Node2", i64 24
  %"f::Node2.inner2_ptr.inner2_ptr2.unbox" = load i64, ptr %"f::Node2.inner2_ptr.inner2_ptr2", align 8
  %6 = add i64 %"f::Node2.inner2_ptr.inner2_ptr2.unbox", %"p::Int64"
  %memoryref_offset9 = shl i64 %6, 3
  %7 = getelementptr i8, ptr %memoryref_data, i64 %memoryref_offset9
  %memoryref_data15 = getelementptr i8, ptr %7, i64 -8
  %8 = load double, ptr %memoryref_data15, align 8
  %9 = fsub double %4, %8
  %10 = fmul double %9, 2.000000e+00

;; load the constant coefficient 100 from the tree
  %"f::Node2.unbox" = load i64, ptr %"f::Node2", align 8
  %11 = sitofp i64 %"f::Node2.unbox" to double

;; second term: locate and load x[i-1] and the constant 1; evaluate x[i-1] - 1
  %"f::Node2.inner2_ptr16" = getelementptr inbounds i8, ptr %"f::Node2", i64 32
  %"f::Node2.inner2_ptr16.unbox" = load i64, ptr %"f::Node2.inner2_ptr16", align 8
  %12 = sub i64 %"p::Int64", %"f::Node2.inner2_ptr16.unbox"
  %"f::Node2.inner2_ptr16.inner2_ptr" = getelementptr inbounds i8, ptr %"f::Node2", i64 40
  %"f::Node2.inner2_ptr16.inner2_ptr.unbox" = load i64, ptr %"f::Node2.inner2_ptr16.inner2_ptr", align 8
  %13 = add i64 %12, %"f::Node2.inner2_ptr16.inner2_ptr.unbox"
  %memoryref_offset23 = shl i64 %13, 3
  %14 = getelementptr i8, ptr %memoryref_data, i64 %memoryref_offset23
  %memoryref_data29 = getelementptr i8, ptr %14, i64 -8
  %15 = load double, ptr %memoryref_data29, align 8
  %"f::Node2.inner2_ptr16.inner2_ptr30" = getelementptr inbounds i8, ptr %"f::Node2", i64 48
  %"f::Node2.inner2_ptr16.inner2_ptr30.unbox" = load i64, ptr %"f::Node2.inner2_ptr16.inner2_ptr30", align 8
  %16 = sitofp i64 %"f::Node2.inner2_ptr16.inner2_ptr30.unbox" to double
  %17 = fsub double %15, %16
  %18 = fmul double %17, 2.000000e+00

;; reverse pass: propagate the seed adjoint down the tree by the chain rule
  %19 = fmul double %11, %"adj::Float64"
  %20 = fmul double %10, %19
  %21 = fmul double %5, %20

;; accumulate the first term's contributions into y[i-1] and y[i]
  %memoryref_data32 = load ptr, ptr %"y::Array", align 8
  %22 = getelementptr i8, ptr %memoryref_data32, i64 %memoryref_offset
  %memoryref_data40 = getelementptr i8, ptr %22, i64 -8
  %23 = load double, ptr %memoryref_data40, align 8
  %24 = fadd double %21, %23
  store double %24, ptr %memoryref_data40, align 8
  %25 = getelementptr i8, ptr %memoryref_data32, i64 %memoryref_offset9
  %memoryref_data66 = getelementptr i8, ptr %25, i64 -8
  %26 = load double, ptr %memoryref_data66, align 8
  %27 = fsub double %26, %20
  store double %27, ptr %memoryref_data66, align 8

;; accumulate the second term's contribution into y[i-1]
  %28 = fmul double %18, %"adj::Float64"
  %29 = getelementptr i8, ptr %memoryref_data32, i64 %memoryref_offset23
  %memoryref_data92 = getelementptr i8, ptr %29, i64 -8
  %30 = load double, ptr %memoryref_data92, align 8
  %31 = fadd double %28, %30
  store double %31, ptr %memoryref_data92, align 8
  ret ptr %"y::Array"
}
\end{lstlisting}

\section{Per-Instance Benchmark Results}
\label{app:tables}

This appendix presents the complete per-instance benchmark data.
For each suite, we report the wall time in seconds for the callbacks \texttt{NLPModels.obj}, \texttt{NLPModels.cons!}, \texttt{NLPModels.grad!}, \texttt{NLPModels.jac\_coord!}, and \texttt{NLPModels.hess\_coord!}.
The backend column names the framework and, in parentheses, the platform of \Cref{tab:hardware} on which it ran.
An asterisk on a platform label marks timings taken in \texttt{Float32} arithmetic.
For each instance, the best time in each column is shown in bold.

{\scriptsize
\setlength{\tabcolsep}{2pt}
\setlength{\LTcapwidth}{\textwidth}
\LTleft=0pt plus 1fil minus 1fill
\LTright=0pt plus 1fil minus 1fill
\begin{longtable}{@{\extracolsep{\fill}}l | rrrr | l | rrrrr}
  \caption{Per-instance callback times (seconds) for the Luk\v{s}an--Vl\v{c}ek benchmark suite.}
  \label{tab:results:lv} \\
  \toprule
  \textbf{instance} & \textbf{\(n\)} & \textbf{\(m\)} & \textbf{nnzj} & \textbf{nnzh} & \textbf{backend} & \textbf{obj} & \textbf{cons} & \textbf{grad} & \textbf{jac} & \textbf{hess} \\
  \midrule
  \endfirsthead
  \caption[]{(continued)} \\
  \toprule
  \textbf{instance} & \textbf{\(n\)} & \textbf{\(m\)} & \textbf{nnzj} & \textbf{nnzh} & \textbf{backend} & \textbf{obj} & \textbf{cons} & \textbf{grad} & \textbf{jac} & \textbf{hess} \\
  \midrule
  \endhead
  augmented\_lagrangian-20 & 20 & 18 & 54 & 106 & CPU (C3) & \textbf{9.00e-08} & \textbf{1.10e-07} & \textbf{1.40e-07} & \textbf{1.60e-07} & \textbf{4.70e-07} \\
   &  &  &  &  & CPU-4T (C3) & 2.30e-07 & 1.80e-07 & 4.40e-07 & 2.80e-07 & 7.70e-07 \\
   &  &  &  &  & CPU-8T (C3) & 2.10e-07 & 1.90e-07 & 3.90e-07 & 2.90e-07 & 7.70e-07 \\
   &  &  &  &  & CPU-16T (C3) & 2.10e-07 & 1.90e-07 & 4.30e-07 & 2.80e-07 & 7.80e-07 \\
   &  &  &  &  & CUDA (N1) & 4.45e-05 & 5.20e-06 & 2.18e-05 & 1.15e-05 & 1.74e-05 \\
   &  &  &  &  & CUDA (N2) & 4.41e-05 & 7.91e-06 & 3.12e-05 & 1.60e-05 & 2.56e-05 \\
   &  &  &  &  & CUDA (N3) & 3.55e-05 & 6.45e-06 & 2.56e-05 & 1.22e-05 & 2.10e-05 \\
   &  &  &  &  & CUDA (N4) & 2.34e-05 & 4.64e-06 & 1.85e-05 & 8.58e-06 & 1.50e-05 \\
\bottomrule
\end{longtable}
}

{\scriptsize
\setlength{\tabcolsep}{2pt}
\setlength{\LTcapwidth}{\textwidth}
\LTleft=0pt plus 1fil minus 1fill
\LTright=0pt plus 1fil minus 1fill
\begin{longtable}{@{\extracolsep{\fill}}l | rrrr | l | rrrrr}
  \caption{Per-instance callback times (seconds) for the COPS benchmark suite.}
  \label{tab:results:cops} \\
  \toprule
  \textbf{instance} & \textbf{\(n\)} & \textbf{\(m\)} & \textbf{nnzj} & \textbf{nnzh} & \textbf{backend} & \textbf{obj} & \textbf{cons} & \textbf{grad} & \textbf{jac} & \textbf{hess} \\
  \midrule
  \endfirsthead
  \caption[]{(continued)} \\
  \toprule
  \textbf{instance} & \textbf{\(n\)} & \textbf{\(m\)} & \textbf{nnzj} & \textbf{nnzh} & \textbf{backend} & \textbf{obj} & \textbf{cons} & \textbf{grad} & \textbf{jac} & \textbf{hess} \\
  \midrule
  \endhead
  bearing-(5, 5) & 49 & 28 & 28 & 360 & CPU (C3) & \textbf{1.84e-06} & \textbf{6.00e-08} & \textbf{2.07e-06} & \textbf{3.00e-08} & \textbf{2.24e-06} \\
   &  &  &  &  & CPU-4T (C3) & 2.19e-06 & 5.80e-07 & 2.44e-06 & 6.40e-07 & 3.30e-06 \\
   &  &  &  &  & CPU-8T (C3) & 2.11e-06 & 5.40e-07 & 2.47e-06 & 6.90e-07 & 3.30e-06 \\
   &  &  &  &  & CPU-16T (C3) & 2.18e-06 & 5.90e-07 & 2.52e-06 & 6.80e-07 & 3.05e-06 \\
   &  &  &  &  & CUDA (N1) & 5.38e-05 & 2.17e-05 & 3.33e-05 & 2.51e-05 & 4.76e-05 \\
   &  &  &  &  & CUDA (N2) & 5.93e-05 & 2.80e-05 & 6.43e-05 & 3.70e-05 & 8.35e-05 \\
   &  &  &  &  & CUDA (N3) & 5.35e-05 & 2.62e-05 & 4.31e-05 & 3.59e-05 & 6.53e-05 \\
   &  &  &  &  & CUDA (N4) & 3.40e-05 & 1.79e-05 & 3.02e-05 & 2.44e-05 & 4.57e-05 \\
\bottomrule
\end{longtable}
}

{\scriptsize
\setlength{\tabcolsep}{2pt}
\setlength{\LTcapwidth}{\textwidth}
\LTleft=0pt plus 1fil minus 1fill
\LTright=0pt plus 1fil minus 1fill
\begin{longtable}{@{\extracolsep{\fill}}l | rrrr | l | rrrrr}
  \caption{Per-instance callback times (seconds) for PGLIB-OPF (polar formulation).}
  \label{tab:results:opf-polar} \\
  \toprule
  \textbf{instance} & \textbf{\(n\)} & \textbf{\(m\)} & \textbf{nnzj} & \textbf{nnzh} & \textbf{backend} & \textbf{obj} & \textbf{cons} & \textbf{grad} & \textbf{jac} & \textbf{hess} \\
  \midrule
  \endfirsthead
  \caption[]{(continued)} \\
  \toprule
  \textbf{instance} & \textbf{\(n\)} & \textbf{\(m\)} & \textbf{nnzj} & \textbf{nnzh} & \textbf{backend} & \textbf{obj} & \textbf{cons} & \textbf{grad} & \textbf{jac} & \textbf{hess} \\
  \midrule
  \endhead
  case3\_lmbd & 24 & 28 & 103 & 141 & CPU (C3) & \textbf{2.00e-08} & \textbf{2.10e-07} & \textbf{2.00e-08} & \textbf{3.40e-07} & \textbf{5.10e-07} \\
   &  &  &  &  & CPU-4T (C3) & 2.40e-07 & 4.18e-05 & 2.30e-07 & 2.18e-05 & 2.29e-05 \\
   &  &  &  &  & CPU-8T (C3) & 1.30e-07 & 3.95e-05 & 2.20e-07 & 2.17e-05 & 2.32e-05 \\
   &  &  &  &  & CPU-16T (C3) & 2.30e-07 & 4.05e-05 & 3.30e-07 & 2.19e-05 & 2.32e-05 \\
   &  &  &  &  & CUDA (N1) & 3.91e-05 & 2.11e-04 & 1.79e-05 & 1.32e-04 & 1.49e-04 \\
   &  &  &  &  & CUDA (N2) & 3.94e-05 & 2.42e-04 & 2.71e-05 & 1.94e-04 & 2.26e-04 \\
   &  &  &  &  & CUDA (N3) & 3.58e-05 & 2.72e-04 & 2.41e-05 & 1.85e-04 & 2.08e-04 \\
   &  &  &  &  & CUDA (N4) & 2.28e-05 & 1.92e-04 & 1.70e-05 & 1.37e-04 & 1.53e-04 \\
\bottomrule
\end{longtable}
}

\vspace{0.1cm}
\begin{flushright}
	\scriptsize \framebox{\parbox{2.5in}{Government License: The
			submitted manuscript has been created by UChicago Argonne,
			LLC, Operator of Argonne National Laboratory (``Argonne").
			Argonne, a U.S. Department of Energy Office of Science
			laboratory, is operated under Contract
			No. DE-AC02-06CH11357.  The U.S. Government retains for
			itself, and others acting on its behalf, a paid-up
			nonexclusive, irrevocable worldwide license in said
			article to reproduce, prepare derivative works, distribute
			copies to the public, and perform publicly and display
			publicly, by or on behalf of the Government. The Department of Energy will provide public access to these results of federally sponsored research in accordance with the DOE Public Access Plan. http://energy.gov/downloads/doe-public-access-plan. }}
	\normalsize
\end{flushright}

\end{document}